# Chapters on Algebraic Surfaces

**Miles Reid**

### Foreword

This is a first graduate course in algebraic geometry. It aims to give the student a lift up into the subject at the research level, with lots of interesting topics taken from the classification of surfaces, and a human-oriented discussion of some of the technical foundations, but with no pretence at an exhaustive treatment. I hope that graduate students can use some of these chapters as a reader through the subject, maybe in parallel with a conventional textbook. The early chapters introduce topics that are useful throughout projective and algebraic geometry, make little demands, and lead to fun calculations. The intermediate chapters introduce elements of the technical language gradually, whereas the later chapters get into the substance of the classification of surfaces.

I have given several sets of lectures on algebraic geometry at the graduate or advanced undergraduate level over the last few years. The 1993 Park City summer school, with the challenge of taking a numerous audience through from a basic level up to the sky, was perhaps the most fun. I have thus accumulated a considerable body of notes, without ever having time to organise the material properly in book form. The format of this series allows me to present lecture notes, which are at the same time (somewhat preliminary) extracts of chapters of a more ambitious book on surfaces. A number of the chapters here are based on the notes and examples sheets handed out in class at Park City, with only minor editing. I hope that the reader will excuse their many obvious shortcomings. One of these is the absence of any pictures, a consequence of the way I have put together the notes as computer files; I would not dream of lecturing any part of this material without drawing scribbles every minute or two, and you must make up your own figures if you want the argument to reach deep down to the parts of your intellect that words cannot reach. Some of the chapters have a serious structural fault, namely the tendency to break off into an incoherent list of further topics just as they seems to be getting somewhere. Chapter 4 on singularities should of course be broken up into an introductory discussion of Du Val singularities and the various games we want to play with them, backed up by numerous exercises and worked examples, and a separate chapter giving a serious treatment of resolutions of surface singularities, with adjunction and cohomology; similar remarks apply to Chapter 3 on K3s. I would be very grateful if you would notify me of other possible improvements, and buy my subsequent masterpiece.

[1]Math Inst., Univ. of Warwick, Coventry CV4 7AL, England
**E-mail address:** `Miles@Maths.Warwick.Ac.UK`

## Thanks

Some of the material of this book has been lectured at graduate courses at Univ. of Warwick, East China Normal University, Shanghai, Univ. of Tokyo, Univ. of Utah, as well as the Park City summer school. In addition to debts of gratitude stretching back 20 years that it would be tedious to enumerate, I am particularly grateful to Koji Ohno for his work with my Tokyo lecture notes. I also thank Gavin Brown for help preparing and proof-reading the Park City notes, and Nick Shepherd-Barron for providing several arguments used in Chapter E.

## Status of this draft

This is almost the final version of notes submitted to the 1993 volume of the IAS/AMS Park City lecture notes series. (Extrapolating past experience suggests that there are probably a few mathematical mistatements, and a further 20–30 spelling and typesetting errors to be discovered; please let me know a.s.a.p. if you find some.)

## Contents

# CHAPTER 1. The cubic surface

There are several ways of treating the cubic surface $X \subset \mathbb{P}^3$ and its 27 lines, for example, in terms of elementary coordinate geometry in $\mathbb{P}^3$, or as a blowup of $\mathbb{P}^2$ in 6 points, as a blowup of $\mathbb{P}^1 \times \mathbb{P}^1$ in 5 points, etc. The difficulty with approaches in this style is that because they start by finding one or more lines of $X$, they do not give the full symmetry of the configuration of lines. The symmetry group of the configuration has order

$$51,840 = 2^7 \cdot 3^4 \cdot 5 = 27 \cdot 5! \cdot 2^4 = 27 \cdot 16 \cdot 10 \cdot 6 \cdot 2 = \cdots$$

The configuration of lines does, however, have an entirely symmetric description in terms of a certain lattice—the divisor class group or middle cohomology $\operatorname{Pic} X = H^2(X, \mathbb{Z})$. The lines $L \subset X$ are the solutions in this lattice of the equations $HL = 1$, $L^2 = -1$ where $H$ is the class of the hyperplane section.

As well as discussing the cubic surface, this section introduces and gives examples of the following notions: $-1$-curve, blowup, conic bundle, $\mathbb{P}^1$-bundle, divisor class group $\operatorname{Pic} X$, intersection numbers, linear system and rational map, del Pezzo surface, normal surface singularity and its resolution.

## Summary

The cubic surface $X_3 \subset \mathbb{P}^3$

1. is rational;
2. is isomorphic to the plane $\mathbb{P}^2$ blown up at 6 points $\Phi = \{P_1, \dots, P_6\}$;
3. is the image of a birational map $\varphi\colon \mathbb{P}^2 \to \mathbb{P}^3$ defined by the cubics through $\Phi$;
4. has a configuration of 27 lines described either in terms of geometry in $\mathbb{P}^3$, or in terms of $\mathbb{P}^2$ blown up in 6 points;
5. has an associated lattice $A(X) \cong \mathbb{Z}^7$ that can be constructed in terms of the configuration of lines, with scalar product that can be diagonalised to $\operatorname{diag}(1, -1, \dots, -1)$; also $A(X) = \operatorname{Pic} X$ and $H^\perp \subset A(X)$ is isomorphic to the root lattice $E_6$;
6. the 27 lines and roots of $E_6$; double-six and simple roots;
7. del Pezzo surfaces in general

## 1.1. Main tricks

Let $X \subset \mathbb{P}^3$ be a nonsingular cubic surface. I'm interested in the lines and the triangles of $X$. Here a *triangle* is a set of 3 distinct coplanar lines $L_1, L_2, L_3 \subset X$, so that $L_1 + L_2 + L_3 = X \cap H$ is a hyperplane section. The following 5 statements allow me to find and organise all the lines of $X$.

0. $X$ contains at least one line $L$. (This is nontrivial, see UAG or Shafarevich, Basic Algebraic Geometry, Chapter I.)
1. Any two intersecting lines determine a triangle. (Obvious.)
2. If $L_1, L_2, L_3$ is a triangle and $M$ a fourth line of $X$, then $M$ meets exactly one of $L_1, L_2, L_3$. (Almost obvious.)
3. Main trick: there are exactly 10 lines meeting $L$, falling into 5 coplanar pairs $L_i, L_i'$; the pairs are disjoint, that is, $(L_i \cup L_i') \cap (L_j \cup L_j') = \emptyset$ for $i \neq j$. (See [**UAG**] or Shafarevich [**Sh**], Chapter IV or [**Beauville**], IV.15.)
4. In particular, there exist two disjoint lines $L$ and $M$.

## 1.2. All the lines of $X$

Here I find and give names to all the lines of $X$. Fix disjoint lines $L$ and $M$. By 1.1.3, $L$ takes part in exactly 5 triangles $L, L_i, L'_i$ for $i = 1, \dots, 5$, and $M$ meets exactly one of $L, L_i, L'_i$ for each $i$. By renumbering, let these be $L_1, \dots, L_5$. Then the 10 lines meeting $M$ are $L''_1, \dots, L''_5$ as in the figure.

**Claim.** *There are* 10 *further lines* $L_{klm}$ *which meet* $L_k, L_l, L_m$ *and not* $L_i, L_j$, *where* $\{i, j, k, l, m\} = \{1, 2, 3, 4, 5\}$.

**Proof.** Notice that if $i \neq j$ then $L'_i$ does not meet $M$ nor $L_j$, so it must intersect the third line of that triangle, that is, $L'_i$ meets $L''_j$ for $i \neq j$.

Choose some $L'_i$, for clarity $L'_1$, and consider the 5 triangles involving $L'_1$. These are $L, L_1, L'_1$, and 4 others determined by the intersecting pairs $L'_1, L''_j$ for $j = 2, \dots, 5$. Consider for example the triangle $L'_1, L''_2, N$. Then by 1.1.3, $N$ is disjoint from the line pairs of the other 4 triangles, in particular $N$ does not intersect $L''_3, L''_4, L''_5$. Therefore $N$ intersects $L_3, L_4$ and $L_5$.

It is easy to see that the lines $L, M, L_i, L'_i, L''_i$ and $L_{klm}$ exhaust all the lines of $X$.

## 1.3. The lattice $A(X)$

I define a lattice $A(X)$, a kind of miniature version of the divisor class group of any nonsingular variety. $A(X)$ is defined as an Abelian group with generators and relations. The generators are the 27 lines of $X$; the relations are simply "triangle = constant". More formally, $A$ is the free Abelian group on the 27 lines, modulo the set of relations $L + L' + L'' = M + M' + M''$ whenever $L, L', L''$ and $M, M', M''$ are triangles.

**Proposition.** $A = \mathbb{Z}^7$*; a basis is* $L_1, \dots, L_4, L'_5, L''_5, L_5$

**Proof.** The triangles containing $L$ are $L + L_i + L'_i$ for $i = 1, \dots, 5$. Thus I get relation $L + L_i + L'_i = L + L_5 + L'_5$ for $i = 1, \dots, 5$, so that

$$L'_i = L_5 + L'_5 - L_i. \tag{1}$$

Arguing in the same way on the triangles containing $M$ gives $L''_i = L_5 + L''_5 - L_i$.

I showed in 1.2 that $L'_i + L''_j + L_{klm}$ is a triangle when $\{i, j, k, l, m\}$ is a permutation of $\{1, 2, 3, 4, 5\}$, and hence

$$L_{klm} = L + L_i - L''_j. \tag{2}$$

Finally, $L_1 + L_{123} + L_{145}$ is another triangle, so that

$$L + L_1 + L'_1 = L_1 + L_{123} + L_{145} = L_1 + 2L + L_4 - L''_5 + L_2 - L''_3,$$

and therefore

$$L = L'_1 - L_4 + L''_5 - L_2 + L''_3 = 2(L_5 + L'_5 + L''_5) - L_1 - L_2 - L_3 - L_4 - L'_5. \tag{3}$$

Therefore, all the lines, hence everything in $A(X)$ can be written as integral combinations of the seven classes $L_1, \dots, L_4, L'_5, L''_5, L_5$. I prove that these elements are linearly independent in $A(X)$ in the next section by introducing a scalar product on $A(X)$.

## 1.4. The scalar product

I claim that there is a scalar product $A(X) \times A(X) \to \mathbb{Z}$ such that:

1. Two distinct lines $L$, $L'$ have $LL' = 0$ or $1$ according as they are disjoint or intersect.
2. $L^2 = -1$ for any line $L$.
3. $L(M + M' + M'') = 1$ for any line $L$ and triangle $M, M', M''$.

Indeed, formally, I can set $LL'$ and $L^2$ to be anything I wish, provided that my choice is compatible with the equivalence relation defining $A$. If $M, M', M''$ is a triangle and $L$ a line distinct from $M, M', M''$ then we know that $L$ meets exactly one of $M, M', M''$, so that (1) implies (3). If on the other hand $L = M$ then $LM' = LM'' = 1$ and $L^2 = -1$ again implies (3). This proves that the scalar product is well defined.

**Proposition.** *The scalar products of the* 7 *elements*

$$e_0 = L_5 + L_5' + L_5'', \quad e_1 = L_1, \dots, e_4 = L_4, \quad e_5 = L_5' \quad \textit{and} \quad e_6 = L_5''$$

*is given by* $e_0^2 = 1$, $e_i^2 = -1$ *for* $i = 1, \dots, 6$ *and* $e_i e_j = 0$.

*In particular,* $e_0, \dots, e_6$ *form a basis of* $A(X)$.

**Proof.** $e_1 = L_1, \dots, e_4 = L_4$, $e_5 = L_5'$ and $e_6 = L_5''$ are 6 disjoint lines, and the first 4 are also disjoint from $L_5$, so that most of the multiplication table follows at once by construction of the scalar product. The only assertions still to prove concerns $e_0$. Thus

$$e_0 e_5 = e_0 L_5' = (L_5 + L_5' + L_5'')L_5' = 1 - 1 + 0 = 0,$$

and similarly $e_0 e_6 = e_0 L_5'' = 0$. Finally

$$e_0^2 = e_0(L_5 + L_5' + L_5'') = e_0 L_5 = (L_5 + L_5' + L_5'')L_5 = -1 + 1 + 1 = 1. \quad \text{Q.E.D.}$$

## 1.5. Symmetric treatment of the lines

**Claim.** *Write* $h$ *for the class of a triangle. Then* $h^2 = 3$ *and* $hx \equiv x^2 \bmod 2$ *for every* $x \in A(X)$.

**Proof.** If $L + L' + L''$ is any triangle then $h = L + L' + L''$ and $h^2 = h(L + L' + L'') = 1 + 1 + 1 = 3$.

For the second part, note that for any lattice, $x^2 \bmod 2$ is a *linear* function $A \to \mathbb{F}_2$, because $(x + y)^2 = x^2 + 2xy + y^2 \equiv x^2 + y^2 \bmod 2$. Therefore it is enough to check that $hx \equiv x^2 \bmod 2$ holds for any set of generators of $A$. But $A(X)$ is generated by lines $L$, and we know $hL = 1$ and $L^2 = -1$.

To conclude: a cubic surface $X$ has an associated lattice $A(X) \cong \mathbb{Z}^7$, with a scalar product that can be diagonalised to $\operatorname{diag}(1, -1, \dots, -1)$, and an element $h \in A(X)$ such that $h^2 = 3$ and $hx \equiv x^2 \bmod 2$ for every $x \in A$. It is an easy result of lattice theory that the pair $A(X)$ and $H$ is uniquely determined up to isomorphism by the stated properties.

Inside $A(X)$, the lines are the solutions of the equations $hL = 1$ and $L^2 = -1$. This is a completely symmetric description of the configuration of lines.

**Exercise**

In the basis $e_0, \dots, e_6$ of 1.3–4, show that

$$h = 3(L_5 + L_5' + L_5'') - L_1 - L_2 - L_3 - L_4 - L_5' - L_5'' = 3e_0 - e_1 - \cdots - e_6.$$

and that the 27 classes of lines are

6 classes $e_i$ for $i = 1, \dots, 6$;
15 classes $e_0 - e_i - e_j$ for $i \neq j = 1, \dots, 6$;
6 classes $2e_0 - e_{i_1} - \cdots - e_{i_5}$ for 5 distinct $i_1, \dots, i_5 = 1, \dots, 6$.

What is $A(X)$? It was derived above in terms of lines, but maybe it has some more natural description in terms of $X$. One answer is that if $k = \mathbb{C}$, then

$$A(X) = H_2(X, \mathbb{Z}) = H^2(X, \mathbb{Z}).$$

It is not hard to see by topological arguments that this is a lattice of rank 7, with a scalar product defined by cup product of signature $1, -6$ and discriminant 1. Moreover, the hyperplane section $H = X \cap \mathbb{P}^2$ defines a cohomology class $h \in H^2(X, \mathbb{Z})$ with $h^2 = 3$ and $hx \equiv x^2 \bmod 2$ for every $x \in H^2(X, \mathbb{Z})$.

Taking a line $L$ into the homology class $[L]$ defines a map $A(X) \to H_2(X, \mathbb{Z})$, well defined because all hyperplane sections of $X$ are homologous. The proof that this is injective and surjective are rather similar to the arguments of the next section.

## 1.6. The divisor class group Pic $X$

I now use the above results on lines on the cubic surface to motivate general ideas of divisors on varieties. A sum of lines $\sum n_i L_i$ is a particular case of a divisor, the relation "all triangles are equal" $L + L' + L'' = M + M' + M''$ a particular case of linear equivalence of divisors, the group $A(X)$ a particular case of the divisor class group Pic $X$, and the scalar product on $A(X)$ given by $L_1 L_2 = 0, 1$ or $-1$ a particular case of intersection numbers.

I run briefly through the theory, which is in use throughout these notes; if you have not seen this material before, you should read for example [**Sh**], Chapter III for more details. Let $X$ be a normal variety (for example, nonsingular). A *prime divisor* of $X$ is an irreducible codimension 1 subvariety $\Gamma \subset X$. A *divisor* $D = \sum n_i \Gamma_i$ is a formal linear combination of prime divisors $\Gamma_i \subset X$ with coefficients $n_i \in \mathbb{Z}$. The group of divisor Div $X$ is the free Abelian group generated by prime divisors. $D$ is *effective*, written $D \geq 0$, if all $n_i \geq 0$.

Recall that if $\Gamma \subset X$ is a prime divisor then $\mathcal{O}_{X,\Gamma}$, the *local ring* at $\Gamma$, is the subring of $f \in k(X)$ regular at some point of $\Gamma$ (therefore regular on a dense open set); $\mathcal{O}_{X,\Gamma}$ is a discrete valuation ring, that is, there is a valuation $v_\Gamma \colon k(X) \setminus 0 \to \mathbb{Z}$ such that $\mathcal{O}_{X,\Gamma} = \{ f \mid v_\Gamma(f) \geq 0 \} \cup \{0\}$. Here for $0 \neq f \in k(X)$, we say that

$v_\Gamma(f) > 0$ if and only if $f$ has a zero along $\Gamma$ of order $v_\Gamma(f)$;
$v_\Gamma(f) < 0$ if and only if $f$ has a pole along $\Gamma$ of order $-v_\Gamma(f)$;
$v_\Gamma(f) = 0$ if and only if $f$ and $f^{-1}$ are regular along $\Gamma$

If $0 \neq f \in k(X)$, the expression $\operatorname{div} f = \sum_\Gamma v_\Gamma(f)\Gamma$ defines the divisor of $f$; the sum runs over all prime divisors, but is finite, since both $f$ and $f^{-1}$ are regular on dense open sets of $X$. A *principal divisor* is a divisor of the form $\operatorname{div} f$. Two divisors $D$ and $D'$ are *linearly equivalent*, written $D \overset{\text{lin}}{\sim} D'$ or $D \sim D'$ if $D - D'$ is a principal divisor, that is,

$$D - D' = \operatorname{div} f \qquad \text{for some } 0 \neq f \in k(X).$$

If $X$ is a nonsingular variety, its *divisor class group* $\operatorname{Pic} X$ is the group of divisors modulo linear equivalence, $\operatorname{Pic} X = \operatorname{Div} X/\sim$.

Beware that there are two or three other equivalence relations on divisors (and many others on algebraic cycles) in common use. To get the idea of linear equivalence, suppose that $D$, $D'$ are effective with no common components and $D - D' = \operatorname{div} f$. This means that $f$ has zeros on $D$ and poles on $D'$. I can view $f$ as a rational map $f\colon X \dashrightarrow \mathbb{P}^1$, and $D = f^{-1}(0)$, $D' = f^{-1}(\infty)$ (the locus of indeterminacy of $f$ is a subset of $X$ of codimension $\geq 2$). Thus $D_t = f^{-1}(t)$ is a divisor that moves with $t \in \mathbb{P}^1$ from $D$ at $t = 0$ to $D'$ at $t = \infty$.

**Example.** Let $X \subset \mathbb{P}^3$ be a nonsingular cubic surface and $L+L'+L''$ and $M+M'+M''$ two triangles cut out on $X$ by the planes $A = 0$ and $B = 0$ where $A$ and $B$ are linear forms on $\mathbb{P}^3$. Then $A/B \in k(X)$ and $\operatorname{div}(A/B) = L+L'+L''-M-M'-M''$, so that any two triangles are linearly equivalent. Since all the defining relations for $A$ were of this form, this implies that there is a well-defined map $\alpha\colon A(X) \to \operatorname{Pic} X$.

I will prove in the next two sections that $\alpha$ is injective and surjective.

## 1.7. Intersection numbers

Here are some basic facts about divisors on a projective nonsingular surface $Y$, which are discussed in more detail in Chapter A below. One can define intersection numbers of divisors $D_1D_2$ such that

1. $D_1D_2$ is bilinear in each factor and symmetric;
2. $D_1D_2$ only depends on $D_1, D_2$ up to linear equivalence: that is,
$$\mathcal{O}_Y(D_1) \cong \mathcal{O}_Y(D_1') \implies D_1D_2 = D_1'D_2;$$
3. if $D_1, D_2 \geq 0$ and have no common components then $D_1D_2 = \sum_P (D_1D_2)_P$, where the sum runs over all $P \in D_1 \cap D_2$ and
$$(D_1D_2)_P = \dim_k \mathcal{O}_{Y,P}/(I_{D_1} + I_{D_2})\mathcal{O}_{Y,P} = \dim_k \mathcal{O}_{Y,P}/(f_1, f_2),$$
where $D_1, D_2$ are locally defined by $f_1, f_2$;
4. if $C$ is an irreducible curve then $CD = \deg_C \mathcal{O}_C(D)$.

The properties (1–3) uniquely define the intersection number as a bilinear pairing $\operatorname{Div} Y \times \operatorname{Div} Y \to \mathbb{Z}$ (see, for example, Hartshorne [**H1**], Chapter V). The point of (2) is that $D_1D_2$ is well defined on $\operatorname{Pic} Y$, not just on $\operatorname{Div} Y$. The bilinear pairing I gave on $A(S)$ satisfies (1–3) so is the intersection pairing on $\operatorname{Pic} S$ (under the map $\alpha$ defined at the end of 1.6).

Notice the useful fact that for irreducible curves $C, C'$, the only way $CC' < 0$ is possible is for $C = C'$. The selfintersection $C^2$ can be interpreted as the degree of the normal bundle $N_{C/Y}$.

## 1.8. Conic bundles and the cubic surface

Let $L \subset S$ be a line and $\varphi = \varphi_L \colon S \dashrightarrow \mathbb{P}^1$ be the projection away from $L$. I can realise this map as follows. Take $M = \mathbb{P}^1 \subset \mathbb{P}^3$ disjoint from $L$. If $P \in S \setminus L$ then $P$ and $L$ span a unique plane $\Pi = \mathbb{P}^2 \subset \mathbb{P}^3$. Say $\Pi \cap M = \pi \in \mathbb{P}^1$ and define $\varphi(P) = \pi$. $\varphi$ is not *a priori* defined at points in $L$ but in fact I have

**Lemma.** *$\varphi\colon S \setminus L \to \mathbb{P}^1$ extends to a morphism $\varphi\colon S \to \mathbb{P}^1$.*

I give two proofs.

**First Proof.** Let $L = (x = y = 0)$ so $S = Ax + By$ where $A$ and $B$ are forms of degree 2 in $\{x, y, z, t\}$ which have no common zero on L (by nonsingularity). So the ratio $x : y = B : A$ is everywhere well defined. This ratio defines a morphism to $\mathbb{P}^1$ and is clearly just the map described above. Q.E.D.

**Second Proof.** This proof is by intersection numbers. Let $H = L + F$ be a hyperplane section of $S$ through $L$. Moving $H$ to a linearly equivalent section I can see that $HL = 1$ and $FL = FH = 2$ so $H^2 = 3$ and $F^2 = F(H - L) = 0$. Since the $F$ are effective and distinct they must be disjoint so $\varphi$ is well defined. Q.E.D.

By construction, the fibres of $\varphi$ are the plane conics residual to $L$. By the Main Trick 1.1.3, five of these are line pairs, the remainder nonsingular conics. Also, $F^2 = 0$, as we saw in the second proof. In particular, if $F = L_1 + L_2$, then $L_1L_2 = 1$, so that $L_i^2 = -1$, which agrees with the choice made in 1.4.2.

**Remark.** $L$ is not a section of the bundle. Since $LF = 2$, it cuts a generic fibre twice.

## 1.9. Other birational models for $S$

The conic bundle structure described in 1.8 can be constructed starting from any line $L$ lying on $S$. I have two disjoint lines $L$ and $M$ lying on $S$, provided by 1.1.4, so I can construct

$$\varphi = \varphi_L \times \varphi_M \colon S \to \mathbb{P}^1 \times \mathbb{P}^1.$$

Since $\varphi_L$ is the conic bundle obtained as the linear projection away from $L$, it is clear from 1.1.3 that it takes the 5 line pairs meeting $L$ to 5 distinct points of $\mathbb{P}^1$. Thus both $\varphi_L$ and $\varphi_M$ are constant on the 5 lines $L_1, \dots, L_5$ meeting both of them, so that $\varphi = \varphi_L \times \varphi_M$ contracts these 5 lines $L_1, \dots, L_5$ to points $Q_1, \dots, Q_5$ of $\mathbb{P}^1 \times \mathbb{P}^1$.

In fact $S$ is $\mathbb{P}^1 \times \mathbb{P}^1$ blown up in five distinct points (see below for blowup).

If I identify $\mathbb{P}^1 \times \mathbb{P}^1$ with the quadric $Q \subset \mathbb{P}^3$, and compose the map $\varphi\colon S \to \mathbb{P}^1 \times \mathbb{P}^1$ with the linear projection $Q \dashrightarrow \mathbb{P}^2$ from $Q_5$ then I obtain a morphism $\psi\colon S \to \mathbb{P}^2$ that contracts $L_1, \dots, L_4, L_5', L_5''$ to 6 points of $\mathbb{P}^2$. See Ex. 1.7 for this map in coordinates.

## 1.10. Blowup

We're going to see that $S$ is isomorphic to $\mathbb{P}^2$ blown up in 6 points "in general position", with the basis of $A(S)$ corresponding to $e_0 =$ a line in $\mathbb{P}^2$ (with $e_0^2 = 1$), and the 27 classes of lines given as:

$e_i$ the inverse image of the blowup of $P_i$;
$e_0 - e_i - e_j$ the line through $P_i, P_j$;
$2e_0 - e_1 - \cdots - e_5$ the conic through $P_1, \dots, P_5$.

For this I need to say what blowup means. (Compare also [**Sh**], Chapter II, §4.)

Let $P \in S$ be a point of a nonsingular surface. Then there exists a nonsingular surface and a morphism $\sigma\colon S_1 \to S$ such that

(i) $\sigma$ is an isomorphism $S_1 \setminus \sigma^{-1}(P) \cong S \setminus P$;

(ii) $\sigma^{-1}(P) \cong \mathbb{P}^1$.

$S_1$ and $\sigma$ are in fact uniquely determined by these properties.

**Construction**

Suppose first that $S = \mathbb{A}^2$ with coordinates $x, y$, and $P = (0,0)$. Define $S_1 \subset \mathbb{A}^2 \times \mathbb{P}^1$ to be the closed graph of the rational map $(x,y) \mapsto (x:y)$. If $(u:v)$ are homogeneous coordinates of $\mathbb{P}^1$ then $S_1$ is defined by $x/u = y/v$, that is $xv = yu$.

So if $(x,y) \neq (0,0)$, the ratio $(u:v)$ is well defined, which gives (i). Also $\sigma^{-1}(0,0)$ is clearly the whole of $\mathbb{P}^1$. Finally, to see that $S_1$ is nonsingular, note that since $(u:v)$ are homogeneous coordinates of $\mathbb{P}^1$, if $u \neq 0$ then I can assume $u = 1$, so that the open set $(u \neq 0) \subset S_1$ is the surface in $\mathbb{A}^3$ with coordinates $(x,y,v)$ defined by $y = xv$.

More generally, if $P \in S$ is any nonsingular point of a surface, let $x, y$ be local coordinates at $P$; by passing to a small enough open set, I can assume that $S$ is affine, and $m_P = (x,y)$ (where $m_P$ is the maximal ideal of functions vanishing at $P$ in the coordinate ring of $S$). The same construction works. The curve $C$ that comes out of blowing up a nonsingular point on a surface satisfies $C \cong \mathbb{P}^1$ and $C^2 = -1$. I call such a curve a $-1$-*curve*.

**Castelnuovo's criterion.** *Let $Y$ be a nonsingular surface and $C$ a curve in $Y$. Then there is a morphism $\varphi\colon Y \to X$ to a nonsingular surface $X$ contracting $C$ to a point and an isomorphism outside $C$ if and only if $C$ is a $-1$-curve.*

One proof of this is given in Contraction Theorem 4.15 below.

**Remark.** A blowup is also called a monoidal transformation, sigma process, etc.; blowups play a star role in birational geometry of surfaces and in resolution of singularities. For example, if $C \subset S = \mathbb{A}^2$ is a curve which has a singularity of multiplicity $m$ at the origin, $\sigma\colon S_1 \to S$ the blowup, then $\sigma^{-1}(C) = L \cup C'$ where $C' \subset S_1$ is "less singular" than $C$.

**Exercise.** Show that $C : (x^2 = y^3) \subset \mathbb{A}^2$ is resolved by a single blowup, whereas $(x^2 = y^5)$ needs two.

## 1.11. The cubic surface as $\mathbb{P}^1$ blown up 6 times

**Theorem.** *Suppose $\Phi = \{P_1, \dots, P_6\} \subset \mathbb{P}^2$ is a set of 6 points with the following "general position" properties: (i) no 2 points coincide; (ii) no 3 collinear; (iii) not all 6 on a conic. Let $S$ be the blowup of $\mathbb{P}^2$ at the 6 points $\Phi$.*

*Then the vector space of cubic forms on $\mathbb{P}^2$ vanishing at $\Phi$ is 4-dimensional, and if $F_1, \dots, F_4$ is a basis then the rational map $\mathbb{P}^2 \dashrightarrow \mathbb{P}^3$ defined by $P \mapsto (F_1(P) : \cdots : F_4(P))$ induces an isomorphism of $S$ with a nonsingular cubic surface $S \cong S_3 \subset \mathbb{P}^3$.*

*This construction is a 2-sided inverse to the map $\psi$ described at the end of 1.9 (see also Ex. 1.7).*

The fact that $S_3(P_1, \dots, P_6)$ is 4-dimensional is proved in [**UAG**], §1. The proof of the theorem can be proved by similar arguments, but I omit it for lack of time.

## 1.12. Final remarks on cubic surface

All the topics below have been extensively studied:

### 1. Weyl group of $E_6$ and Galois theory of cubic surfaces

If $S$ is defined over an algebraically nonclosed field $K$, the 27 lines will not in general be defined over $K$, but will usually require a field extension $L/K$. Then $\mathrm{Gal}(L/K)$ acts by permuting the 27 lines, so defines a homomorphism $\varphi\colon \mathrm{Gal}(L/K) \to W(E_6)$ to the symmetry of the configuration of lines, which by what I said in 1.5 is the Weyl group $W(E_6$ of the root system in the lattice $h^\perp \subset A(X)$ (compare Ex. 1.4). For example if $K = \mathbb{R}$ and $L = \mathbb{C}$ then complex conjugation interchanges some of the lines; the lines of the real locus $S(\mathbb{R}) \subset S(\mathbb{C})$ are just the fixed lines.

The image of $\varphi$ is a measure of how complicated the arithmetic properties of $S$ are. The question of whether $S$ is rational over the given field $K$ (that is, does there exist birational map $S \dashrightarrow \mathbb{P}^2$ defined over $K$?) can be discussed in terms of $\operatorname{im}\varphi$.

### 2. Monodromy

If $\{S_t\}_{t \in T}$ is a family of cubic surfaces, say with $S_t$ nonsingular for every $t \in T_0 \subset T$ then the 27 lines of $S_t$ define a 27-to-1 cover $L \to T_0$. Moving round a closed loop in $T_0$ permutes the lines among themselves, which defines a homomorphism $\psi\colon \pi_1(T_0, t_0) \to W(E_6)$.

### 3. Singular cubic surfaces

If $S$ is a cubic surface with an ordinary double point $P \in S$ as its only singularity, it turns out that it still has 15 lines not through $P$, and 6 lines through $P$ (which should be "counted with multiplicity 2"). $S$ can still be obtained from a linear system of cubics in $\mathbb{P}^2$ through 6 points, but these must be in "special position", for example 3 collinear: then the line joining these 3 is contracted to the singularity of $S$ (see Ex. 1.12).

### 4. Cubic surfaces and del Pezzo surfaces in general

Cubic surfaces are the best known of the series of del Pezzo surfaces. A *del Pezzo surface* is a nonsingular projective surface $S$ with $-K_S$ ample. In addition to the cubic $S_3$, examples include the complete intersection of two quadrics $S_4 = Q_1 \cap Q_2 \subset \mathbb{P}^4$, the double plane $S_2 \to \mathbb{P}^2$ with quartic branch curve, and $S_6 \subset \mathbb{P}^6$, the section of $\mathbb{P}^1 \times \mathbb{P}^1 \times \mathbb{P}^1 \subset \mathbb{P}^7$ by a general hyperplane $\mathbb{P}^6$. Just as the cubic is obtained by blowing up $\mathbb{P}^2$ in 6 general points (then making the anticanonical embedding), the other del Pezzo surfaces are obtained by blowing up $\mathbb{P}^2$ in $k$ general points for some $k = 0, \dots, 8$ (except $\mathbb{P}^1 \times \mathbb{P}^1$, which is also a del Pezzo).

Del Pezzo surfaces have been studied in many contexts of geometry, singularity theory and number theory (see for example [**Manin**], Manin and Tsafsman [**M–Ts**]). They have recently acquired extra significance: in Mori theory (of minimal models, compare Chapter D), a principal roles is played by *Mori fibre spaces* $X \to Y$, that is, morphisms whose fibres $S$ are varieties with ample $-K_S$. Thus del Pezzo surfaces and 1-parameter families of del Pezzo surfaces appear in a fundamental way in the classification of surfaces and 3-folds.

## Exercises to Chapter 1

**1.** Determine the entire configuration of lines in terms of the notation $L$, $M$, $L_i$, $L'_i$, $L''_i$, $L_{klm}$. Compare [**UAG**], §7.

**2.** A *double six* is two sextuples of disjoint lines $l_1, \dots, l_6$ and $m_1, \dots, m_6$ such that

$$l_i \cap m_j = \begin{cases} \text{one point} & \text{if } i \neq j, \\ \emptyset & \text{otherwise.} \end{cases}$$

Let $l_1 = L_1, \dots, l_4 = L_4, l_5 = L'_5, l_6 = L''_5$. Find $m_1, \dots, m_6$ forming a double six.

**3.** If $l_1, \dots, l_6, m_1, \dots, m_6$ form a double six, and have classes $l_i = e_i$ and $m_i = 2e_0 -$ (the five the other $e_j$), show how to index the remaining 15 lines so that $n_{ij}$ meets $l_i, l_j, m_i, m_j$. When does $n_{ij}$ meet $n_{kl}$?

**4.** Prove that the six classes $f_1 = L_1 - L_2, f_2 = L_2 - L_3, f_3 = L_3 - L_4, f_4 = L_4 - L'_5, f_5 = L'_5 - L''_5$ and $f_6 = L_4 - L_5$ base the orthogonal complement of $H$ in $A(X)$, where $H$ is the hyperplane section (that is, the class of a triangle), and that $f_i^2 = -2$ and $f_i f_j = 0$ or $1$ for $i \neq j$. (If we define a graph with nodes $\{f_i\}$ and an edge joining $f_i$ and $f_j$ if and only if $f_i f_j = 1$ then we obtain the graph $E_6$, a Dynkin diagram which appears in many areas of math.)

**5.** This is an exercise on $\mathbb{P}^1 \times \mathbb{P}^1$. Show that $\operatorname{Pic}(\mathbb{P}^1 \times \mathbb{P}^1) = \mathbb{Z} \oplus \mathbb{Z}$. Writing an element of $\operatorname{Pic}(\mathbb{P}^1 \times \mathbb{P}^1)$ as $(a, b)$, work out an intersection formula for two curves lying on $\mathbb{P}^1 \times \mathbb{P}^1$.

**6.** Show that $\{\frac{x}{y}, \frac{z}{t}, 1\}$ is a basis of $V_{|F|}$. Write out in coordinates the map to projective space that this defines.

**7.** Suppose that $L : (x = y = 0)$, $M : (z = t = 0)$ and $L_5 : (y = t = 0)$ lie on a nonsingular cubic $X \subset \mathbb{P}^3$. Prove that

$$\alpha : X \dashrightarrow \mathbb{P}^2, \quad \text{defined by } (x, y, z, t) \mapsto (xt : yz : yt)$$

extends to a morphism. [Hint: $(x : y)$ is well defined at any point of $L$ and $(y : t)$ is well defined at any point of $L_5$.]

**8.** Find three lines $M_1, M_2, M_3$ such that

$$L_1 + L_2 + M_3 = L_2 + L_3 + M_1 = L_1 + L_3 + M_2 = 2H - L - M - M_5$$

in $A(X)$. Prove that

$$(L_1 + L_2 + M_3) \cap (L_2 + L_3 + M_1) \cap (L_1 + L_3 + M_2) = \emptyset,$$

and use this to give an alternative proof of Ex. 1.7.

**9.** Show that $\alpha$ in Ex. 1.7 contracts the six lines $L_1, \dots, L_4, L_5', L_5''$ to points and is an isomorphism outside these.

**10.** Let $L$ be a line in $\mathbb{P}^2$ and $P_1, P_2, P_3$ noncollinear points. Describe the map to projective space determined by the linear system $|2L - P_1 - P_2 - P_3|$ and its inverse in terms of blowing up points and contracting $-1$-curves. This map is called the *elementary quadratic transformation.*

**11.** Let $S_m$ be a nonsingular surface of degree $m$ in $\mathbb{P}^3$. Generalise the proof of 1.8 to show that if $L = \mathbb{P}^1$ lies on $S_m$ then $L^2 = 2 - m$. Compare this to what you know on $S_2$ and $S_3$. The surface $S_4$ is an example of a K3 surface, and is studied in greater detail later.

**12.** Suppose that $\Phi = \{P_1, \dots, P_6\} \subset \mathbb{P}^2$ is a set 6 points which is in general position except that $P_1, P_2, P_3$ are collinear, on a line $L$. Prove that cubics through $\Phi$ define a birational map $\varphi \colon \mathbb{P}^2 \dashrightarrow S \subset \mathbb{P}^3$ to a singular cubic surface. It blows up the 6 points of $\Phi$, then contracts the line $L$ to a singular point $P \in S$, and is otherwise an isomorphism.

**13.** Let $S \subset \mathbb{P}^3$ be a cubic surface with an ordinary double point $P \in S$. Study the projection map $\pi \colon S \dashrightarrow \mathbb{P}^2$ from $P$, and show that it contracts 6 lines of $S$ through $P$ to a set of 6 points $\Phi \subset \mathbb{P}^2$ lying on a conic. Study the inverse map $\varphi \colon \mathbb{P}^2 \to S$ of $\pi$ as in the preceding question.

# CHAPTER 2. Rational scrolls

This chapter describes scrolls, especially the rational normal scrolls. These varieties occur throughout projective and algebraic geometry, and the student will never regret the investment of time studying them. One reason for presenting them here is that they can be discussed with very little background, and can be used to illustrate many constructions of algebraic geometry with substantial examples. I use them here to give simple examples of rational surfaces, K3 surface, elliptic surfaces, and surfaces with pencils of curves of genus 2, 3, 4, etc.

In intrinsic terms, a scroll is a $\mathbb{P}^{n-1}$-bundle $F \to C$ over a curve $C$, that is, an algebraic fibre bundle, isomorphic to $U_i \times \mathbb{P}^{n-1}$ over small Zariski open sets $U_i \subset C$, and glued by transition functions given by morphisms $U_i \cap U_j \to \mathrm{PGL}(n)$. It can be written as the projectivisation $F = \mathbb{P}(E)$ of a vector bundle $E$, and the study of general scrolls is essentially equivalent to that of vector bundles over curves. In the case $C = \mathbb{P}^1$ everything is much simpler, because the base curve $\mathbb{P}^1$ is a very explicit object, and every vector bundle is a direct sum of line bundles. Thus any question about scrolls can be solved in very explicit terms. I give some examples in the text, and many more in the exercises (see for example Ex. 2.6–2.9).

As well as discussing scrolls, this section introduces and gives examples of the following notions: linear system, free linear system, very ample linear system and projective embedding, quadrics of rank 3 and 4, determinantal variety, base locus of linear system, divisor class group $\operatorname{Pic} X$, intersection numbers, Veronese surface and cones over it, vector bundles over curves, projectivised bundle $\mathbb{P}_C(E)$, Chern numbers, Harder–Narasimhan filtration, K3 surface, elliptic surface, Weierstrass normal form, surface with a pencil of curves of genus $g = 1, 2, 3, 4, \ldots$.

## Summary

1. $\mathbb{F} = \mathbb{F}(a_1, \ldots, a_n)$ is defined as a quotient of the $(n+2)$-dimensional space $(\mathbb{A}^2 \setminus 0) \times (\mathbb{A}^n \setminus 0)$ by an action of two copies of the multiplicative group $\mathbb{G}_m$. There is a projection morphism $\pi\colon \mathbb{F} \to \mathbb{P}^1$ making $\mathbb{F}$ into a $\mathbb{P}^{n-1}$ fibre bundle.
2. Rational functions on $\mathbb{F}$ are defined as ratios of bihomogeneous polynomials. $\mathbb{F}$ has an embedding into $\mathbb{P}^N$ with the fibres $\mathbb{P}^{n-1}$ of $\pi$ mapping to $(n-1)$-planes of $\mathbb{P}^N$.
3. The divisor class group of $\mathbb{F}$ can be generated by two elements, $\operatorname{Pic} \mathbb{F} = \mathbb{Z}L \oplus \mathbb{Z}M$, where $L$ is a fibre of $\pi$ and $M$ is a relative hyperplane.
4. $\mathbb{F}$ contains negative subscrolls $B_c = \mathbb{F}(a_k, \ldots, a_n)$ corresponding to the "unstable" filtration of the integers $a_1, \ldots, a_n$, that is, when $a_1, \ldots, a_{k-1} \ge c > a_k, \ldots, a_n$.
5. The base locus of the linear system $|aL + bM|$ is determined in terms of negative subscrolls $B_c \subset \mathbb{F}$; nonsingularity conditions on the general $D \in |aL+bM|$ impose combinatorial conditions on the numerical data, and often lead to finite lists.
6. Applications of scrolls: varieties in $\mathbb{P}^n$ of small degree, del Pezzo's theorem, Castelnuovo varieties.
7. Fibred surfaces in scrolls.

## 2.1. Reminder: $\mathbb{P}^{l-1}\times\mathbb{P}^{m-1}$

I start by recalling the product of two projective spaces as treated in elementary textbooks (compare, for example, [**Sh**], Chapter I or [**UAG**], §5), which is a very useful analogy for scrolls. $\mathbb{P}^{l-1}\times\mathbb{P}^{m-1}$ can be defined as the quotient of $(\mathbb{C}^l\setminus 0)\times(\mathbb{C}^m\setminus 0)$ by the action of two copies of the multiplicative group $\mathbb{C}^*\times\mathbb{C}^*$ acting separately on the two factors:

$$(x_1,\dots,x_l;y_1,\dots,y_m)\mapsto(\lambda x_1,\dots,\lambda x_l;\mu y_1,\dots,\mu y_m)\qquad\text{for }(\lambda,\mu)\in\mathbb{C}^*\times\mathbb{C}^*.$$

Subvarieties of $\mathbb{P}^{l-1}\times\mathbb{P}^{m-1}$ are defined by *bihomogeneous* polynomials, that is, polynomials that are homogeneous separately in $x_1,\dots,x_l$ and $y_1,\dots,y_m$, and rational functions on $\mathbb{P}^{l-1}\times\mathbb{P}^{m-1}$ as quotients of two bihomogeneous polynomial of the same bidegree. Next $\mathbb{P}^{l-1}\times\mathbb{P}^{m-1}$ has the Segre embedding into usual projective space,

$$\mathbb{P}^{l-1}\times\mathbb{P}^{m-1}\hookrightarrow\mathbb{P}^{lm-1}$$

defined by bilinear forms

$$(x_1,\dots,x_l;y_1,\dots,y_m)\mapsto\big(u_{ij}=x_iy_j\big)_{\substack{i=1\dots l\\ j=1\dots m}}.$$

The image is defined by equations $\operatorname{rank}(u_{ij})\le 1$.

**Remark.** As an algebraic geometer, I should say that $\mathbb{P}^{l-1}\times\mathbb{P}^{m-1}$ is the quotient of the variety $(\mathbb{A}^l\setminus 0)\times(\mathbb{A}^m\setminus 0)$ by the action of the algebraic group $\mathbb{G}_m\times\mathbb{G}_m$. Here $\mathbb{A}^l$ is affine space, the variety corresponding to $k^l$ for a field $k$, and $\mathbb{G}_m$ is the algebraic group corresponding to the multiplicative group $k^*$. If it makes life simpler, you can replace $\mathbb{A}^l$ by $\mathbb{C}^l$ and $\mathbb{G}_m$ by $\mathbb{C}^*$ throughout.

## 2.2. Definition of $\mathbb{F}(a_1,\dots,a_n)$

Let $a_1,\dots,a_n$ be integers. I define the scroll $\mathbb{F}=\mathbb{F}(a_1,\dots,a_n)$ as the quotient of $(\mathbb{A}^2\setminus 0)\times(\mathbb{A}^n\setminus 0)$ by an action of $\mathbb{G}_m\times\mathbb{G}_m$, the product of two multiplicative groups. Write $t_1,t_2$ for coordinates on $\mathbb{A}^2$ and $x_1,\dots,x_n$ on $\mathbb{A}^n$, and $\lambda$ and $\mu$ for elements of the two factors of $\mathbb{G}_m\times\mathbb{G}_m$, that is, $(\lambda,\mu)\in\mathbb{G}_m\times\mathbb{G}_m$. The action is given as follows:

$$\begin{aligned}(\lambda,1)\colon(t_1,t_2;x_1,\dots,x_n)&\mapsto(\lambda t_1,\lambda t_2;\lambda^{-a_1}x_1,\dots,\lambda^{-a_n}x_n);\\(1,\mu)\colon(t_1,t_2;x_1,\dots,x_n)&\mapsto(t_1,t_2;\mu x_1,\dots,\mu x_n).\end{aligned}$$

Note first that the ratio $t_1:t_2$ is preserved by the action of $\mathbb{G}_m\times\mathbb{G}_m$, so that the projection to the first factor defines a morphism $\pi\colon\mathbb{F}\to\mathbb{P}^1$:

$$\begin{array}{ccc}(\mathbb{A}^2\setminus 0)\times(\mathbb{A}^n\setminus 0)&\longrightarrow&\mathbb{F}(a_1,\dots,a_n)\\ {\scriptstyle p_1}\big\downarrow&&\big\downarrow{\scriptstyle\pi}\\ (\mathbb{A}^2\setminus 0)&\longrightarrow&\mathbb{P}^1\end{array}$$

**Remark.** Compared to 2.1, I have restricted to the case $l=2$ (so that the first factor is $\mathbb{P}^1$), and generalised the group action to allow it to mix up the two factors, so that $\mathbb{F}\to\mathbb{P}^1$ can be a nontrivial $\mathbb{P}^{n-1}$ fibre bundle. The material of 2.2–7 follows exactly the remaining steps of 2.1.

## 2.3. $\mathbb{F}(a_1,\dots,a_n)$ as a fibre bundle

Above any given ratio $(t_1 : t_2) \in \mathbb{P}^1$, I can normalise to fix the values of $t_1, t_2$ with the given ratio, and this takes care of the action of the first factor $\mathbb{G}_m$; after this, the fibre of $\pi$ over $(t_1 : t_2)$ consists of the set of ratios $(x_1 : \cdots : x_n)$, forming a copy of $\mathbb{P}^{n-1}$. Thus every fibre of the projection map $\pi\colon \mathbb{F} \to \mathbb{P}^1$ is isomorphic to $\mathbb{P}^{n-1}$. As I show in Theorem 2.5 below, a good 19th century way of understanding $\mathbb{F}$ is to embed it in a projective space so that the fibres of $\pi$ are linearly embedded $(n-1)$-planes.

**Remark.** $\pi\colon \mathbb{F} \to \mathbb{P}^1$ is an example of a *fibre bundle* with fibre $\mathbb{P}^{n-1}$ and structure group the diagonal subgroup of $\mathrm{PGL}(n)$. More explicitly, on the affine piece $U_0 = (t_2 \neq 0) \subset \mathbb{P}^1$, I set $t_2 = 1$, so that $\pi^{-1}(U_0) = \mathbb{A}^1 \times \mathbb{P}^{n-1}$, with $t_1 = t_1/t_2$ the affine coordinate in the first factor and $(x_1 : \cdots : x_n)$ homogeneous coordinates in the second. Similarly, over $U_\infty = (t_1 \neq 0)$ I get $\pi^{-1}(U_\infty) = U_\infty \times \mathbb{P}^{n-1}$. The affine coordinates $t_1 = t_1/t_2$ and $t_2 = t_2/t_1$ on $U_0$ and $U_\infty$ are related on the overlap $U_0 \cap U_\infty$ in the usual way by $t_1 = 1/t_2$, and the two open sets $\pi^{-1}(U_0) = U_0 \times \mathbb{P}^{n-1}$ and $\pi^{-1}(U_\infty) = U_\infty \times \mathbb{P}^{n-1} \subset \mathbb{F}$ are glued together as follows:

$$t_1 : 1; x_1^{(0)} : \cdots : x_n^{(0)} \xrightarrow{t_1^{-1}} 1, 1/t_1; t_1^{a_1} x_1^{(0)} : \cdots : t_1^{a_n} x_n(0) = 1, t_2; x_1^{(\infty)} : \cdots : x_n^{(\infty)}.$$

(The arrow is the action of $\lambda = t_1^{-1} \in \mathbb{G}_m$). In brief, $\mathbb{F}$ is the union of two copies of $\mathbb{A}^1 \times \mathbb{P}^{n-1}$ glued together by $t_1 \mapsto t_1^{-1}$ in the first factor and $\mathrm{diag}(t_1^{a_1}, \dots, t_1^{a_n})$ in the second.

## 2.4. Bihomogenous polynomials

Rational functions on $\mathbb{F}$ are defined as ratios of bihomogenous polynomials, that is, eigenfunctions of the action of $\mathbb{G}_m \times \mathbb{G}_m$. I write down some vector spaces of bihomogeneous functions. I've already given one, the space $\langle t_1, t_2 \rangle$ of homogeneous polynomials of degree 1 in $t_1, t_2$ and degree 0 in $x_1, \dots, x_n$; the ratio $t_1 : t_2$ defines the projection $\pi\colon \mathbb{F} \to \mathbb{P}^1$.

Next, consider the functions that are linear in $x_1, \dots, x_n$. It's clear that this means that the second factor of $\mathbb{G}_m \times \mathbb{G}_m$ acts by $(1, \mu)\colon h \mapsto \mu h$. Consider polynomials that are also invariant under the action of the first factor, that is, $(\lambda, 1)\colon h \mapsto h$. Obviously, to cancel the group action $x_i \mapsto \lambda^{-a_i} x_i$, the linear term $x_i$ must be accompanied by a monomial $t_1^b t_2^c$ with $b + c = a_i$, and hence the vector space of $\mu$-invariant polynomials is based by

$$S^{a_1}(t_1, t_2)x_1, \dots, S^{a_n}(t_1, t_2)x_n,$$

where $S^a(t_1, t_2) = \{t_1^a, t_1^{a-1}t_2, \dots, t_2^a\}$ is the set of monomials of degree $a$ in $t_1, t_2$. Of course, $S^a = \emptyset$ if $a < 0$ and $S^0 = \{1\}$. The notation $S^a$ stands for the $a$th symmetric tensor power: if $\langle \Sigma \rangle$ denotes the vector space spanned by a set $\Sigma$ and $\mathrm{Sym}^a$ the $a$th symmetric tensor power of a vector space then $\langle S^{a_1}(t_1, t_2) \rangle = \mathrm{Sym}^a \langle t_1, t_2 \rangle$.

In the same way, the space of polynomials that are linear in $x_1, \dots, x_n$ and in the $\lambda^e$ eigenspace of the first factor is based by

$$S^{a_1+e}(t_1, t_2)x_1, \dots, S^{a_n+e}(t_1, t_2)x_n;$$

its dimension is $\sum_{i=1,\dots,n}^{+}(a_i + d + 1)$, where $\sum^{+}$ means you only take the sum of the terms that are $\geq 0$.

There is a similar description of the bihomogeneous polynomials of any bidegree, that is, of degree $d$ in $x_1, \dots, x_n$ and *extra degree* $e$ in the $t_i$; this is the vector space based by the monomials $t_1^{e_1} t_2^{e_2} x_1^{d_1} \cdots x_n^{d_n}$ with $\sum d_i = d$, and $e_1 + e_2 = \sum_{i=1}^{n} d_i a_i + e$. I discuss this in more detail later.

**2.5. Theorem** (Linear embeddings $\mathbb{F} \hookrightarrow \mathbb{P}^N$). *Suppose that $a_1, \dots, a_n > 0$; then the ratios between the bihomogeneous polynomials*

$$S^{a_1}(t_1, t_2)x_1, \dots, S^{a_n}(t_1, t_2)x_n \tag{2.5.1}$$

*define an embedding $\varphi\colon \mathbb{F}(a_1, \dots, a_n) \hookrightarrow \mathbb{P}^N$ (where $N = \sum_{i=1}^{n}(a_i + 1) - 1$) in such a way that every fibre $\mathbb{P}^{n-1}$ of $\pi$ goes into a linearly embedded $(n-1)$-plane.*

*The image is the subvariety of $\mathbb{P}^N$ defined by the determinantal equations*

$$\operatorname{rank}\begin{pmatrix} u_1 & u_2 & \dots & u_{a_1} & u_{a_1+2} & \dots & u_{a_1+a_2+1} & \dots & u_N \\ u_2 & u_3 & \dots & u_{a_1+1} & u_{a_1+3} & \dots & u_{a_1+a_2+2} & \dots & u_N + 1 \end{pmatrix} \leq 1.$$

The matrix here has $n$ blocks of size $2 \times a_i$; in each block, the $(1, j)$th entry for $j \geq 2$ repeats the $(2, j-1)$st entry. The meaning of the determinantal equation is that if the monomials (2.5.1) are listed as $t_1^{a_1}x_1, t_1^{a_1-1}t_2x_1, \dots, t_1^{a_2}x_2, \dots$ then the ratio $t_1 : t_2$ equals the ratio between the first and second rows, that is

$$\frac{t_1}{t_2} = \frac{t_1^{a_1}x_1}{t_1^{a_1-1}t_2x_1} = \frac{t_1^{a_1-1}t_2x_1}{t_1^{a_1-2}t_2^2x_1} = \dots = \frac{t_1^{a_2}x_2}{t_1^{a_2-1}t_2x_2} = \dots = \frac{t_1t_2^{a_n-1}x_n}{t_2^{a_n}x_n}.$$

**Proof.** The proof is very similar to that for the Segre embedding $\mathbb{P}^{l-1} \times \mathbb{P}^{m-1} \hookrightarrow \mathbb{P}^{lm-1}$ (see [**Sh**], Chapter I).

$\mathbb{F}$ is covered by a number of open sets $U_{ij} : (t_i \neq 0, x_j \neq 0)$ for $i = 1, 2$ and $j = 1, \dots n$, each isomorphic to $\mathbb{A}^n$. The piece $U_{11}$ is typical. The $n$ affine coordinates on it are $t_2/t_1$ and $t_1^{a_i - a_1}x_i/x_1$ for $i = 2, \dots, n$.

The set of monomials include $t_1^{a_1}x_1, t_1^{a_1-1}t_2x_1, t_1^{a_2}x_2, \dots, t_n^{a_n}x_n$. The first of these is nonzero everywhere on $U_{11}$, so that the ratio is well defined there. The $n$ affine coordinates of $U_{11}$ are precisely given by the ratios between $t_1^{a_1}x_1$ and the $n$ succeeding monomials, so that these embed $U_{11}$.

In the given determinantal equations, clearly if $u_1 = 1$ then all the remaining $u_i$ are determined by $u_2$ and $u_{a_1+2}, u_{a_1+a_2+3}, \dots$ corresponding to the $n$ affine coordinates of $U_{11}$. Q.E.D.

**Remarks.** (a) "Linear generation" of scrolls. The image variety $\mathbb{F}(a_1, \dots, a_n) \subset \mathbb{P}^N$ has the following description in projective geometry. Consider a fixed copy of $\mathbb{P}^1$ with homogeneous coordinates $t_1, t_2$, and $n$ embeddings $v_{a_i}\colon \mathbb{P}^1 \hookrightarrow \mathbb{P}^{a_1}$ defined by

$$(t_1 : t_2) \mapsto (t_1^{a_i} : t_1^{a_i-1}t_2 : \dots : t_2^{a_i});$$

this is the $a_i$th *Veronese embedding*, and the image $\Gamma_i = v_{a_i}(\mathbb{P}^1)$ is called the *rational normal curve* of degree $a_i$. Embed all the projective spaces $\mathbb{P}^{a_i} \hookrightarrow \mathbb{P}^N$ as linearly independent subspaces of a common $\mathbb{P}^N$ with $N = \sum_{i=1}^{n}(a_i + 1) - 1$. Now the curves $\Gamma_i$ are all identified with $\mathbb{P}^1$, so that it makes sense to take the linear span of corresponding points. This is $\mathbb{F}(a_1, \dots, a_n) \subset \mathbb{P}^N$; prove this using the determinantal equations as an exercise. (See Ex. 2.14.)

(b) If the assumption of the theorem is weakened to $a_i \geq 0$ then the ratio of the monomials (2.5.1) still defines a morphism $\varphi\colon \mathbb{F} \to \mathbb{P}^N$ that embeds each fibre of $\pi$ as an $(n-1)$-plane, but if some $a_i = 0$ then $x_i$ appears only in a single monomial $S^{a_i}(t_1, t_2)x_i = \{x_i\}$, and $\varphi(\mathbb{F})$ is a cone. The determinantal equations still make perfectly good sense, but the coordinate corresponding to $x_i$ does not appear in the equations. (See Ex. 2.15.)

(c) "Linear equations of scrolls". There is a classical description of the determinantal equations of Theorem 2.5 as $c = \operatorname{codim}(\mathbb{F} \subset \mathbb{P}^N)$ quadrics through a $(N-2)$-plane. (See Ex. 2.16.)

## 2.6. Particular cases

$\mathbb{F}(1,0)$ is a surface scroll, $\varphi\colon \mathbb{F}(1,0) \to \mathbb{P}^2$ is the blowup of a point. More generally $\varphi\colon \mathbb{F}(1,0,\dots,0) \to \mathbb{P}^n$ is the pencil of hyperplanes through a given codimension 2 linear subspace.

$\mathbb{F}(2,1) \subset \mathbb{P}^4$ is the cubic scroll.

$\mathbb{F}(1,1) \cong \mathbb{P}^1 \times \mathbb{P}^1 \cong Q \subset \mathbb{P}^3$ is the nonsingular quadric surface with a choice of projection.

$\mathbb{F}(2,0) \to Q' \subset \mathbb{P}^3$ is the standard resolution of the ordinary quadric cone (blowup).

More generally $\mathbb{F}(1,1,\underbrace{0,\dots,0}_{n-2}) \to Q_4 \subset \mathbb{P}^{n+1}$ is a resolution of a quadric of rank 4 associated with a chosen family of generators, and $\mathbb{F}(2,\underbrace{0,\dots,0}_{n-1}) \to Q_3 \subset \mathbb{P}^{n+1}$ is the standard resolution of a quadric of rank 3.

$\varphi\colon \mathbb{F}(a,0) \to \overline{\mathbb{F}}_a \subset \mathbb{P}^a$ is the blowup of the cone over a rational normal curve. The surface scroll $\mathbb{F}(a,0) \cong \mathbb{F}(a+b,b)$ for any $b \in \mathbb{Z}$ is usually called $\mathbb{F}_a$. As I discuss below, the exceptional curve of the resolution $B = \varphi^{-1}(0) \subset \mathbb{F}_a$ is a section of $\pi\colon \mathbb{F}_a \to \mathbb{P}^1$ with $B^2 = -1$.

$\mathbb{F}(0,\dots,0) \cong \mathbb{P}^1 \times \mathbb{P}^{n-1}$ with $\varphi\colon \mathbb{F}(0,\dots,0) \to \mathbb{P}^{n-1}$ the second projection. This is the only case with all $a_i \geq 0$ for which $\varphi$ is not birational.

**2.7. Lemma.** *The divisor class group of the scroll $\mathbb{F}$ is the free Abelian group*

$$\operatorname{Pic}\mathbb{F} = \mathbb{Z}L \oplus \mathbb{Z}M,$$

*with generators the following two divisor classes: $L$ is the class of a fibre of $\pi$, and $M$ the class of any monomial $t_1^b t_2^c x_i$ with $b + c = a_i$. (If all the $a_i > 0$, then $M$ is the divisor class of the hyperplane section under the embedding $\mathbb{F} \subset \mathbb{P}^N$ of Theorem 2.5.)*

**Proof.** First note that any two fibres of $\pi\colon \mathbb{F} \to \mathbb{P}^1$ are linearly equivalent: because a fraction $\alpha(t_1,t_2)/\beta(t_1,t_2)$, where $\alpha, \beta$ are linear forms, is a rational function on $\mathbb{F}$ with divisor the difference of two fibres. Thus the divisor class $L$ of a fibre is well defined.

To see $M$ more clearly, let $F_i \subset \mathbb{F}$ be the locus defined by $x_i = 0$; this is clearly the subscroll $F_i = \mathbb{F}(a_1,\dots,\widehat{a_i},\dots,a_n)$. Then the divisors $a_iL + F_i$ are all linearly equivalent, and define the divisor class $M$. Indeed, the fraction $t_1^{a_i}x_i/t_1^{a_j}x_j$ is a rational function on $\mathbb{F}$ with divisor $(a_iL+F_i)-(a_jL+F_j)$.

$L$ and $M$ are linearly independent in $\operatorname{Pic}\mathbb{F}$, since if $aL+bM \overset{\text{lin}}{\sim} 0$ then restricting to any fibre $\mathbb{P}^{n-1}$ of $\pi$ gives $b=0$, and then clearly $a=0$. Finally, I have to prove that every divisor of $\mathbb{F}$ is linearly equivalent to $aL+bM$ for some $a,b\in\mathbb{Z}$. Indeed, any irreducible codimension 1 subvariety $C\subset\mathbb{F}$ is defined by a single bihomogeneous polynomial equation in the sense of 2.4; to see this, take the inverse image in $\mathbb{A}^2\times\mathbb{A}^n$, and argue as in the case of usual projective space. If $C$ is defined by $f$ with given bidegree, it is obvious how to fix up a monomial $t_1^{a_1d+e}x_1^d$ with the same bidegree, so that $f/t_1^{a_1d+e}x_1^d$ is a rational function, and $C \overset{\text{lin}}{\sim} eL+dM$. Q.E.D.

**Remark.** In this notation, the canonical class of $\mathbb{F}(a_1,\dots,a_n)$ is given by

$$K_{\mathbb{F}} \overset{\text{lin}}{\sim} -2L-\sum F_i \overset{\text{lin}}{\sim} (-2+\sum a_i)L-nM.$$

See A.10 and Ex. A.13 for details.

## 2.8. Negative subscrolls $B_b \subset \mathbb{F}$ and the base locus of linear systems

Linear systems on general varieties are discussed below. Here I treat from an elementary point of view the linear system $|eL+dM|$ on $\mathbb{F}$, the family of divisors of $\mathbb{F}$ parametrised by the vector space of bihomogeneous polynomials of degree $d$ in the $x_i$ and extra degree $e$ in the $t_i$. I assume $d\ge 1$.

**Definition.** The *subscroll* corresponding to a subset $\{a_{i_1},\dots,a_{i_m}\}\subset\{a_1,\dots,a_n\}$ is the subvariety $\mathbb{F}(a_{i_1},\dots,a_{i_m})\subset\mathbb{F}(a_1,\dots,a_n)$ defined by the equations $x_j=0$ for $j\notin\{a_{i_1},\dots,a_{i_m}\}$. It is clearly a scroll in its own right with bihomogeneous coordinates $t_1:t_2;x_{1_1}:\dots:x_{1_m}$.

For any $b$, define the *negative subscroll* $B_b\subset\mathbb{F}(a_1,\dots,a_n)$ to be the subscroll corresponding to the subset $\{a_i \mid a_i\le b\}$. Suppose now for convenience that $a_1\le a_2\le\dots\le a_n$. Then

$$B_b=\mathbb{F}(a_1,\dots,a_m)\subset\mathbb{F}(a_1,\dots,a_n),$$

where $m$ is determined by $a_1\le\dots\le a_m\le b<a_{m+1}\le\dots\le a_n$.

As shown by (1) of the following proposition, the point of the definition is that the $B_b$ have a tendency to be base locus of linear systems.

**Theorem.** *(1) The base locus of $|-(b+1)L+M|$ is exactly $B_b$.*

*(2) Suppose that $b=a_m$. Then $B_b$ is contained with multiplicity $<\mu$ in the base locus of $|eL+dM|$ if and only if*

$$e+bd+(a_n-b)(\mu-1)\ge 0. \tag{2.8.1}$$

**Proof.** (1) An element of the linear system $|eL+M|$ is a hypersurface in $\mathbb{F}$ defined by a form $f$ which is a sum of monomials $S^{a_i+e}(t_1,t_2)x_i$ for $i=1,\dots,n$. Clearly, if $a_i\le b<-e$ then $x_i$ doesn't appear in any such monomial; therefore $f$ vanishes identically on the locus $x_{m+1}=\dots=x_n$. I told you so!

(2) The proof of (2) will make more sense after thinking about the worked examples 2.10–11 and drawing the corresponding Newton polygons.

An element of $|eL+dM|$ is defined by a bihomogeneous polynomial of bidegree $d, e$. Monomials having degree $\geq \mu$ in $x_{m+1}, \dots, x_n$ vanish with multiplicity $\mu$ along $B_b$. Thus the assertion is that there exists a monomials of bidegree $d, e$ of degree $< \mu$ in $x_{m+1}, \dots, x_n$. if and only if (2.8.1) holds. To make a monomial of extra degree $e$, the term $x_m^{d-\mu+1} x_n^{\mu-1}$ must be accompanied by a monomial of degree

$$e + a_m(d - \mu + 1) + a_n(\mu - 1),$$

which equals the left-hand side of (2.8.1), since $b = a_m$. This is obviously the highest accompanying degree of any of the allowed monomials. Q.E.D.

## 2.9. Special case: the surface scroll $\mathbb{F}_a = \mathbb{F}(0, a)$

Here $a \geq 0$, and $\mathbb{F}_a = \mathbb{F}(0, a)$ is the surface scroll. I adopt the notation $B = B_0 : (x_2 = 0)$ for the negative section (the point $(1, 0)$ on every fibre $\mathbb{P}^1$), and $A = L$ for the fibre. Then $\operatorname{Pic} \mathbb{F}_a = \mathbb{Z}A \oplus \mathbb{Z}B$, and the intersection numbers are

$$A^2 = 0, AB = 1, B^2 = -a.$$

**Proof.** From Proposition 2.5, the morphism $\varphi\colon \mathbb{F}_a \to \overline{\mathbb{F}}_a \subset \mathbb{P}^a$ defined by $|M| = |aA + B|$ is the natural resolution of the cone over the rational normal curve of degree $a$, with $B$ contracting to the vertex. The curve $(x_1 = 0) : M \overset{\text{lin}}{\sim} aA + B$ maps to a hyperplane section, and is disjoint from $B$. Hence $B(aA+B) = 0$. Since $A^2 = 0$ and $AB = 1$ are obvious, this completes the proof. Q.E.D.

In this case $|eL+dM| = |(e+ad)A+dB|$, so that the conclusion of Theorem 2.8, (1) have very simple interpretations:

$$B \text{ is fixed in } |(e + ad)A + dB| \iff e < 0 \iff B\bigl((e + ad)A + dB\bigr) < 0$$

and

$$\begin{aligned} \mu B \text{ is fixed in } |(e + ad)A + dB| &\iff e + a(\mu - 1) < 0 \\ &\iff B\bigl((e + ad)A + (d - \mu + 1)B\bigr) < 0. \end{aligned}$$

## 2.10. Worked example: The Maroni invariant of a trigonal curve

A curve $C$ (of genus $g \geq 3$, assumed to be nonhyperelliptic) is *trigonal* if it has a 3-to-1 map $C \to \mathbb{P}^1$, or equivalently, if it has a $g_3^1$, a free linear system $|D|$ with $\dim |D| = 1$ and $\deg D = 3$. Consider the canonical model $C \subset \mathbb{P}^{g-1}$. Then geometric RR says that 3 points $P_1 + P_2 + P_3$ on $C$ move in a $g_3^1$ if and only if they are collinear in $\mathbb{P}^{g-1}$. (Compare 3.2 and [**4 authors**], Chapter III, §3.) It follows at once that the canonical model of a trigonal curve is contained in a rational normal surface scroll $C \subset \mathbb{F}(a_1, a_2) \subset \mathbb{P}^{g-1}$ where $g = a_1 + a_2 + 2$, and the pencil $|A|$ on $\mathbb{F}$ cuts out the $g_3^1$ on $C$.

Order the $a_i$ as $a_1 \leq a_2$, set $a = a_2 - a_1$, and, as before, write $A$ for the fibre of $\mathbb{F}_a \to \mathbb{P}^1$ and $B \subset \mathbb{F}_a$ for the negative section. Then $\mathbb{F}(a_1, a_2) \subset \mathbb{P}^{g-1}$ is $\mathbb{F}_a$ embedded by $a_2A + B$. The canonical curve $C \subset \mathbb{F}_a$ is linearly equivalent to $\alpha A + 3B$ for some $\alpha$, and computing degree gives $\alpha = a + a_2 + 2$. By Theorem 2.8, the surface scroll $\mathbb{F}_a$ contains a nonsingular curve $C \subset \mathbb{F}_a$ linearly equivalent to $\alpha A + 3B$ if and only if $\alpha - 3a = B(\alpha A + 3B) \geq 0$. Therefore $\alpha \geq 3a$, that is, $a + a_2 + 2 \geq 3a$, which works out finally as

$$3a \leq g + 2, \quad \text{or} \quad 3a_2 \leq 2g - 2, \quad \text{or} \quad 3a_1 \geq g - 4.$$

Thus the quantity $a$ is a further invariant of trigonal curves, the *Maroni invariant.* The final inequality says in particular that $\mathbb{F}(a_1, a_2)$ can only be a cone if $g = 4$.

## 2.11. Worked example: Elliptic surfaces $X \subset \mathbb{F} = \mathbb{F}(a_1, a_2, a_3)$ and Weierstrass fibration

I give a typical application of Theorem 2.8. Let $\mathbb{F} = \mathbb{F}(a_1, a_2, a_3) \to \mathbb{P}^1$ be a 3-fold scroll with $a_1 \leq a_2 \leq a_3$ and $X \subset \mathbb{F}$ a surface meeting the general fibre of $\mathbb{F} \to \mathbb{P}^1$ in a nonsingular cubic curve. Then $X \in |(k + 2 - \sum a_i)L + 3M|$ for some $k \in \mathbb{Z}$; note that the class of $X$ is unchanged if I change $(a_1, a_2, a_3) \mapsto (a_1 - \nu, a_2 - \nu, a_3 - \nu)$ for some $\nu \in \mathbb{Z}$ and $M \mapsto M - \nu L$. To tidy up the calculation, I will assume later that $a_1 = 0$. (The class of $X$ is arranged to that the canonical class of $X$ is $K_X = kL_{|X}$, by the adjunction formula, compare A.11 below.)

I write out the equation of $X$ as a relative cubic

$$\sum_{i+j+k=3} c_{ijk}(t_1, t_2) x_1^i x_2^j x_3^k,$$

and keep track of the degrees $\deg c_{ijk} = (k + 2 - \sum a_i) + ia_1 + ja_2 + ka_3$ of the accompanying homogeneous terms in the Newton polygon:

$$\begin{array}{ccccccc} & & & k+2+2a_1-a_2-a_3 & & & \\ & & k+2+a_1-a_3 & & k+2+a_1-a_2 & & \\ & k+2+a_2-a_3 & & k+2 & & k+2-a_2+a_3 & \\ k+2-a_1+2a_2-a_3 & & k+2-a_1+a_2 & & k+2-a_1+a_3 & & k+2-a_1-a_2+2a_3 \end{array} \qquad (*)$$

In order for $X$ to be nonsingular at the general fibre of $\mathbb{F} \to \mathbb{P}^1$, its base locus is restricted by two conditions: $B_{a_2} \not\subset X$ and $2B_{a_1} \not\subset X$. These conditions just say that the general fibre of $X \to \mathbb{P}^1$ does not break up as the line $x_3 = 0$ plus a conic (so at least one of the degrees on the left-hand side of the Newton polygon is $\geq 0$, that is, $k + 2 + 2a_2 = -a_3 \geq 0$), and does not have $(1, 0, 0)$ as a double point (so that at least one of the degrees in the top corner is $\geq 0$, that is, $k + 2 - a_2 \geq 0$).

The criterion of Theorem 2.8 is of the form

$$(k + 2 - \sum a_i) + 3b + (a_n - a_m)(\mu - 1) \geq 0,$$

which works out as follows:

**Condition $B_{a_2} \not\subset X$:**

$b = a_2$, $\mu = 1$, so that $(k + 2 - \sum a_i) + 3a_2 + 0(a_3 - a_2) \geq 0$; setting $a_1 = 0$ and repeating the usual assumptions on the $a_i$ gives

$$a_2 + k + 2 \geq a_3 - a_2 \geq 0. \qquad (1)$$

**Condition $2B_{a_1} \not\subset X$:**

$b = a_1$, $\mu = 2$, so that $k + 2 - \sum a_i + 3a_1 + (a_3 - a_1) \ge 0$; setting $a_1 = 0$ gives

$$k + 2 \ge a_2 \ge 0. \tag{2}$$

For fixed value of $k$, (1) and (2) have solutions

$$a_2 = 0, \dots, k+2;$$
$$a_3 = a_2, \dots, 2a_2 + k + 2.$$

It is fun to consider the extreme cases of these inequalities. Referring to the Newton polygon, one sees that:

1. If $k + 2 < a_2 + a_3$ then the curve $B_{a_1} \subset \operatorname{Bs}|(k + 2 - \sum a_i)L + 3M|$.
2. If $k + 2 < a_3$ then every $X \in |(k + 2 - \sum a_i)L + 3M|$ contains $B_{a_1}$ and is tangent along it to $B_{a_2}$. In this case, the general $X$ has singularities on $B_{a_1}$ at the $k + 2 - a_2$ zeros of $c_{201}(t_1, t_2)$.
3. If $k+2 < a_3 - a_2$ then every $X \in |(k+2-\sum a_i)L + 3M|$ contains $B_{a_1}$ and has a flex along $B_{a_2}$.

The extreme case of the inequalities (1–2) are $a_2 = k + 2$, $a_3 = 3(k + 2)$. In this case the critical coefficients of $x_2^3$ and $x_1^2 x^3$ are homogeneous forms in $t_1, t_2$ of degree zero, that is constants, so that $X$ has equation

$$1 \cdot x_2^3 + 1 \cdot x_1^2 x_3 + \text{other terms.}$$

In other words, $X \subset \mathbb{F}$ is a nonsingular surface, with every fibre the Weierstrass normal form of an elliptic curve.

## 2.12. Final Remarks on Scrolls

### Minimal degree

Scrolls occur throughout projective algebraic geometry as projective varieties of minimal degree: del Pezzo's theorem (from the early 1880s) says that an irreducible $d$-dimensional variety $V$ spanning $\mathbb{P}^n$ has degree $\ge n - d + 1$, and if equality holds then $V$ is a linearly embedded scroll (as in Theorem 2.5), a cone over a scroll (as in Remark 2.5, (b)), or one of the sporadic cases: $\mathbb{P}^n$ itself, a quadric hypersurface $Q \subset \mathbb{P}^n$, the Veronese surface $W \subset \mathbb{P}^5$ or a cone over $W$. See [**Bertini**] or Eisenbud and Harris [**E–H**] for proofs of different vintages, or do it for yourself (Ex. 2.19).

Hypersurfaces in scrolls play a similar role in the study of curves whose degree is small compared to the genus, or surfaces of general type with $K^2$ small compared to $p_g$. Compare Ex. 2.24 or [**Harris**] or [**4 authors**], Chapter 3 (including the exercises).

### Surfaces with a pencil of curves

Many surfaces come with a natural pencil of curves of small genus; for example, Castelnuovo and Horikawa showed that surfaces with $K^2 = 3p_g - 7$ for which the canonical map $\varphi_K$ is birational (and $p_g \ge 7$) are naturally relative quartics in a scroll $\varphi_K(X) \subset \mathbb{F} \subset \mathbb{P}^{p_g - 1}$. These surfaces can therefore be studied as hypersurfaces in an explicit rational 3-fold. Compare Ex. 2.24–25.

### Scrolls over curves of genus $\geq 1$

An $n$-dimensional scroll can more generally be defined as a $\mathbb{P}^{n-1}$-bundle $F \to C$ over any curve $C$, that is, a fibre bundle with fibre $\mathbb{P}^{n-1}$ and structure group $\mathrm{PGL}(n)$. It can be proved that every scroll $F$ is the projectivisation $F = \mathbb{P}(\mathcal{E})$ of a rank $n$ vector bundle $\mathcal{E}$ over $C$ (compare Tsen's theorem in C.4 below). The assumption that the base curve is $\mathbb{P}^1$ is a major simplifying feature, which makes it possible to give a completely elementary self-contained treatment: in this case there is no effort in saying what the base curve $\mathbb{P}^1$ and the vector bundle $\mathcal{E}$ is: every vector bundle over $\mathbb{P}^1$ is a direct sum of $\mathcal{O}_{\mathbb{P}^1}(a_i)$. (This is a famous theorem, traditionally attributed to Grothendieck, Atiyah, Birkhoff, Hilbert, Gauss, Euler, Archimedes, ... ).

However, for the knowledgeable reader, essentially each part of the discussion here carries over to the more general case. This was one of the prime motivations of the theory of algebraic vector bundles over curves in the 1950s. The positive subsheaves $\bigoplus_{a_i \geq c} \mathcal{O}_{\mathbb{P}^1}(a_i)$ that correspond to the negative sections of the scrolls generalise to the Harder–Narasimhan filtration of $\mathcal{E}$, etc. See for example [**H1**], Chapter V, §2.

## Exercises to Chapter 2

**1.** Prove that $\mathbb{F}(a) \cong \mathbb{P}^1$ for any $a \in \mathbb{Z}$.

**2.** Prove that $\mathbb{F}(0,0) \cong \mathbb{P}^1 \times \mathbb{P}^1$. Generalise to $\mathbb{F}(0,\dots,0)$ (with $n$ zeros).

**3.** Show how to cover $\mathbb{F}(a_1,\dots,a_n)$ by $2n$ standard affine pieces isomorphic to $\mathbb{A}^n$, and write down the transition functions glueing any two pieces.

**4.** Prove that $\mathbb{F}(3,1) \cong \mathbb{F}(2,0)$ by comparing coordinates patches, and that $\mathbb{F}(2,0) \cong \mathbb{F}(1,-1)$.

**5.** Prove that in general

$$\mathbb{F}(a_1,\dots,a_n) \cong \mathbb{F}(a_1+b,\dots,a_n+b)$$

for any $b \in \mathbb{Z}$. [Hint: Every element of the group $\mathbb{C}^* \times \mathbb{C}^*$ can be written as a product of $(\lambda,1)$ and $(1,\mu)$ (for suitable $\lambda,\mu$), or alternatively as a product of $(\lambda,\lambda^b)$ and $(1,\mu')$ (for suitable $\lambda,\mu'$). In other words, the two actions of $\mathbb{C}^* \times \mathbb{C}^*$ only differ by an automorphism.] Deduce that the assumption $a_1 = 0$ is harmless if you're only interested in $\mathbb{F}$ up to isomorphism.

How is $\mathbb{F}(1,1) \cong \mathbb{F}(-1,-1)$ reconciled with Theorem 2.5?

**6.** Use the description of $\operatorname{Pic}\mathbb{F}$ and Theorem 2.8, (1) to prove that

$$\mathbb{F}(a_1,\dots,a_n) \cong \mathbb{F}(b_1,\dots,b_n) \iff \{a_1,\dots,a_n\} = \{b_1+c,\dots,b_n+c\} \text{ for some } c \in \mathbb{Z}.$$

**7.** Which of the following are rational functions on the named scrolls?

1. $x_1$ on $\mathbb{F}(0)$.
2. $t_1x_2/x_1$, $t_1x_1/x_2$ and $t_2x_1/x_2$ on $\mathbb{F}(1,0)$.
3. $(x_1^2+x_2)/t_1x_3$ on $\mathbb{F}(1,2,3)$.
4. $(x_1^3+t_1t_2x_1x_2^2)/(t_2x_1^2x_2+t_1^9x_3^3)$ on $\mathbb{F}(1,2,4)$.

Decide which of the following are bihomogeneous polynomials on $\mathbb{F}(0,3,5)$:

$$x_1^2 + x_2, \quad t_1^3 x_2 + x_1^2, \quad t_1^3 x_2 + x_1, \quad x_3, \quad t_1 t_2 x_3 + x_2.$$

**8.** Just as for projective space, a nonzero bihomogeneous polynomial $g$ of bidegree $(d,e) \neq (0,0)$ is not a well-defined function on a scroll $\mathbb{F}$. Prove that for $P \in \mathbb{F}$ the condition $g(P) = 0$ is well defined, so that $g$ defines a hypersurface in $\mathbb{F}$. [Hint: If you don't see this, prove it by lifting $P \in \mathbb{F}$ to different representatives $\widetilde{P} \in (\mathbb{C}^2 - 0) \times (\mathbb{C}^n - 0)$, and evaluating $g$ at these points.]

**9.** Convince yourself that any two curves of the same degree in $\mathbb{P}^2$ are linearly equivalent. Now prove the same for any two hypersurfaces $X_{d,e} \subset \mathbb{F}$ of the same bidegree on a scroll $\mathbb{F}$. [Hint: Because the ratio of their equations is a rational function.]

Let $\mathbb{F}_a = \mathbb{F}(0,a)$ be the surface scroll, and $D_1, D_2$ the sections defined by $(x_1 = 0)$ and $(x_2 = 0)$. Find all divisors linearly equivalent to $D_1$ and containing $D_2$.

**10.** Let $Q_3 \subset \mathbb{P}^3$ be a quadric of rank 3 and $\mathbb{F}_2 \to Q_3$ its natural resolution (see 2.6). Study curves in $Q_3$ in terms of $\mathbb{F}_2$.

**11.** (a) Using bihomogeneous polynomials of bidegree $(2,-4)$ on $\mathbb{F}(1,2,3)$, write down a nonsingular surface $X_{2,-4} \subset \mathbb{F}(1,2,3)$ of bidegree $(2,-4)$. How many singular fibres does the conic bundle $X_{2,-4} \to \mathbb{P}^1$ have?

(b) The same question for $X_{2,-3} \subset \mathbb{F}(1,2,3)$.

(c) Describe (in terms of its fibres) the locus $(t_1^2 x_2^3 + t_1 t_2 x_1^2 x_3 = 0) \subset \mathbb{F}(1,2,4)$.

**12.** Let $X = X_{2,e} \subset \mathbb{F}(a_1,a_2,a_3)$ be a surface of bidegree $(2,e)$. The fibres of $X \to \mathbb{P}^1$ are plane conics. Prove that, if $X$ is nonsingular, then every fibre is either a nonsingular conic or line pair. [Hint: You have to show that a double line leads to a singularity of $X$.]

Deduce a formula for the number of line pairs. [Hint: Singular conics are detected by a determinant, and you have to find its degree in $t_1, t_2$. Compare [**UAG**], proof of Proposition 7.3.]

**13.** Consider $\mathbb{F}(a_1,\dots,a_n)$ with some of the $a_i < 0$. When is the rational map of Theorem 2.5 defined? When is it the constant map? What is the dimension of its image? When is it in fact a morphism? Compare this with Theorem 2.8.

**14.** Prove the statement on linear generation of scrolls given in Remark 2.5, (a). [Hint: Write down the equation of all the rational normal curves of degree $a_i$, then the condition that corresponding points are joined up. Compare with the equations in Theorem 2.5.]

**15.** Generalise Theorem 2.5 to the case that some of the $a_i = 0$; compare Remark 2.5, (b).

**16.** Let $\Pi = \mathbb{P}^{N-2} : (x_0 = x_1 = 0) \subset \mathbb{P}^N$ be a codimension 2 linear subspace, and let $Q_1, \dots, Q_{N-n}$ be linearly independent quadrics containing $\Pi$. Prove that the intersection $\bigcap Q_i$ consists of $\Pi$ together with an $n$-dimensional variety $\overline{\mathbb{F}}$ that is the image of a scroll under a linear embedding. [Hint: The ratio $(x_0 : x_1)$ defines a rational map $\mathbb{P}^N \to \mathbb{P}^1$, whose fibres are the pencil of hyperplanes through $\Pi$. Fibre-by-fibre, each $Q_i$ defines a hyperplane. If the statement is true for some $c$ then $Q_{c+1}$ is a divisor in the scroll for $c$.]

**17.** Let $\mathbb{F}_a = \mathbb{F}(a, 0)$ be the surface scroll as in 2.9, and consider the linear system $|D| = |\alpha A + 4B|$ for suitable $\alpha \in \mathbb{Z}$. Prove that $|D|$ is very ample for $\alpha > 4a$, and find out what happens when $\alpha = 4a$. Use the Newton polygon argument to prove that for $3a \le \alpha < 4a$, the general element of $|\alpha A + 4B|$ is of the form $B + C$ where $B$ is the negative section and $C$ is a nonsingular curve having $\alpha - 3a$ transverse points of intersection with $B$.

If $\alpha$ is even, study the elliptic surface obtained as double cover of $\mathbb{F}_a$ branched in a general element of $|\alpha A + 4B|$.

**18.** Suppose that $a_1 < \cdots < a_n$. Show that an automorphism of $\mathbb{F}(a_1, \dots, a_n)$ compatible with the projection $\mathbb{F} \to \mathbb{P}^1$ is of the form

$$\begin{pmatrix} x_1 \\ \vdots \\ x_n \end{pmatrix} \mapsto M \begin{pmatrix} x_1 \\ \vdots \\ x_n \end{pmatrix},$$

where $M = \{m_{ij}\}$ is an uppertriangular matrix with entries $m_{ij}(t_1, t_2)$ homogeneous polynomials of degree $a_j - a_i$.

What happens if $a_1 \le \cdots \le a_n$, with some equalities allowed?

**19.** Prove del Pezzo's theorem: an irreducible surface spanning $\mathbb{P}^n$ and of degree $n - 1$ is either a scroll $\mathbb{F}(a_1, a_2)$ with $a_1 + a_2 = n - 1$, or a cone $\overline{\mathbb{F}}(n-1, 0)$, or $\mathbb{P}^2$ if $n = 2$ or the Veronese surface if $n = 5$.

**20.** Suppose that $a_1 \le \cdots \le a_n$ and $b_1 \le \cdots \le b_m$. Prove that there exists a surjective homomorphism

$$\mathcal{O}(a_1) \oplus \cdots \oplus \mathcal{O}(a_n) \twoheadrightarrow \mathcal{O}(b_1) \oplus \cdots \oplus \mathcal{O}(b_m)$$

if and only if $m \le n$ and for every $i$,

$$a_i \le b_i, \quad \text{and if } (a_1, \dots, a_i) \ne (b_1, \dots, b_i) \text{ then also } b_{i+1} \le a_i.$$

If $0 < a_1$, deduce necessary and sufficient conditions for $\mathbb{F}(b_1, \dots, b_{n-1})$ to be a hyperplane section of $\mathbb{F}(a_1, \dots, a_n)$.

**21.** Problem: find necessary and sufficient conditions for the existence of a short exact sequence

$$0 \to \mathcal{O}(c_1) \oplus \cdots \oplus \mathcal{O}(c_{n-m}) \to \mathcal{O}(a_1) \oplus \cdots \oplus \mathcal{O}(a_n) \to \mathcal{O}(b_1) \oplus \cdots \oplus \mathcal{O}(b_m) \to 0.$$

**22.** If $a_1 \le a_2$ and $a'_1 \le a'_2$, prove that $\mathcal{O}(a_1) \oplus \mathcal{O}(a_2)$ has a small deformation isomorphic to $\mathcal{O}(a'_1) \oplus \mathcal{O}(a'_2)$ if and only if $a_1 + a_2 = a'_1 + a'_2$ and $a_1 \le a'_1 \le a'_2 \le a_1$. [Hint: You can find small deformations of $\mathbb{F}(a_1, a_2)$ by taking it as a "special" hyperplane section of a 3-fold $\mathbb{F}(b_1, b_2, b_3)$, then varying the hyperplane.]

**23.** Problem: find necessary and sufficient conditions for $\mathcal{O}(a_1) \oplus \cdots \oplus \mathcal{O}(a_n)$ to have a small deformation isomorphic to $\mathcal{O}(b_1) \oplus \cdots \oplus \mathcal{O}(b_m)$.

**24.** By analogy with the relative cubics of 2.10, consider the scroll $\mathbb{F}(a_1, a_2, a_3)$ with $0 \le a_1 \le a_2 \le a_3$, and relative quartic surfaces $X \in |(2 - \sum a_i)A + 4M|$. (These are the surfaces of general type on the Castelnuovo–Horikawa line $K^2 = 3p_g - 7$.)

**25.** By analogy with the elliptic surface of Ex. 2.17 obtained as double cover of the surface scroll $\mathbb{F}_a = \mathbb{F}(a, 0)$ branched in a curve $|D| = |2\alpha A + 4B|$, study the linear system $|D| = |2\alpha A + 6B|$ on $\mathbb{F}_a$ and the double cover branched in $|D| = |2\alpha A + 6B|$. (These are the surfaces of general type on the Max Noether–Horikawa line $K^2 = 2p_g - 4$.)

## CHAPTER A. Curves on surfaces and intersection numbers

This chapter discusses intersection numbers, and gives practice at calculating them. At its most primitive, the idea is to count the number of points of intersection of two distinct irreducible curves $C$ and $C'$ on a surface:

$$CC' = \#\{C \cap C'\}.$$

For example, Bézout's theorem (due to Sir Isaac Newton, as everyone knows) states that two plane curves of degree $n$ and $m$ intersect in $nm$ points counted with multiplicities. I give many simple examples of calculating intersection numbers. I prove that curves on a surface that can be contracted to a point by a nonconstant morphism have negative selfintersection. I also discuss the canonical class of an algebraic surface and the formula for the genus of a curve.

Formal treatments of the material of this section can be found in [**Sh**], Chapter IV, [**H1**], Chapter V or in much more rigour, generality and detail in Fulton [**F2**–**F3**]. Fulton's introductory book [**F1**] contains an extremely concrete discussion of all the issues discussed here in the particular case of plane curves.

### Summary

1. Intersection numbers $D_1D_2$ are defined and their properties listed. Examples.
2. The selfintersection $C^2$ of a curve is defined, and is related to the degree of the normal bundle.
3. The selfintersection $C^2$ may be negative. The intersection matrix of exceptional curves is negative definite; that of components of a fibre is negative semidefinite.
4. Intersection numbers can be constructed rather concretely in terms of divisor on curves, or more abstractly in terms of coherent cohomology.
5. Bézout's theorem and the Euler characteristic in coherent cohomology.
6. The canonical class $K_X$; canonical class of $\mathbb{P}^n$ and of scrolls $\mathbb{F}$. The adjunction formula and the genus formula $2g - 2 = (K_X + C)C$.

### A.1. The formal statement

I start with the formal statement from 1.7, the ingredients of which are discussed later. $X$ is a nonsingular projective surface throughout this section.

**Theorem** (see [**H1**], Chapter V or [**Sh**], Chapter IV). *One can define an intersection pairing* $\operatorname{Div} X \times \operatorname{Div} X \to \mathbb{Z}$*, written* $(D_1, D_2) \mapsto D_1D_2$*, with the following four properties:*

1. $D_1D_2$ *is bilinear in each factor and symmetric.*
2. $D_1D_2$ *depends on* $D_1, D_2$ *only up to linear equivalence: that is,*

$$D_1 \overset{\text{lin}}{\sim} D_1' \implies D_1D_2 = D_1'D_2.$$

3. *If* $D_1, D_2 \geq 0$ *are effective divisors and have no common components then* $D_1D_2 = \sum_P (D_1D_2)_P$*, where the sum runs over all points of intersection*

*$P \in D_1 \cap D_2$, and the local intersection number $(D_1 D_2)_P$ is defined as follows:*

$$(D_1 D_2)_P = \dim_k \mathcal{O}_{X,P}/(f_1, f_2) = \dim_k \mathcal{O}_{X,P}/(I_{D_1} + I_{D_2});$$

*here $\mathcal{O}_{X,P}$ is the local ring of $P \in X$, and $D_1, D_2$ are defined near $P$ by local equations $f_1, f_2$, or (equivalently) by ideals $I_{D_1}$, $I_{D_2} \subset \mathcal{O}_{X,P}$.*

4. *If $C$ is a nonsingular curve then $C^2 = \deg_C N_{X|C}$ is the degree of the normal line bundle $N_{X|C}$ to $X$ along $C$. More generally, if $C$ is irreducible, then $C^2 = \deg_C \mathcal{O}_C(C)$, where $\mathcal{O}_C(C)$ is the invertible sheaf on $C$ (normal sheaf) discussed in the next section.*

The properties (1–3), or (1, 3–4), uniquely define the intersection number as a bilinear pairing $\operatorname{Div} X \times \operatorname{Div} X \to \mathbb{Z}$. The point of (2) is that $D_1 D_2$ is well defined on the divisor class group $\operatorname{Pic} X = \operatorname{Div} X / \overset{\text{lin}}{\sim}$, not just on $\operatorname{Div} X$. I say a few words in A.8–9 about how the pairing can be constructed and the theorem proved.

## A.2. Discussion

Notice that if $D_1, D_2 \geq 0$ are effective divisors with no common components, then each local contribution $(D_1 D_2)_P$ is the dimension of a vector space, therefore $\geq 0$. In fact $(D_1 D_2)_P > 0$ if $P \in D_1 \cap D_2$, since then $f_1, f_2 \in m_P$. Therefore $D_1 D_2 \geq 0$, and $D_1 D_2 = 0$ implies that $D_1$ and $D_2$ are disjoint.

For an irreducible curve, $C^2$ is called its *selfintersection*. If $C$ moves in a linear system then $C \overset{\text{lin}}{\sim} C'$, so that $C^2 = CC' \geq 0$. Or you can argue geometrically that if $C$ moves (more generally in an algebraic family) then the normal bundle $N_{X|C}$ has a section, so again $C^2 = \deg N_{X|C} \geq 0$.

However, as I said in 1.4, a line on a cubic surface $L \subset X$ has $L^2 = -1$, and we will see presently that $C^2 < 0$ happens in lots of interesting cases (see Theorem A.7). A curve having negative selfintersection $C^2 < 0$ may be the cause of psychological discomfort to students, so it's worth underlining the following point, which follows immediately from the above discussion.

**Proposition.** *If $D_1, D_2 \geq 0$ are effective divisors and $D_1 D_2 < 0$ then $D_1$ and $D_2$ have at least one common component $C$ that is a curve with $C^2 < 0$.*

## A.3. First example, the surface scroll $\mathbb{F}_a = \mathbb{F}(0, a)$

Here $\mathbb{F}_a = \mathbb{F}(0, a)$ is the surface scroll, with $a \geq 0$. Write $D_1 : (x_1 = 0)$ and $D_2 : (x_2 = 0)$ for the sections defined by $x_1, x_2$; then $D_1$ and $D_2$ are disjoint, since $(x_1, x_2) \neq (0, 0)$ at every point of $\mathbb{F}_a$. Therefore $D_1 D_2 = 0$.

Also, $x_1$ and $S^a(t_1, t_2) x_2$ are both bihomogeneous of the same bidegree $1, 0$, it follows that $D_1 \overset{\text{lin}}{\sim} D_2 + aL$ where $L$ is the divisor class of a fibre (compare Lemma 2.7). Obviously $D_1 L = D_2 L = 1$ (a single point transversally), so that

$$D_1^2 = D_1(D_2 + aL) = a \qquad \text{and} \qquad D_2^2 = D_2(D_1 - aL) = -a.$$

I write $A = L$ and $B = D_2$ in this example, which gives the traditional basis for $\operatorname{Pic} \mathbb{F}$, with $A^2 = 0$, $AB = 1$ and $B^2 = -a$. Notice that I can now write down the intersection number of any two curves on $\mathbb{F}_a$ by using this basis of $\operatorname{Pic} \mathbb{F}$.

By Remark 2.5, (b), $x_1$ and $S^a(t_1,t_2)x_2$ define a morphism $\varphi\colon \mathbb{F}_a \to \overline{\mathbb{F}}_a \subset \mathbb{P}^a$ such that $\varphi$ contracts $B$ to a point $Q$, and is an isomorphism $\mathbb{F}_a \setminus B \cong \overline{\mathbb{F}}_a \setminus Q$. The equations defining $\overline{\mathbb{F}}_a$ are the relations holding between the $a+1$ monomial $S^a(t_1,t_2)x_2$, that is, the determinantal equations

$$\operatorname{rank}\begin{pmatrix} u_0 & u_1 & \dots & u_{a-1} \\ u_1 & u_2 & \dots & u_a \end{pmatrix} \le 1.$$

(I've renumbered the $u_i$ slightly. $u_{-1} = x_1$ is a linear form on $\mathbb{P}^a$ that doesn't appear in the equations.) Thus $\overline{\mathbb{F}}_a$ is the cone over the rational normal curve of degree $a$. The singularity $Q \in \overline{\mathbb{F}}_a$ (given by the same equations in affine variables $u_0, \dots, u_a$) is the cyclic quotient singularity $\mathbb{C}^2/(\mathbb{Z}/a)$, where the cyclic group $\mathbb{Z}_a$ acts on $\mathbb{C}^2$ by $(x,y) \mapsto (\varepsilon x, \varepsilon y)$ for a primitive $a$th root $\varepsilon$.

Thus $\mathbb{F}_a$ is a resolution of the singularity of $\overline{\mathbb{F}}_a$, and the curve $B$ with $B^2 = -a$ is the exceptional locus.

## A.4. Intersection numbers on the $n$-fold rational scroll $\mathbb{F}(a_1, \dots, a_n)$

As in Chapter 2, I write $M$ for the divisor class of any polynomial in the vector space spanned by

$$S^{a_1}(t_1,t_2)x_1, \dots, S^{a_n}(t_1,t_2)x_n.$$

and $L$ for the class of a fibre. $L^2 = 0$, since any fibre is linearly equivalent to a disjoint fibre.

**Claim.** *$M^{n-1}L = 1$ and $M^n = \sum a_i$.*

Indeed, let $D_i : (x_i = 0)$, so that $D_i \overset{\text{lin}}{\sim} -a_iL + M$. On each fibre $\mathbb{P}^{n-1}$, the $D_i$ for $i = 1, \dots, n$ are the coordinate hyperplanes $x_i = 0$, and, as before, $\bigcap D_i = \emptyset$; therefore $D_1D_2\cdots D_n = 0$. Also, any $n-1$ of the $D_i$ (say $D_2, \dots, D_n$) intersect transversally in a curve section of $\mathbb{F} \to \mathbb{P}^1$ (the point $(1,0,\dots,0)$ in each fibre). Therefore

$$1 = LD_2\cdots D_n = L(-a_2L + M)\cdots(-a_nL + M) = LM^{n-1}$$

(because $L^2 = 0$), and

$$\begin{aligned} 0 = D_1D_2\cdots D_n &= (-a_1L + M)(-a_2L + M)\cdots(-a_nL + M) \\ &= -\sum a_iLM^{n-1} + M^n. \end{aligned}$$

**Remark.** The argument here is essentially Grothendieck's construction of Chern classes: if $E$ is a vector bundle of rank $n$ over any variety $X$, $\mathbb{P}(E)$ the corresponding $\mathbb{P}^{n-1}$ bundle, and $M$ the tautological hyperplane divisor class, then $M$ satisfies the relation $M^n - c_1(E)M^{n-1} + \cdots = 0$. In our case, this just says $M^n = \sum a_i = c_1(\mathcal{O}(a_1) \oplus \cdots \oplus \mathcal{O}(a_n))$.

## A.5. Fibred surfaces

Let $X = X_{2,e} \subset \mathbb{F}(a_1,a_2,a_3)$ be a nonsingular surface which is a divisor of bidegree $(2,e)$ in a 3-fold scroll. It can be shown (see Ex. 2.12) that every fibre of $X \to \mathbb{P}^1$ is either a nonsingular conic or a line pair. Any two distinct fibres $F$ and $F'$ of $X \to \mathbb{P}^1$ are disjoint, so that $CF = CF' = 0$ for every component $C \subset F$.

If $F = L_1 + L_2$ is a line pair, I have

$$0 = L_1F = L_1^2 + L_1L_2 = L_1^2 + 1.$$

Therefore $L_1^2 = -1$, and similarly, $L_2^2 = -1$.

More generally, if $X = X_{d,e} \subset \mathbb{F}(a_1, a_2, a_3)$ is a nonsingular surface fibred in curves of degree $d$, and one of the fibres of $X \to \mathbb{P}^1$ breaks up as a line $L$ plus a curve $C$ of degree $d-1$, then

$$0 = LF = L(L + C) = L^2 + d - 1, \qquad \text{so that } L^2 = -(d-1).$$

## A.6. Another example

I discussed in 2.10 the example of the nonsingular surface

$$X = X_{3,-3k-6} \subset \mathbb{F}(0, k+2, 3k+6)$$

turning up as the extreme case of the nonsingularity inequalities $a_2 \le k+2$, $a_3 \le 2a_2 + k + 2$. It has inflectional tangency along the curve $C : (x_2 = x_3 = 0)$ to the hyperplane $D_3 : (x_3 = 0)$. Therefore $D_3$ restricts on $X$ to $3C$.

To calculate the selfintersection $C^2$ of $C$ on $X$, note first that $D_3 + (3k+6)L \overset{\text{lin}}{\sim} D_1$ on the scroll $\mathbb{F}$, because the difference is the divisor of $f_{3k+6}(t_1, t_2)x_3/x_1$, with $f$ homogeneous of degree $3k+6$. The rational function $f_{3k+6}(t_1, t_2)x_3/x_1$ does not have $X$ as zero or pole, so it restricts to a rational function on $X$ having as zeros $3C$ plus $3k+6$ fibres $L_X$ of $X \to \mathbb{P}^1$, and as pole a divisor $D_{1,X}$ supported on the locus $(x_1 = 0) \cap X$. Therefore also on $X$ I have $3C + (3k+6)L_X \overset{\text{lin}}{\sim} D_{1,X}$. However, $D_1$ is disjoint from $C$ (since $C$ meets each fibre in $(1,0,0)$ and $D_1$ is defined by $x_1 = 0$, so that $C(3C + (3k+6)L_{|X}) = CD_{1,X} = 0$, and therefore $C^2 = -k-2$. (Ex. A.8 provides an independent check of this.)

## A.7. Contracted divisors are negative

In A.3, I discussed an example of a curve $B$ with negative selfintersection $B^2 = -a$ which is contracted to a point by a birational morphism. Also, in A.5, I gave an example of a curve $L_1$ with negative selfintersection which is a component of a fibre of a morphism $X \to C$. This is typical. Part (1) of the following theorem applies directly to resolutions of surface singularities.

**Theorem** (Negativity of exceptional locus).

1. *Let $\{\Gamma_i\}_{i=1}^k$ be a bunch of curves on a nonsingular surface $Y$ contracted to points by a projective birational morphism $f\colon Y \to X$. Then for every $(n_1, \dots, n_k) \in \mathbb{Z}^k \setminus 0$*

$$q(n_1, \dots, n_k) = \left(\sum n_i\Gamma_i\right)^2 = \sum_{i=1}^k \sum_{j=1}^k n_in_j\Gamma_i\Gamma_j < 0.$$

*In other words, the quadratic form $q$ represented by the symmetric matrix $(\Gamma_i\Gamma_j)$ is negative definite.*

2. *Let* $\{\Gamma_i\}_{i=1}^k$ *be a bunch of* $k$ *curves such that* $\Sigma = \bigcup \Gamma_i$ *is connected, and suppose that there is a surjective morphism* $f\colon X \to C$ *to a nonsingular curve* $C$ *which contracts* $\Sigma$ *to a point* $Q \in C$. *Then every* $(n_1, \dots, n_k) \in \mathbb{Z}^k \setminus 0$ *satisfies*

$$q(n_1, \dots, n_k) = \left(\sum n_i \Gamma_i\right)^2 \le 0.$$

*In other words,* $q$ *is negative semidefinite. Moreover,* $q(n_1, \dots, n_k) = 0$ *holds if and only if* $\sum n_i \Gamma_i$ *is proportional to the fibre. More precisely, if* $t \in m_Q \subset \mathcal{O}_{C,Q}$ *is a local parameter at* $Q$*, and* $g = f^*(t) \in k(X)$ *is the rational function on* $X$ *obtained as the pullback of* $t$*, then* $\sum n_i \Gamma_i$ *and the connected component of* $\operatorname{div} t$ *at* $f^{-1}Q$ *are rational multiples of one another.*

**Proof.** For simplicity, I concentrate on proving (1) under the additional assumption that the bunch of curves $\{\Gamma_i\}_{i=1}^k$ has connected union $\Sigma = \bigcup \Gamma_i$, and that $f\colon Y \to X$ contracts $\Sigma$ to a point $Q$, and is an isomorphism $Y \setminus \Sigma \to X \setminus Q$.

The proof breaks up into two steps; the first is geometric and the second an argument in the algebra of quadratic forms. (There is an alternative proof using the Hodge algebraic index theorem, see Ex. A.16.)

**Claim.** *There is an effective divisor* $D = \sum a_i \Gamma_i$ *such that* $D\Gamma_i \le 0$ *for every* $i$ *(that is,* $-D$ *is nef on the* $\Gamma_i$*), and* $D\Gamma_i < 0$ *for at least one* $\Gamma_i$.

**Proof.** Let $g \in m_{X,Q} \subset \mathcal{O}_{X,Q}$ be a function that is regular at $Q = f(\Gamma_i)$ and vanishes there, and write $C = \operatorname{div} g \subset X$ for its divisor in $X$; then the function $g$ on $Y$ (in other words, $f^*(g) \in \Gamma(Y, \mathcal{O}_Y)$) is regular and its divisor $\operatorname{div}_Y g = f^*C$ on $Y$ is of the form $D + C'$, where $D = \sum a_i \Gamma_i$ is exceptional, all $a_i > 0$, and $C'$ is the birational transform of $C$. Now $(\operatorname{div} g)\Gamma_i = 0$, so the first part of the claim follows: $C'$ has no exceptional components, so that $C'\Gamma_i \ge 0$, therefore $D\Gamma_i \le 0$.

There is a tricky bit to the last part of the claim which I explain below. To get the result with the minimum of pain, choose first any curve $C_0' \subset Y$ which is not exceptional, but intersects at least one of the $\Gamma_i$ (for example, take $C_0'$ to be a hyperplane section of $Y$). Now assume that $g \in \mathcal{O}_{X,Q}$ vanishes on $C_0 = f(C_0')$, and write as before $C = \operatorname{div} g$ and $f^*C = C' + D$. Then $C'\Gamma_i \ge 0$ for every $i$, and $> 0$ for at least on $i$ (because $C_0' \subset C'$). Then the claim follows by $(f^*C)\Gamma_i = 0$. Q.E.D.

**Remark.** In fact, it is true in general that for any $g \in m_{X,Q} \subset \mathcal{O}_{X,Q}$, the birational transform $C'$ of the curve $C = \operatorname{div} g$ intersects every connected component of $f^{-1}Q$; however, this depends on a hard theorem of Zariski: first replace $X$ by its normalisation $X' = \operatorname{Spec}_X(f_*\mathcal{O}_Y)$; the resolution $f$ factors through a morphism $f'\colon Y \to X'$ to a normal variety $X'$. Then Zariski's connectedness theorem ([**EGA III**$_1$], 4.3.1 or [**H1**], Corollary III.11.4) asserts that every fibre of $f'$ is connected (so that the birational transform of a connected divisor $C$ must intersect every component of the exceptional locus): the basic reason is that if $f'^{-1}P$ has several connected components then $(f'_*\mathcal{O}_Y)_P$ has a function that vanishes on one and not on the other.

**Proof of Claim $\implies$ Theorem A.7.** I get a contradiction from $E^2 \ge 0$, where $E = \sum b_i \Gamma_i$ with some $b_i \neq 0$. First, $D$ contains every curve $\Gamma_i$ (that is, every $a_i > 0$): this is obvious if $D\Gamma_i < 0$ (see A.2), and also if $D\Gamma_i = 0$ and $D$ contains some curve $\Gamma_j$ intersecting $\Gamma_i$; by assumption $D$ contains a curve in every connected component of $\bigcup \Gamma_i$, which gives what I said. Next, I can assume that $E > 0$. For if $E = E_1 - E_2$ with $E_1, E_2 > 0$ and $E_1, E_2$ having no common components then $0 \le E^2 = E_1^2 + E_2^2 - 2E_1E_2 \le E_1^2 + E_2^2$, so that, say $E_1^2 \ge 0$.

If some curves $\Gamma_i$ do not appear in $E$ (that is, $b_i = 0$), I can just delete them from the bunch $\Gamma_i$: by doing this, I decrease $D$ at every deleted curve, and therefore the assumption that $D\Gamma_i < 0$ for at least one $\Gamma_i$ in every connected component of $\bigcup \Gamma_i$ is preserved. The contradiction is now by induction on the number of components $\Gamma_i$ appearing in $E$ with $b_i > 0$. Write $c = \min\{b_i/a_i\}$, and set $E' = E - cD \ge 0$; then by construction, at least one $\Gamma_i$ has coefficient $b'_i = b_i - ca_i = 0$ in $E'$. Now $-D$ nef implies

$$E'^2 = (E - cD)E' \ge EE' = (E - cD)E \ge E^2 \ge 0.$$

By induction, I have a contradiction unless $E' = 0$, but then $E = cD$ and $E^2 = cED < 0$, because $D$ is negative on at least one curve of $E$. Q.E.D.

**Sketch proof of (2).** First, a brief discussion of the background: a surface fibration $f : X \to C$ is a projective morphism of a surface onto a curve with connected general fibre. A particular fibre may be reducible or have multiple components: write $f^{-1}P = \bigcup \Gamma_i$ for the set theoretic fibre and $f^*P = F = \sum a_i\Gamma_i$ for the scheme theoretic fibre over $P$.

Let $t \in m_Q \subset \mathcal{O}_{C,Q}$ be a local parameter on $C$ at $P$ and $g = f^*(t) \in k(X)$ be as in the statement of (2), and set $Z = \operatorname{div} t$. Then

$$Z = \sum a_i\Gamma_i - \text{part disjoint from } \Sigma = F - F'$$

satisfies $a_i > 0$ and

$$Z\Gamma_i = 0 \qquad \text{for all } \Gamma_i. \tag{$**$}$$

This is the geometric part of the proof of (2), with $F$ the analog of $D$ in the proof of (1).

I omit the algebraic part (but see Ex. A.17). If $E^2 > 0$ then the same kind of argument gives a contradiction. The idea for treating the case of equality $E > 0$, $E^2 = 0$ is to rework the preceding proof with only weak inequalities $\le 0$. This is just a matter of being careful about the logic.

## A.8. How to construct the pairing

There are two different strategies to the proof of Theorem A.1. One is the rather straightforward idea of defining $D_1D_2$ by the property (3). First, the local intersection multiplicity $(D_1D_2)_P$ is bilinear whenever it makes sense (that is, for effective divisor without common components at $P$). For this, you have to prove that if $f_1, f_2$ and $g \in \mathcal{O}_{X,P}$ are such that $f_1f_2$ and $g$ have no common factors, then the finite dimensional $k$ vector space $\mathcal{O}_{X,P}/(f_1f_2, g)$ can be made up as an extension of $\mathcal{O}_{X,P}/(f_1, g)$ and $\mathcal{O}_{X,P}/(f_2, g)$. More precisely, you can prove the algebraic fact that multiplication by $f_2$ takes $\mathcal{O}_{X,P}/(f_1, g)$ isomorphically to the kernel of the reduction map $\mathcal{O}_{X,P}/(f_1f_2, g) \to \mathcal{O}_{X,P}/(f_2, g)$. In other words, there is a short exact sequence

$$0 \to \mathcal{O}_{X,P}/(f_1, g) \to \mathcal{O}_{X,P}/(f_1f_2, g) \to \mathcal{O}_{X,P}/(f_2, g) \to 0.$$

The definition can then be extended to noneffective divisors by bilinearity. In order to extend the definition to divisors with common components, one must prove that every curve $C$ moves under linear equivalence to a divisor (not necessarily effective) not involving $C$. This is easy. However, in order that this definition be well defined, you must to prove (2). This boils down eventually to arguing on the number of zeros and poles of a rational function on a nonsingular curve. This is the construction given in [**Sh**], Chapter IV.

## A.9. Bézout's theorem

The other strategy for proving Theorem A.1 is more abstract and cohomological in nature, but is extremely convenient to use if you know how: you just define the pairing by

$$D_1 D_2 = \chi(\mathcal{O}_X(H-D_1-D_2)) - \chi(\mathcal{O}_X(H-D_1)) - \chi(\mathcal{O}_X(H-D_2)) + \chi(\mathcal{O}_X(H)) \quad (!)$$

where $H$ is any divisor, for example a large multiple of a hyperplane section (I discuss the Euler–Poincaré characteristic $\chi$ in B.9, viii.) This is right because $\chi(D)$ depends on $D$ as an inhomogeneous quadratic polynomial (by RR on a surface, see B.9, ix), and the formula (!) is just the associated homogeneous bilinear form obtained from the leading term (see Ex. A.15). The properties of Theorem A.1 then follow from easy cohomological manipulations (see [**H1**], Chapter V or [**Beauville**], I.4).

It is illuminating to see how this definition works in the particular case of the projective plane.

**Bézout's theorem.** *In $\mathbb{P}^2$, two curves of degree $n$ and $m$ with no common components meet in $nm$ points, counted with the multiplicity of intersection in Theorem A.1, (3).*

**Sketch proof** (see [**F1**]). Suppose that $C : (f = 0)$ and $D : (g = 0)$, where $f$ and $g$ are homogeneous polynomials of degree $n$ and $m$. Then since the polynomial ring $k[x, y, z]$ and the local rings $\mathcal{O}_{X,P}$ are UFDs and $f$ and $g$ have no common factors, it is easy to see that the sequence of sheaves

$$0 \to \mathcal{O}(-n-m) \xrightarrow{-g,f} \mathcal{O}(-n) \oplus \mathcal{O}(-m) \xrightarrow[g]{f} \mathcal{O} \to \bigoplus_P \mathcal{O}_P/(f,g) \to 0$$

is exact. Straightforward manipulations from this show that for any large $d$ we have the exact sequence of finite dimensional vector spaces

$$0 \to V_{d-n-m} \xrightarrow{-g,f} V_{d-n} \oplus V_{d-m} \xrightarrow[g]{f} V_d \to \bigoplus_P \mathcal{O}_P/(f,g) \to 0$$

where $V_d = k[x, y, z]_d$ is the space of homogeneous polynomials in $k[x, y, z]$ of degree $d$. Now we know the dimension of all the first 3 terms, since $V_d = \binom{d+2}{2}$. Therefore (applying the alternating sum formula to the preceding display, with $d \mapsto d-2$) I get the triumphant conclusion

$$\dim_k \bigoplus_P \mathcal{O}_P/(f,g) = \binom{d-n-m}{2} - \binom{d-n}{2} - \binom{d-m}{2} + \binom{d}{2} = nm. \quad \text{Q.E.D.}$$

## A.10. Enter the Hero: The canonical divisor class $K_X$

Let $X$ be a nonsingular $n$-dimensional algebraic variety or complex manifold. If $z_1, \dots, z_n$ are local algebraic or complex analytic coordinates then I want to think of $\mathrm{d}z_1 \wedge \cdots \wedge \mathrm{d}z_n$ as a *(complex) volume element.* Now, in the algebraic case, take any elements $f_1, \dots, f_n \in k(X)$ that form a (separable) transcendence basis, that is, such that $k(f_1, \dots, f_n) \subset k(X)$ is a finite algebraic extension. In the complex case, choose global meromorphic functions $f_1, \dots, f_n$ that are algebraically independent.

Also pick any $0 \neq g \in k(X)$. I write formally $s = g\mathrm{d}f_1 \wedge \cdots \wedge \mathrm{d}f_n$, and call it a *rational (meromorphic) n-form.* I'm not going to worry too much about what it *is.* The point is just that locally, I can compare it to the local volume element $\mathrm{d}z_1 \wedge \cdots \wedge \mathrm{d}z_n$ by means of the Jacobian determinant:

$$s = g\mathrm{d}f_1 \wedge \cdots \wedge \mathrm{d}f_n = J \cdot g\mathrm{d}z_1 \wedge \cdots \wedge \mathrm{d}z_n, \quad \text{where} \quad J = J\left(\frac{f_1, \dots, f_n}{z_1, \dots, z_n}\right) = \det\left|\frac{\partial f_i}{\partial z_j}\right|.$$

(The reader is invited to provide his or her own interpretation of $\partial f_i / \partial z_j$ in the algebraic case, as an exercise.) Note that by the chain rule, a different choice of local coordinates $z_1, \dots, z_n$ multiplies $J$ by an invertible function, so that the zeros and poles of $J$ are well defined. This makes it possible to determine the zeros and poles of $s$: namely, define the valuation of $s$ at a prime divisor $\Gamma$ by

$$v_\Gamma(s) := v_\Gamma(J \cdot g), \qquad \text{and set} \quad \operatorname{div} s = \sum_\Gamma v_\Gamma(s)\Gamma$$

Note that this is *not the same* as the valuation of a function: we don't know $s = g\mathrm{d}f_1 \wedge \cdots \wedge \mathrm{d}f_n$ as a function, and so we can't compare it with the unit function 1. It is a rational $n$-form, so we can compare it with the volume element $\mathrm{d}z_1 \wedge \cdots \wedge \mathrm{d}z_n$.

**Definition.** The *canonical class* of $X$ is the divisor class $K_X = \operatorname{div} s$ where $s$ is a rational $n$-form. It is a well-defined divisor class, because two different $n$-forms $s$ and $s'$ are related by $s = hs'$ with $0 \neq h \in k(X)$ a rational function, and then, plainly, $\operatorname{div} s = \operatorname{div} s' + \operatorname{div} h$.

## A.11. The adjunction formula

The canonical class $K_X$ is the great hero of classification of varieties. There are many important things to say about it, and I return to this in several of the subsequent chapters. Here I just restrict myself to the single point that it is easy to calculate if you know $X$ well enough. Here and in the exercises, we're going to see lots of cases where it can be calculated by a single trick.

**Adjunction formula.** *Let $X \subset Y$ be a nonsingular hypersurface in a nonsingular variety. Then $K_X = (K_Y + X)_X$.*

Here the restriction of divisor classes means that I first take a divisor $D$ on $Y$ with $D \overset{\text{lin}}{\sim} K_Y + X$ and such that $D$ does not contain $X$, and then intersect $D$ with $X$ to get the divisor class $D_X = (K_Y + X)_X$. The construction of the restriction will be treated more formally in 3.1, where I give the traditional proof of the adjunction formula. (See also Ex. 3.26 for a more high-brow proof in terms of Serre duality.)

## A.12. The genus formula

It is known that for a nonsingular curve $C$ of genus $g$, the canonical class $K_C$ has degree $2g-2$. It thus follow from the adjunction formula that for a nonsingular curve $C \subset X$ on a nonsingular surface,

$$2g - 2 = (K_X + C)C.$$

This formula can be checked against almost every selfintersection number calculated in this section (see Ex. A.1–8). For this, I first need a few direct calculations as a starting point.

**A.13. Proposition.** *The canonical class of projective space $\mathbb{P}^n$ is given by*

$$K_{\mathbb{P}^n} = -(n+1)H,$$

*where $H$ is the class of a hyperplane $\mathbb{P}^{n-1} \subset \mathbb{P}^n$. The canonical class of the scroll $\mathbb{F}(a_1, \dots, a_n)$ is given by*

$$K_{\mathbb{F}} = -2L - \sum D_i \overset{\text{lin}}{\sim} -\Big(\sum a_i + 2\Big)L - nM \overset{\text{lin}}{\sim} -\Big(\sum (a_i + a_1) + 2\Big)L - nD_1;$$

*here, as in 2.7, $L$ is the class of a fibre, $D_i : (x_i = 0)$, and $M$ the divisor class of $a_i L + D_i$.*

Ex. A.9–14 outline several alternative proofs of the proposition. The crudest statement and proof for $\mathbb{P}^n$ is that if $x_1, \dots, x_n$ are affine coordinates on the 0th coordinate piece $\mathbb{A}^n \subset \mathbb{P}^n$ then $s = \mathrm{d}x_1 \wedge \cdots \wedge \mathrm{d}x_n$ has a pole of order $n+1$ along the hypersurface at infinity $x_0$. This can be checked by a direct change of variables argument (see Ex. A.11).

A slightly more sophisticated approach coming from toric geometry is to write down the $n$-forms

$$s^{(i)} = (-1)^i \frac{\mathrm{d}x_0^{(i)}}{x_0^{(i)}} \wedge \cdots \wedge \begin{pmatrix} \text{omit } i\text{th} \\ \text{factor} \end{pmatrix} \wedge \cdots \wedge \frac{\mathrm{d}x_n^{(i)}}{x_n^{(i)}}$$

on the $i$th piece $\mathbb{A}^n_{(i)}$ (here $x_j^{(i)} = x_j/x_i$ are the usual affine coordinates). This visibly has simple poles along the $n$ coordinate hyperplanes of $\mathbb{A}^n_{(i)}$. However, one can check easily that all the $s^{(i)}$ are equal, so define a single rational $n$-form $s_{\text{best}}$ on $\mathbb{P}^n$ with simple poles along all $n+1$ coordinate hyperplanes (see Ex. A.12).

Implicit in this argument is the observation that the complement of the $n+1$ coordinate hyperplanes in $\mathbb{P}^n$ is a big torus $\mathbb{T} = \mathbb{G}_m^n$ (where $\mathbb{G}_m$ is the algebraic group corresponding to the multiplicative group of the field $k^*$), and that any $n$ of the $n+1$ vector fields $x_i \frac{\partial}{\partial x_i}$ (they sum to zero on $\mathbb{P}^n$ by the Euler relation) is a basis for the tangent bundle of $\mathbb{T}$, with logarithmic zeros along the coordinate hyperplanes.

Similar remarks apply to the scroll $\mathbb{F}(a_1, \dots, a_n)$ (see Ex. A.13).

## A.14. Problem

The $(n+1)$-form on $\mathbb{A}^{n+1} \setminus 0$

$$t = \frac{\mathrm{d}x_0}{x_0} \wedge \cdots \wedge \frac{\mathrm{d}x_n}{x_n}$$

has simple poles along each of the coordinate hyperplanes, and is invariant under the action of the multiplicative group $\mathbb{G}_m$. It should be possible to prove that $s_{\text{best}}$ is in an intrinsic way a contraction of $t$ against the Euler vector field

$$\sum_{i=0}^{n} x_i \frac{\partial}{\partial x_i}.$$

## Exercises to Chapter A

**1.** Let $X = X_d \subset \mathbb{P}^3$ be a nonsingular surface of degree $d$, and $L \subset X$ a line. Calculate $L^2$ in $X$. Check that this is what you know for $d = 1, 2, 3$.

**2.** $X$ and $L$ as in Ex. A.1. Prove that planes through $L$ cut out a pencil $|F|$ of curves with $F^2 = 0$. Deduce that the linear projection $\mathbb{P}^3 \to \mathbb{P}^1$ away from $L$ extends to a morphism $X \to \mathbb{P}^1$.

**3.** Let $X = X_d \subset \mathbb{P}^3$ be a nonsingular surface of degree $d$, and suppose that $X$ has a plane section $\Pi$ that decomposes as a union of two curves $\Pi \cap X = A + B$ of degrees $a$ and $b$. Prove that $(AB)_{\mathbb{P}^2} = (AB)_X$ (that is, $AB$ calculated in the plane $\Pi$ or in $X$ is the same thing). Calculate $A^2$ and $B^2$ in $X$.

**4.** Let $X = X_{3,e} \subset \mathbb{F}(a_1, a_2, a_3)$ be a nonsingular surface in the scroll of the indicated bidegree. Every fibre $F$ of $X \to \mathbb{P}^1$ is a plane cubic curve. Prove that if $F$ is a union of a line $L$ and conic $Q$ then $L^2 = Q^2 = -2$. If $F$ is a union of 3 lines $L_1, L_2, L_3$, calculate $L_i^2$.

**5.** Use Proposition A.13 and the adjunction formula A.11 to find the canonical class of a surface $X_d \subset \mathbb{P}^3$ of degree $d$. Do the same for the hypersurface in a 3-fold scroll $X_{d,e} \subset \mathbb{F}(a_1, \dots, a_3)$ of the indicated bidegree.

**6.** Let $L \subset X_d \subset \mathbb{P}^3$ be a line on a nonsingular surface, as in Ex. A.1. Check your answers to Ex. A.1 and Ex. A.5 against the genus formula A.12.

**7.** A nonsingular plane curve $A$ of degree $a$ has genus $\binom{a-1}{2}$. If $A \subset X_d$ is as in Ex. A.3, check your answers to Ex. A.3 and Ex. A.5 against the genus formula.

**8.** Consider the example $X = X_{3,-3k-6} \subset \mathbb{F}(0, k+2, 3k+6)$ treated in 2.12 and A.6. Show that $K_X = kL$ and use the genus formula to calculate the selfintersection of any section of $X \to \mathbb{P}^1$. Compare with the result of the calculation in A.6.

**9.** Consider affine coordinates $x_1, \dots, x_n$ on $\mathbb{A}^n \subset \mathbb{P}^n$ where $x_0 \neq 0$, and verify that, as stated in A.13, $s = \mathrm{d}x_1 \wedge \cdots \wedge \mathrm{d}x_n$ has a pole of order $n + 1$ along the hypersurface at infinity $x_0$. [Hint: Use coordinates $y_0, y_2, \dots, y_n$ with $y_0 = 1/x_1$ and $y_i = x_i/x_1$.]

**10.** Check that $s^{(0)} = s^{(1)}$ in the discussion after Proposition A.13.

**11.** Let $\mathbb{F}(a_1, \dots, a_n)$ be the scroll. On the $1,1$th coordinate patch $\mathbb{A}^n_{11}$ (see Ex. 2.3) $t_1 = 1, x_1 = 1$, write down the $n$-form

$$s = \mathrm{d}t_2 \wedge \mathrm{d}x_2 \wedge \cdots \wedge \mathrm{d}x_n.$$

Check that, viewed as a rational $n$-form on $\mathbb{F}$, this has a pole of order $n$ along the "horizontal" fibre at infinity $D : (x_1 = 0)$, and a pole of order $2 + \sum(a_i + a_1)$ along the "vertical" fibre at infinity $L_\infty : (t_1 = 0)$. [Hint: For the first statement, work on $\mathbb{A}^n_{12}$ given by $t_1 = 1, x_2 = 1$; then $t_1$ plays no role in the calculation, which is the same as in Ex. A.9. For the second statement, work on $\mathbb{A}^n_{21}$ given by $t_2 = 1, x_1 = 1$; the change of variable $t_2 = 1/t_1$ gives the factor $t_1^{-2}$ and the change $y_j = x_j t_1^{-a_j}$ does the rest.]

**12.** By analogy with the more sophisticated argument given after Proposition A.13, write down a rational $n$-form $s^{(ij)}$ on each of the $2n$ affine pieces of $\mathbb{F}(a_1, \dots, a_n)$ with simple poles along each coordinate hyperplane, and check that they coincide on overlaps. This defines a rational $n$-form on $\mathbb{F}$ with simple poles along every coordinate hyperplane $x_i = 0$ and $t_i = 0$, proving the statement in Proposition A.13.

**13.** $\mathbb{F}(a_1, \dots, a_n)$ is the quotient of $\mathbb{C}^2 \setminus 0 \times \mathbb{C}^n \setminus 0$ by a certain action of $\mathbb{C}^* \times \mathbb{C}^*$. Write down the two Euler relations holding between the $n+2$ vector fields

$$t_i \frac{\partial}{\partial t_i}, \quad \text{and} \quad x_j \frac{\partial}{\partial x_j}.$$

**14.** Solve Problem A.14, and generalise to the scroll $\mathbb{F}(a_1, \dots, a_n)$. Try to get your solution published. [Hint: Choose a journal that doesn't use me as a referee.]

**15.** Let $V$ be a $k$-vector space, $\varphi\colon V \to k$ a quadratic polynomial map (possibly inhomogeneous), and $h \in V$ any element. Define $\psi$ by

$$\psi(x, y) = \tfrac{1}{2}\bigl\{\varphi(h) - \varphi(h - x) - \varphi(h - y) + \varphi(h - x - y))\bigr\}.$$

Prove that $\psi(x, y)$ is independent of $h$, that it is a symmetric bilinear form, and that $\psi(x, x)$ is the leading term of $\varphi$.

**16.** Give an alternative proof of negative definiteness (Theorem A.7) using the Hodge algebraic index theorem (Corollary D.2.2.). [Hint: Assume $X$ and $Y$ are projective, and let $H \in \operatorname{Pic} Y$ be the pullback of an ample divisor on $X$. Then $H^2 > 0$ and $H$ is orthogonal to all the exceptional curves.]

**17.** Prove negative semidefiniteness for a fibre $F = \sum a_i \Gamma_i$ of a surface fibration (Theorem A.7, (2)). [Hint: In more detail, assume that $\bigcup \Gamma_i$ is connected. Then $F\Gamma_i = 0$; if $\operatorname{hcf}\{a_i\} = m$ then $F = mF'$, where $F' = \sum a_i' \Gamma_i$, and the $a_i' = a_i/m$ have no common factor. Now adapt the argument of A.7 to prove that $E = \sum b_i \Gamma_i$ has $E^2 \le 0$, with equality if and only if $E$ is an integer multiple of $F'$. Compare also [**Beauville**], Corollary VIII.4.]

**18.** Prove the same statement as in the preceding exercise using the Hodge algebraic index theorem (Corollary D.2.2.).

## CHAPTER B. Sheaves and coherent cohomology

Sheaf theory is a language to treat geometric data (functions, vector fields, etc.) on a space $X$ in terms of the same kind of data on open sets $U \subset X$. You probably know that in Molière's play *Le Bourgeois Gentilhomme*, M. Jourdain was amazed to learn from his grammar teacher that he had been speaking all of his life *in prose*! In the same way, I hope to persuade you that you have been using some of the ideas of sheaf theory without knowing it ever since your first calculus course.

It is of course out of the question to give a reasonable treatment of sheaves and cohomology in a single lecture. Instead, I give examples of some of the main classes of sheaves occurring in algebraic geometry, and discuss their role in the foundational crises of the subject through the ages. I also try to explain the definition of coherent sheaves, and to highlight the specific features of coherent cohomology which make it different from other cohomology theories used in topology, differential geometry and algebraic geometry.

If all the definitions in this section are intimidating to the younger reader, I assure you that with coherent cohomology, the whole of this type of algebraic geometry becomes a game with fixed rules and just a few standard gambits. This section concludes with a list of "Rules of the Game" for coherent cohomology, which can be taken as axioms, or read up in several references. Remember the Zariski quote[2] we heard a few nights ago: before Serre's [**FAC**], just a few *maestri* who had spent all their lives contemplating the intricacies of the black arts could say when some restriction map was surjective, and all you could do was to believe them; after [**FAC**], any idiot could write down exact sequences and deduce any number of such statements.

### Summary

1. You already know lots of sheaves. Continuous functions on $\mathbb{R}$, vector fields on a manifold
2. Definition of sheaf: sections $\Gamma(U, \mathcal{F})$ over an open and "glueing conditions"
3. The structure sheaf $\mathcal{O}_X$, intrinsic definition of a variety
4. Other sheaves on algebraic varieties, $\mathcal{O}_X(D)$, $\Omega^1_X$, $\mathcal{O}_X(K_X) = \Omega^n_X$
5. Subsheaves, stalks and quotient sheaves. Surjective must be defined on stalks, that's where cohomology comes in
6. Coherent sheaves. Homomorphisms between locally free sheaves
7. Rules of coherent cohomology
8. The Nakai–Moishezon criterion for surfaces

## B.1. You already know lots of sheaves

When we talk of functions on a space $X$, we often actually mean functions defined on some open subset $U \subset X$. For example, in elementary calculus, a function $f(x)$ of a single real variable might be defined on some interval $(a, b) \subset \mathbb{R}$. If $f$ is defined on some big interval $(A, B)$, the property that $f$ is continuous (or differentiable, or real analytic) is defined locally at every point $P \in (A, B)$, and for $P$, only depends

[2]Compare [**Parikh**], p. 138; Carol Parikh gave an interesting evening talk on Zariski's life and opinions at the Park City summer school.

on the behaviour of $f$ on any smaller interval $(a,b) \ni P$. *Voilà, M. Jourdain,* a sheaf!

In the same way, in a first course on manifolds, we learn that a vector field on a manifold $X$ is a thing that can be written $\sum u_i \frac{\partial}{\partial x_i}$ in terms of local coordinates $x_1, \dots, x_n$ on a coordinate patch $U \subset X$. It is natural in discussing manifolds to work locally over a small open sets (for example, coordinate patches), and very awkward to insist that the only vector fields we use are global vector fields on the whole of $X$. For any open subset $V \subset X$, write

$$\Gamma(V, \mathcal{T}_X) = \bigl\{\text{vector fields on } V\bigr\}$$

for the set of all vector fields on $V$. The *tangent sheaf* $\mathcal{T}_X$ is this data: all possible vector fields defined on all possible open subsets $V \subset X$. Or in other words, the assignment $V \mapsto \Gamma(V, \mathcal{T}_X)$.

## B.2. The structure sheaf $\mathcal{O}_X$ of a variety

Let $X$ be an irreducible quasiprojective variety with its Zariski topology, and $k(X)$ its field of rational functions. For $P \in X$, we know what it means for a rational function $f \in k(X)$ to be *regular* at $P$ (see for example [**UAG**], §5), and the set of rational functions that are regular at $P$ is the *local ring* $\mathcal{O}_{X,P} \subset k(X)$. For any Zariski open $U \subset X$, define

$$\Gamma(U, \mathcal{O}_X) = \bigl\{f \in k(X) \bigm| f \text{ is regular at all } P \in U\bigr\} = \bigcap_{P \in U} \mathcal{O}_{X,P} \subset k(X).$$

The *structure sheaf* of $\mathcal{O}_X$ is the assignment $U \mapsto \Gamma(U, \mathcal{O}_X)$. In this case everything is very simple, because all the $\Gamma(U, \mathcal{O}_X)$ are subrings of the fixed function field $k(X)$.

For $\emptyset \neq V \subset U$ a smaller Zariski open set, regular functions on $U$ obviously restrict to regular functions on $V$, defining a inclusion $\Gamma(U, \mathcal{O}_X) \subset \Gamma(V, \mathcal{O}_X)$. Also, since every function regular at $P$ is regular on some open neighbourhood of $P$, it follows that

$$\mathcal{O}_{X,P} = \bigcup_{U \ni P} \Gamma(U, \mathcal{O}_X) \subset k(X).$$

In sheaf theory, the inclusion $\Gamma(U, \mathcal{O}_X) \subset \Gamma(V, \mathcal{O}_X)$ is called a *restriction map* $\mathrm{Res}_{U,V}$ (or sometimes $\rho_{U,V}$), and $\mathcal{O}_{X,P}$ the *stalk* of $\mathcal{O}_X$ at $P$.

## B.3. The definition of a sheaf

In general, a *presheaf* $\mathcal{F}$ on a topological space $X$ is a way of assigning to every open subset $U \subset X$ a set $\Gamma(U, \mathcal{F})$ (or ring, or group, or object of any category; the symbol $\Gamma(U, \mathcal{F})$ is pronounced "the sections of $\mathcal{F}$ over $U$"), and to every inclusion $V \subset U$ a restriction map $\mathrm{Res}_{U,V} \colon \Gamma(U, \mathcal{F}) \to \Gamma(V, \mathcal{F})$ (or ring homomorphism, or group homomorphism, or morphism of any category). Restriction is supposed to be transitive, in the obvious sense that $\mathrm{Res}_{U,V} \circ \mathrm{Res}_{V,W} = \mathrm{Res}_{U,W}$ whenever $W \subset V \subset U$.

A *sheaf* is a presheaf with "glueing conditions": a section $s \in \Gamma(U, \mathcal{F})$ is uniquely determined by its restrictions $\mathrm{Res}_{U,U_\alpha}(s)$ to open sets $U_\alpha$ that cover $U$; and given an open cover, and sections $s_\alpha \in \Gamma(U_\alpha, \mathcal{F})$ having equal restrictions to all overlaps $U_\alpha \cap U_\beta$, the $s_\alpha$ come by restricting a section $s \in \Gamma(U, \mathcal{F})$.

If $\mathcal{F}$ is a presheaf whose $\Gamma(U, \mathcal{F})$ are sets of maps defined on $U$ (say, with values in some set $\Sigma$), then the first glueing condition holds trivially, and the second holds only if the restrictions placed on maps $U \to \Sigma$ are "local in nature". For example, continuous or differentiable functions on a $C^\infty$ manifold form a sheaf, since continuity and differentiability are local properties of a function. In the same way, sections of a vector bundle $\pi\colon E \to X$, or more generally, the continuous sections of any map $\pi\colon F \to X$ form a sheaf; see Ex. B.4–5 for details.

## B.4. The definition of a variety

The first key service that the language of sheaves performs is to give a satisfactory definition of variety: irreducible affine algebraic sets and their sheaf of regular functions were discussed in B.2. A *variety* $X, \mathcal{O}_X$ is a topological space $X$ with a sheaf $\mathcal{O}_X$ of rings of functions $U \to k$ such that $X, \mathcal{O}_X$ is locally isomorphic to an irreducible affine algebraic set with its sheaf of regular functions. (This can be easily generalised to allow reducible varieties, or indeed general $k$-schemes.)

Thus sheaves solve an important foundational problem of algebraic geometry, the intrinsic definition of variety. You know from a first course in algebraic geometry (for example, [**UAG**] or [**Sh**]) that quasiprojective varieties are quite convenient to define and treat using the tricks of homogeneous coordinates. However, it's unsatisfactory to take this as the formal definition of variety, because you only define $X$ together with extrinsic data of an embedding $X \subset \mathbb{P}^N$.

Compare how the topologist, in defining a manifold, gets away without using sheaves: a topological manifold $M$ is a topological space locally homeomorphic to a ball in $\mathbb{R}^n$; continuity of functions is completely determined by the topology (of course), so that there is no need to specify the sheaf of continuous functions. For a differentiable manifold, you require in addition that the local charts $\varphi_i\colon U_i \to$ (ball in $\mathbb{R}^n$) satisfy $\varphi_i \circ \varphi_j^{-1}$ is differentiable (wherever defined). Because of this assumption, the condition that a function on an open subset $V \subset M$ is differentiable is well defined. An equivalent definition would be to specify the sheaf $\mathcal{E}_{M,\mathrm{diff}}$ of differentiable functions on $M$, and assume that the local charts are isomorphisms of ringed spaces.

In case anyone still hasn't got the idea of what all the fuss is about, I repeat: you can't just define an algebraic variety $X$ to be a point set, or a set with a Zariski topology, because all plane curves would be homeomorphic. If you define it as embedded in a space (affine or projective), you get the notion of rational function $f \in k(X)$ and regular function, and hence isomorphism of varieties, but you also get extrinsic stuff which may have less to do with $X$ than with the embedding. The sheaf $\mathcal{O}_X$ specifies the regular functions on every open set. A variety is a space $X$ together with a notion of regular function on opens of $X$. By giving $\mathcal{O}_X$, you give every possible isomorphism of opens of $X$ with subvarieties of affine or projective space.

## B.5. Other sheaves on algebraic varieties

I continue the theme that you already know lots of sheaves. If $X$ is a normal variety and $D = \sum n_\Gamma \Gamma$ an effective divisor on $X$, rational functions on $X$ with poles at worst $D$ form a sheaf $\mathcal{O}_X(D)$, the *divisorial sheaf* of $D$. Recall that zeros and poles of rational functions are interpreted in terms of discrete valuations $v_\Gamma$ on the local rings $\mathcal{O}_{X,\Gamma}$ at every prime divisor $\Gamma$ (compare 1.6). The condition that $0 \neq f \in k(X)$ has a pole of order at worst $n$ along $\Gamma$ reads $v_\Gamma(f) \geq -n$. To define the sheaf $\mathcal{O}_X(D)$, I tell you what its sections are on every open set $U \subset X$:

$$\Gamma(U, \mathcal{O}_X(D)) = \big\{ f \in k(X) \bigm| v_\Gamma(f) \geq -n_\Gamma \text{ for every } \Gamma \text{ such that } \Gamma \cap U \neq \emptyset \big\}.$$

These are local conditions, so that $\mathcal{O}_X(D)$ is a sheaf. Since $\operatorname{div} f = \sum_\Gamma v_\Gamma(f)\Gamma$, the condition on the right can also be written $\operatorname{div} f + D \geq 0$ on $U$.

**Historical discussion.** In the case $U = X$ a projective variety (for example a nonsingular projective curve), the space

$$\mathcal{L}(D) = \Gamma(X, \mathcal{O}_X(D))$$

is traditionally called the *RR space* of $D$. The RR theorem on a curve gives its dimension:

$$\ell(D) = \dim \mathcal{L}(D) \geq 1 - g + \deg D \qquad \text{with equality if } \deg D > 2g - 2, \qquad \text{(RR)}$$

where $g = g(X)$ is the genus of $X$. This result was the subject of an earlier foundational crisis in algebraic geometry. It follows from the maximum principle that the only global holomorphic functions on a compact Riemann surface are the constants. To find some maps from $X$ to $\mathbb{P}^1 = \mathbb{C} \cup \infty$, Riemann allowed functions which are everywhere meromorphic, with a finite number of poles, which he wrote as a divisor $D$. He then claimed a proof of (RR) based on a Dirichlet principle from electrostatics: imagine that your Riemann surface is made of beaten copper, and the poles as point electric charges; then physical intuition says that the potential equations for the electric field must have a solution. Riemann used (RR) to prove that a compact Riemann surface can be embedded in a projective space $\mathbb{P}^N_\mathbb{C}$, and is hence an algebraic curve (the *Riemann existence theorem*). Unfortunately, Weierstrass pointed out at once that Riemann's Dirichlet principle was false as stated; in fact, 30 or 40 years later Hilbert gave a revised statement and claimed proof of the Dirichlet principle that was also erroneous, although, by all accounts, Hilbert was so famous by then that nobody dared tell him. In any case, Clebsch, Max Noether and Brill proved (RR) for algebraic curves by purely algebraic means shortly after Riemann.

On an irreducible variety, $\mathcal{O}_X(D)$ is defined as a subsheaf of the constant sheaf $k(X)$. If $\operatorname{div} f = D - D'$ then multiplying by $f$ in $k(X)$ clearly takes $\mathcal{O}_X(D)$ into $\mathcal{O}_X(D')$. Thus linearly equivalent divisors $D$ give rise to isomorphic divisorial sheaves $\mathcal{O}_X(D)$. See Ex. B.13 for more details.

In A.10 I introduced the canonical divisor class $K_X$ as the divisor $\operatorname{div} s$ of a rational $n$-form $s = g\mathrm{d}f_1 \wedge \cdots \wedge \mathrm{d}f_n$. And $K_X$ is defined up to linear equivalence, because any other rational $n$-form $s'$ is of the form $s' = hs$ for some $h \in k(X)$, so that $\operatorname{div}(s') = \operatorname{div} h + \operatorname{div} s$. Therefore, there is a well-defined divisorial sheaf $\mathcal{O}_X(K_X)$ on $X$. In fact it is easier and more intrinsic to introduce the sheaf $\omega_X = \Omega^n_X$ first, prove that it is divisorial (on a nonsingular $X$ it is *locally free of rank* 1, or *invertible*), and then to define $K_X$ as any divisor such that $\mathcal{O}_X(K_X) = \Omega^n_X$. See Ex. B.16. In fact, although the divisor $K_X$ is only a divisor class, the sheaf $\omega_X = \Omega^n_X = \mathcal{O}_X(K_X)$ is canonical.

## B.6. Subsheaves, stalks and quotient sheaves

There are two very different kinds of definitions and arguments in sheaf theory, those that take place at the level of the spaces of sections $\Gamma(U, \mathcal{F})$, and those that are local at every point $P \in X$ and involve stalks. I start with a few instances of the first. A *homomorphism* of sheaves $\varphi\colon \mathcal{F} \to \mathcal{G}$ is a collection of morphisms $\varphi_U : \Gamma(U, \mathcal{F}) \to \Gamma(U, \mathcal{G})$ for all opens $U$ that commutes with restrictions $\mathrm{Res}_{U,V}$ of $\mathcal{F}$ and $\mathcal{G}$:

$$\begin{array}{ccc} \Gamma(U,\mathcal{F}) & \to & \Gamma(U,\mathcal{G}) \\ \downarrow & \copyright & \downarrow \\ \Gamma(V,\mathcal{F}) & \to & \Gamma(V,\mathcal{G}) \end{array}$$

Similarly, a *subsheaf* $\mathcal{F} \subset \mathcal{G}$ is a collection of subobjects $\Gamma(U, \mathcal{F}) \subset \Gamma(U, \mathcal{G})$ that themselves form a sheaf under the same restriction maps. If $f\colon X \to Y$ is a continuous map of topological spaces, and $\mathcal{F}$ a sheaf on $X$, then it is immediate to check that $U \mapsto \Gamma(f^{-1}U, \mathcal{F})$ for open sets $U \subset Y$ gives a sheaf on $Y$, the sheaf-theoretic *pushforward* $f_*\mathcal{F}$. All of these definitions and constructions can be made for presheaves just as well as for sheaves.

The *stalk* of a sheaf $\mathcal{F}$ at a point $P \in X$ is the direct limit $\mathcal{F}_P = \varinjlim_{U \ni P} \Gamma(U, \mathcal{F})$. The limit looks intimidating, but this is just another case of M. Jourdain's prose. Namely, the direct limit is the set of all sections $s \in \Gamma(U, \mathcal{F})$ over all open sets containing $U$, modulo the equivalence relation $s = s'$ if they coincide on some smaller neighbourhood of $P$; in other words, $\mathcal{F}_P$ consists of *germs* of sections at $P$. For example, if $\mathcal{O}_{\mathrm{an}}$ is the sheaf of holomorphic functions on $\mathbb{C}$ then the stalk $\mathcal{O}_{\mathrm{an},0}$ consists of all power series with positive radius of convergence, and a germ is an analytic function on some neighbourhood of 0; different germs are defined on different neighbourhoods. The stalks of the structure sheaf $\mathcal{O}_X$ of a variety $X$ are the local rings $\mathcal{O}_{X,P}$, and in this case the direct limit is simply a union, as mentioned in B.2. Another example: a common definition of a tangent vector to a manifold $M$ at a point $P$ is as a derivation of functions defined near $P$. The derivation acts on germs of smooth functions: it looks at the function only in an arbitrarily small neighbourhood of $P$.

Now we say that a homomorphism of sheaves $\varphi\colon \mathcal{F} \to \mathcal{G}$ is *surjective* if it induces surjective maps $\varphi_P : \mathcal{F}_P \to \mathcal{G}_P$ on each stalk. Surjectivity onto $\Gamma(U, \mathcal{G})$ for all $U$ is the wrong requirement. Another way of saying the same thing is as follows: if $s \in \Gamma(U, \mathcal{G})$, I don't require that $s$ itself comes from some $t \in \Gamma(U, \mathcal{F})$, but only that this holds in a small neighbourhood of any $P \in U$. I now give a baby example (see also Ex. B.10).

Let $P_1, P_2, P_3 \in \mathbb{P}^n$ be 3 distinct points, not in the hyperplane $(x_0 = 0)$. On $\mathbb{P}^n$, consider the sheaf $\mathcal{O}_{\mathbb{P}^n}(1)$ of linear forms, which is the sheaf defined by

$$\Gamma(U, \mathcal{O}_{\mathbb{P}^n}(1)) = \left\{ \frac{f}{g} \in k(x_0, \dots, x_n) \;\middle|\; \begin{array}{l} f, g \in k[x_0, \dots, x_n] \text{ homog. of degree} \\ d+1 \text{ resp. } d, \text{ and } g(P) \neq 0 \text{ at } P \in U \end{array} \right\}.$$

Now I can find a linear form not vanishing at $P_i$, so that the evaluation map $\mathcal{O}_{\mathbb{P}^n}(1) \to k_P$ defined by $s \mapsto (s/x_0)(P)$ is surjective. Here $k_P$ is the "skyscraper sheaf" with sections over $U$ either zero if $P \notin U$ or a copy of $k$ if $P \in U$. The kernel is the sheaf of linear forms vanishing at $P$, that is, $m_P \cdot \mathcal{O}_{\mathbb{P}^n}(1)$. Now consider the evaluation map at all 3 points at once:

$$0 \to \mathcal{I} \cdot \mathcal{O}_{\mathbb{P}^n}(1) \to \mathcal{O}_{\mathbb{P}^n}(1) \to k_{P_1} \oplus k_{P_2} \oplus k_{P_3} \to 0, \qquad (*)$$

where for brevity I write $\mathcal{I} = m_{P_1} \cdot m_{P_2} \cdot m_{P_3}$ for the ideal sheaf of $\{P_1, P_2, P_3\}$, that is, the subsheaf of $\mathcal{O}_{\mathbb{P}^n}$ consisting of regular functions on an open $U$ vanishing at $\{P_1, P_2, P_3\} \cap U$. Now in sheaf theory we say that the evaluation map in $(*)$ is surjective, because it is surjective locally at every point. Is it surjective on global sections? The global map evaluates linear forms in $x_0, \dots, x_n$ on 3 points, which is surjective if $P_1, P_2, P_3$ span a plane in $\mathbb{P}^n$, and not surjective if they are collinear. I thus get the exact sequence

$$
\begin{aligned}
0 &\to \text{linear forms vanishing at } P_1, P_2, P_3 \to \langle x_0, \dots, x_n \rangle \to k_{P_1} \oplus k_{P_2} \oplus k_{P_3} \to \\
&\to \text{linear dependences among } P_1, P_2, P_3 \to 0
\end{aligned}
$$

In other words, the homomorphism of sheaves is surjective, but it gives a homomorphism on global sections which is not necessarily surjective. In more general language, the last display is written

$$
\begin{aligned}
0 &\to H^0(\mathbb{P}^n, \mathcal{I} \cdot \mathcal{O}_{\mathbb{P}^n}(1)) \to H^0(\mathbb{P}^n, \mathcal{O}_{\mathbb{P}^n}(1)) \to H^0(\mathbb{P}^n, k_{P_1} \oplus k_{P_2} \oplus k_{P_3}) \to \\
&\to H^1(\mathbb{P}^n, \mathcal{I} \cdot \mathcal{O}_{\mathbb{P}^n}(1)) \to 0,
\end{aligned}
$$

where $H^1(\mathbb{P}^n, \mathcal{I} \cdot \mathcal{O}_{\mathbb{P}^n}(1))$ is the first cohomology group of the sheaf $\mathcal{I} \cdot \mathcal{O}_{\mathbb{P}^n}(1)$.

I mention two other constructions of sheaf theory that use stalks: the *quotient sheaf* $\mathcal{G}/\mathcal{F}$ is defined to have stalk $\mathcal{G}_P/\mathcal{F}_P$ at every $P \in X$; and the *sheaf-theoretic pullback* $f^{-1}\mathcal{F}$ of a sheaf on $Y$ by a morphism $f\colon X \to Y$ is defined to have stalk $\mathcal{F}_{f(P)}$ at $P \in X$. I omit the definition, which involves the notion of the associated sheaf of a presheaf. (This material is completed in Ex. B.27–29.)

## B.7. Coherent sheaves

The sheaves of algebraic geometry $\mathcal{O}_X$, $\mathcal{O}_X(D)$, $\mathcal{O}_{\mathbb{P}^n}(r)$ we have met so far are all coherent sheaves; so is the ideal sheaf $\mathcal{I}_Y \subset \mathcal{O}_X$ of a subvariety $Y \subset X$. The adjective coherent means that they are sheaves of modules over $\mathcal{O}_X$, with a finiteness condition, and closely related to the structure sheaf $\mathcal{O}_X$.

The general progression is presheaf, sheaf, sheaf of $\mathcal{O}_X$-modules, (quasi-) coherent sheaf, locally free sheaf. I have not been through the general definitions particularly carefully; it should be clear what the definition of sheaf of $\mathcal{O}_X$-modules is. If you have trouble see, for example, [**FAC**] or [**H1**], Chapter II. The definition of (quasi-)coherent involves tension between the requirements of generality and explicitness: namely, the definition is that $\mathcal{F}$ should be a sheaf of $\mathcal{O}_X$-modules, and $\mathcal{F}$ should be locally isomorphic to the cokernel of a homomorphism between free sheaves. In other words, on local pieces $U$, there should exist a resolution

$$
\mathcal{O}_U^{\oplus N} \to \mathcal{O}_U^{\oplus M} \to \mathcal{F}_{|U} \to 0 \tag{1}
$$

(that is, an exact sequence of sheaves of $\mathcal{O}_X$-modules). The tension comes because to have an intrinsic definition you want the condition for all sufficiently small neighbourhoods $U$, but to have an explicit construction you want only that there exists a cover of $X$ by opens $U$ with the property.

The tension is solved in the best possible way: thanks to cohomology, we can have it both ways! In other words, if I have an open cover $X = \bigcup U_i$ of $X$ by affine sets $U_i$ such that (1) holds for each $U_i$, then the same condition (1) holds for every affine open set $U \subset X$. This is the content of Rule ii below. It is proved in [**FAC**], [**H1**], Chapter II and [**Sh**], Chapter VII.

Now *quasicoherent* is condition (1) for all $U$ of some (any) open cover with arbitrary cardinals $N$ and $M$. *Coherent* is the same with finite $N$ and $M$. Notice that for a fixed affine $U$, (1) and the vanishing of $H^1$ (see Rule ii) gives an sequence

$$k[U]^N \to k[U]^M \to \Gamma(U,\mathcal{F}) \to 0, \tag{2}$$

where $k[U] = \Gamma(U, \mathcal{O}_X)$ is the affine coordinate ring of $U$. In other words, $\Gamma(U,\mathcal{F}) = F$ is just an arbitrary $k[U]$-module; moreover, (1) implies that $\mathcal{F}$ is determined on $U$ by $F$ and localisation. This is the construction of the sheaf $\widetilde{F}$ on $U$ from a $k[U]$-module $F$.

It is instructive to compare the condition in (1) with the topologist's notion of a map of vector bundles. (1) is a homomorphism of free sheaves $\mathcal{O}_U^{\oplus N} \to \mathcal{O}_U^{\oplus M}$, and so for a choice of bases is determined by a $N \times M$ matrix $A$ with coefficients in $\Gamma(U, \mathcal{O}_X)$. In algebraic geometry, we must allow the rank of the matrix $A$ to vary from point to point. It is upper semicontinuous in any case, since $\operatorname{rank} A \le r$ is a closed condition.

## B.8. Examples

**Example 1.** If $Y \subset X$ is a subvariety (subscheme) of an affine variety defined by $f_1 = \cdots = f_n = 0$ then the structure sheaf $\mathcal{O}_Y$ is determined by the exact sequence

$$\mathcal{O}_X^n \xrightarrow{F} \mathcal{O}_X \to \mathcal{O}_Y \to 0, \qquad \text{where } F = (f_1, \dots, f_n).$$

Obvious $F$ has rank 1 outside $Y$ and rank 0 along $Y$.

**Example 2.** Now let $X = \mathbb{A}^n$ and $Y = \mathbb{A}^{n-1} : (x_n = 0) \subset X$, and consider the surjective homomorphism $p\colon \mathcal{O}_X^r \to \mathcal{O}_Y$ determined by $(f_1, \dots, f_r) \mapsto f_{r|Y}$. Obviously the kernel of $p$ is the subsheaf $\mathcal{O}_X^{r-1} \oplus \mathcal{I}_Y \subset \mathcal{O}_X^r$. Since $\mathcal{I}_Y = x_n \cdot \mathcal{O}_X$ is itself a free sheaf, $\ker p \cong \mathcal{O}_X^n$. Using this isomorphism I get a short exact sequence

$$0 \to \mathcal{O}_X^n \xrightarrow{A} \mathcal{O}_X^n \to \mathcal{O}_Y \to 0, \qquad \text{where } A = \operatorname{diag}(1, \dots, 1, x_n). \tag{3}$$

The construction $\ker p$ passes from a locally free sheaf $\mathcal{E}$ to $\mathcal{E}' = \ker p \subset \mathcal{E}$, where $p$ is the composite of restricting to a divisor $Y \subset X$ and a projection of $\mathcal{E}_{|Y}$ to a quotient bundle. This construction is well known as a standard elementary transformation of vector bundles.

Notice that the homomorphism $A$ of sheaves is injective even at points of $Y$ where $A$ drops rank. The point is that the map of sheaves only looks at sections over opens, or stalks, and does not look at the *fibre* of the vector bundle. The stalk looks like $\mathcal{O}_{X,P}^n$, a free module over the local ring, whereas the fibre looks like the quotient $\mathcal{E}_P / m_P \mathcal{E}_P$, which is a $k$-vector space.

**Example 3.** A similar example. Suppose that $X = \mathbb{A}^2$ and that $Y \subset X$ is the subscheme defined by $f = g = 0$; suppose for simplicity that $Y$ only lives at one point, that is $V(f,g) = \{P\}$. The subscheme $Y$ is the point $P$ with structure sheaf the finite dimensional ring $\mathcal{O}_Y = \mathcal{O}_{X,P}/(f,g)$. Then $\mathcal{O}_{X,P}$ is a UFD, so it's easy to check that the following sequence

$$0 \to \mathcal{O}_{X,P} \xrightarrow{-g,f} \mathcal{O}_{X,P} \oplus \mathcal{O}_{X,P} \xrightarrow{\substack{f\\ g}} \mathcal{O}_{X,P} \to \mathcal{O}_Y \to 0$$

is exact. It's called the *Koszul complex* of $f, g$; its construction only depends on the fact that $f, g$ forms a regular sequence in $\mathcal{O}_{X,P}$.

Now the point of this example is that every section of a locally free sheaf of rank 2 with only zeros in codimension 2 looks like this.

As a rule, traditional topologists have only allowed maps of constant rank between vector bundles, which is equivalent to saying that the kernel, image and cokernel are locally direct summands. As we have seen in Examples 1–3, the more general notion of sheaf homomorphism between locally free sheaves is very useful in algebraic geometry.

## B.9. Rules of coherent cohomology

This table of rules states the main useful results of coherent cohomology at a fairly simple level of generality. I will take them as axioms throughout. For the proofs, see [**FAC**]. Anyone complaining that the paper is in French will receive a blast of unpleasant sarcasm.

Actually, the hard thing is not to get used to these rules, but to understand what a coherent sheaf is.

### Data 1

For any variety $X$ over $k$ and any (quasi-) coherent sheaf $\mathcal{F}$ on $X$ there is a $k$-vector space $H^i(X, \mathcal{F})$, that is functorial in $\mathcal{F}$. In other words a homomorphism of sheaves of $\mathcal{O}_X$-modules $a\colon \mathcal{F} \to \mathcal{G}$ gives rise to a linear map $a_*\colon H^i(X, \mathcal{F}) \to H^i(X, \mathcal{G})$, with obvious compatibilities. (To answer the topologist's immediate question: there is no functoriality for morphisms of varieties $X$ as yet. Sheaf cohomology is a property of the category of sheaves over a fixed $X$.)

### Data 2

If $0 \to \mathcal{F}' \to \mathcal{F} \to \mathcal{F}'' \to 0$ is a short exact sequence of quasicoherent sheaves on $X$ then there is a coboundary map

$$d_i\colon H^i(X, \mathcal{F}'') \to H^{i+1}(X, \mathcal{F}'),$$

again functorial in exact sequences.

So far, $H^*(X, \text{blank})$ is a *cohomological $\delta$-functor*, if you like that kind of thing. This data satisfies the following conditions:

i. **Sections** $H^0$

$$H^0(X, \mathcal{F}) = \Gamma(X, \mathcal{F}),$$

the space of sections of a sheaf, as in the definition of a sheaf.

ii. **Affine varieties** $X$

If $X$ is affine then

$$H^i(X, \mathcal{F}) = 0 \qquad \text{for all } i > 0.$$

Moreover, $H^0(X, \mathcal{F})$ is a sufficiently big module over the affine coordinate ring $k[X] = H^0(X, \mathcal{O}_X)$, so that the following localisation works:

$$H^0(U, \mathcal{F}) = H^0(X, \mathcal{F}) \otimes_{k[X]} H^0(U, \mathcal{O}_X) \qquad \text{for every open } U \subset X;$$
$$\mathcal{F}_P = H^0(X, \mathcal{F}) \otimes_{k[X]} \mathcal{O}_{X,P} \qquad \text{for every point } P \in X.$$

Actually, you have to prove all this before the notion of coherent sheaf is reasonably intrinsic (compare B.7).

iii. **Dimension**

$$H^i(X, \mathcal{F}) = 0 \qquad \text{for all } i > \dim X$$

The topologist who at last finds some mild satisfaction should beware that I mean the dimension of $X$ as an algebraic variety, e.g., an algebraic curve has dimension 1 (although over the complexes it's a Riemann surface).

iv. **Long exact sequence**

If $0 \to \mathcal{F}' \to \mathcal{F} \to \mathcal{F}'' \to 0$ is a short exact sequence of quasicoherent sheaves on $X$ then the functoriality homomorphisms of Data 1 and the coboundary homomorphisms of Data 2 give a cohomology long exact sequence

$$\cdots \to H^i(X, \mathcal{F}') \to H^i(X, \mathcal{F}) \to H^i(X, \mathcal{F}'') \to$$
$$\to H^{i+1}(X, \mathcal{F}') \to \cdots$$

v. **Finite dimensionality**

If $\mathcal{F}$ is coherent and $X$ is proper (for example, projective) then

$$H^i(X, \mathcal{F}) \qquad \text{is finite dimensional over } k \text{ for any } i.$$

One traditionally writes $h^i(X, \mathcal{F}) = \dim_k H^i(X, \mathcal{F})$.

vi. **Ample line bundle, Serre vanishing**

Suppose that $X \subset \mathbb{P}^n$ is a closed subvariety. Let $\mathcal{O}_X(1) = \mathcal{O}_X \otimes \mathcal{O}_{\mathbb{P}^n}(1)$ be the invertible sheaf obtained by restricting $\mathcal{O}_{\mathbb{P}^n}(1)$ to $X$; this is the sheaf having the homogeneous coordinates $x_0, \dots, x_n$ as sections (compare B.6 for a definition of $\mathcal{O}_{\mathbb{P}^n}(1)$). Write $\mathcal{O}_X(r)$ for the $r$ times tensor product of $\mathcal{O}_X(1)$, and for any quasicoherent sheaf $\mathcal{F}$ on $X$, write $\mathcal{F}(r) = \mathcal{F} \otimes \mathcal{O}_X(r)$.

Then given a coherent sheaf $\mathcal{F}$, there exists $N$ such that all the following hold for all $r \geq N$: the space of global sections $H^0(X, \mathcal{F}(r))$ is big enough so that

$$\operatorname{Im}\{H^0(X, \mathcal{F}(r)) \to H^0(U, \mathcal{F}(r))\} \qquad \text{generates } H^0(U, \mathcal{F}(r))$$

as a module over $H^0(U, \mathcal{O}_X)$ for every open $U \subset X$, and

$$\operatorname{Im}\{H^0(X, \mathcal{F}(r)) \to \mathcal{F}(r)_P\} \qquad \text{generates } \mathcal{F}(r)_P$$

as a module over $\mathcal{O}_{X,P}$ for every point $P \in X$. In other words, $\mathcal{F}(r)$ is *generated by its* $H^0$. Moreover,

$$H^i(X, \mathcal{F}(r)) = 0 \qquad \text{for all } i > 0.$$

This is called *Serre vanishing.*

Actually everything in (vi), apart from the language, is a trivial consequence of (ii) applied to the affine cone $C_X$ over $X$. A coherent sheaf $\mathcal{F}$ on $X$ corresponds to a finitely generated graded module $\bigoplus_{r\geq 0} H^0(X, \mathcal{F}(r))$ over the usual homogeneous coordinate ring of $X \subset \mathbb{P}^n$, the affine coordinate ring of $C_X$. Thus the language is just a formal way of saying the usual correspondence between homogeneous polynomials and functions on $X \subset \mathbb{P}^n$.

### vii. **Serre duality**

Let $X$ be a nonsingular projective $n$-fold and $K_X$ its canonical divisor class, so that $\mathcal{O}_X(K_X) = \Omega^n_X = \bigwedge^n \Omega^1_X$. Then

$$H^n(X, \mathcal{O}_X(K_X)) \qquad \text{is a 1-dimensional vector space} \cong k.$$

Assuming nobody objects, I pick a generator and write $\ = k$.

For any invertible sheaf $\mathcal{L} = \mathcal{O}_X(D)$, write $\mathcal{L}^{-1} = \mathcal{H}om_{\mathcal{O}_X}(\mathcal{L}, \mathcal{O}_X) = \mathcal{O}_X(-D)$. Then there is a canonical pairing

$$H^i(X, \mathcal{L}) \times H^{n-i}(X, \mathcal{L}^{-1} \otimes_{\mathcal{O}_X} \mathcal{O}_X(K_X)) \to k,$$

or

$$H^i(X, \mathcal{O}_X(D)) \times H^{n-i}(X, \mathcal{O}_X(K_X - D)) \to k,$$

which establishes a duality between the two groups.

**Remark.** You can ignore this remark *en première lecture.* Of course Serre duality can be generalised to singular $X$ and arbitrary coherent sheaves $\mathcal{F}$. However, in the same way that Poincaré duality for singular cohomology requires a manifold, or at least a space satisfying a suitable local duality, the general form is a bit complicated. If $X$ is Cohen–Macaulay and $\dim X = n$ then there exists a sheaf $\omega_X$, the *Grothendieck dualising sheaf,* such that $H^n(X, \omega_X) = k$, and for any coherent sheaf $\mathcal{F}$ there is a canonical pairing

$$H^i(X, \mathcal{F}) \times \operatorname{Ext}^{n-i}_{\mathcal{O}_X}(\mathcal{F}, \omega_X) \to H^n(X, \omega_X) = k$$

which establishes a duality between the two groups.

If $X$ is not Cohen–Macaulay, for example if it has components of different dimension, then you can't expect a duality that works in a single dimension ($i$ against $n - i$), and $\omega_X$ is replaced by a complex.[3]

[3] I have colloquial lecture notes on this topic which I may include in a later edition. See also, for example, [**R1**], App. to §2.

### viii. Euler–Poincaré characteristic $\chi(X, \mathcal{F})$ and Hilbert polynomial

Whenever the dimensions are finite, I write $h^i(X, \mathcal{F}) = \dim_k H^i(X, \mathcal{F})$. Define the Euler–Poincaré characteristic of $\mathcal{F}$ by $\chi(X, \mathcal{F}) = \sum_{i=0}^{\infty}(-1)^i h^i(X, \mathcal{F})$. Although its definition involves all the cohomology groups $H^i(X, \mathcal{F})$, this alternating sum is in fact a much more elementary quantity. For example, from the cohomology long exact sequence (iv), it follows at once that

$$\chi(X, \mathcal{F}) = \chi(X, \mathcal{F}') + \chi(X, \mathcal{F}'').$$

For $X \subset \mathbb{P}^n$, the numerical function given by $r \mapsto \chi(X, \mathcal{F}(r))$ is a polynomial, called the *Hilbert polynomial* of $\mathcal{F}$.

### ix. Riemann–Roch

For a divisor $D$ on an algebraic curve $C$

$$\chi(\mathcal{O}_C(D)) = h^0(\mathcal{O}_C(D)) - h^1(\mathcal{O}_C(D)) = \chi(\mathcal{O}_C) + \deg D,$$
$$\text{and} \quad \chi(\mathcal{O}_C) = 1 - g, \qquad \text{where } g = g(C) \text{ is the genus.}$$

For a divisor $D$ on an algebraic surface $X$

$$\chi(\mathcal{O}_X(D)) = h^0(\mathcal{O}_X(D)) - h^1(\mathcal{O}_X(D)) + h^2(\mathcal{O}_X(D)) = \chi(\mathcal{O}_X) + \tfrac{1}{2}D(D - K_X),$$
$$\text{and} \quad \chi(\mathcal{O}_X) = \frac{1}{12}(c_1^2 + c_2) = \frac{1}{12}(K_X^2 + e(X)),$$

where $K_X$ is the canonical class of $X$ and $e(X) =$ the topological Euler number of $X$, the alternating sum of Betti numbers.

For an arbitrary sheaf on a projective variety $X$,

$$\chi(X, \mathcal{F}) = \int \operatorname{ch} \mathcal{F} \cdot \operatorname{Td}_X,$$

where $\operatorname{ch} \mathcal{F}$ and $\operatorname{Td}_X$ are certain characteristic classes of the sheaf $\mathcal{F}$ and the tangent sheaf of $X$, and the integral sign (also pronounced "evaluate on the fundamental class of $X$") means in practice that you take the sum of homogeneous terms of degree $\dim X$, interpret them as a zero dimensional cycle on $X$, and then as an integer.

You can't be a grown-up algebraic geometer until you have memorised these formulas. Eventually you have to learn what they mean, and how to calculate with them as well.

## Exercises to Chapter B

**1.** Prove that the constant sheaf $\mathbb{Z}$ cannot be made into a sheaf of $\mathcal{O}_X$-modules.

**2.** Recall that the stalk $\mathcal{F}_P$ of a sheaf $\mathcal{F}$ at a point $P$ is defined as the direct limit of the sets of sections $\Gamma(U, \mathcal{F})$ taken over all $U \ni P$. If $\varphi\colon \mathcal{F} \to \mathcal{G}$ is a morphism of sheaves, show how to define $\varphi_P\colon \mathcal{F}_P \to \mathcal{G}_P$, and check that it is well defined.

**3.** If $X$ is a variety with structure sheaf $\mathcal{O}_X$ and $P \in X$ a point, prove that the stalk $\mathcal{O}_{X,P}$ is a local ring.

**4.** Let $\pi\colon F \to X$ be a continuous map of topological spaces. A section of $\pi$ over an open set $U \subset X$ is a map $s\colon U \to F$ such that $\pi \circ s = \mathrm{id}_U$. Prove that sections of $F$ form a sheaf $\mathcal{F}$.

**5.** Let $\pi\colon F \to X$ be a vector bundle of rank $r$ (in the continuous, differentiable, complex analytic or algebraic categories). Prove that $\mathcal{F}$ constructed in the preceding exercise is a locally free sheaf of rank $r$ over the appropriate structure sheaf of $X$.

**6** (harder). Prove that there is an equivalence of categories between vector bundles and locally free sheaves. You'll need to choose one of the continuous, differentiable, complex analytic or algebraic categories, and be careful to ensure that the two sides of your equivalence have the same morphisms; one (boring) possibility is to allow only isomorphisms as morphisms.

**7.** An *affine structure* on an $n$-dimensional manifold $M$ is an atlas consisting of charts $\varphi_i\colon U_i \xrightarrow{\sim}$ ball in $\mathbb{R}^n$ such that the glueing maps $\varphi_i \circ \varphi_j^{-1}$ are affine linear transformations $x \mapsto Ax + B$. Show how to introduce a sheaf of affine linear functions on $M$, and to give an alternative definition of manifold with affine structure based on an affine linear structure sheaf.

**8.** If $\varphi\colon \mathcal{F} \to \mathcal{G}$ is a morphism of sheaves, define $\operatorname{Im}\varphi \subset \mathcal{G}$ and prove that $\varphi\colon \mathcal{F} \to \operatorname{Im}\varphi$ is surjective. [Hint: As in B.6, use stalks $\mathcal{F}_P$.]

**9.** Let $X$ be a projective variety, $D$ a divisor on $X$, and $\mathcal{O}_X(D)$. Choose a basis $f_0, \dots, f_n$ of the RR space $\mathcal{L}(D)$, that is, rational functions $f_i \in k(X)$ with $\operatorname{div} f_i + D \ge 0$; write $\varphi_{|D|}\colon X \to \mathbb{P}^n$ for the map defined by the ratio $f_0 : \cdots : f_n$. On the other hand, there is a map $\varphi_{\mathcal{L}}\colon X \to \mathbb{P}^n$ defined by the global section of $\mathcal{L}$. Think through the definitions, and show that these two maps are identical.

**10.** On the complex plane $\mathbb{C}$ (with the complex topology), let $\mathcal{O}_{\text{an}}$ be the sheaf of holomorphic functions, and $\mathcal{O}^*_{\text{an}}$ the sheaf of invertible holomorphic functions; check you have mastered the language by writing down displayed formulas with $\{\,|\,\}$ for the sections of $\mathcal{O}_{\text{an}}$ and $\mathcal{O}^*_{\text{an}}$ over an appropriate domain.

Show that the exponential map $f \mapsto \exp(f)$ defines a morphism of sheaves $\exp\colon \mathcal{O}_{\text{an}} \to \mathcal{O}^*_{\text{an}}$, and that it is surjective. Prove that the kernel is the constant sheaf $2\pi i\mathbb{Z}$. Consider the exact sequence

$$0 \to 2\pi i\mathbb{Z} \to \mathcal{O}_{\text{an}} \xrightarrow{\exp} \mathcal{O}^*_{\text{an}} \to 0.$$

We agreed earlier that exp is surjective as a morphism of sheaves. Show that if $U \subset \mathbb{C}$ is the annular region $0 < a < |z| < b$ then $\Gamma(U, \mathcal{O}_{\text{an}}) \to \Gamma(U, \mathcal{O}^*_{\text{an}})$ is not surjective.

Find a necessary and sufficient condition on an open set $U \subset \mathbb{C}$ such that $\Gamma(U, \mathcal{O}_{\text{an}}) \to \Gamma(U, \mathcal{O}^*_{\text{an}})$ is surjective.

(In this question, the sheaves are not coherent algebraic sheaves. Only $\mathcal{O}_{\text{an}}$ is a coherent analytic sheaf. The two sheaves $2\pi i\mathbb{Z}$ and $\mathcal{O}^*_{\text{an}}$ are sheaves of Abelian groups, but obviously cannot be made into $\mathcal{O}_{\text{an}}$-modules.)

**11.** Show that

$$0 \to \mathcal{I}_Y \to \mathcal{O}_X \to \mathcal{O}_Y \to 0$$

is an exact sequence of sheaves on $X$ whenever $Y$ is a subvariety.

**12.** (a) Let $C$ be a projective curve. Prove that

$$C \cong \mathbb{P}^1 \iff g(C) = 0.$$

[Hint: Use RR.]

(b) Let $C \subset Y$ be an irreducible curve in a nonsingular surface. Show that if $K_Y C = 0$ and $C^2 = -2$ then $C \cong \mathbb{P}^1$ so that $C$ is a $-2$-curve.

**13.** Prove that any locally free sheaf of $\mathcal{O}_X$-modules of rank 1 $\mathcal{L}$ (invertible sheaf) on a nonsingular variety $X$ is of the form $\mathcal{L} \cong \mathcal{O}_X(D)$. Prove that $\mathcal{O}_X(D) \cong \mathcal{O}_X(D')$ as sheaves of $\mathcal{O}_X$-modules if and only if $D \overset{\text{lin}}{\sim} D'$ (linear equivalence was defined in 1.6).

**14.** Let $X$ be a nonsingular $n$-fold. Define a *rational* 1*-form* to be an expression $\sum f_i \mathrm{d}g_i$ with $f_i, g_i \in k(X)$ modulo the Leibnitz rules $\mathrm{d}a = 0$ for $a \in k$ and $\mathrm{d}(fg) = f\mathrm{d}g + g\mathrm{d}f$. Write $\Omega^1_{k(X)/k}$ for the set of rational 1-forms. Prove that it is an $n$-dimensional vector space over $k(X)$ with basis $\mathrm{d}g_1, \dots, \mathrm{d}g_n$, where $g_1, \dots, g_n$ is any (separable) transcendence basis of $k(X)/k$.

**15.** The sheaf $\Omega^1_X$ of regular 1-forms is defined by imposing regularity conditions on rational 1-forms; in other words, if $s \in \Omega^1_{k(X)/k}$, then $s$ is regular at a point $P \in X$ if and only if it can be written $\sum f_i \mathrm{d}g_i$ with $f_i, g_i \in \mathcal{O}_{X,P}$. Prove that if $z_1, \dots, z_n$ are local coordinates at a point $P \in X$ then $\mathrm{d}z_1, \dots, \mathrm{d}z_n$ are local generators of $\Omega^1_X$ in a neighbourhood of $P$.

If you're happy with the tangent sheaf $\mathcal{T}_X$ or tangent bundle $T_X$ of $X$, show that $\Omega^1_X$ can be identified with the sheaf of linear forms on $\mathcal{T}_X$ or $T_X$. That is, $\Omega^1_X = \mathcal{H}om_{\mathcal{O}_X}(\mathcal{T}_X, \mathcal{O}_X)$, (the sheaf Hom, defined by setting the stalk at $P$ equal to $\Omega^1_{X,P} = \mathrm{Hom}_{\mathcal{O}_{X,P}}(\mathcal{T}_{X,P}, \mathcal{O}_{X,P})$) or $\Gamma(U, \Omega^1_X)$ is the set of morphisms $T_{X|U} \to k \times U$ that commute with the projection to $U$ and are linear in each fibre.

**16.** Set $\Omega^n_X = \bigwedge^n \Omega^1_X$, the sheaf of regular $n$-forms. Prove that $\Omega^n_X = \mathcal{O}_X(K_X)$.

**17.** Show that $\mathcal{O}_{\mathbb{P}^n}(H) \cong \mathcal{O}_{\mathbb{P}^n}(1)$ where $H$ is any hyperplane. Extend to $\mathcal{O}_{\mathbb{P}^n}(kH)$ for any $k \in \mathbb{Z}$.

**18.** Give a definition of $\mathcal{O}_{\mathbb{P}^n}(r)$ in terms of ratios $f/g$ of homogeneous polynomials in $x_1, \dots, x_n$ of degree $d + r$ and $d$ respectively (compare B.6). Do the same for $\mathcal{O}_{\mathbb{F}}(eL + dM)$ on the scrolls in terms of bihomogeneous polynomials. Notice that here the space is constructed in terms of a group action, its structure sheaf in terms of invariant rational functions, and the other eigenspaces (character spaces) of rational functions correspond to divisorial sheaves (locally free sheaves of rank 1).

**19.** In the notation of Chapter 2, let $\mathbb{F} = \mathbb{F}(a_1, \dots, a_n)$ be the scroll, and $M$ the divisor class linearly equivalent to $D_i + a_i L$. Prove that the pushforward of $\mathcal{O}_X(M)$ is a sheaf of $\mathcal{O}_{\mathbb{P}^1}$-modules isomorphic to $\mathcal{O}_{\mathbb{P}^1}(a_1) \oplus \cdots \oplus \mathcal{O}_{\mathbb{P}^1}(a_n)$.

**20.** Construct an example of an invertible sheaf $\mathcal{L}$ on a variety $X$ generated by its $H^0$, but whose sections do not separate points. Construct an example so that the sections separate points but do not separate tangent vectors.

**21.** State and prove Bézout's theorem on the surface scroll $\mathbb{F}_a = \mathbb{F}(a, 0)$; in other words, if $C$ and $D$ are curves of bidegree $d, e$ and $d', e'$ with no common components, state a formula for the number of points of $C \cap D$ counted with multiplicities, and prove it by the argument sketched in A.9 of the notes. You'll need to figure out the dimension of the space of forms of bidegree $d, e$.

**22.** From now on $C$ is a nonsingular projective curve. I assume known that $\deg \operatorname{div} f = 0$ for any rational function $f \in k(C)$, that is, a rational function has the same number of zeros and poles (counted with multiplicities).

(a) Prove that if $\deg D < 0$ then $H^0(C, \mathcal{O}_C(D)) = 0$.

(b) Prove that $h^0(\mathcal{O}_C(D-P)) = h^0(\mathcal{O}_C(D))$ or $h^0(\mathcal{O}_C(D)) - 1$.

(c) Prove that $h^0(\mathcal{O}_C(D)) \le \deg D + 1$, and that equality holds, with $\deg D \ge 1$, if and only if $C \cong \mathbb{P}^1$.

**23.** Given that $\chi(\mathcal{O}_C) = 1 - g$, use induction and short exact sequences of the form

$$0 \to \mathcal{O}_C(D-P) \to \mathcal{O}_C(D) \to k_P \to 0$$

to prove RR. [Hint: If $0 \to V_0 \to \cdots \to V_n \to 0$ is an exact sequence of finite dimensional vector spaces then $\sum_{i=0}^n (-1)^i \dim V_i = 0$. This exercise is carried out in [**H1**], Chapter IV.]

**24.** Use RR and Serre duality to prove that $\deg K_C = 2g - 2$ and $h^0(K_C) = g$.

To prove that the number $g$ appearing in RR is the same as the number $g$ in the famous picture of the surface with $g$ holes, you have to use that $T_C$ is the dual of $K_C$, so has degree $2-2g$, and some form of Gauss–Bonnet: the number of zeros of a regular vector field, counted with their indexes, equals the Euler characteristic. Any argument involving coherent cohomology (polynomials) on one side and topology on the other is automatically deeper than anything purely in algebraic geometry or purely in topology.

**25.** For points $P_1, \cdots, P_k \in \mathbb{P}^2$, write $h^0(\mathbb{P}^2, \mathcal{I}_{P_1+\cdots+P_k} \cdot \mathcal{O}(2))$ for the vector space of conics through $P_1, \cdots, P_k$, and $h^1(\mathbb{P}^2, \mathcal{I}_{P_1+\cdots+P_k} \cdot \mathcal{O}(2))$ for the space of linear dependence relations between the conditions $P_1, \dots, P_k$ impose on conics (compare the example in B.6). State and prove the results of [**UAG**], §1 on the dimension of the space conics through points $P_1, \cdots, P_k$ in terms of coherent cohomology groups

$$h^1(\mathbb{P}^2, \mathcal{I}_{P_1+\cdots+P_k} \cdot \mathcal{O}(2)).$$

**26.** Let $P_1, \dots, P_9 \in C \subset \mathbb{P}^2$ be 9 distinct points contained in a nonsingular cubic curve. Suppose that $h^1(\mathbb{P}^2, \mathcal{I}_{P_1+\cdots+P_9} \cdot \mathcal{O}(3)) \ne 0$. Prove that the surface $S = \mathrm{Bl}_{P_1,\dots,P_9} \mathbb{P}^2$ obtained by blowing up $P_1, \dots, P_9$ has an elliptic fibration $S \to \mathbb{P}^1$.

**27.** The sheafication $\mathrm{sh}(\mathcal{F})$ of a presheaf. If $\mathcal{F}$ is a presheaf, there is an *associated sheaf* or *sheafication* $\mathrm{sh}(\mathcal{F})$ which satisfies the universal mapping property for homomorphisms from $\mathcal{F}$ to a sheaf. Construct $\mathrm{sh}(\mathcal{F})$ and prove the universal mapping property. The idea is to consider the stalks $\mathcal{F}_P$, and set

$$\Gamma(U, \mathcal{F}) = \text{good maps } P \mapsto s_P \in \mathcal{F}_P \text{ for all } P \in U,$$

where "good" means that all the $s_Q$ for $Q$ in some small neighbourhood $V_P$ of $P$ are the restrictions of some $s \in \Gamma(V_P, \mathcal{F})$. If you have trouble with this question, refer to [**H1**], Chap. II or one of the books on sheaf theory.

**28.** If $\mathcal{F} \subset \mathcal{G}$ is a subsheaf, construct the quotient sheaf $\mathcal{G}/\mathcal{F}$ as the associated sheaf of the presheaf $U \mapsto \Gamma(U, \mathcal{G})/\Gamma(U, \mathcal{F})$, and prove that it has the universal mapping property for maps from $\mathcal{G}$ to a sheaf killing $\mathcal{F}$. Prove also that its stalks are $\mathcal{G}_P/\mathcal{F}_P$, so that the sequence $0 \to \mathcal{F} \to \mathcal{G} \to \mathcal{G}/\mathcal{F} \to 0$ is exact.

**29.** If $f\colon X \to Y$ is a continuous map of topological spaces and $\mathcal{F}$ is a sheaf on $Y$, construct the sheaf theoretic pullback $f^{-1}\mathcal{F}$, whose stalk at $P \in X$ is $\mathcal{F}_{f(P)}$. Prove that it has the universal mapping property for sheaves $\mathcal{G}$ on $X$ such that there exists a sheaf homomorphism $\mathcal{F} \to f_*\mathcal{G}$.

Incidentally, you mustn't write $f^*$ for $f^{-1}$, because $f^*$ is usually reserved for the pull back of sheaves of $\mathcal{O}_X$-module, given by $f^{-1}\mathcal{F} \otimes_{f^{-1}\mathcal{O}_Y} \mathcal{O}_X$.

# CHAPTER C. Guide to the classification of surfaces

The classification of surfaces goes back to Castelnuovo and Enriques in the early years of the 20th century. It divides algebraic surfaces into 4 big classes. There will be lots to say about surfaces in each of the 4 classes, and the characteristic methods and results for each class are different. The subject matter thus divides naturally into the 4 separate classes, together with the logical division into cases, or the *proof of classification.*

This chapter introduces the classification of surfaces in overall terms. It is primarily intended as a guide to results that lots of people need to use, and most of the results are given without proof. However, to give a veneer of modernity I cast the classification within the logical structure of Mori theory: the dichotomy $K_X$ nef or otherwise, and the numerical dimension of $K_X$ form the main logical framework, taking precedence over criteria for rationality and ruledness and Kodaira dimension. The underlying material is of course essentially the same as that of Castelnuovo, Enriques and Kodaira. I return to the proof of the classification in Chapters D–E below.

I assume that $X$ is a nonsingular projective surface over $\mathbb{C}$ to make life easier for myself and the reader. It is of course often reasonable to weaken these conditions: a topologist or analyst might want compact complex surfaces without the assumption of algebraicity, a number theorist might want singular projective surfaces over a finite field or a number field, students interested in higher dimensional geometry will need to know how the results for surfaces work in several other contexts. Getting the right level of generality is also important for the internal development of the subject.

## Summary

1. Numerical invariants
2. Birational versus biregular classification: blowups
3. Results of minimal model theory
4. Ruled and rational surfaces, Tsen's theorem
5. Classification of surfaces with $K_X$ nef by $\nu$
6. Kodaira dimension, the statement $\kappa = \nu$
7. More precise description of $\kappa = 0$

## C.1. Invariants

### C.1.1. Numerical invariants

The main numerical invariants are

1. the Betti numbers $B_i$ with the Euler number $e(X)$;
2. the signature $(B_2^+, B_2^-)$ of the quadratic form $Q_X$ on $H^2(X,\mathbb{Z})/\operatorname{Tors}$;
3. the Hodge numbers $h^{p,q}$ with the Poincaré characteristic $\chi(\mathcal{O}_X)$;
4. the Chern classes $c_1 = [-K_X] \in H^2(X,\mathbb{Z})$ and $c_2 \in H^4(X,\mathbb{Z}) = \mathbb{Z}$, and the Chern numbers $c_1^2$ and $c_2$.

Fortunately there are lots of relations between these numbers: only two or three of them are independent. There are also three other discrete invariants that are not simply numbers, but are closely related to the numerical invariants: the quadratic form $Q_X$ on $H^2(X,\mathbb{Z})/\operatorname{Tors} X$ (up to isomorphism), the torsion subgroup $\operatorname{Tors} X = \operatorname{Tors} H^2(X,\mathbb{Z})$ and the fundamental group $\pi_1(X)$. The intersection form on $H^2(X,\mathbb{Z})$ is completely determined from a knowledge of the numerical invariants (including $c_1 \in H^2(X,\mathbb{Z})$), and in turn, for simply connected surfaces, it is known that this determines the homotopy type and the homeomorphism type of $X$.

### C.1.2. Topological invariants

A complex surface is a compact 4-manifold with a chosen orientation. Accordingly, its Betti numbers $B_i = \operatorname{rank} H^i(X,\mathbb{Z})$ satisfy Poincaré duality

$$B_0 = B_4 = 1, \qquad B_1 = \operatorname{rank} H^1(X,\mathbb{Z}) = B_3.$$

Hence the Euler number is given by

$$e(X) = \sum(-1)^i B_i = 2 - 2B_1 + B_2.$$

(Hirzebruch's convention is to write $e(X)$ for the Euler characteristic, rather than something like $\chi_{\text{top}}(X,\mathbb{Z})$, to distinguish it from Euler–Poincaré characteristics such as $\chi(\mathcal{O}_X)$ in sheaf cohomology; I follow this throughout, and urge you to do likewise.) In addition, Poincaré duality says that the symmetric bilinear form defined by cup product

$$\cup\colon H^2(X,\mathbb{Z}) \times H^2(X,\mathbb{Z}) \to H^4(X,\mathbb{Z}) = \mathbb{Z}$$

is a perfect pairing on $H^2(X,\mathbb{Z})/\operatorname{Tors} X = \operatorname{Hom}(H_2(X,\mathbb{Z}),\mathbb{Z})$. This means, equivalently, that $\cup$ induces an isomorphism $H^2(X,\mathbb{Z})/\operatorname{Tors} X \to \operatorname{Hom}_{\mathbb{Z}}(H^2(X,\mathbb{Z}),\mathbb{Z})$, or that in any $\mathbb{Z}$-basis of $H^2(X,\mathbb{Z})/\operatorname{Tors} X$, the symmetric matrix representing the pairing is unimodular, that is, has determinant $\pm 1$.

### C.1.3. Analytic invariants

Most of the analytic invariants are defined in terms of the sheaves of differentials: recall that $T_X$ is a rank 2 complex vector bundle over $X$, whose sections are vector fields on (open subsets of) $X$; the dual bundle corresponds to the sheaf of Kähler differentials $\Omega^1_X = \operatorname{Hom}_{\mathcal{O}_X}(T_X,\mathcal{O}_X)$. This is a locally free sheaf of $\mathcal{O}_X$ modules of rank 2, based locally by $\mathrm{d}x, \mathrm{d}y$ where $x, y$ are local coordinates. The sheaf of canonical differentials $\Omega^2_X = \mathcal{O}_X(K_X)$ is the sheaf of holomorphic 2-forms, defined as $\bigwedge^2 \Omega^1_X$, and is the locally free sheaf of $\mathcal{O}_X$ modules of rank 1, based locally by $\mathrm{d}x \wedge \mathrm{d}y$.

Global sections of these bundles give us important invariants: the *geometric genus*

$$p_g(X) = h^0(X,\Omega^2_X) = \dim \Gamma(X,\Omega^2_X),$$

that is, the number of global holomorphic canonical differentials, and the *irregularity* $q = h^0(X,\Omega^1_X)$. *Hodge theory* allows us to express some of the other invariants of $X$ in terms of these: firstly, it defines a decomposition

$$H^1(X,\mathbb{Z}) \otimes_{\mathbb{Z}} \mathbb{C} = H^1(X,\mathbb{C}) = H^{0,1} \oplus H^{1,0},$$

where $H^{0,1} = H^1(X,\mathcal{O}_X)$ and $H^{1,0} = H^0(X,\Omega^1_X)$; in this decomposition $H^{0,1}$ is the complex conjugate of $H^{1,0}$, so that is particular

$$h^{0,1} = h^{1,0} = q \quad \text{and} \quad B_1 = \operatorname{rank} H^1(X,\mathbb{Z}) = \dim_{\mathbb{C}} H^1(X,\mathbb{C}) = 2q.$$

Next, it tells us the signature of the pairing on $H^2(X,\mathbb{R}) = H^2(X,\mathbb{Z}) \otimes \mathbb{R}$: diagonalising this quadratic form over $\mathbb{R}$ gives the pairing as a diagonal matrix with $B_2^+$ entries $+1$ and $B_2^-$ entries $-1$; the pair $(B_2^+, B_2^-)$ is the Sylvester inertia index of the form. Now Hodge theory tells us that $B_2^+ = 2p_g + 1$. In particular, $p_g$ is an invariant of the (oriented) homotopy type of $X$.

More precisely, the Hodge decomposition of $H^2(X,\mathbb{C})$ is

$$H^2(X,\mathbb{C}) = H^{0,2} \oplus H^{1,1} \oplus H^{2,0},$$

where $H^{2,0} = H^0(X,\Omega_X^2)$ and $H^{0,2}$ is the complex conjugate of $H^{2,0}$, so that the vector subspace

$$V_{\mathbb{C}}^+ = H^{2,0} \oplus H^{0,2} \subset H^2(X,\mathbb{C})$$

is invariant under complex conjugation, and corresponds to a $2p_g$-dimensional real subspace $V^+ \subset H^2(X,\mathbb{R})$. Then if $L \in H^2(X,\mathbb{Z})$ is the class of a hyperplane section of $X$ (the first Chern class of the line bundle $\mathcal{O}_X(1)$ for some embedding $X \subset \mathbb{P}^N$), by Hodge theory $L \in H^{1,1}$, and the Hodge index theorem says that the intersection form is positive definite on $V^+ \oplus \mathbb{R} \cdot L$, and negative definite on the "primitive part" $H_0^{1,1}$, the orthogonal complement of $L$ in $H^{1,1}$.

**Remark.** More generally, the set-up of Hodge theory on a Kähler manifold is

$$H^n(X,\mathbb{C}) = \oplus_{p+q=n} H^{p,q}(X),$$

where $H^{p,q}(X) = H^q(X,\Omega_X^p)$. One proves that $H^{p,q}$ and $H^{q,p}$ are complex conjugates, and determines the signature of the intersection pairing on the "primitive" cohomology, etc. What I said above is enough for complex projective surfaces.

### C.1.4. Chern classes

A crude consequence of the presence of the complex structure is that the tangent bundle $T_X$ has the structure of a complex vector bundle with fibre $\mathbb{C}^2$. This vector bundle has characteristic classes $c_1(T_X) \in H^2(X,\mathbb{Z})$ and $c_2(T_X) \in H^4(X,\mathbb{Z})$. Each of these objects has an alternative interpretation: since $\Omega_X^1$ is the dual of $T_X$ and $\Omega_X^2$ is the determinant bundle $\bigwedge^2 \Omega_X^1$,

$$c_1(T_X) = -c_1(\Omega_X^1) = -c_1(\Omega_X^2) = -[K_X] \in H^2(X,\mathbb{Z}).$$

Also, a general result is that $c_2(T_X) = e(X)$: morally speaking, the top Chern class of $T_X$ should be the number of zeros of a generic section, but the number of zeros of a generic vector field is the Euler number.

The final relation between the numerical invariants is Noether's formula:

$$\chi(\mathcal{O}_X) = 1 - q + p_g = \frac{1}{12}(c_1^2 + c_2) = \frac{1}{12}(K_X^2 + e(X));$$

here I interpret $c_1^2$ and $c_2 \in H^4(X,\mathbb{Z}) = \mathbb{Z}$ as integers. The formula is part of the set-up of Hirzebruch–RR. Since $c_2(X) = e(X) = 2 - 4q + B_2$, and since $q$ and $p_g$ are homotopy invariants, so is $c_1^2$.

### C.1.5. Exercises

1. For regular surfaces (that is, surfaces with $q = 0$), show how to express all the numerical invariants in terms of $K_X^2$ and $\chi = \chi(\mathcal{O}_X)$, and verify the Hirzebruch signature formula
$$B_2^+ - B_2^- = 4\chi - e(X) = \frac{1}{3}(c_1^2 - 2c_2).$$
2. Assume that $c_1^2 \geq 0$, $c_2 \geq 0$ and $c_1^2 \leq 3c_2$. (Lots of surfaces satisfy these conditions.) Prove that $B_2^- = 0$ is only possible if $X$ satisfies $q = 0$, $B_2 = 1$, $c_1^2 = 9$ and $c_2 = 3$ (that is, $X$ has the numerical invariants of $\mathbb{P}^2$). This proves that the quadratic form on $H^2$ is either $\mathbf{1}$ or is indefinite.

### Solutions

1. $q = 0$, $p_g = \chi - 1$, $e(X) = 12\chi - K_X^2$,
$B_2 = e(X) - 2$, $B_2^+ = 2p_g + 1 = 2\chi - 1$, $B_2^- = e(X) - 2\chi - 1$.
2. $B_2^- = B_2 - 2p_g - 1$
$= e(X) - 2 + 4q - 2p_g - 1$
$= c_2 - \frac{1}{6}(c_1^2 + c_2) + 2q - 1$
$= \frac{1}{3}c_2 + \frac{1}{6}(3c_2 - c_1^2) + 2q - 1.$
So $B_2^- = 0$ gives $q = 0$ and $2c_2 + (3c_2 - c_1^2) = 6$, which has only one solution compatible with $c_1^2 + c_2 \equiv 0 \bmod 12$.

### C.1.6. Quadratic form on $H^2(X, \mathbb{Z})$

Consider a quadratic form $(H, Q)$ over $\mathbb{Z}$, that is, a free $\mathbb{Z}$-module (lattice) $H \cong \mathbb{Z}^n$, together with a quadratic map $Q\colon H \to \mathbb{Z}$. There is a whole theory of isomorphism classes of quadratic forms over $\mathbb{Z}$; let me start by going through the easy bits of this theory; a nondegenerate quadratic form $(H, Q)$ has a *signature* $(n^+, n^-)$, a *discriminant* $\det Q$, and a *parity*. Parity just means the following: if $Q(x) \in \mathbb{Z}$ is even for all $x \in H$ then the quadratic form is *even*, otherwise it's *odd*.

For $H^2(X, \mathbb{Z})$ of a complex surface, all these invariants are under control: the signature $(B_2^+, B_2^-)$ is discussed in C.1.3 above; the quadratic form $Q_X$ is unimodular (that is, $\det Q_X = 1$) by Poincaré duality; and the parity is determined by $c_1(X)$ or $K_X$:
$$Q_X(x) \equiv c_1 \cdot x \bmod 2 \quad \text{for all } x \in H^2(X, \mathbb{Z}).$$
Hence $Q_X$ is even if and only if $c_1(X)$ maps to $0 \in H^2(X, \mathbb{Z}/2)$, or equivalently, $K_X$ is divisible by 2 as a divisor class. A surface having this property is *even* or is a *spin surface.*

A unimodular quadratic form $(H, Q)$, assumed to be indefinite, is determined up to isomorphism by its signature $(n^+, n^-)$ and its parity. In fact if $Q$ is odd then it is isomorphic to the diagonal form:
$$Q \text{ is odd} \implies (H, Q) \cong \mathbf{1}^{n^+} \oplus (-\mathbf{1})^{n^-},$$
and there is a similar standard form for an even indefinite form
$$Q \text{ is even} \implies (H, Q) \cong \begin{pmatrix} 0 & 1 \\ 1 & 0 \end{pmatrix}^a \oplus (\pm E_8)^b,$$

where $\begin{pmatrix} 0 & 1 \\ 1 & 0 \end{pmatrix}$ is the standard hyperbolic lattice, $E_8$ is the lattice based by 8 vector $e_i$ in bijection with the vertexes of the Dynkin diagram $E_8$ and quadratic form given by $e_i^2 = 2$, $e_ie_j = 0$ or $-1$ according as the corresponding vertexes are unjoined or joined, and $-E_8$ the same with the opposite signs $e_i^2 = -2$, $e_ie_j = 0$ or $+1$. That is, $\pm E_8 = \mathbb{Z}^8$ with the quadratic form

$$\pm\begin{pmatrix} -2 & 0 & 0 & 1 & 0 & 0 & 0 & 0 \\ 0 & -2 & 1 & 0 & 0 & 0 & 0 & 0 \\ 0 & 1 & -2 & 1 & 0 & 0 & 0 & 0 \\ 1 & 0 & 1 & -2 & 1 & 0 & 0 & 0 \\ 0 & 0 & 0 & 1 & -2 & 1 & 0 & 0 \\ 0 & 0 & 0 & 0 & 1 & -2 & 1 & 0 \\ 0 & 0 & 0 & 0 & 0 & 1 & -2 & 1 \\ 0 & 0 & 0 & 0 & 0 & 0 & 1 & -2 \end{pmatrix}$$

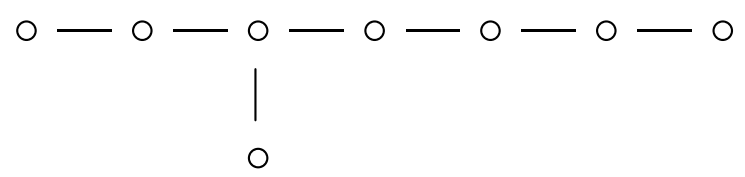

Definite quadratic forms over $\mathbb{Z}$ are of course much more complicated, but fortunately, we never have to deal with them (compare Exercise 1.5 above).

### C.1.7. Exercise

All these invariants are usually easy enough to calculate if your surface is given fairly explicitly. For example, suppose that $X = X_d \subset \mathbb{P}^3_{\mathbb{C}}$ is a nonsingular hypersurface of degree $d$. It is simply connected by Lefschetz theory, and in particular $q = 0$. The adjunction formula gives $K_X = (-4+d)H = \mathcal{O}_X(-4+d)$, so that

$$p_g = \begin{cases} 0 & \text{for } d \le 3 \\ \binom{d-1}{3} & \text{for } d \ge 4 \end{cases} \quad \text{and } K_X^2 = d(d-4)^2.$$

Finally, $K_X$ is divisible by 2 in $\operatorname{Pic} X$ if and only if $d-4$ is even, that is, $d$ is even. It's an easy exercise to determine all the invariants of $X$ from this.

## C.2. Birational versus biregular classification: blowups

There are several different points of view on what may be meant by the classification of algebraic varieties, among them the following three:

1. birational, that is, up to birational equivalence;
2. biregular, that is, up to isomorphism;
3. projective, that is, up to projective equivalence.

To say that $X$ is a projective surface means that $X$ is capable of being embedded as a closed subvariety of $\mathbb{P}^N$ for some $N$. The gap between (2) and (3) is the problem of finding a projective embedding of an abstract surface $X$, or separating off intrinsic properties of $X \subset \mathbb{P}^N$ from the properties of the ambient space and of the embedding. I leave this for the moment. Most of this section and the next is concerned with the gap between (1) and (2).

Notice that for nonsingular projective curves, (1) = (2): for $C$ and $C'$ any two nonsingular projective curves, any birational map from $C$ to $C'$ is an isomorphism. The blowup (see 1.10) shows that this is not the case for surfaces. The main point of this section is that blowups account for the whole of the gap between (1) and (2).

**Remark.** For surfaces having a minimal model with $K_X$ nef, it turns out that the minimal model is unique, which restores (1) = (2). A birational map between two surfaces with $K_X$ nef is necessarily an isomorphism.

**C.2.1. Proposition.** (a) Blowups. *Given a surface $X$ and a point $P \in X$, there exists a surface $X_1$ and a morphism $\sigma\colon X_1 \to X$ such that $\sigma^{-1}P = E$ is a curve, and $\sigma$ restricts to an isomorphism $X_1 \setminus E \to X \setminus P$.*

*Here $E \cong \mathbb{P}^1$ and $E^2 = -1$, so that $E$ is a $-1$-curve; if $C$ is a curve of $X$ then $\sigma^*C = C' + mE$, where $m = \operatorname{mult}_P C$, $C'$ is the birational transform of $C$ (that is, the closure in $X_1$ of $\sigma^{-1}(C \setminus P)$), and the points of $C'$ over $P$ correspond one-to-one to the distinct tangent lines to $C$.*

(b) Castelnuovo's contractibility criterion. *Given a surface $X$ and a $-1$-curve $E \subset X$, there exists a contraction $\sigma\colon X \to Y$ such that $E$ maps to a point $P$ of a (non-singular projective) surface $Y$, and $\sigma$ is the blowup of $P \in Y$.*

For (a), see 1.10. (b) is proved in Chapter 4.

**C.2.2. Theorem.** (a) Resolution of indeterminacies. *Let $f\colon X \dashrightarrow \mathbb{P}^N$ be a rational map. Then there exists a commutative diagram*

$$\begin{array}{ccc} & Y & \\ {}^{g}\swarrow & & \searrow^{h} \\ X & \dashrightarrow & \mathbb{P}^N \end{array}$$

*where $g$ is a composite $Y = X_n \to \cdots \to X_1 \xrightarrow{\sigma_1} X$ of blowups, and $h$ is a morphism.*

(b) Factorisation of birational morphisms. *Let $h\colon X \to Y$ be a birational morphism between nonsingular projective surfaces. Then $h$ is the composite of a chain of blowups $X \xrightarrow{\sigma_1} Y_1 \to \cdots \to Y_N = Y$.*

### C.2.3. Topological view of a blowup

For surfaces over $\mathbb{C}$, a blowup $\sigma\colon X_1 \to X$ corresponds to taking a connected sum $X\#Q$ of $X$ with the 4-manifold $Q = \overline{\mathbb{P}^2_{\mathbb{C}}}$, that is, $\mathbb{P}^2_{\mathbb{C}}$ with the opposite orientation.

The idea is that it replaces a small disc around $P$ by a tubular neighbourhood of the $-1$-curve $E$; it is not hard to see that this is diffeomorphic to $\mathbb{P}^2_{\mathbb{C}} \setminus$ (disc), but with the opposite orientation (both of them are disc bundles over $\mathbb{P}^1_{\mathbb{C}} \cong S^2$, but the "core" copy of $\mathbb{P}^1_{\mathbb{C}}$ has selfintersection $\pm1$ in the two cases).

### C.2.4. Fundamental asymmetry

It's important to note the asymmetry here: taking connected sum with $Q$ is a blowup, so is an operation of algebraic geometry. But taking connected sum with $\mathbb{P}^2_{\mathbb{C}}$ itself cannot be an operation of algebraic geometry: it can never happen that $X$ and $X\#\mathbb{P}^2_{\mathbb{C}}$ are both diffeomorphic to algebraic surfaces (with oriented connected sum and diffeomorphism preserving the orientation), if only because $B_2^+$ passes from odd to even.

At a much deeper level, the contrast can be seen in a very striking way from the following two results:

**Theorem** (Moishezon–Mandelbaum). *For most naturally occurring large classes of algebraic surfaces, the connected sum $M\#\mathbb{P}^2_{\mathbb{C}}$ has a decomposition*

$$M\#\mathbb{P}^2_{\mathbb{C}} \overset{\text{diffeo}}{\sim} (\mathbb{P}^2_{\mathbb{C}})^{\#n}\#Q^{\#m}.$$

**Theorem** (Donaldson). *If $M$ is the 4-manifold underlying an algebraic surface, there does not exist any connected sum decomposition*

$$M \overset{\text{diffeo}}{\sim} M_1\#M_2$$

*with each of $M_1, M_2$ having $B_2^+ > 0$.*

## C.3. Results of minimal model theory

The following is a modern statement of minimal model theory. It is proved in Chapter D. In the classical theory, it would be a corollary of a whole chain of results. Recall that the *canonical divisor* $K_X$ is the divisor corresponding to the invertible sheaf $\Omega^2_X = \mathcal{O}_X(K_X)$ and that a divisor $D$ on a surface $X$ is *nef* ("numerically eventually free") if $DC \geq 0$ for every curve $C \subset X$.

**C.3.1. Theorem.** *Let $Y$ be any nonsingular projective surface. Then there is a chain*

$$Y \xrightarrow{\sigma_1} Y_1 \to \cdots \to Y_N = X$$

*such that each $\sigma_i\colon Y_{i-1} \to Y_i$ is the contraction of a single $-1$-curve $E_i$, and $X$ satisfies*

*either (i) $K_X$ is nef;*

*or (ii) $X \cong \mathbb{P}^2$ or a $\mathbb{P}^1$-bundle over a curve.*

More crudely, every surface is birational to a surface with $K_X$ nef, or to a $\mathbb{P}^1$-bundle over a curve or $\mathbb{P}^2$. This form of the result is easy to remember, and is very useful. On the other hand, the information it contains can be analysed to give results that are more complicated, but much more general. Notice that as stated, the theorem refers to three different situations:

1. a contraction $\sigma\colon X \to X_1$ of a $-1$-curve $E$, with $E \cong \mathbb{P}^1$, $E^2 = -1$, $K_XE = -1$;
2. a $\mathbb{P}^1$-bundle $X \to C$ with fibre $F \cong \mathbb{P}^1$, $F^2 = 0$, $K_XF = -2$;
3. $X \cong \mathbb{P}^2$; here I can write the constant morphism $\mathbb{P}^2 \to$ pt., and consider a line $L \subset \mathbb{P}^2$ with $L \cong \mathbb{P}^1$, $L^2 = 1$, $K_XL = -3$.

All three of these situations can be classed together as a morphism $\varphi\colon X \to Y$ such that $-K_XC > 0$ for every curve $C$ in a fibre of $f$. Thus Theorem C.3.1 is closely related to the following statement.

**C.3.2. Theorem.** *Let $X$ be a surface, and suppose that $K_X$ is not nef. Then there exists a morphism $\varphi\colon X \to Y$ (with $\dim Y = 0, 1$ or $2$) contracting at least one curve of $X$ to a point, and such that $-K_XC > 0$ for every curve $C$ in a fibre* *of $\varphi$. There is no loss of generality in assuming $\varphi_*\mathcal{O}_X = \mathcal{O}_Y$, which means that $Y$ is normal and the general fibre of $\varphi$ is connected.*

Theorem C.3.2 is more primitive than Theorem C.3.1, but much more general in scope: it holds for surfaces over any field $k$, for a large class of singular surfaces, and for nonsingular varieties of any dimension over $\mathbb{C}$. The modern strategy for proving Theorem C.3.1 is first to prove Theorem C.3.2, then to analyse the morphism $\varphi\colon X \to Y$ for which $-K_X$ is ample.

**C.3.3. Proposition.** *Let $\varphi\colon X \to Y$ be as in the conclusion of Theorem 2. Then one of the following 3 cases holds.*

$\dim Y = 2$. *Then $Y$ is a nonsingular surface, and $\varphi\colon X \to Y$ is the blowdown of a number of disjoint $-1$-curves;*

$\dim Y = 1$. *Then $Y$ is a nonsingular curve, and $\varphi\colon X \to Y$ is a conic bundle; the general fibre of $\varphi$ is isomorphic to $\mathbb{P}^1$, whereas a special fibre is isomorphic to a line pair in $\mathbb{P}^2$;*

$Y = \mathrm{pt}$. *In this case, the constant morphism $\varphi\colon X \to Y$ does not give any information, but $-K_X$ is an ample divisor.*

Surfaces in the final case are called *del Pezzo surfaces*, and they are similar in most of their properties to the cubic surfaces $X_3 \subset \mathbb{P}^3$ discussed in Chapter 1. The main step in deducing Theorem C.3.1 from Theorem C.3.2 is to prove that if $X$ is a del Pezzo surface then either $B_2(X) \geq 2$ and $X$ admits another contraction $X \to X'$ of one of the other two types; or $B_2(X) = 1$ and $X \cong \mathbb{P}^2$.

## C.4. Ruled and rational surfaces, Tsen's theorem

I define a *ruled surface* to mean a nonsingular projective surface $X$ together with a fixed morphism $\varphi\colon X \to B$ to a base curve $B$, whose fibre $F_{k(B)}$ over the generic point is a nonsingular curve of genus 0. The main theorem of this section, proved below, is Tsen's theorem: $F_{k(B)}$ is isomorphic over $k(B)$ to $\mathbb{P}^1_{k(B)}$. To discuss the results of minimal model theory, I assume this for the moment. That is, I assume that $\varphi\colon X \to B$ is a given morphism with generic fibre isomorphic to $\mathbb{P}^1_{k(B)}$.

There are many alternatives, going back to the ancients: you could require only that $\varphi$ exists, without taking the responsibility for specifying it, or you could allow $\varphi$ to be a rational map $\varphi\colon X \dashrightarrow B$. For the purposes of this introduction, call this "birationally ruled". Thus the plane $\mathbb{P}^2$ or any rational surface is birationally ruled. The $\kappa = -\infty$ part of the classification of surfaces (due to Castelnuovo and Enriques) can be stated in the form

$$P_{12}(X) = 0 \iff \kappa = -\infty \iff X \text{ is birationally ruled.}$$

Compare Corollary E.1.

### C.4.1. Models

According to the definition, $X$ has a morphism $X \to B$ with generic fibre $\mathbb{P}^1_{k(B)}$. Choosing an isomorphism of the generic fibre defines a birational map

$$\begin{array}{ccc} X & \dashrightarrow & B \times \mathbb{P}^1 \\ & \searrow \quad \swarrow & \\ & B & \end{array}$$

Since the definition of ruled surface only fixes the generic fibre, the general picture is as follows: $X$ is isomorphic to $B \times \mathbb{P}^1$ over a Zariski open set of $B$, but over finitely many points $\{P\} \in B$ the fibre can have lots of components, say $\varphi^{-1}P = \bigcup \Gamma_i$. We know that the intersection matrix $(\Gamma_i\Gamma_j)_{ij}$ is negative semidefinite (see Theorem A.7 and Ex. A.17). It's an exercise to see from first principles that if the fibre is not $\mathbb{P}^1$ then at least one of the $\Gamma_i$ is a $-1$-curve.

There are other conditions that can be added to a ruled surface that restrict the types of singular fibres. These conditions appear in a natural way in classical work, and more especially in Mori theory:

1. $-K_X$ ample characterises conic bundles. Every singular fibre of $\varphi$ is a line pair, a union of two $-1$-curves $L_1 \cup L_2$ meeting transversally at one point. In fact $-K_X$ is relatively very ample, and embeds $X$ into a $\mathbb{P}^2$-bundle over $B$.
2. $-K_X$ ample and $\rho(\varphi) = 1$ characterise $\mathbb{P}^1$-bundle. The condition $\rho(\varphi) = 1$ means that every curve in a fibre is numerically a multiple of the fibre. This implies that every fibre is nonsingular. Using Tsen's theorem, it's not hard to see that then $X = \mathbb{P}(E)$ is the projectivisation of a rank 2 vector bundle over $B$.

The case of conic bundle is important over an algebraically nonclosed field $k$: it can happen that $X \to B$ has $-K_X$ ample, so is geometrically a conic bundle, and has $\rho(\varphi) = 1$, so that there are no curves defined over $k$ in the fibres of $X \to B$ that can be contracted out by a birational morphism; but $X \to B$ can still have singular fibres, line pairs $L_1 \cup L_2$ with $L_1$ and $L_2$ conjugate over $k$.

**C.4.2. Theorem** (Tsen's theorem). *Let $k$ be an algebraically closed field and $k \subset F$ a function field in one variable over $k$ (this means that $F = k(B)$, where $B$ is an irreducible algebraic curve over $k$). Let $X_F$ be a nonsingular projective curve over $F$, and assume that $X_F$ is absolutely irreducible, in the sense that it remains irreducible over the algebraic closure of $F$. Then the following conditions are equivalent:*

1. *The curve $X_F$ has $g = 0$ in RR;*
2. *$X_F$ is isomorphic over $F$ to a nonsingular plane conic;*
3. *$X_F$ becomes isomorphic to $\mathbb{P}^1$ after some field extension $F \subset F'$, that is, $X_F \otimes_F F' \cong \mathbb{P}^1_{F'}$;*
4. *$X_F$ is isomorphic over $F$ to $\mathbb{P}^1_F$.*

**Proof.** $2 \implies 3$ and $4 \implies 1$ are obvious. The implication $1 \implies 2$ is standard use of RR: the anticanonical line bundle $-K_{X_F}$ has degree 2, and is very ample by the usual criterion: $\ell(-K_{X_K}) = 3$, and $\ell(-K_{X_K} - D) = 1$ for any divisor $D$ of degree 2. All this works over an arbitrarily field $F$.

The point is thus to prove $2 \implies 4$. I explain how to use the assumptions that $F = k(B)$ is a function field in one variable over an algebraically closed field. The point is that $k(B)$ is the field of fractions of a ring $k[B]$ with "linear growth", namely the coordinate ring of an affine model of $B$. In more detail, $k[B] = \bigcup k[B]_d$ (polynomials of degree $\le d$) with multiplication satisfying $k[B]_{d_1} \cdot k[B]_{d_2} \subset k[B]_{d_1+d_2}$. Moreover, $k[B]_d$ is finite dimensional over $k$, with dimension growing as a linear function: for some constants $a$, $b_1, b_2$

$$ad - b_1 \le \dim k[B]_d \le ad + b_2 \quad \text{(same } a!\text{)}. \tag{1}$$

There are various ways of proving this, using RR on the curve $B$, or the commutative algebra definition of $\dim k[B] = 1$ in terms of Hilbert–Samuel functions, but in any case the fact that $B$ is 1-dimensional is crucial.

Now $X_F$ is isomorphic over $k(F)$ to a plane conic

$$X_F : q(x,y,z) = 0 \subset \mathbb{P}^2_F \tag{2}$$

defined by a quadratic form $q$ in $x, y, z$ with coefficients in $F$. Notice that since $X_F$ is a plane conic, to prove that it is isomorphic to $\mathbb{P}^1_F$, it is enough to prove that is has a point defined over $F$ (the argument is given in [**UAG**], (1.7) and Ex. 1.5).

The point is thus very simple: under the assumption that $F$ is a function field in 1 variable, I can always solve (2). Multiplying through by denominators, I can assume that the coefficients of $q$ are in $k[B]$, and even in $k[B]_c$ for some choice of $c$. The plan of the proof is to try to solve (2) with $x, y, z \in k[B]_d$ subject to $q(x,y,z) = 0 \in k[B]_{2d+c}$.

The punch-line is that by (1) there are $\geq 3(ad - b_1)$ free variables in the choice of $x, y, z$, and $\leq a(2d+c) + b_2$ polynomial conditions on them in the equation $q(x,y,z) = 0 \in k[B]_{2d+c}$. To explain this gently, $x, y, z$ are any elements of the vector space $k[B]_d$; if $e_1, \dots, e_N$ is a basis then

$$x = \sum u_i e_i, \quad y = \sum v_i e_i, \quad z = \sum w_i e_i$$

give $3N$ variable $u_i, v_i, w_i$. In a similar way, the equation $q(x,y,z) = 0$ can be written out as $\dim k[B]_{2d+c}$ quadratic equations in these variables.

Since for $d \gg 0$, there are more free variables than equations; over the algebraically closed $k$, it follows that $u_i, v_i, w_i$ can be given values in $k$ so that $x, y, z$ satisfy $q(x,y,z) = 0$.

### C.4.3. Birational maps

If $X$ and $Y$ are rational surfaces, the set of birational maps $X \dashrightarrow Y$, or the group $\operatorname{Bir} X$ of birational selfmaps is very big. This set (resp. group) is the obstruction to the uniqueness of the model of $X$ as a $\mathbb{P}^1$-bundle over a curve or $\mathbb{P}^2$. For $X = C \times \mathbb{P}^1$, $\operatorname{Bir} X$ contains $\mathrm{PGL}(2, k(C))$. For $X = \mathbb{P}^2$, the group $\operatorname{Bir} X$ is the Cremona group, and is qualitatively something like the free group on a continuum of generators.

## C.5. Classification assuming $K_X$ nef

### Aim

Assuming $X$ a surface with $K_X$ nef, classify into 3 cases according to the numerical properties of $K_X$ (the 3 possibilities $\nu = 0, 1, 2$ for the numerical dimension of $K_X$), and according to the analytic properties of $K_X$ (the 3 possibilities $\kappa = 0, 1, 2$, for the Kodaira dimension of $K_X$); statement of $\kappa = \nu$.

### C.5.1. General fact. $D$ nef $\implies D^2 \geq 0$

Recall that by definition, $D$ is nef if and only if $DC \geq 0$ for all curves $C \subset X$, so that $D \cdot (\sum n_i C_i) \geq 0$ for all effective divisors $\sum n_i C_i$ with $n_i \geq 0$. See D.2 for a proof.

Morally, the reason C.5.1 holds is that

$$D \text{ nef} \implies D \text{ is “close to being effective”} \implies D^2 \geq 0.$$

### C.5.2. Definition of $\overset{\text{num}}{\sim}$

Two divisors $D_1$ and $D_2$ on a surface $X$ are numerically equivalent (written $D_1 \overset{\text{num}}{\sim} D_2$) if $D_1 C = D_2 C$ for all curves $C \subset X$. Notice that

$$D_1 \overset{\text{lin}}{\sim} D_2 \iff \mathcal{O}_X(D_1) \cong \mathcal{O}_X(D_2) \implies D_1 \overset{\text{num}}{\sim} D_2$$

because $DC = \deg_C \mathcal{O}_X(D)$.

### C.5.3. Division into cases according to $\nu$

Assume $K_X$ is nef, so that $K_X^2 \geq 0$; then the following 3 cases are all-inclusive and mutually exclusive:

$\nu = 0$. $K_X \overset{\text{num}}{\sim} 0$, that is, $K_X C = 0$ for all $C \subset X$;

$\nu = 1$. $K_X C > 0$ for some curve $C \subset X$, but $K_X^2 = 0$;

$\nu = 2$. $K_X^2 > 0$.

You should think of $\nu$ as being defined by this division into cases. Formally, one can write down the definition of $\nu$ as a formula as follows:

$$\nu = \max\bigl\{k \bigm| K_X^k \overset{\text{num}}{\not\sim} 0\bigr\},$$

although this has no other meaning than the case division just given.

**C.5.4. Theorem** (Weak form of main theorem). *There exists a morphism $\varphi\colon X \to Y$ with $\dim Y = \nu(X)$ and $\varphi_*\mathcal{O}_X = \mathcal{O}_Y$ (that is, $Y$ is normal and the generic fibre of $\varphi$ is connected), such that*

$$K_X C = 0 \iff \varphi(C) = pt.$$

*There are 3 separate statements here:*

*$\nu = 0$. The statement is vacuous: $\varphi$ maps $X$ to a point.*

*$\nu = 1$. This is the most substantial case; here $\varphi\colon X \to Y$ is a fibration of $X$ over a curve, with $K_X(\text{fibre}) = 0$. A general fibre $E$ is a nonsingular curve with $E^2 = 0$, so that $K_X E = 0$ implies that $g(E) = 1$, and thus $\varphi$ is an elliptic fibration.*

*$\nu = 2$. In this case the morphism $\varphi\colon X \to Y$ is birational, and contracts at most a finite number of connected configurations of $-2$-curves to Du Val singularities; also, $\mathcal{O}_Y(K_Y)$ is an ample line bundle, and $K_X = \varphi^* K_Y$.*

### C.5.5. Definition of Kodaira dimension

The vector spaces $\Gamma(X, (\Omega_X^2)^m) = \Gamma(\mathcal{O}_X(mK_X))$ provide further invariants of $X$ and $K_X$. Define the *plurigenera* $P_m$ of $X$ by $P_m(X) = \dim \Gamma(X, \mathcal{O}_X(mK_X))$. It is formal to see that the $P_m$ must satisfy one of the following:

Case $\kappa = -\infty$. $P_m = 0$ for all $m > 0$;

Case $\kappa = 0$. $P_m \leq 1$ for all $m > 0$, and $P_m = 1$ for some $m$;

Case $\kappa = 1$ or 2. There exist constants $a, b > 0$ such that $am^\kappa < P_m < bm^\kappa$ for all sufficiently large $m$.

It can be shown that if $\kappa \geq 0$ then $\kappa = \dim \varphi_{mK_X}(X)$, where if $P_m \neq 0$ then $\varphi_{mK_X}\colon X \dashrightarrow \mathbb{P}^{P_m - 1}$ is the rational map defined by $\mathcal{O}_X(mK_X)$.

Note that in contrast to $\nu$, the invariants $P_m$ and $\kappa$ are analytic invariants of $X$: $P_m$ is the number of global holomorphic $m$-times canonical differentials on $X$.

### C.5.6. Statement of $\kappa = \nu$

The statement $\kappa = \nu$ is true, and is one of the central structural results of the theory. It means the following:

Case $\nu = 0$. $K_X \overset{\text{num}}{\sim} 0 \implies \exists m$ such that $mK_X \overset{\text{lin}}{\sim} 0$.

Case $\nu = 1$. $K_X \overset{\text{num}}{\not\sim} 0$ but $K_X^2 = 0 \implies \exists m$ such that $P_m \geq 2$.

Case $\nu = 2$. $K_X^2 > 0 \implies P_m$ grows quadratically in $m$. (This is easy.)

There is a more precise description of case $\nu = 0$.

**C.5.7. Theorem.** *If $\nu$=0 then one of the following 4 cases holds:*

1. $K_X = 0$, $H^1(\mathcal{O}_X) = 0$, *that is, $X$ is a K3 surface (see Chapter 3);*
2. $K_X = 0$, $h^1(\mathcal{O}_X) = 2$, *and $X$ is an Abelian surface.*
3. $K_X \overset{\text{lin}}{\not\sim} 0$ *but* $2K_X \overset{\text{lin}}{\sim} 0$ *and* $H^1(\mathcal{O}_X) = 0$*; these conditions are the definition of an Enriques surface.*
4. $K_X \overset{\text{lin}}{\not\sim} 0$ *but* $mK_X \overset{\text{lin}}{\sim} 0$ *for* $m = 2, 3, 4$ *or* $6$, *and* $h^1(\mathcal{O}_X) = 1$. *Then $X$ is a bielliptic surface, that is, $X = E_1 \times E_2/(\mathbb{Z}/m)$ is a quotient of a product of two elliptic curves by the cyclic group $\mathbb{Z}/m$ acting by translations on one factor and a group automorphism in the other.*

The results sketched in this chapter are taken up again with complete proofs in Chapters D–E below.

# CHAPTER 3. K3s

Let $X$ be a nonsingular projective surface. The formal definition of K3 surface is $K_X = 0$ and $H^1(\mathcal{O}_X) = 0$. As I discuss presently, the historic geometric definition of Enriques and Fano is closely related: $X$ is a K3 if $X \subset \mathbb{P}^g$ has a hyperplane section that is a canonical curve $C \subset \mathbb{P}^{g-1}$. Linear systems on K3s are closely related to the geometry of canonical curves. For example, the dichotomy between hyperelliptic curves and canonically embedded curves carries over to K3s. There is a general principle that a "sufficiently special" linear system $g^r_d$ on a curve $C \subset X$ is cut out by something on the surface. Exceptional classes of K3s (monogonal, hyperelliptic, trigonal, tetragonal, etc.) have natural models in terms of divisors in scrolls.

K3s occupy a special place in the classification of surfaces: a surface $X$ with $K_X = 0$ is either a K3, a special kind of Enriques surface (in characteristic 2 only), or an Abelian surface (that is, $X$ is a commutative algebraic group, with trivial tangent bundle $T_X \cong \mathbb{C}^2$, and hence $h^1(\mathcal{O}_X) = 2$); thus you could make your own original definition of K3 by taking $K_X = 0$, together with any other conditions that exclude Abelian surfaces (such as $\pi_1(X) = 0$).

K3s also occupy a special place in the curriculum for anyone trying to master algebraic surfaces. Although considerably simpler than surfaces of general type, they are a marvellous testing ground for your understanding of linear systems, cohomology, vanishing theorems, the structure sheaf $\mathcal{O}_D$ of a nonreduced divisor $D$, the relation between geometry of linear systems and the algebra of graded rings, singularities, intersection numbers of curves and quadratic forms, Hodge structures, moduli, and many other things.

Any treatment of K3s would be incomplete without a discussion of their moduli, both algebraic and analytic, and their Hodge theory and period map. Although outside the scope of these notes, I mention these topics briefly.

## Summary

1. Canonical class, adjunction
2. Canonical curves and RR for curves, geometric form of RR
3. Historic definition of K3
4. RR for surfaces and its "proof", how it simplifies for a K3, proof of Hodge algebraic index theorem, arithmetic genus of a divisor $p_a D$
5. Easy properties of linear systems on a K3, nef and free, elliptic pencils and monogonal divisors
6. Numerically connected divisors, Ramanujam vanishing
7. Projective embeddings of K3s according to Saint-Donat: hyperelliptic, trigonal, etc., special K3s as divisors in scrolls
8. Special linear systems on K3 sections
9. Analytic theory, moduli

## 3.1. Restriction and adjunction

There are two memorable exact sequences associated with a nonsingular $n$-fold $X$ and a codimension 1 subvariety $Y \subset X$. First

$$0 \to \mathcal{I}_{X,Y} \to \mathcal{O}_X \to \mathcal{O}_Y \to 0, \tag{1}$$

where $\mathcal{I}_Y = \mathcal{I}_{X,Y} = \mathcal{O}_X(-Y)$ is the ideal sheaf defining $Y$; it is an invertible sheaf (line bundle) on $X$. The conormal sheaf $(\mathcal{N}_{X|Y})^\vee = \mathcal{I}_Y/\mathcal{I}_Y^2$ is the restriction of this invertible sheaf to $Y$, also written $\mathcal{I}_{Y|Y} = \mathcal{I}_Y \otimes_{\mathcal{O}_X} \mathcal{O}_Y$ or $\mathcal{O}_Y(-Y)$.

Tensoring (1) with $\mathcal{O}_X(Y)$ gives

$$0 \to \mathcal{O}_X \to \mathcal{O}_X(Y) \to \mathcal{O}_Y(Y) \to 0, \tag{2}$$

where $\mathcal{O}_Y(Y) = \mathcal{O}_X(Y)_{|Y} = N_{X|Y}$ is the normal bundle.

Second,

$$0 \to T_Y \to T_{X|Y} \to N_{X|Y} \to 0. \tag{3}$$

Taking determinants gives

$$\det T_X = \det T_Y \otimes N_{X|Y}, \qquad \text{that is,} \qquad K_Y = (K_X + Y)_{|Y}. \tag{4}$$

(4) is called the *adjunction formula.* Note that the formula

$$2g(C) - 2 = (K_X + C)C$$

for the genus of a nonsingular curve $C \subset X$ on a surface is a particular case. Alternative treatments of adjunction are given in Ex. 3.25–26.

**Example 1.** $K_{\mathbb{P}^n} = \mathcal{O}(-n-1)$. Therefore if $X \subset \mathbb{P}^n$ is a hypersurface of degree $d$ then $K_{\mathbb{P}^n} = \mathcal{O}_{\mathbb{P}^n}(-(n+1))$ and $K_X = \mathcal{O}_X(d-(n+1))$. If $X$ is the complete intersection of two hypersurfaces of degrees $d_1, d_2$ then $K_X = \mathcal{O}_X(d_1+d_2-(n+1))$, etc.

**Example 2.** Let $X = X_4 \subset \mathbb{P}^3$ be a nonsingular quartic surface. Then $K_X = 0$. Also, since $H^i(\mathbb{P}^n, \mathcal{O}_{\mathbb{P}^n}(k)) = 0$ for all $0 < i < n$, it follows from the cohomology exact sequence

$$0 \to \mathcal{O}_{\mathbb{P}^3}(-4) \to \mathcal{O}_{\mathbb{P}^3} \to \mathcal{O}_X \to 0$$

that $H^1(\mathcal{O}_X) = 0$. Thus $X$ is a K3. Similarly the complete intersections $X_{2,3} \subset \mathbb{P}^4$ and $X_{2,2,2} \subset \mathbb{P}^5$ are K3s. See Ex.3.1–8 for a more proof and more examples.

**Example 3.** Let $X = X_{2,3} \subset \mathbb{P}^1 \times \mathbb{P}^2$ be a nonsingular divisor with the indicated bihomogeneity. Then it is easy to see that $K_{\mathbb{P}^1\times\mathbb{P}^2} = \mathcal{O}(-X)$ so that $K_X = 0$; and moreover, $H^1(\mathcal{O}_{\mathbb{P}^1\times\mathbb{P}^2}) = H^2(\mathcal{O}_{\mathbb{P}^1\times\mathbb{P}^2}(-X)) = 0$ so that $X$ is again a K3.

## 3.2. RR for curves and canonical curves

Let $C$ be a nonsingular projective curve. Recall that a divisor $D = \sum n_i P_i$ is a formal sum of points with integer coefficients. The degree of $D$ is $\deg D = \sum n_i$. The RR space or associated vector space $\mathcal{L}(D) = H^0(\mathcal{O}_C(D))$ has dimension $\ell(D) = \dim \mathcal{L}(D) = h^0(\mathcal{O}_X(D))$. The RR theorem states that $C$ has a genus $g \geq 0$ and a canonical divisor class $K_C$ such that for every divisor $D$ on $C$,

$$\ell(D) - \ell(K_C - D) = 1 - g + \deg D. \tag{1}$$

It follows easily from (1) that $\deg K_C = 2g-2$ and $\ell(K_C) = g$. Here $\ell(D) = h^0(\mathcal{O}_C(D))$ and $\ell(K_C - D) = h^0(\mathcal{O}_C(K_C - D)) = h^1(\mathcal{O}_C(D))$, so that (1) is a formula for $\chi(\mathcal{O}_C(D))$. In this form, the theorem is a fairly trivial consequence of the theorems of coherent cohomology (see Ex. B.23–24).

A divisor class $D$ on $C$ is *special* if $\ell(K_C - D) \neq 0$, or equivalently, $\ell(D) > 1 - g + \deg D$. Thus $K_C$ is the biggest special divisor class on $C$: a divisor $D$ is special if and only if (up to linear equivalence) $D \leq K_C$.

Hyperelliptic curves are counterexamples to many results in the theory of curves. A curve $C$ is *hyperelliptic* if it is birational to $y^2 = f_{2g+2}(x)$, or, equivalently, has a $g^1_2$, a divisor class $D$ with $\deg D = 2$, $\ell(D) = 2$. It is not hard to see that the canonical class $K_C$ is then composed of the $g^1_2$, that is, $K_C = (g-1)D$ and $\mathcal{L}(K_C) = S^{g-1}\mathcal{L}(D)$ (polynomials of degree $g-1$); in other words, if $x_1, x_2$ is a basis of $\mathcal{L}(D)$ then $x_1^{g-1}, x_1^{g-2}x_2, \dots, x_2^{g-1}$ is a basis of $K_C$. The canonical map $\varphi_{K_C}\colon C \to \mathbb{P}^{g-1}$ of a hyperelliptic curve is obtained by composing the double cover $C \to \mathbb{P}^1$ defined by the $g^1_2$ with the $(g-1)$st Veronese embedding $\mathbb{P}^1 \hookrightarrow \mathbb{P}^{g-1}$. A curve of genus 2 is automatically hyperelliptic, since $K_C$ is a $g^1_2$.

The following result is classical:

**Theorem** (Max Noether). *If $C$ is nonhyperelliptic of genus $\geq 3$ then $\varphi_{K_C}\colon C \to \mathbb{P}^{g-1}$ is an embedding such that $K_C$ is the hyperplane section divisor. Moreover, $C \subset \mathbb{P}^{g-1}$ is projectively normal, that is, $H^0(\mathbb{P}^{g-1}, \mathcal{O}(k)) \twoheadrightarrow H^0(C, \mathcal{O}_C(kK_C))$ for all $k \geq 1$.*

The image curve $C = \varphi_{K_C}(C) \subset \mathbb{P}^{g-1}$ is called a *canonical curve.*

**Idea of proof.** See [**4 authors**] for the proof. The idea is that $\ell(K_C) = g$; and if $P_1, P_2 \in C$ are any two points then $\ell(K_C - P_1 - P_2) = g-2$, for otherwise RR gives $\ell(P_1 + P_2) = 2$, so $P_1 + P_2$ is a $g^1_2$, and $C$ is hyperelliptic. Therefore functions in $\ell(K_C)$ distinguish $P_1$ and $P_2$.

The final sentence of the theorem is the first step in the Petri analysis [**4 authors**]. If $P_1, \dots, P_g \in C$ are $g$ "sufficiently general" points, I choose coordinates $x_1, \dots, x_g$ of $\mathbb{P}^{g-1}$ such that $P_i = (0, \dots, 1, \dots, 0)$. Then $x_1, x_2 \in \mathcal{L}(K_C - P_3 - \dots - P_g)$ span a free pencil. It then follows from the "free pencil trick" that

$$\left.\begin{aligned} x_1^2, x_1x_2, x_1x_3, \dots, x_1x_g \\ x_2^2, x_2x_3, \dots, x_2x_g \\ x_3^2, \dots, x_g^2 \end{aligned}\right\} \in \mathcal{L}(2K_C)$$

are $3g-3$ linearly independent elements. Passing up from $kK_C$ for $k \geq 3$ is similar (actually, a little easier).

The canonical embedding $C \subset \mathbb{P}^{g-1}$ of a nonhyperelliptic curve allows the following restatement of RR as a result in linear projective geometry. For $D$ an effective divisor on $C$, write $\langle D\rangle \subset \mathbb{P}^{g-1}$ for the smallest linear subspace of $\mathbb{P}^{g-1}$ containing $D$ as a subscheme of $C$. If $D = P_1 + \dots + P_d$ consists of distinct points, this is just their ordinary linear span $\langle P_1 + \dots + P_d\rangle$. Then $|K_C - D|$ is the linear system cut out residually to $D$ by hyperplanes of $\mathbb{P}^{g-1}$ through $\langle D\rangle$, and therefore $\ell(K_C - D) = g - 1 - \dim\langle D\rangle$.

**Theorem** (Geometric form of RR). *Let $C \subset \mathbb{P}^{g-1}$ be a canonical curve, and $D$ an effective divisor on $X$. Then $\ell(D) = \deg D - \dim\langle D\rangle$, in other words, $\ell(D) - 1 = \dim|D| =$ the number of linear dependence relations between the points of $D$.*

## 3.3. The historic definition of K3

A variety $X \subset \mathbb{P}^n$ is *linearly normal* or *embedded by a complete linear system* if $H^0(\mathbb{P}^n, \mathcal{O}(1)) \to H^0(X, \mathcal{O}(1))$ is surjective; in other words, $X$ is not the linear projection of a variety spanning a higher dimensional space $X \subset \mathbb{P}^{n+1}$.

**Theorem.** *Let $X \subset \mathbb{P}^n$ be a nonsingular surface. Then $X$ is a K3 embedded by a complete linear system if and only if one (every) nonsingular hyperplane section is a canonical curve.*

**Proof.** Let $C$ be a nonsingular hyperplane section. Then the adjunction formula and the assumption $K_X = 0$ give $K_C = \mathcal{O}_C(C)$. Also $H^0(X, \mathcal{O}_X(C)) \to H^0(C, K_C)$ is surjective from the cohomology exact sequence of

$$0 \to \mathcal{O}_X \to \mathcal{O}_X(C) \to \mathcal{O}_C(C) \to 0$$

and the assumption $H^1(\mathcal{O}_X) = 0$. Therefore $C$ is a canonical curve.

Conversely, suppose that $C$ is a canonical curve. Then $n = g$, and by Theorem 3.2, $H^0(\mathbb{P}^{g-1}, \mathcal{O}(k)) \to H^0(C, \mathcal{O}_C(kK_C))$ is surjective for all $k \geq 0$. But since $C \subset X \subset \mathbb{P}^g$, it follows from the commutative diagram of restriction maps

$$\begin{array}{ccc} H^0(\mathbb{P}^g, \mathcal{O}(k)) & \to & H^0(X, \mathcal{O}_X(k)) \\ \downarrow & & \downarrow \\ H^0(\mathbb{P}^{g-1}, \mathcal{O}(k)) & \to & H^0(C, \mathcal{O}_C(kK_C)) \end{array}$$

that $H^0(X, \mathcal{O}_X(kC)) \to H^0(C, \mathcal{O}_C(kK_C))$ is surjective for all $k \geq 0$.

Consider the cohomology exact sequence

$$\begin{aligned} \cdots \to H^0(\mathcal{O}_X(kC)) \to H^0(\mathcal{O}_C(kK_C)) \to \\ H^1(\mathcal{O}_X((k-1)C)) \to H^1(\mathcal{O}_X(kC)) \to H^1(\mathcal{O}_C(kK_C)) \to \\ H^2(\mathcal{O}_X((k-1)C)) \to \cdots \end{aligned}$$

Now since $H^0(X, \mathcal{O}_X(kC)) \to H^0(C, \mathcal{O}_C(kK_C))$ is surjective, it follows that

$$H^1(\mathcal{O}_X) \hookrightarrow H^1(\mathcal{O}_X(C)) \hookrightarrow \cdots \hookrightarrow H^1(\mathcal{O}_X(kC))$$

for all $k$. But $H^1(\mathcal{O}_X(kC)) = 0$ for $k \gg 0$ by Serre vanishing. Therefore $H^1(\mathcal{O}_X) = 0$.

Also, when $k = 1$, I get that $0 \neq H^1(\mathcal{O}_C(K_C)) \hookrightarrow H^2(\mathcal{O}_X)$. Hence $H^2(\mathcal{O}_X) \neq 0$, so that by duality $H^0(\mathcal{O}_X(K_X)) \neq 0$. Thus $K_X$ is linearly equivalent to an effective divisor $D \geq 0$. But the adjunction formula $K_C = (K_X + C)_{|C}$ together with the assumption $K_C = C_{|C}$ implies that $D \cap C = \emptyset$, so that $D = 0$. Thus $K_X = 0$. Q.E.D.

Linear systems on K3s are closely related to the geometry of canonical curves. If $C \subset X$ and $A = P_1 + \cdots + P_d$ is a $g^r_d$ on $C$, with $d$ fairly small, then by the geometric form of RR,

$$\mathbb{P}^g \supset \langle P_1, \dots, P_d \rangle = \Pi = \mathbb{P}^{d-1-r}.$$

This property of $\{P_1, \dots, P_d\}$ does not depend on $C$ through $\Pi$, so that all curves $C$ through $\{P_1, \dots, P_d\}$ also have a $g^r_d$. Thus, very roughly, one expects that all curves in a given linear system on $X$ have the same "very special" linear systems.

## 3.4. RR on a surface

If $X$ is a surface and $D$ a divisor then RR takes the form

$$\chi(\mathcal{O}_X(D)) = \chi(\mathcal{O}_X) + \tfrac{1}{2}D(D - K_X). \tag{1}$$

**"Proof".** The idea is to use induction on the components of $D$, as in the proof of RR for curves (see Ex. B.23). I can get from 0 to $D = \sum n_i C_i$ by successively adding or subtracting an irreducible curve $C$; and (1) holds for $D = 0$, so that (1) for all $D$ will follow by induction if I can prove the formula

$$\begin{aligned}\chi(\mathcal{O}_X(D)) - \chi(\mathcal{O}_X(D - C)) &= \big(\chi(\mathcal{O}_X) + \tfrac{1}{2}D(D - K_X)\big) \\ &\quad - \big(\chi(\mathcal{O}_X) + \tfrac{1}{2}(D - C)(D - C - K_X)\big) \\ &= -\tfrac{1}{2}(K_X + C)C + DC.\end{aligned} \tag{2}$$

Suppose for simplicity that $C$ is a nonsingular curve. Then the cohomology long exact sequence of $0 \to \mathcal{O}_X(D - C) \to \mathcal{O}_X(D) \to \mathcal{O}_C(D) \to 0$ gives

$$\chi(\mathcal{O}_X(D)) - \chi(\mathcal{O}_X(D - C)) = \chi(\mathcal{O}_C(D)). \tag{3}$$

Now the genus of $C$ is given by the adjunction formula: $2g(C) - 2 = (K_X + C)C$ and $\mathcal{O}_C(D)$ is a line bundle on $C$ of degree $DC$, so that $\chi(\mathcal{O}_C(D))$ by RR on $C$:

$$\chi(\mathcal{O}_C(D)) = 1 - g(C) + \deg(\mathcal{O}_C(D)) = -\tfrac{1}{2}(K_X + C)C + DC.$$

It is not too difficult to get around the assumption that $C$ is nonsingular. This "proves" (2) and hence (1).

Recall from A.9 that one way of defining intersections numbers $D_1D_2$ is based on knowing that $\chi(\mathcal{O}_X(D))$ is a quadratic function of $D$. If this is your definition, then the proof just given is, on the face of it, circular. Rather than as a logical proof, it is better to think of it as a compatibility between all the ingredients in the formula: intersection numbers, the genus of a curve, the degree of a divisor on a curve, etc. I suggest accepting RR as an axiom for the present until you have time to learn the general form of Hirzebruch RR or Grothendieck RR all at one go.

**3.5. Corollary** (the Hodge algebraic index theorem). *If $H$ is ample on $X$ then $HD = 0$ implies $D^2 \le 0$; moreover, if $D^2 = 0$ then $D \overset{\text{num}}{\sim} 0$.* Here *numerical equivalence* $D \overset{\text{num}}{\sim} 0$ means that $D\Gamma = 0$ for every curve $\Gamma \subset X$.

Another way of stating this, which is useful in calculations, is that if $D_1, D_2$ are divisors and $(\lambda D_1 + \mu D_2)^2 > 0$ for some $\lambda, \mu \in \mathbb{R}$ then the determinant

$$\det \begin{vmatrix} D_1^2 & D_1D_2 \\ D_1D_2 & D_2^2 \end{vmatrix} \le 0,$$

with equality if and only if some nonzero rational linear combination is numerically equivalent to zero, that is $\alpha D_1 + \beta D_2 \overset{\text{num}}{\sim} 0$.

I leave the proof as an exercise; or see Corollary D.2.2 below.

## 3.6. The arithmetic genus of a curve

I use RR to extend the definition of genus to singular curves, and even arbitrary effective divisors. If $C \subset X$ is an irreducible curve, it has a nonsingular model $\nu\colon C^\nu \to C$ given by normalisation, and the genus of $C^\nu$ is traditionally called the *geometric genus* of $C$; a better term might be *birational genus*.

**Definition.** Let $D$ be an effective divisor on a surface $X$. The *arithmetic genus* $p_aD$ of $D$ is defined by

$$2p_aD - 2 = (K_X + D)D. \tag{1}$$

If $D = C$ is a nonsingular curve then the definition is just the adjunction formula, so that $p_aC = g(C)$. (See the discussion in 4.10–11 for more information.)

Consider the quotient sheaf $\mathcal{O}_D = \mathcal{O}_X/\mathcal{O}_X(-D)$, where $\mathcal{O}_X(-D)$ is the ideal sheaf of regular functions vanishing along $D$. This is a coherent sheaf on $X$ by construction, and can be viewed as the structure sheaf of the subscheme $D \subset X$; see 3.10 below for a more detailed discussion. Be that as it may, the cohomology long exact sequence of $0 \to \mathcal{O}_X(-D) \to \mathcal{O}_X \to \mathcal{O}_D \to 0$ gives

$$\chi(\mathcal{O}_D) = \chi(\mathcal{O}_X) - \chi(\mathcal{O}_X(-D)),$$

and plugging in RR gives

$$\chi(\mathcal{O}_D) = 1 - p_aD. \tag{2}$$

**Proposition.** *Let $C$ be an irreducible curve, possibly singular, on a surface $X$, and $\nu\colon C^\nu \to C$ its normalisation. Then*

$$p_aC = g(C^\nu) + \sum_{P_i \in \operatorname{Sing} C} \delta(P_i),$$

*where $\delta(P_i) > 0$ are numerical invariants of the singularities of $C$.*

*In particular, if $p_aC = 0$ then $C \cong \mathbb{P}^1$.*

**Proof.** There is an exact sequence of sheaves on $C$

$$0 \to \mathcal{O}_C \to \nu_*(\mathcal{O}_{C^\nu}) \to \mathcal{N} \to 0, \tag{3}$$

where $\nu_*$ is the sheaf theoretic pushforward (see B.6), and $\mathcal{N}$ the cokernel. Now I claim that $\mathcal{N}$ consists of finite dimensional vector spaces $\mathcal{N}_{P_i}$ of dimension $\delta(P_i)$ supported at the singular points $P_1, \dots, P_k$ of $C$.

To see this, note that $\nu_*(\mathcal{O}_C^\nu)$ is a coherent sheaf of $\mathcal{O}_C$-modules by finiteness of normalisation; hence so is $\mathcal{N}$. Moreover, outside the singular points, $\nu$ is an isomorphism and $\mathcal{O}_C = \nu_*(\mathcal{O}_C^\nu)$, so that $\mathcal{N}$ is supported at $\{P_1, \dots, P_k\} \subset C$. If $n_i \subset \mathcal{O}_{X,P_i}$ is ideal annihilating $\mathcal{N}_{P_i}$ then $V(n_i) = \{P_i\}$. By the Nullstellensatz, $n_i$ contains a power $m_{P_i}^k$ of the maximal ideal at $P_i$, so that $\mathcal{N}_{P_i}$ is a finitely generated module over $\mathcal{O}_{X,P_i}/m_{P_i}^k$ and is finite dimensional.

Now the cohomology long exact sequence of (3) gives $\chi(\nu_*(\mathcal{O}_{C^\nu})) = \chi(\mathcal{O}_C) + \sum \delta(P_i)$, that is,

$$p_aC = g(C^\nu) + \sum \delta(P_i), \tag{4}$$

where $\delta(P_i) = \dim_k \mathcal{N}_{P_i}$. Here I have used that $H^i(C, \nu_*(\mathcal{O}_{C^\nu})) = H^i(C^\nu, \mathcal{O}_{C^\nu})$; for $i = 0$ this comes from the definition of pushforward $\nu_*$, for $i = 1$ because $\nu$ has only zero dimensional fibres. If $P_i$ is singular then $C^\nu \not\cong C$ at $P_i$, so that $\delta(P_i) > 0$.

Notice that (4) expresses $p_a C$ as a sum of terms which are $\geq 0$. If $p_a C = 0$ it follows that $\sum \delta(P_i) = 0$, so that there are no $P_i$, and $C = C^\nu$ has $g(C) = 0$, that is, $C \cong \mathbb{P}^1$. Q.E.D.

## 3.7. RR on a K3

If $X$ is a K3 then $\chi(\mathcal{O}_X) = 2$ and $K_X = 0$, so that RR takes the simpler form

$$\chi(\mathcal{O}_X(D)) = h^0(D) - h^1(D) + h^0(-D) = 2 + \tfrac{1}{2}D^2.$$

**Corollary.** *(i) $D^2 \geq -2$ implies $H^0(D) \neq 0$ or $H^0(-D) \neq 0$.*

*(ii) $D^2 \geq 0$ implies $D \overset{\text{lin}}{\sim} 0$, or $h^0(D) \geq 2$, or $h^0(-D) \geq 2$.*

*(iii) If $D$ is an effective divisor on $X$ with $h^0(D) = 1$ then $D'^2 \leq -2$ for every divisor $D'$ with $0 < D' \leq D$, and in particular $D$ is a sum of $-2$-curves with $D^2 \leq -2$.*

If $C \subset X$ is an irreducible curve on a K3 then $C^2$ is even, and the adjunction formula gives $2p_a C - 2 = C^2 \geq -2$, and $C^2 = -2$ implies that $C \cong \mathbb{P}^1$ by Proposition 3.6. A curve $C \cong \mathbb{P}^1$ with $C^2 = -2$ is called a $-2$-*curve*; these curves are important throughout the classification of surfaces. Note that any curve on a K3 has $C^2$ even, and $C^2 < 0$ only for $-2$-curves.

The bunch of $-2$-curves occurring as components of $D$ in (iii) is quite restricted, but it's not necessarily true that it has negative definite intersection matrix. See Ex. 3.12.

## 3.8. Easy properties of linear systems

If $D$ is a divisor on a variety $X$ with $\mathcal{L}(D) \neq 0$, the projective space $|D| = \mathbb{P}_*(\mathcal{L}(D))$ parametrises effective divisors $D' \geq 0$ linearly equivalent to $D$. This is called a *complete linear system.* (Example: a pencil $\lambda F + \mu G$ or net $\lambda_1 F_1 + \lambda_2 F_2 + \lambda_3 F_3$ of plane curves of degree $d$ is an example of a linear system in the plane; the complete linear system $|\mathcal{O}_{\mathbb{P}^2}(d)|$ is the system of all curves of degree $d$, of dimension $\binom{d+2}{2} - 1$). In general, $|D|$ may have a *base locus* $\operatorname{Bs}|D|$, that is, a subscheme $\Sigma$ contained in every divisor $D' \in |D|$. It is traditional to treat the codimension 1 part and the codimension $\geq 2$ part of $\operatorname{Bs}|D|$ separately. The *fixed part* of $D$ is the biggest divisor $F$ such that $F \leq D'$ for every $D' \in |D|$; or in other words, the gcd of all $D' \in |D|$.

On a general variety it may be hard to predict what $\operatorname{Bs}|D|$ looks like. However, on a K3 things are very nice. Recall that a divisor $D$ is *nef* (numerically eventually free) if $D\Gamma \geq 0$ for every curve $\Gamma \subset X$.

**Theorem.** *(a) If $X$ is a K3 and $D$ any effective divisor on $X$ then I can subtract of an effective sum of $-2$-curves $F = \sum n_i \Gamma_i$ to get $M = D - F$ such that $M$ is effective and nef (possibly zero), $M^2 \geq D^2$ and $H^0(X, \mathcal{O}_X(M)) = H^0(X, \mathcal{O}_X(D))$.*

*(b) If $D > 0$ is nef and $D^2 = 0$ then $D = aE$, where $|E|$ is a free pencil.*

*(c) If $D$ is nef and $D^2 > 0$ (that is, $D$ is* nef and big*) then $H^1(D) = 0$, so that $h^0(D) = 2 + \frac{1}{2}D^2$.*

*(d)* The first dichotomy: *if $D$ is nef and big then either $|D|$ has no fixed part, or $D = aE + \Gamma$, where $|E|$ is a free pencil and $\Gamma$ an irreducible $-2$-curve such that $E\Gamma = 1$. In this case $D$ is* monogonal.

**Confusion.** I clear up some possible sources of confusion. The last condition in (a) means that $F \subset \operatorname{Bs}|D|$, in particular $H^0(F) = 0$, but $F$ does not have to be the whole base locus: an example is the monogonal linear system in (d). The decomposition $D = M + F$ is not the same as the Zariski decomposition, which looks like $D = P + N$ as $\mathbb{Q}$-divisors, with $P$ nef and $N$ orthogonal to $P$; here $F \geq N$, and $F > N$ sometimes happens (see Ex. 3.13).

**Proof.** (a) If $D$ is nef, I'm home setting $F = 0$. An easy point that is crucial in all work on surfaces: if $D$ is effective and $D\Gamma < 0$ then $\Gamma$ is a component of $D$ and $\Gamma^2 < 0$. Because, I can certainly write $D = a\Gamma + D'$ with $a > 0$, $D' > 0$ and $\Gamma \not\subset D'$; then $D\Gamma = a\Gamma^2 + D'\Gamma$. Since $D'\Gamma \geq 0$ and $a \geq 0$, the only way this can be negative is if $\Gamma^2 < 0$ and $a > 0$.

If $D$ is not nef then there is a $\Gamma$ with $D\Gamma < 0$, necessarily in the fixed part of $|D|$. Set $D_1 = D - \Gamma$, so that $D_1^2 = D^2 - 2D\Gamma + \Gamma^2$. Now $\Gamma^2 = -2$, and $D\Gamma < 0$, and therefore $D_1^2 \geq D^2$. Also, obviously $H^0(D_1) = H^0(D)$. Continuing by induction proves (a).

(b) $D^2 = 0$ gives $h^0(D) \geq 2$, so that $D$ moves in a nontrivial linear system $|D|$. Write $|D| = |M| + F$ with $F$ the fixed part; then $M$ is mobile, hence also nef. As in (a), this is an easy point that occurs again and again: for any curve $\Gamma \subset X$, since $M$ is mobile, there is an effective divisor in $|M|$ not containing $\Gamma$, so $M\Gamma \geq 0$.

Thus $0 = D^2 = DM + DF$, so that $DM = DF = 0$; next, $M^2 + MF = 0$ so that $M^2 = MF = F^2 = 0$. Now if $F \neq 0$, I get $h^0(F) \geq 2$ from RR, contradicting $F$ fixed. Therefore $|D| = |M|$ has no fixed part. But $D^2 = 0$ implies that $D$ is free; because $D_1, D_2 \in |D|$ with no common components have $D_1D_2 = 0$, therefore $D_1 \cap D_2 = \emptyset$.

I claim that every element of $|D|$ is made up of components of fibres of a morphism $f\colon X \to \mathbb{P}^1$ with connected general fibre. Indeed, $D^2 = 0$ implies that the morphism $\varphi_D\colon X \to \mathbb{P}^{h^0(D)-1} = \mathbb{P}^*(H^0(D))$ defined by $|D|$ has image a curve $C$. If $X \to \widetilde{C} \to C$ is the Stein factorisation of $\varphi_D$ then a general fibre $E$ of $X \to \widetilde{C}$ is connected. In fact necessarily $\widetilde{C}$ is isomorphic to $\mathbb{P}^1$, since $h^1(\mathcal{O}_X) = 0$, and $\widetilde{C} = C$ since $\widetilde{C} \to \mathbb{P}^{h^0(D)-1}$ is defined by a complete linear system. Now a general fibre of $\varphi$ is an irreducible curve $E$ (by the first Bertini theorem, see [**Sh**], Chapter II, Theorem 6.1). Also $E$ is reduced, since if $E = nE'$ then also ${E'}^2 = 0$ and $E'$ also moves in a linear system.

Then it is easy to see that $D \overset{\text{lin}}{\sim} aE$, where $h^0(D) = a+1$: because some element in $|D|$ certainly intersects $E$, and therefore contains it, and then I can apply the argument to $D - E$, $D - 2E$, etc. This proves (b).

(c) The statement that $D$ nef and big implies $H^1(D) = 0$ is proved in 3.11–12 below, and I assume it for the moment.

(d) Write $|D| = |M| + F$, with $|M|$ mobile and $F$ fixed; note that $M$ is also effective and nef. There are two cases:

**Case $M^2 > 0$.** I prove that then $F = 0$. Indeed, $M$ is also nef and big, so that $H^1(M) = 0$, and so $h^0(M) = h^0(D)$ together with RR implies $M^2 = D^2$. Now

$$D^2 = D(M + F) \geq DM = (M + F)M \geq M^2,$$

and equality implies that $DF = MF = 0$; hence also $F^2 = 0$. If $F \neq 0$ then $H^0(F) \geq 2$, contradicting that $F$ is fixed.

**Case $M^2 = 0$.** Then $M = aE$ by (b), so that $h^0(M) = a + 1$, and RR gives $h^0(D) \geq 2 + \frac{1}{2}D^2$. Thus

$$\begin{aligned} a + 1 = h^0(M) = h^0(D) &\geq 2 + \tfrac{1}{2}D^2 \\ &= 2 + \tfrac{1}{2}(aE + F)^2 = 2 + aEF + \tfrac{1}{2}F^2 \\ &= 2 + \tfrac{1}{2}(aEF + F^2) + \tfrac{1}{2}aEF \geq 2 + \tfrac{1}{2}aEF. \end{aligned}$$

Now $EF \geq 2$ contradicts this inequality. Clearly $EF = 0$ is impossible, since then $DE = 0$, which would contradicts the index theorem. Therefore $EF = 1$ and $F^2 \leq -2$.

There is a unique irreducible component $\Gamma$ of $F$ with $E\Gamma = 1$, and $F = \Gamma + F'$, say, with $E\Gamma = 0$. Then $D = aE + \Gamma + F'$; but it's easy to check that $D' = aE + \Gamma$ is also nef, and $D'^2 = 2a - 2 = D^2$. It follows as above that $(F')^2 = 0$, so that $F' = 0$. Q.E.D.

## 3.9. Numerically connected divisors on surfaces

Let $X$ be a nonsingular projective surface. For an effective divisor $D$, the quotient sheaf $\mathcal{O}_D = \mathcal{O}_X/\mathcal{O}_X(-D)$ is a sheaf of rings, which can be viewed as the structure sheaf of a subscheme $D \subset X$. For a nonreduced divisor, $\mathcal{O}_D$ involves the nilpotents of scheme theory, but in a fairly mild way. The key to controlling the global sections $H^0(\mathcal{O}_D)$, and to the proof of Theorem 3.8, (c), is the notion of numerically connected divisor due to Franchetta and C. P. Ramanujam.

**Definition.** Let $k \in \mathbb{Z}$ (usually with $k \geq 0$). An effective divisor $D$ is *numerically $k$-connected* if $D_1D_2 \geq k$ for every effective decomposition $D = D_1 + D_2$ with $D_1, D_2 > 0$. (Of course, an irreducible curve is $k$-connected for all $k$, since the condition is vacuous.)

**Example.** If $E$ is irreducible with $E^2 = 0$, then $2E$ is numerically 0-connected but numerically 1-disconnected. If $|E|$ moves in a free pencil then it is the fibre of a morphism $f\colon X \to B$ to a base curve $B$, say $E = f^{-1}P$. Then the structure sheaf of the subscheme $aE \subset X$ is a module over $\mathcal{O}_{B,P}/m_P^a \cong k[\varepsilon]/\varepsilon^a$, and $H^0(\mathcal{O}_{aE}) = k[\varepsilon]/\varepsilon^a$; you can think of the nilpotent sections as (dual to) a normal vector field to $X$ along $E$ that goes out $a - 1$ infinitesimal steps. See Ex. 3.20.

**Example.** If $\sigma\colon Y \to X$ is the blowup of a nonsingular point $P \in X$ and $E$ the exceptional curve then $E$ is a $-1$-curve, that is $E \cong \mathbb{P}^1$ and $E^2 = -1$. Obviously $2E$ is not 0-connected. It is easy to see that $H^0(\mathcal{O}_{2E}) = \mathcal{O}_{X,P}/m_P^2 \cong k[x,y]/(x,y)^2$. See Ex. 3.21.

**3.10. Lemma** (*Dévissage* of $\mathcal{O}_D$). *If $D = D' + D''$ with $D', D'' > 0$ then there is an exact sequence*

$$0 \to \mathcal{O}_{D''}(-D') \to \mathcal{O}_D \to \mathcal{O}_{D'} \to 0;$$

*here $\mathcal{O}_{D''}(-D') = \mathcal{O}_X(-D') \otimes \mathcal{O}_{D''}$.*

**Proof.** The exact sequence defining $\mathcal{O}_{D''}$ is

$$0 \to \mathcal{O}_X(-D'') \to \mathcal{O}_X \to \mathcal{O}_{D''} \to 0;$$

Tensoring by $\mathcal{O}_X(-D')$ gives

$$0 \to \mathcal{O}_X(-D) \to \mathcal{O}_X(-D') \to \mathcal{O}_{D''}(-D') \to 0.$$

The ideal sheaves satisfy $\mathcal{O}_X(-D) \subset \mathcal{O}_X(-D')$, so that $\mathcal{O}_D \to \mathcal{O}_{D'}$ is surjective, and the kernel can be calculated by the snake lemma applied to the commutative diagram

$$\begin{array}{ccccccccc}
0 & \to & \mathcal{O}_X(-D) & \to & \mathcal{O}_X & \to & \mathcal{O}_D & \to & 0 \\
 & & \downarrow & & \| & & \downarrow & & \\
0 & \to & \mathcal{O}_X(-D') & \to & \mathcal{O}_X & \to & \mathcal{O}_{D'} & \to & 0 \\
 & & \downarrow & & & & & & \\
 & & \mathcal{O}_{D''}(-D') & & & & & & \\
 & & \downarrow & & & & & & \\
 & & 0. & & & & & &
\end{array}$$

Q.E.D.

**3.11. Lemma.** *Let $X$ be a nonsingular projective surface and $D$ an effective divisor.*

*(i) $D$ nef and $D^2 > 0$ implies that $D$ is numerically* 1*-connected.*

*(ii) $D$ numerically* 1*-connected implies that $H^0(\mathcal{O}_D) =$* constants.

*(iii) Suppose that $D$ is numerically* 1*-connected, and let $\mathcal{L}$ be a line bundle on $D$ with $\mathcal{L}\Gamma_i \le 0$ for every component $\Gamma_i$ of $D$. Then*

$$H^0(\mathcal{L}) \neq 0 \implies \mathcal{L} \cong \mathcal{O}_D.$$

**Proof.** (i) Let $D = D_1 + D_2$ with $D_1, D_2 \ge 0$. Then $D$ nef gives

$$D_1^2 + D_1D_2 = DD_1 \ge 0$$
$$D_1D_2 + D_2^2 = DD_2 \ge 0$$

Now if $D_1D_2 \le 0$ then $D_1^2D_2^2 \ge (D_1D_2)^2 \ge 0$, which contradicts the Hodge algebraic index theorem (Corollary 3.5) unless $D_1$ or $D_2 = 0$.

(ii) $H^0(\mathcal{O}_D)$ is a finite dimensional algebra over $k$, so that if $h^0(\mathcal{O}_D) \neq 1$ then it contains either an idempotent $e \neq 0, 1$, or a nilpotent element $\varepsilon \neq 0$ with $\varepsilon^2 = 0$. If $e \neq 0, 1$ is an idempotent then $e \cdot \mathcal{O}_D$ and $(e-1) \cdot \mathcal{O}_D$ are ideals defining disjoint subschemes $D_1, D_2 \subset D$, so that $D$ is disconnected.

If $0 \neq \varepsilon \in H^0(\mathcal{O}_D)$ is a nilpotent element and $D_1 \subset D$ is the greatest divisor on which it vanishes then the ideal sheaf $\mathcal{I}_{D_1} \subset \mathcal{O}_D$ is generated by $\varepsilon$ outside a set of codimension $\geq 1$. By Lemma 3.10, this kernel is $\mathcal{O}_{D_2}(-D_1)$; thus $\mathcal{O}_{D_2}(-D_1)$ is generated by $\varepsilon$ at the generic point of each component, and it follows that $D_1\Gamma \leq 0$ for every component $\Gamma \subset D_2$. In particular $D_1D_2 \leq 0$, contradicting 1-connected.

I give a second proof which is more elementary: suppose that $D_1 \subset D$ is a divisor for which $h^0(\mathcal{O}_{D_1}) = 1$, for example an irreducible component of $D$. By the 1-connectedness assumption, if $D_1 < D$ then $D_1(D - D_1) \geq 1$, so there is a component $\Gamma$ of $D - D_1$ such that $D_1\Gamma \geq 1$. Then I claim that also $h^0(\mathcal{O}_{D_1+\Gamma}) = 1$. Indeed, the sequence

$$H^0(\mathcal{O}_\Gamma(-D_1)) \to H^0(\mathcal{O}_{D_1+\Gamma}) \to H^0(\mathcal{O}_{D_1})$$

is exact, the first term is 0, and the third term 1-dimensional. Then $h^0(\mathcal{O}_D) = 1$ follows by induction.

(iii) is an exercise (see Ex. 3.27). Q.E.D.

## 3.12. Ramanujam vanishing

I start by completing the proof of Theorem 3.8, (c) on K3s. $H^0(D) \neq 0$, so that I can assume that $D$ is effective. By assumption it is nef and big, so that $H^0(\mathcal{O}_D) =$ constants by Lemma 3.11. Now in the exact sequence

$$\begin{aligned} 0 = H^0(\mathcal{O}_X(-D)) \to H^0(\mathcal{O}_X) \to H^0(\mathcal{O}_D) \to \\ \to H^1(\mathcal{O}_X(-D)) \to H^1(\mathcal{O}_X), \end{aligned}$$

$H^1(\mathcal{O}_X) = 0$ by assumption ($X$ is a K3), and $H^0(\mathcal{O}_X)$ also equals the constant functions, and hence maps onto $H^0(\mathcal{O}_D)$. Therefore $H^1(-D) = 0$, which is the dual of $H^1(D)$. Q.E.D.

On a surface $X$ with $H^1(\mathcal{O}_X) = 0$, the same argument proves the vanishing $H^1(\mathcal{O}_X(-D)) = 0$ if $D$ is effective and 1-connected, in particular if $D$ is effective, nef and big. Dually, $H^1(\mathcal{O}_X(K_X + D)) = 0$. This is a weak form of Kodaira vanishing.

**Discussion.** Kodaira vanishing states that on a nonsingular $n$-fold over a field of characteristic zero,

$$H \text{ ample} \quad \Longrightarrow \quad H^i(\mathcal{O}_X(K_X + H)) = 0 \quad \text{for all } i > 0,$$

or dually, $H^j(\mathcal{O}_X(-H)) = 0$ for all $j < n$. Although the statement is purely in terms of coherent cohomology, the result is deeper than algebraic geometry in characteristic zero. There is no proof purely within coherent cohomology, and the result is false in characteristic $p$, already for surfaces. The two known proofs in all dimensions are Kodaira's (representing a cohomology class by a harmonic form, and integrating to give something which must be zero by Stokes' theorem, and at the same time strictly positive if the cohomology class is nonzero), and Deligne and Illusie's (by reducing modulo a sufficiently large prime $p$, and analysing the splitting of the characteristic $p$ de Rham complex by Frobenius and Cartier operators). Generalisations of Kodaira's proof in a number of directions have been given, and this is a key ingredient in Mori theory and higher dimensional classification; compare the discussion in D.3 below.

A modern view of what vanishing *really is* has been given by Kollár: vanishing comes about when a coherent cohomology group admits a topological interpretation. Thus the argument given above works because the coherent cohomology group $H^0(\mathcal{O}_D)$ has an interpretation in terms of connectedness of $D$.

There are two different proofs of vanishing for algebraic surfaces in characteristic zero (without the assumption $H^1(\mathcal{O}_X) = 0$): the Bogomolov–Mumford proof based on rank 2 vector bundles, their Chern numbers and the inequality $c_1^2 < 4c_2$, and instability,[4] and C. P. Ramanujam's, which I now sketch.

**3.13. Theorem.** *If $D$ is effective, nef and big then $H^1(\mathcal{O}_X(-D)) = 0$.*

**Sketch proof.** We have already seen that $D$ nef and big implies $D$ 1-connected, and $H^0(\mathcal{O}_D) =$ constants. Then cohomology gives an exact sequence

$$0 \to H^1(\mathcal{O}_X(-D)) \to H^1(\mathcal{O}_X) \xrightarrow{\alpha(D)} H^1(\mathcal{O}_D).$$

Now C. P. Ramanujam proves that $D$ effective and big also implies the injectivity $\alpha(D)\colon H^1(\mathcal{O}_X) \hookrightarrow H^1(\mathcal{O}_D)$.

**Key lemma.** *For any effective divisor $D$, write $\alpha(D)\colon H^1(\mathcal{O}_X) \to H^1(\mathcal{O}_D)$. Then, in characteristic zero, $\ker\alpha(D)$ depends only on the reduced divisor $D_{\text{red}}$.*

**Corollary.** *If $D$ is effective and big then $\alpha(D)$ is injective.*

**Proof of Corollary, Step 1.** The condition *$D$ big* means that $|nD|$ defines a birational map $\varphi_{|nD|}\colon X \dashrightarrow X' \subset \mathbb{P}^N$ for some $n > 0$; this is also popularly expressed in terms of the *Iitaka dimension* of $D$ as $\kappa(X, D) = 2$. Now if $|nD|$ is a big linear system, I can write $nD = D' + F$ (equality of divisors) with $D', F \geq 0$, and $|D'|$ a big linear system without fixed part. The map $\alpha(D')$ factors as

$$H^1(\mathcal{O}_X) \to H^1(\mathcal{O}_{nD}) \to H^1(\mathcal{O}_{D'}),$$

so that obviously $\alpha(D')$ injective implies $\alpha(nD)$ injective, and by the Key Lemma, $\alpha(nD)$ injective implies $\alpha(D)$ injective. Thus it is enough to prove the corollary for $D'$. Hence from now on, replacing $D$ by $D'$, I assume $D$ is nef and big.

**Step 2.** It follows from big that for any ample $H$, there exists an $m \gg 0$ such that $h^0(mD - H) \neq 0$, so that $mD \overset{\text{lin}}{\sim} H + E$ with $E$ effective; indeed, just take a form in $\mathbb{P}^N$ vanishing on $\varphi(H)$. Taking a still larger multiple gives $lmD \overset{\text{lin}}{\sim} lH + lE = H' + lE$, where $H' \overset{\text{lin}}{\sim} lH$ is a very ample divisor. Then $H' + (lE)_{\text{red}} \overset{\text{lin}}{\sim} lH + E_{\text{red}}$, which is as ample as I like, so that, by Serre vanishing, $H^1(\mathcal{O}_X(-H' - E_{\text{red}})) = 0$.

**Step 3.** Therefore, in particular, $\alpha(H' + E_{\text{red}})$ is injective, so by using the Key Lemma, also $\ker\alpha(H' + lE) = 0$.

**Step 4.** Now $H' + lE \overset{\text{lin}}{\sim} lmD$ is nef and big, so that by Lemma 3.11, $h^0(\mathcal{O}_{H'+lE}) = 1$, and together with $\ker\alpha(H' + lE) = 0$ this implies that $H^1(\mathcal{O}_X(-H' - lE)) = 0$. Now $lmD \overset{\text{lin}}{\sim} H' + lE$ gives $\mathcal{O}_X(-lmD) \cong \mathcal{O}_X(-H' - lE)$, so also $H^1(\mathcal{O}_X(-lmD)) = 0$.

[4] This proof uses essentially only the leading term in RR, and an argument on inseparable sections of scrolls. I hope to make this into a later chapter in these notes.

**Step 5.** Therefore $\ker \alpha(mD) = 0$. Hence by using the lemma again, since $D$ and $mD$ have the same reduced divisor, $\ker \alpha(D) = 0$.

## 3.14. Proof of the Key Lemma 3.13

This is where we could hope to see Kollár's principle in action: we're talking about the map $H^1(\mathcal{O}_X) \to H^1(\mathcal{O}_D)$. If the group $H^1(\mathcal{O}_D)$ appearing here was topological in nature (say, homotopy invariant), then the group for $D$ and for $D_{\text{red}}$ should be the same, since they are the same topological space. Unfortunately, this is not the case.

C. P. Ramanujam's idea is to view $H^1(\mathcal{O}_X)$ and $H^1(\mathcal{O}_D)$ as the Zariski tangent spaces to the Picard schemes $\operatorname{Pic}^0 X$ and $\operatorname{Pic}^0 D$. Assume we are over $\mathbb{C}$ and give $X$ the complex topology. Then the Picard group $\operatorname{Pic} X = H^1(\mathcal{O}_X^*)$ fits into the exact sequence of sheaf cohomology

$$\begin{aligned} 0 \to H^1(X, \mathbb{Z}) \to H^1(X, \mathcal{O}_X) \to H^1(\mathcal{O}_X^*) \to \\ \to H^2(X, \mathbb{Z}) \end{aligned}$$

Thus $\operatorname{Pic}^0 X = \ker\{H^1(\mathcal{O}_X^*) \to H^2(X, \mathbb{Z})\}$ is the moduli space of topologically trivial line bundles. By Hodge theory, $H^1(X, \mathbb{Z}) \subset H^1(X, \mathcal{O}_X)$ is a full lattice (that is, $\mathbb{Z}^{2q} \subset \mathbb{C}^q$ as a discrete cocompact subgroup), so that the quotient $\operatorname{Pic}^0 X = H^1(X, \mathcal{O}_X)/H^1(X, \mathbb{Z})$ is naturally an Abelian variety.

Now $\operatorname{Pic}^0 D =$ the identity component of $H^1(\mathcal{O}_D^*)$ fits into a similar exact sequence, and can be made into a Lie group $\operatorname{Pic}^0 D = H^1(D, \mathcal{O}_D)/H^1(D, \mathbb{Z})$ in the same way. However, it can shown that $H^1(D, \mathbb{Z}) \subset H^1(D, \mathcal{O}_X)$ is a full lattice if and only if $D$ is a nonsingular curve. Indeed, if $\nu \colon C \to D_{\text{red}} \subset D$ is the resolution of singularities of the reduced subscheme of $D$, then there is a natural homomorphism $\operatorname{Pic}^0 D \to \operatorname{Pic}^0 C$ (the product of the Jacobians of the components of $C$), and the kernel comes from the nonreduced structure of $D$ or the singularities of $D_{\text{red}}$. The nonreduced structure contributes only additive subgroups (a direct sum of copies of $\mathbb{C}^+$), and the singularities of $D_{\text{red}}$ contribute either additive or multiplicative groups (a direct sum of copies of $\mathbb{C}^+$ and $\mathbb{C}^*$).

Now a morphisms from an Abelian variety $A$ to an additive or multiplicative group $G$ is necessarily zero, for example because $A$ is compact and $G$ is affine. Therefore, there is no nonconstant map from $\operatorname{Pic}^0 X$ to $\ker\{\operatorname{Pic}^0 D \to \operatorname{Pic}^0 C\}$, and hence

$$\ker\{\operatorname{Pic}^0 X \to \operatorname{Pic}^0 D\} = \ker\{\operatorname{Pic}^0 X \to \operatorname{Pic}^0 D_{\text{red}}\} = \ker\{\operatorname{Pic}^0 X \to \operatorname{Pic}^0 C\}.$$

Since $\alpha(D)$ and $\alpha(D_{\text{red}})$ are the derivatives at 0 of these maps, this concludes the sketch proof of the Key Lemma. Q.E.D.

## 3.15. Final remarks on K3s

### The 2nd dichotomy

Saint-Donat's theory of linear systems on K3s continues with a criterion for $\varphi_{|D|}$ to define an embedding modulo $-2$-curves: if $D$ is nef and big, and not monogonal, then $|D|$ is free, and $\varphi_{|D|}$ is either 2-to-1 to its image (then $D$ is *hyperelliptic*), or $\varphi_{|D|}$ is birational to a normal surface $\overline{X}$, and is an embedding, except that it contracts $-2$-curves $\Gamma$ with $D\Gamma = 0$ to Du Val singularities. In the first case $X \to D$ is a double cover of $\mathbb{P}^2$, or the Veronese surface, or a rational normal scroll. The hyperelliptic case is characterised by the fact that $X$ has $\varphi_{|D|}$ an elliptic pencil $|E|$ with $DE = 2$, or (in the Veronese case) a free linear system $|B|$ with $B^2 = 2$ and $D = 2B$.

### The 3rd dichotomy

If $D$ is nef and big, and $\varphi_D \colon X \to \overline{X} \subset \mathbb{P}^g$ is birational, then either $\overline{X}$ is contained in a 3-fold $W$ contained in the intersection of all quadrics through $\overline{X}$, or $\overline{X}$ is an ideal theoretic intersection of quadrics. In the first case $X \to D$ is contained in the 3-fold $W$, which is either $\mathbb{P}^3$ or $Q \subset \mathbb{P}^4$ or the cone over the Veronese surface, or a rational normal scroll. In the final case $D$ is *trigonal*. The trigonal case is characterised by the fact that $X$ has $\varphi_{|D|}$ an elliptic pencil $|E|$ with $DE = 3$, or a free linear system $|B|$ with $B^2 = 2$ and $DB = 5$.

These dichotomies certainly continue for a while. In general, for the reason mentioned briefly at the end of 3.3, if one curve $C \in |D|$ on $X$ has a very special linear system, then one expects all the other curves $C \in |D|$ to have a closely related linear system. For example, if $C$ has a $g^1_d$ and $d$ is small compared to $D^2$ then $C$ is cut out on $X$ by a linear system. Lazarsfeld has shown that if $C$ has a $g^r_d$ which is special in the sense of Brill–Noether theory then $\operatorname{Pic} X$ must be strictly bigger than $\mathbb{Z} \cdot D$. Green and Lazarsfeld have proved that all $C \in |D|$ have the same Clifford index. It seems to be known that with a single well-known exception involving $g^1_4$s and $g^2_6$s, all $C \in |D|$ have the same gonality. I believe there are open research problems in this area.

### Moduli and periods

If you fix a primitive sublattice $L$ of rank $\rho \le 20$ and signature $(1, \rho - 1)$ of the K3 lattice $H^2(X, \mathbb{Z}) \cong 2\begin{pmatrix} 0 & 1 \\ 1 & 0 \end{pmatrix} \oplus 3E_8$ then K3s having $\operatorname{Pic} X \supset L$ form a nonempty moduli space of dimension $20 - \rho$. In particular, there are countably many moduli spaces of algebraic K3s, each depending on 19 moduli, with (generically) $\operatorname{Pic} X = \mathbb{Z} \cdot D$ with $D^2 = 2g - 2$.

Analytic K3s depend on a single irreducible 20-dimensional moduli space. The Teichmüller space (parametrising K3s plus a basis of $H^2(X, \mathbb{Z})$) is just one half of a 20-dimensional quadric. To study K3s from this point of view, the essential result is the Torelli theorem, which says that a polarised K3 surface is uniquely determined by its Hodge structure. Every complex K3 has a Kähler metric.

A more leisurely discussion of these topics can be found in [**3 authors**].

## Exercises to Chapter 3

**1.** Prove that $H^1(\mathbb{P}^1, \mathcal{O}(k)) = 0$ for all $k \ge -1$. [Hint: It's true if $k \gg 0$ by Serre vanishing; use the cohomology of $0 \to \mathcal{O}(k-1) \to \mathcal{O}(k) \to k_P \to 0$ for $P \in \mathbb{P}^1$, and induction.]

**2.** Prove that

$$H^0(\mathbb{P}^n, \mathcal{O}(k)) = 0 \qquad \text{for all } k < 0;$$
$$H^n(\mathbb{P}^n, \mathcal{O}(k)) = 0 \qquad \text{for all } k > -n-1;$$
$$H^i(\mathbb{P}^n, \mathcal{O}(k)) = 0 \qquad \text{for all } 0 < i < n \text{ and all } k \in \mathbb{Z}.$$

[Hint: As in the preceding exercise, argue by induction on $n$ using the cohomology of $0 \to \mathcal{O}_{\mathbb{P}^n}(k-1) \to \mathcal{O}_{\mathbb{P}^n}(k) \to \mathcal{O}_{\mathbb{P}^{n-1}}(k) \to 0$ for $\mathbb{P}^{n-1} \subset \mathbb{P}^n$.]

**3.** If $X = X_d \subset \mathbb{P}^n$ is a nonsingular hypersurface, prove that

$$H^0(X, \mathcal{O}(k)) = 0 \qquad \text{for all } k < 0;$$
$$H^{n-1}(X, \mathcal{O}(k)) = 0 \qquad \text{for all } k > d-n-1;$$
$$H^i(X, \mathcal{O}(k)) = 0 \qquad \text{for all } 0 < i < n-1 \text{ and all } k \in \mathbb{Z}.$$

[Hint: Use the previous exercises, and argue on $0 \to \mathcal{O}_{\mathbb{P}^n}(k-d) \to \mathcal{O}_{\mathbb{P}^n}(k) \to \mathcal{O}_X(k) \to 0$.]

**4.** Let $X$ be nonsingular and $Y \subset X$ a nonsingular subvariety of codimension $c$. Define the normal sheaf $N_{X|Y}$ as the dual of $\mathcal{I}_Y/\mathcal{I}_Y^2$ and prove that it is locally free of rank $c$. Prove that it fits into an exact sequence with $T_Y$ and the restriction $T_{X|Y}$. State and prove the adjunction formula giving $K_Y$ in terms of $K_X$.

**5.** Prove that every nonsingular complete intersection $X_4 \subset \mathbb{P}^3$, $X_{2,3} \subset \mathbb{P}^4$ or $X_{2,2,2} \subset \mathbb{P}^5$ is a K3.

**6.** Show that the list of K3 complete intersections in $\mathbb{P}^n$ (of hypersurfaces of degree $\geq 2$) of the preceding question is complete.

**7.** Find all values of $a_1, a_2, a_3$ and $e$ for which the general hypersurface $X_{3,e} \subset \mathbb{F}(a_1, a_2, a_3)$ is a nonsingular K3. [Hint: The canonical class of $\mathbb{F}$ is given in Proposition A.9. The criterion for $X \in |eL + 3M|$ to be a nonsingular cubic curve generically over $\mathbb{P}^1$ is a Newton polygon argument (referred to in A.6). Note that the result Worked Example 2.10 does not give nonsingularity in general: when $|eL + 3M|$ has a base locus, you have to check for isolated singularities along the base locus.]

**8.** Find all values of $a_1, \ldots, a_4$ and $d_1, d_2$ for which the general codimension 2 complete intersection $Q_{d_1,2} \cap Q'_{d_2,2} \subset \mathbb{F}(a_1, \ldots, a_4)$ is a nonsingular K3.

**9.** Prove the Hodge algebraic index theorem, Corollary 3.5. [Hint: If $D^2 > 0$, deduce from RR that either $nD$ or $-nD$ is equivalent to a effective divisor for $n \gg 0$ (if $K_X - nD$ grows quadratically, you have to invent a restriction argument to show that also $-nD$ grows quadratically). Prove then that $HD \neq 0$. Now for the case of equality, if $HD = 0$, and $D^2 = 0$ but $D\Gamma \neq 0$ for some curve $\Gamma$, find a linear combination $D' = \alpha D + \beta A + \gamma H$ with $(D')^2 > 0$ but $DH = 0$.]

**10.** Prove the det $\begin{vmatrix} D_1^2 & D_1D_2 \\ D_1D_2 & D_2^2 \end{vmatrix} \leq 0$ form of the Hodge algebraic index theorem given in 3.5.

**11.** Show that $C^2$ is even and $\geq -2$ for every irreducible curve $C$ on a K3, in particular $X$ does not contain any $-1$-curves. Prove Corollary 3.7 (the easy corollaries of RR on a K3).

**12.** Let $C \subset \Pi = \mathbb{P}^2 \subset \mathbb{P}^3$ be a nonsingular plane conic, and $X = X_4 \subset \mathbb{P}^3$ a quartic surface tangent to $\Pi$ along $C$, but otherwise general. Prove that $X$ has six nodes at points of $C$. [Hint: Take $C : (x_0 = q(x_1, x_2, x_3) = 0) \subset \mathbb{P}^3$ and $X : q^2 + x_0 f_3$.]

Let $f \colon Y \to X$ be the blowup of the 6 nodes, and $H = f^*\mathcal{O}_X(1)$, so that $\varphi_{|H|} \colon Y \to \mathbb{P}^3$ defines $f$. Show that $|H|$ contains a divisor of the form $2C + \sum_{i=1}^6 E_i$, where $C$ and the $E_i$ are $-2$-curves, and $h^0(X, \mathcal{O}_X(C + \sum_{i=1}^6 E_i)) = 1$. Show that the intersection matrix of $C$ and the $E_i$ is not negative definite.

**13.** Let $D$ be an effective divisor on a surface $X$. The Zariski decomposition of $D$ is an expression $D = P + N$ (pronounced "positive plus negative") with $P$ a nef $\mathbb{Q}$-divisor, $N = \sum q_i\Gamma_i$ with $q_i \in \mathbb{Q}$, $q_i > 0$, the intersection matrix $\Gamma_i\Gamma_j$ is negative definite, and $P\Gamma_i = 0$ for all $\Gamma_i$. Show that this exists and is unique. [Hint: You can get uniqueness by a simple argument in quadratic forms. For existence: if $D$ is nef there is nothing to prove. Otherwise $D\Gamma_1 < 0$, so $\Gamma_1^2 < 0$, and subtracting off a suitable multiple of $\Gamma_1$ I get $D_1 = D - a_1\Gamma_1$ with $D_1\Gamma_1 = 0$. Proceed in the same way, keeping $D_i$ orthogonal to $\Gamma_1, \dots, \Gamma_i$.]

Let $D$ be the divisor in Ex. 3.12. Find its Zariski decomposition $D = N + P$ and its fixed part $F$, and compare with the remark in Theorem 3.8.

**14.** Let $C$ and $D$ be irreducible curves. When is $nC$ numerically $n$-connected? When is $nC + D$ numerically $n$-connected?

**15.** Suppose that $D = \sum_{i=1}^k n_i\Gamma_i$, with all $\Gamma_i^2 = 0$. When does $D$ fail to be numerically 1-connected? For any $n$, give an example of such a divisor which is numerically $n$-connected but not numerically $n+1$-connected.

**16.** Suppose that $C_1, C_2, C_3$ form a triangle $D = C_1 + C_2 + C_3$, and each $C_i^2 < 0$. List all cases of $nD$ numerically 2-connected.

**17.** Work out the exact sequence of Lemma 3.10 for a reduced divisor $D$. Understand why $h^0(\mathcal{O}_D) = 1$ if $D$ is connected.

**18.** Ample certainly implies nef and big for a divisor on a surface. Deduce from Lemma 3.11 that a hyperplane section of an irreducible variety is a connected set.

**19.** Let $C : (x = 0)$ be the $y$-axis in $\mathbb{C}^2$ and let $D = 2C$. Write out explicitly the dévissage of the sheaf $\mathcal{O}_D$. The same question for $D = (x^2y^2 = 0)$.

**20.** Let $E = f^{-1}P$ be a nonsingular reduced fibre of a morphism $f \colon X \to B$ of a surface to a base curve. Prove that $H^0(\mathcal{O}_{aE}) = k[t]/(t^a)$. [Hint: Start with $a = 2$. You can do this in two different ways. First by carrying out the dévissage of Lemma 3.10 systematically, using the fact that the conormal bundle is trivial, $\mathcal{O}_E(-E) \cong \mathcal{O}_E$. Secondly, you can embed $k[t]/(t^a) = \mathcal{O}_{B,P}/m_P^a$ into $H^0(\mathcal{O}_{aE})$ using $f^*$ and the fact that $t^a$ defines exactly the subscheme $aE$; then show that $h^0(\mathcal{O}_{aE}) \le a$ by restricting to a section of $f$.]

**21.** If $E = \sigma^{-1}P$ is a $-1$-curve, the exceptional curve of a blowup $\sigma \colon X \to Y$, prove that $H^0(\mathcal{O}_{aE}) = k[x, y]/(x, y)^a$.

**22.** If $E$ is a $-2$-curve and $Q \in Y$ is the ordinary double point $xz = y^2$ then $H^0(\mathcal{O}_{aE}) = \big(k[x, y, z]/(xz - y^2)\big)/(x, y, z)^a$.

**23.** Let $D = 3L \subset \mathbb{P}^2$ be a plane cubic that happens to be a triple line. Calculate $\chi(\mathcal{O}_D)$ using dévissage, and compare with the formula for the arithmetic genus of $D$ as a plane cubic.

**24 (harder).** Use the C. P. Ramanujam method to prove the following result (due essentially to P. Francia). Let $X$ be a surface $D$ an effective divisor on $X$ and $P \in \operatorname{Sing} D$; write $\sigma\colon Y \to X$ for the blowup of $P$, $E$ for the exceptional curve, $D' = \sigma^*D - E$ and $D'' = \sigma^*D - 2E$. Then $P$ is a base point of $|D+K|$ if and only if $\sigma^*\colon H^0(\mathcal{O}_D) \to H^0(\mathcal{O}_{D''})$ is not surjective. [Hint: By cohomology, $P \in \operatorname{Bs}|D+K|$ if and only if $H^1(m_P \cdot \mathcal{O}_X(D+K)) \to H^1(\mathcal{O}_X(D+K))$ is not injective. Write this on $Y$ in terms of $H^1(\mathcal{O}_Y(\sigma^*(D+K)-E))$, and use Serre duality to express in terms of $H^1(\mathcal{O}_Y) \to H^1(\mathcal{O}'_D)$ and $H^1(\mathcal{O}''_D)$. Then use the C. P. Ramunujam trick together with the fact that $\ker\colon \operatorname{Pic}^0 D' \to \operatorname{Pic}^0 D''$ does not contain any Abelian variety.]

**25.** The adjunction formula via Poincaré residue. $X$ is a nonsingular variety and $Y \subset X$ a nonsingular codimension 1 subvariety. Suppose that $z_0, \dots, z_n$ are local coordinates on $X$, and $Y$ is locally defined $f(z_0, \dots, z_n) = 0$, with $\partial f/\partial z_0 \neq 0$. Prove that there is a sheaf homomorphism

$$\operatorname{Res}\colon \Omega_X^{n+1}(Y) \to \Omega_Y^n \quad \text{taking the local basis} \quad \frac{\mathrm{d}z_0 \wedge \cdots \wedge \mathrm{d}z_n}{f} \mapsto \frac{\mathrm{d}z_1 \wedge \cdots \wedge \mathrm{d}z_n}{\partial f/\partial z_0}.$$

(The exercise is to show that the homomorphism is intrinsic, that is, independent of the choice of local coordinates and $f$). This is called *Poincaré residue*. It fits into an exact sequence of sheaves $0 \to \Omega_X^{n+1} \to \Omega_X^{n+1}(Y) \to \Omega_Y^n \to 0$.

**26.** The adjunction formula via Serre duality. The full form of Serre duality was mentioned in Chapter B, Rule vii: for $X$ projective and Cohen–Macaulay (say nonsingular), there exists a duality pairing

$$H^i(X, \mathcal{F}) \times \operatorname{Ext}^{n-i}_{\mathcal{O}_X}(\mathcal{F}, \omega_X) \to H^n(X, \omega_X) = k$$

for any coherent sheaf $\mathcal{F}$. Let $Y \subset X$ be a codimension 1 subvariety which is a Cartier divisor. Now coherent sheaves on $Y$ are particular cases of coherent sheaves on $X$ via $\mathcal{O}_X \twoheadrightarrow \mathcal{O}_Y$, and it can be checked that $\operatorname{Ext}^{n-i}_{\mathcal{O}_X}(\mathcal{F}, \omega_X) = \operatorname{Ext}^{n-i-1}_{\mathcal{O}_Y}(\mathcal{F}, \omega_Y)$, where $\omega_Y = \mathcal{E}xt^1_{\mathcal{O}_X}(\mathcal{O}_Y, \omega_X)$.

Apply the cohomological $\delta$-functor $\operatorname{Ext}^*_{\mathcal{O}_X}$ to the exact sequence

$$0 \to \mathcal{O}_X(-Y) \to \mathcal{O}_X \to \mathcal{O}_Y \to 0$$

and prove that $\omega_Y$ fits in an exact sequence $0 \to \omega_X \to \omega_X(Y) \to \omega_Y \to 0$. [Hint: $\omega_X(Y) = \mathcal{H}om_{\mathcal{O}_X}(\mathcal{O}_X(-Y), \omega_X)$ and $\mathcal{E}xt^1_{\mathcal{O}_X}(\mathcal{O}_X(-Y), \omega_X) = 0$.]

**27.** Prove Lemma 3.11, (iii). [Hint: For $s \in H^0(\mathcal{L})$, decompose $D = A + B$ so that $A$ is the biggest divisor on which $s$ vanishes. Then consider the exact sequence

$$0 \to \mathcal{L} \otimes \mathcal{O}_B(-A) \to \mathcal{L} \to \mathcal{L}_{|B} \to 0,$$

and make a numerical consequence from the fact that $\mathcal{L}(-A)$ is generically generated by the section $s$.]

# CHAPTER 4. Singularities and surfaces

This chapter turns to surface singularities. The main object of study is a normal (isolated) surface singularity $P \in X$ and a resolution $f\colon Y \to X$, with $f^{-1}P = \bigcup \Gamma_i$ a bunch of curves on a nonsingular surface $Y$. I introduce examples of three classes of singularities, the Du Val singularities, more general rational singularities and elliptic Gorenstein singularities; the latter two classes can be viewed as generalisations of the Du Val singularities. These singularities play an important role all over the classification of surfaces and 3-folds. At the same time, I discuss invariants of singularities, mainly concerned with quantifying the difference between $X$ and $Y$ in the resolution $f\colon Y \to X$. The results here are mostly taken from Artin [**A1**–**2**], Laufer [**L1**–**3**] and my unpublished manuscript [**R**], although much of it was known in some form to Du Val in the 1930s. Any area of math can be given a superior treatment using several decades' hindsight (and not least, Mori's notes [**Mori**]), and in this spirit, this chapter contains new proofs of the foundational results on rational and elliptic Gorenstein surface singularities.

Even for the reader only interested in nonsingular varieties or complex manifolds, there are several reasons for studying singularities. A philosophical point is that the study of varieties and the study of singularities each contains the other many times over. The techniques of this chapter mostly involve a bunch of curves $\{\Gamma_i\}$ on a surface $Y$, and don't really need the singularity $P \in X$ as such; they apply in other situations, notably to a fibre $f^{-1}P$ of a surface fibred over a curve $f\colon X \to C$. The main technical methods used here are those of the previous chapters, divisors on a nonsingular surface, their cohomology, and so on. As I said for K3s at the start of Chapter 3, this material provides excellent practice for calculations involving intersection numbers, coherent cohomology, and results based on these.

## Summary

1. Examples: the ordinary double point and how it arises, the remaining Du Val singularities $A_n$, $D_n$, $E_6$, $E_7$, $E_8$ and their resolutions
2. Numerical cycle and multiplicity; there is a unique minimal exceptional divisor $0 \neq Z_{\text{num}}$, called the *numerical cycle*, with $Z_{\text{num}}\Gamma \le 0$ for every exceptional curve $\Gamma$. In simple cases, properties of $P \in X$ can be expressed in terms of $Z_{\text{num}}$
3. $R^1 f_* \mathcal{O}_Y$, and how you calculate it; there is a unique minimal divisor $Z_1$ on $Y$, the *cohomological cycle*, that carries the cohomology $H^1(\mathcal{O}_Y)$
4. Characterisations of rational singularities
5. Contractibility of a bunch of curves and application to Castelnuovo's criterion and to Du Val singularities. Minimal resolution, minimal models of surfaces
6. Gorenstein condition: $P \in X$ is Gorenstein if and only if the canonical class of $Y$ is represented by an exceptional divisor $Z_K$, the *canonical cycle*
7. Rational and Du Val singularities; how these relate to canonical models of surfaces, projective models of K3s, elliptic pencils
8. Elliptic Gorenstein surface singularities, their numerical characterisation
9. Graded rings on divisors, multiplicity and embedding dimension, proof of the main theorems on rational and elliptic Gorenstein surface singularities

## 4.1. Example: the ordinary double point

As a first example, consider the ordinary quadratic cone

$$P = (0,0) \in X : (xz = y^2) \subset \mathbb{A}^3.$$

This singularity occurs throughout the theory of algebraic surfaces, and can be used to illustrate a whole catalogue of arguments. Because $X$ is a cone with vertex $P$, it has the standard "cylinder" resolution $Y \to X$; that is, the cone is a union of generating lines through $P$, and $Y$ is the disjoint union of the generating lines. In other words, $Y$ is the blowup, or the correspondence between the cone and the base of the cone. This singularity and its resolution have already appeared in Chapter 2, where the surface scroll $\mathbb{F}_2 = \mathbb{F}(0,2)$ has a morphism to $\overline{\mathbb{F}}(0,2) \subset \mathbb{P}^3$ (the ordinary quadratic cone) contracting the negative section.

The exceptional curve $f^{-1}P = \Gamma$ of the resolution is a $-2$-*curve*, satisfying

$$\Gamma \cong \mathbb{P}^1, \quad \Gamma^2 = -2 \quad \text{(or equivalently, } K_Y\Gamma = 0).$$

A particular case of Contraction Theorem 4.15 is that every $-2$-curve is obtained in this way by resolving an ordinary double point (that is, a double point with nondegenerate tangent cone). On minimal surfaces (where there are no $-1$-curves), $-2$-curves are characterised as the irreducible curves with $K_Y\Gamma \le 0$ and $\Gamma^2 < 0$, and they arise naturally for this reason in many important contexts of the classification of surfaces.

The ordinary double point also appears in Chapter 1 on cubic surfaces, and in the theory of projective embeddings of K3s in Chapter 3. As we saw in Chapter 1, a nonsingular cubic surface is isomorphic to $\mathbb{P}^2$ blown up in 6 "general" points $P_1, \dots, P_6$, embedded by the linear system of cubics through $\{P_1, \dots, P_6\}$. If $P_1, P_2, P_3$ are collinear (but the $\{P_i\}$ are otherwise "general") then on the blowup $X = \mathrm{Bl}_6(\mathbb{P}^2)$, the line through $P_1, P_2, P_3$ becomes a $-2$-curve, and is contracted by $X \to \overline{X} \subset \mathbb{P}^3$ to an ordinary double point of a cubic surface. Similarly, as described in 3.15, if $X$ is a K3 and $|D|$ a free nonhyperelliptic linear system on $X$, the birational morphism $\varphi_D \colon X \to \overline{X} \subset \mathbb{P}^g$ contracts exactly curves $\Gamma$ with $D\Gamma = 0$; any such curve is a $-2$-curve (because $K_X\Gamma = 0$ and $\Gamma^2 < 0$), and if there is exactly one then $\varphi_D(\Gamma) = P \in \overline{X}$ is an ordinary double point.

## 4.2. The Du Val singularities

**Example:** $D_4$**.** Consider the singularity

$$P = (0,0,0) \in X : (g = x^2 + y^3 + z^3 = 0) \subset \mathbb{A}^3.$$

The blow-up $X_1 \to X$ is covered by 3 affine pieces, of which I only write down one: consider $\mathbb{A}^3$ with coordinates $x_1, y_1, z$, and the morphism $\sigma \colon \mathbb{A}^3 \to \mathbb{A}^3$ defined by $x = x_1 z, y = y_1 z, z = z$. The inverse image of $X$ under $\sigma$ is defined by

$$g(x_1 z, y_1 z, z) = x_1^2 z^2 + y_1^3 z^3 + z^3 = z^2 g_1, \qquad \text{where } g_1 = x_1^2 + (y_1^3 + 1)z.$$

Here the factor $z^2$ vanishes on the exceptional $(x_1, y_1)$-plane $\mathbb{A}^2 = \sigma^{-1}P \subset \mathbb{A}^3$, and the residual component $X_1 : (g_1 = 0) \subset \mathbb{A}^3$ is the birational transform of $X$. Now clearly the inverse image of $P$ under $\sigma\colon X_1 \to X$ is the $y_1$-axis, and $X_1$ has ordinary double points at the 3 points where $y_1^3+1=0$. (Please check for yourselves that the other affine pieces of the blow-up have no further singular points.) The resolution $Y$ is obtained on blowing up these three points $Y \to X_1 \to X$.

I claim that $f^{-1}P$ consists of $-2$-curves $\Gamma_0, \Gamma_1, \Gamma_2, \Gamma_3$ meeting as follows:

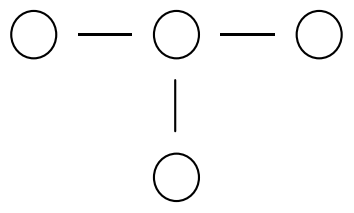

the configuration is the Dynkin diagram $D_4$. To prove this, it is clear that $\Gamma_1, \Gamma_2, \Gamma_3$ are $-2$-curves, since they arise from the blowup of ordinary double points. Also the fact that $\Gamma_0 \cong \mathbb{P}^1$, and $\Gamma_0$ meets each of $\Gamma_i$ transversally in 1 point, can be verified directly from the coordinate description of $Y$. Finally, to see that $\Gamma_0^2 = -2$, note that $y$ is a regular function on $Y$ whose divisor is $\operatorname{div} y = 2\Gamma_0 + \Gamma_1 + \Gamma_2 + \Gamma_3 + C$, where $C$ is the curve $y_1 = 0$ in $Y$, which also meets $\Gamma_0$ transversally in 1 point. Thus

$$0 = (\operatorname{div} y)\Gamma_0 = 2\Gamma_0^2 + \Gamma_0(\Gamma_1+\Gamma_2+\Gamma_3+C) = 2\Gamma_0^2 + 4, \qquad \text{so that } \Gamma_0^2 = -2.$$

The remaining Du Val singularities can be resolved in a similar way (see Ex. 4.3). The equations are

$$\begin{aligned}
A_n:&\qquad x^2+y^2+z^{n+1}=0,\\
D_n:&\qquad x^2+y^2z+z^{n-1}=0,\\
E_6:&\qquad x^2+y^3+z^4=0,\\
E_7:&\qquad x^2+y^3+yz^3=0,\\
E_8:&\qquad x^2+y^3+z^5=0,
\end{aligned}$$

and the resolution is a bunch of $-2$-curves whose configuration is given by the corresponding Dynkin diagram; for example, $E_7$ is the following:

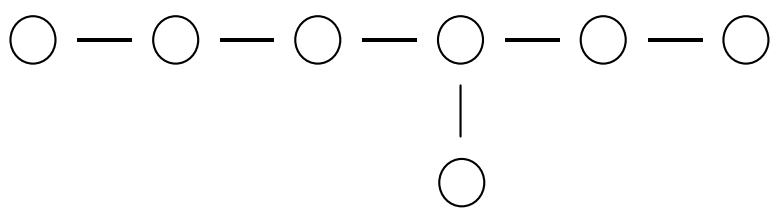

## 4.3. The general set-up

Let $X$ be an affine surface and $P \in X$ a singular point of $X$. I want to talk about isolated surface singularities $P \in X$; old-fashioned singularity theorists always write $(X, P)$, but my notation is shorter and self-documenting. It is reasonable to restrict attention to normal singularities, for reasons I will explain. So $X$ is a surface, that is, a 2-dimensional algebraic variety (over an algebraically closed field $k$, for example, $k = \mathbb{C}$), and $P \in X$ an isolated singular point. I'm usually only interested in a neighbourhood of $P \in X$, so without loss of generality, I can choose $X$ affine and $P \in X$ the only singularity. However, $X$ is contained in a projective variety $\overline{X}$, and for some purposes it may be more convenient to think of $X$ and $Y$ as projective.

A *resolution* of $X$ is a birational projective morphism $f\colon Y \to X$ from a nonsingular surface $Y$. It is an important and difficult theorem that a resolution always exists. I can assume that $f$ gives $Y \setminus f^{-1}P \xrightarrow{\sim} X \setminus P$. By definition, $f$ is *projective* if it factors as $Y \hookrightarrow X \times \mathbb{P}^N \to X$ with the first map a closed embedding. If $X$ is projective, the condition on $f$ is equivalent to $Y$ also projective. In general, $f^{-1}P$ is a connected bunch of projective curves $\bigcup_{i=1}^n \Gamma_i$; it is traditional to say that $f$ is a *good* resolution (or, more recently, a *log* resolution, meaning resolution in the category of *log varieties*) if $f^{-1}P$ is a union of nonsingular curves meeting transversally, but I do not necessarily assume this.

The surface $X$ is normal if and only if $\mathcal{O}_X = f_*\mathcal{O}_Y$. Here $\mathcal{O}_Y$ is the structure sheaf of $Y$. The point is just that as a topological space, the normal variety $X$ is the quotient space $Y/\sim$ where $\sim$ is the equivalence relation that collapses $f^{-1}P$ to a point; and $X$ has structure sheaf

$$f_*\mathcal{O}_Y = \big\{\text{regular functions on } Y \text{ that are constant on } f^{-1}P\big\}.$$

Any other variety structure $X'$ on this topological space (having a morphism $Y \to X'$) is obtained by taking a subring $\mathcal{O}_{X'} \subset \mathcal{O}_X$ of finite codimension, that is, $\mathcal{O}_{X'}$ is a subring containing $m_P^n$ for some $n$, where $m_P \subset \mathcal{O}_{X,P}$ is the maximal ideal. Normal implies nonsingular in codimension 1, and the converse holds for hypersurfaces or complete intersections by the Serre criterion, compare [**YPG**], App. to §3. If my display of erudition confuses you, just concentrate on the case of isolated hypersurface singularities, which are automatically normal, and lots of fun in their own right.

## 4.4. Changing to the analytic category

When studying singularities, I do not spend much time worrying about which category of geometry I work in. A natural alternative to the language used here is to work up to local analytic equivalence: any ordinary double point is locally isomorphic to the singularity $y^2 = xz$ of 2.1 by a suitable analytic coordinate change. This does not make sense in the category of algebraic varieties: ordinary double points appear on singular cubic surfaces, or on (contracted) K3 surfaces, etc., and then, of course, no Zariski neighbourhood of $P \in X$ is isomorphic to the cone $y^2 = xz$ in the category of algebraic varieties. My justification for being sloppy here is that everything I need to know about an isolated singularity $P \in X$ of a variety depends not on the local ring $\mathcal{O}_{X,P}$ itself, but only on a quotient $\mathcal{O}_{X,P}/m_P^N$ for some large $N$.

To be able to use the ideas of singularity theory, especially local normal forms, I could assume that the base field is $\mathbb{C}$, and use local complex geometry (typically, the inverse function theorem). Very roughly, since I am mostly concerned with the invariants of coherent cohomology, affine corresponds to Stein. On the other hand, it is not difficult to work over any field using formal completions or the language of the etale topology. Beware that the normal forms of the singularities used here may not work at all in characteristic $p$, especially when $p$ is small.

## 4.5. The numerical cycle of a resolution $Z_{\rm num}$

Given a resolution $f\colon Y \to X$ of $P \in X$, any nonzero element $g \in m_P \subset \mathcal{O}_{X,P}$ vanishes on the set theoretic fibre $f^{-1}P$, so that $\operatorname{div} g = \sum a_i\Gamma_i + C'$ is as in the proof of Theorem A.7. As we saw there, the exceptional divisor $D = \sum a_i\Gamma_i$ satisfies $D\Gamma_i \le 0$ for every contracted curve $\Gamma_i$.

**Proposition.** *There exists a unique effective exceptional divisor* $Z = Z_{\text{num}} = \sum b_i\Gamma_i > 0$ *such that* $Z > 0$,

$$Z\Gamma_i \leq 0 \qquad \textit{for every curve } \Gamma_i, \tag{1}$$

*and* $Z$ *is minimal with this property.*

**Proof.** There certainly exist some cycle $Z$ with this property, as I just said, so the only thing to prove is that if $Z_1$ and $Z_2$ both satisfy (1) then so does $Z = \gcd(Z_1, Z_2)$. This is more or less obvious: I can write $Z_1 = Z + A_1$ and $Z_2 = Z + A_2$, where $A_1, A_2 \geq 0$ have no common components; then any $\Gamma \not\subset A_1$ has $Z\Gamma = Z_1\Gamma - A_1\Gamma \leq Z_1\Gamma \leq 0$, and similarly for any $\Gamma \not\subset A_2$. Q.E.D.

**Definition.** $Z_{\text{num}}$ is called the *numerical cycle* (or *numerical fundamental cycle*) of the bunch $\{\Gamma_i\}$. The discussion preceding the proposition shows that every $g \in m_P$ vanishes at least along $Z_{\text{num}}$, that is, $\operatorname{div} g \geq Z_{\text{num}}$.

The numerical cycle is often called the *fundamental cycle* in the literature, but lots of other cycles are equally fundamental. In 4.6, I introduce the *fibre cycle* $Z_{\text{f}}$, which is closely related to the scheme theoretic fibre, and later in this chapter the *cohomological cycle* $Z_1$, which carries the cohomology $H^1(Y, \mathcal{O}_Y)$ (see 4.8), and the *canonical cycle* $Z_K$ for a Gorenstein singularity $P \in X$, which carries the canonical class of $Y$ in a neighbourhood of $f^{-1}P$ (see 4.19).

Note that $Z_{\text{num}}$ is very easy to calculate from the intersection matrix $\Gamma_i\Gamma_j$ (see Ex. 4.4).

## 4.6. Fibre cycle $Z_{\text{f}}$ and multiplicity

I define the *fibre cycle* $Z_{\text{f}}$ to be the greatest divisor contained in every $\operatorname{div} g$, that is,

$$Z_{\text{f}} = \gcd\bigl\{\operatorname{div} g \bigm| \text{for all } g \in m_P \subset \mathcal{O}_{X,P}\bigr\}.$$

In other words, consider the ideal sheaf $\mathcal{I} = m_P\mathcal{O}_Y$ of $\mathcal{O}_Y$ generated by all $g \in m_P$. The subscheme of $Y$ defined by $\mathcal{I}$ is the scheme theoretic fibre of the morphism $f$, that is, the subscheme of $Y$ defined by $g = 0$ taken over all $g \in m_P$. By definition the fibre cycle $Z_{\text{f}}$ is the maximal effective divisor contained in the scheme theoretic fibre, that is, such that $\mathcal{I} \subset \mathcal{O}_Y(-Z_{\text{f}})$. In more geometric terms, take a finite dimensional vector subspace $V \subset \mathcal{O}_{X,P}$ which generates $m_P$ as a $\mathcal{O}_{X,P}$-module, and for $g \in V$, consider the divisor $H_g$, a variable hyperplane through $P$. Then for fixed $g$, $\operatorname{div} g = f^*(H_g) = L + \sum b_i\Gamma_i$, where $L$ has no components in common with $\{\Gamma_i\}$ and $\sum b_i\Gamma_i$ is an effective divisor. If $g \in V$ is "fairly general" then $\sum b_i\Gamma_i$ takes on a certain minimal value, which is $Z_{\text{f}}$.

It follows from the remarks in 4.5 that $Z_{\text{f}} \geq Z_{\text{num}}$; equality holds in some important cases: see Theorem 4.17 for rational singularities and Theorem 4.23 for elliptic Gorenstein singularities. However, in more complicated cases, for example, when some of the curves $\Gamma_i$ have large genus, $Z_{\text{f}}$ is not determined by simple numerical properties of the configuration $\{\Gamma_i\}$, and it is not hard to find examples with $Z_{\text{f}} > Z_{\text{num}}$ (see Ex. 4.13).

**Definition.** I take the *multiplicity* of a normal surface singularity $P \in X$ in the pedestrian sense: $\operatorname{mult}_P X$ is the dimension over $k$ of $\mathcal{O}_{X,P}/(h_1, h_2)$ (also called its *length*), where $h_1, h_2 \in m_P$ are sufficiently general elements. In other words, this is the local intersection number at $P$ of two general hyperplane sections through $P$. If $X \subset \mathbb{P}^N$ is normal and of degree $d$ then any two hyperplane sections $H_1, H_2$ with no common components meet in $d$ points counted with multiplicity, so that two general hyperplanes through $P$ intersect in $P$ counted $\operatorname{mult}_P X$ times plus $d - \operatorname{mult}_P X$ free intersection points outside $P$. Thus the multiplicity of $P$ is the number of intersection points of $H_1, H_2$ absorbed[5] into $P$.

**Lemma.** *(i) Suppose that* $m_P\mathcal{O}_Y = \mathcal{O}_Y(-Z_\mathrm{f})$ *with* $Z_\mathrm{f}$ *the fibre cycle. Then* $\operatorname{mult}_P X = -Z_\mathrm{f}^2$.

*(ii) More generally, if* $m_P\mathcal{O}_Y = \mathcal{I}\mathcal{O}_Y(-Z_\mathrm{f})$ *with* $\mathcal{I}$ *an ideal defining a zero dimensional subscheme of* $Y$*, then* $\operatorname{mult}_P X \ge -Z_\mathrm{f}^2 + \dim_k \mathcal{O}_Y/\mathcal{I}$.

**Remark.** $m_P\mathcal{O}_Y = \mathcal{O}_Y(-Z_\mathrm{f})$ holds if and only if the resolution dominates the ordinary blowup of $P \in X$ (that is, the blowup of the maximal ideal $m_P \subset \mathcal{O}_X$). Because either condition says that the ratio between generators of $m_P$ is a regular function on $Y$. Thus, by appropriate choice of resolution, I can always arrange that case (i) holds.

**Proof.** What's going on is the following: Think of $X$ as a projective surface $X \subset \mathbb{P}^N$, and consider the linear system $|m_P\mathcal{O}_X(1)|$ of hyperplane sections $H \ni P$. Taking $f^*$ of these divisors defines a linear system $f^*|m_P\mathcal{O}_X(1)|$ on $Y$ with some base locus, which, following tradition, I break up as the divisorial part $Z_\mathrm{f}$ and the codimension $\ge 2$ part. The intersection of two elements $H_1, H_2 \in |m_P\mathcal{O}_X(1)|$ consists of $P$ with multiplicity $\operatorname{mult}_P X$ plus a free part; I can calculate the free part just as well on $Y$, and I get the intersection outside $f^{-1}P$ of two elements of $|f^*H - Z_\mathrm{f}|$.

In case (i), $|f^*H - Z_\mathrm{f}|$ is free. If $h_1, h_2 \in \mathcal{O}_{X,P}$ are two general elements, and $H_1, H_2 \subset X$ the corresponding divisors, then $f^*(H_1) = M_1 + Z_\mathrm{f}$ and $f^*(H_2) = M_2 + Z_\mathrm{f}$ with $M_1$ and $M_2$ disjoint over $P$, that is, $M_1 \cap M_2 \cap f^{-1}P = \emptyset$. Then $M_1M_2$ is the free intersection of $H_1, H_2$ outside $P$.

Then since $(f^*H)Z_\mathrm{f} = 0$, I get

$$\begin{aligned} H^2 = f^*(H_1)f^*(H_2) &= f^*(H_1)(M_2 + Z_\mathrm{f}) = f^*(H_1)M_2 \\ &= (M_1 + Z_\mathrm{f})M_2 = M_1M_2 + Z_\mathrm{f}(f^*(H_2) - Z_\mathrm{f}) \\ &= M_1M_2 - Z_\mathrm{f}^2. \end{aligned}$$

That is, the base locus $Z_\mathrm{f}$ of $f^*|m_P\mathcal{O}_X(1)|$ absorbs $-Z_\mathrm{f}^2$ points of intersection, as required.

(ii) is similar; you have to prove that two general curves through a zero dimensional base scheme at $P \in Y$ have local intersection number $\ge \dim_k \mathcal{O}_{Y,P}/\mathcal{I}_P$. It is in any case obvious that any points of a zero dimensional base locus make a strictly positive contribution. Q.E.D.

[5] This is proved exactly as the Bézout theorem in A.9. However, note that it is essential for this that $P \in X$ is normal, so that sufficiently general $h_1, h_2 \in m_P$ form a regular sequence, and the usual Koszul sequence $0 \to \mathcal{O}_{X,P} \to \mathcal{O}_{X,P} \oplus \mathcal{O}_{X,P} \to (h_1, h_2) \subset \mathcal{O}_{X,P}$ defined by $h_1, h_2$ is exact. In general, the right definition of the multiplicity $\operatorname{mult}_P X$ is in terms of the leading term of the Hilbert–Samuel function: $h(k) = \dim m_P^k/m_P^{k+1} = (\operatorname{mult}_P X)k + \text{const}$. For a nonnormal isolated surface singularity, the Koszul complex for any $h_1, h_2$ is *never* exact in the middle, and $\dim \mathcal{O}_{X,P}/(h_1, h_2)$ is strictly bigger than $\operatorname{mult}_P X$. See [**YPG**], App. to §3 for more discussion.

## 4.7. $R^1 f_*\mathcal{O}_Y$, and how to calculate it

Here I'm interested in the sheaf $R^1 f_*\mathcal{O}_Y$ for a resolution $f\colon Y \to X$. Assuming that $X$ is affine, we will see that this is the finite dimensional vector space $H^1(Y, \mathcal{O}_Y)$ at $P$.

**Proposition.** *Let $P \in X$ be a surface singularity, $f\colon Y \to X$ its resolution and $P \in U \subset X$ an affine neighbourhood of $P$. Then for any coherent sheaf $\mathcal{F}$ on $Y$,*

1. *$H^p(f^{-1}U, \mathcal{F}) = 0$ for $p \geq 2$, and*
2. *$H^1(f^{-1}U, \mathcal{F})$ is independent of $U$.*

**Definition.** $R^1 f_*\mathcal{F}$ is defined to be $H^1(f^{-1}U, \mathcal{F})$, viewed as a sheaf supported at the point $P \in X$. Note that $R^1 f_*\mathcal{F} = \varinjlim_U H^1(f^{-1}U, \mathcal{F})$, because by (2), the injective limit is constant for affine $U$. Since $H^1(f^{-1}U, \mathcal{F})$ is a $\Gamma(U, \mathcal{O}_X)$-module, it is clear that $R^1 f_*\mathcal{F}$ is a module over $\mathcal{O}_{X,P}$. I will prove shortly that it is a finite dimensional vector space.

**Proof.** $f^{-1}U$ is covered by 2 affine open sets $V_1, V_2$, so that the first assertion follows at once by calculating $H^p(f^{-1}U, \mathcal{F})$ by the Čech complex. If $U = \operatorname{Spec} A$ then arbitrarily small neighbourhoods of $P$ are given by $U_g = \operatorname{Spec} A[g^{-1}]$ for $g \in A$; the Čech complex of the corresponding cover $V_i \cap f^{-1}(U_g)$ is obtained by applying the exact functor $\otimes_A A[g^{-1}]$, and therefore

$$H^1(f^{-1}(U_g), \mathcal{F}) = H^1(f^{-1}U, \mathcal{F}) \otimes_A A[g^{-1}] = H^1(f^{-1}U, \mathcal{F}). \quad \text{Q.E.D.}$$

**Discussion.** My treatment here is self-contained, but I give a brief discussion of some background in sheaf theory for the interested reader. Quite generally, for a morphism $f\colon Y \to X$ and sheaf $\mathcal{F}$ on $Y$, the sheaf $R^p f_*\mathcal{F}$ is defined as the associated sheaf of the presheaf $U \mapsto H^p(f^{-1}U, \mathcal{F})$. This means that $R^p f_*\mathcal{F}$ is the sheaf on $X$ whose stalk at a point $P \in X$ is $(R^p f_*\mathcal{F})_P = \varinjlim_U H^p(f^{-1}U, \mathcal{F})$, where the $\varinjlim$ runs over opens $U \ni P$.

Note that $R^p f_*\mathcal{F}$ is defined in terms of cohomology of thinner and thinner open neighbourhoods of the fibre $f^{-1}P$. In contrast, the holomorphic functions theorem of Zariski and Grothendieck interprets $(R^p f_*\mathcal{F})_P$ for a coherent sheaf $\mathcal{F}$ in terms of the fibre $f^{-1}P$ and sheaves supported on it:

$$R^p f_*\mathcal{F} = \varprojlim_Z H^p(f^{-1}P, \mathcal{F} \otimes_{\mathcal{O}_Y} \mathcal{O}_Z),$$

where the limit is taken over subschemes $Z$ with support in $f^{-1}P$, which in practice means a limit over fatter and fatter infinitesimal neighbourhoods of the fibre. I prove this in the simple case I need in the next section (the general case is done in [**EGA III**$_1$], 4.1.5).

The etymology comes from the idea that an element of the left-hand side is something like a Taylor series expansion of a cohomology class around the fibre, and the point is to prove that it converges to a "holomorphic" cohomology class on a neighbourhood of the fibre.

## 4.8. Cohomological cycle

Let $f^{-1}P = \bigcup_{i=1}^n \Gamma_i$; as in 3.10, given a divisor $D = \sum n_i\Gamma_i$ with $n_i \in \mathbb{Z}$, $n_i \geq 0$, write $\mathcal{O}_Y(-D) \subset \mathcal{O}_Y$ for the ideal sheaf of regular functions on $Y$ vanishing along $D$, and $\mathcal{O}_D = \mathcal{O}_Y/\mathcal{O}_Y(-D)$. This is a sheaf of rings on $f^{-1}P$, the structure sheaf of a subscheme $D \subset Y$. If $D_1 \le D$ is another divisor supported on $\bigcup \Gamma_i$ then $D = D_1 + D_2$ with $D_1, D_2 > 0$, so by Lemma 3.10, the inclusion $\mathcal{O}_Y(-D) \subset \mathcal{O}_Y(-D_1)$ induces a short exact sequence

$$0 \to \mathcal{O}_{D_2}(-D_1) \to \mathcal{O}_D \to \mathcal{O}_{D_1} \to 0 \qquad (2)$$

**Theorem.** *(a) For every $D_1 \le D$ the induced map $H^1(\mathcal{O}_D) \to H^1(\mathcal{O}_{D_1})$ is surjective; the dimension $h^1(\mathcal{O}_D)$ is bounded above.*

*(b) There is a unique smallest $Z_1$ such that $h^1(\mathcal{O}_{Z_1})$ takes the maximal value. Then for every $D \ge Z_1$, the map of (a) is an isomorphism $H^1(\mathcal{O}_D) = H^1(\mathcal{O}_{Z_1})$, and $h^1(\mathcal{O}_D) < h^1(\mathcal{O}_{Z_1})$ if $D \not\ge Z_1$.*

*(c) The exact sequence of sheaves $0 \to \mathcal{O}_Y(-D) \to \mathcal{O}_Y \to \mathcal{O}_D \to 0$ induces a surjective homomorphism $R^1 f_*\mathcal{O}_Y \to H^1(\mathcal{O}_D)$, which is an isomorphism for any $D \ge Z_1$.*

**Definition.** The divisor $Z_1$ is called the *cohomological cycle* of $P \in X$.

The point of the theorem is that $R^1 f_*\mathcal{O}_Y$, defined in 4.7 in terms of open neighbourhoods $f^{-1}U$ of the fibre, is equal to the cohomology group $H^1(\mathcal{O}_{Z_1})$ of an exceptional divisor. This is the same thing as $\varprojlim_D H^1(f^{-1}P, \mathcal{O}_D)$, where the limit runs over fatter and fatter exceptional divisors $D$, as discussed in 4.7. Moreover, the limit is achieved at $Z_1$. The theorem is due to M. Artin in modern form, although it seems to have been known in substance to P. Du Val some 30 years earlier.

**Proof.** (a) The surjectivity of $H^1(\mathcal{O}_D) \to H^1(\mathcal{O}_{D_1})$ follows from (2) since by dimension

$$H^2(\mathcal{O}_{D_2}(-D_1)) = 0.$$

To prove boundedness, note that there exists an effective exceptional divisor $Z = Z_{\text{ample}} = \sum a_i\Gamma_i$ such that $-Z\Gamma_i \ge 0$ for every $\Gamma_i$. This follows from negative definiteness. In other words, $-Z$ is relatively ample for $f$.

It follows that there exists some exceptional divisor $D$ such that $-D\Gamma_i \ge K_Y\Gamma_i$ for every $i$ (for example, take $D$ to be a large multiple of $Z_{\text{ample}}$). Then for every $D' > 0$,

$$H^1(\mathcal{O}_{D'}(-D)) \overset{\mathrm{d}}{=} H^0(\mathcal{O}'_D(K_Y + D' + D))$$

(by Serre duality on $D'$, see 4.10). But you can see that $H^0(\mathcal{O}_{D'}(D') \otimes \mathcal{L}) = 0$ if $\deg_{\Gamma_i} \mathcal{L} < 0$ for all $i$ as an easy exercise (compare Ex. 4.14). Therefore $H^1(\mathcal{O}_{D+D'}) \to H^1(\mathcal{O}_D)$ is also injective for every $D' \ge 0$.

(b) I need to prove that if $D'$ and $D''$ are divisors for which $h^1(\mathcal{O}_{D'}) = h^1(\mathcal{O}_{D''})$ both take the maximum value, then $A = \gcd(D', D'')$, has the same property, that is, $h^1(\mathcal{O}_A) = h^1(\mathcal{O}_{D'}) = h^1(\mathcal{O}_{D''})$.

So let $D' = A + B'$ and $D'' = A + B''$, where $B'$ and $B''$ have no common components, and $C = A + B' + B'' = \operatorname{lcm}(D', D'')$. Then I claim that there exists an exact sequence

$$H^1(\mathcal{O}_C) \to H^1(\mathcal{O}_{D'}) \oplus H^1(\mathcal{O}_{D''}) \to H^1(\mathcal{O}_A) \to 0. \qquad (3)$$

This proves what I want, since $C, D', D''$ all have the same value of $h^1$, so

$$h^1(\mathcal{O}_A) \ge h^1(\mathcal{O}_{D'}) + h^1(\mathcal{O}_{D''}) - h^1(\mathcal{O}_C).$$

(3) comes by chasing the diagram

$$\begin{array}{ccccccc}
 & & H^1(\mathcal{O}_{B''}(-D')) & \to & H^1(\mathcal{O}_{B''}(-A)) & \to & 0 \\
 & & \downarrow & & \downarrow & & \\
H^1(\mathcal{O}_{B'}(-D'')) & \to & H^1(\mathcal{O}_C) & \to & H^1(\mathcal{O}_{D''}) & \to & 0 \\
\downarrow & & \downarrow & & \downarrow & & \\
H^1(\mathcal{O}_{B'}(-A)) & \to & H^1(\mathcal{O}_{D'}) & \to & H^1(\mathcal{O}_A) & \to & 0 \\
\downarrow & & \downarrow & & \downarrow & & \\
0 & & 0 & & 0 & &
\end{array}$$

here $B'$ and $B''$ are without common components, so that $\mathcal{O}_{B'\cap B''}$ has support the finite set $B'\cap B''$, and $H^1(\mathcal{O}_{B'\cap B''}) = 0$ together with the exact sequence

$$0 \to \mathcal{O}_{B''}(-D') \to \mathcal{O}_{B''}(-A) \to \mathcal{O}_{B'\cap B''} \to 0$$

provides the exactness of the top row and the left-hand column. I leave you the pleasure of the diagram chase.

(c) The surjectivity of $H^1(\mathcal{O}_Y) \to H^1(\mathcal{O}_D)$ follows from the assertion $H^2 = 0$ of Proposition 4.7. The point is to prove that it is injective for some $D$, so that in particular $H^1(\mathcal{O}_Y)$ is finite dimensional. If $Z = Z_{\text{ample}}$ is the divisor constructed in (a) then the sheaf $\mathcal{O}_Y(-Z)$ is relatively ample on $Y$ (see Ex. 4.16, essentially by the Nakai–Moishezon ampleness criterion[6]). From this, $H^1(\mathcal{O}_Y(-nZ)) = 0$ for $n \gg 0$ by Serre vanishing, and the result follows. Q.E.D.

**4.9. Corollary.**

1. $R^1f_*\mathcal{O}_Y = 0$ *if and only if* $H^1(\mathcal{O}_D) = 0$ *for every* $D$.
2. $R^1f_*\mathcal{O}_Y$ *is a* 1*-dimensional vector space if and only if* $H^1(\mathcal{O}_D) \le 1$ *for every* $D$ *and* $= 1$ *for some* $D$.

If $P \in X$ is an affine neighbourhood, then $R^1f_*\mathcal{O}_Y$ is the sheaf consisting of $H^1(Y,\mathcal{O}_Y)$ supported at $P$, so that the conditions in (1) and (2) characterise respectively $H^1(Y,\mathcal{O}_Y) = 0$ and $h^1(Y,\mathcal{O}_Y) = 1$.

## 4.10. Serre duality

Let $Y$ be a projective nonsingular surface and $D$ an effective divisor on $Y$. I want to use the fact that the line bundle $\mathcal{O}_D(K_Y + D)$ has the properties of a dualising sheaf for $D$ (for example, in the above proof of Theorem 4.8). I don't need the full strength.

**Proposition.** *If* $L$ *is a divisor on* $Y$ *then there is a natural duality*

$$H^i(\mathcal{O}_D(L)) \overset{\mathrm{d}}{=} H^{1-i}(\mathcal{O}_D(K_Y + D - L)) \qquad \textit{for } i = 0, 1.$$

[6]This is the central point of the proof. Nakai–Moishezon says that if $H$ is a Cartier divisor on an $n$-dimensional complete scheme, then $H^rZ > 0$ for every $r$-dimensional irreducible subvariety $Z$ implies that $H$ is ample. This is completely elementary if you know coherent cohomology: by a dévissage argument, you can see that the numerical condition implies that $H^i(\mathcal{F}\otimes\mathcal{O}(nH)) = 0$ for every $i > 0$, for every coherent sheaf $\mathcal{F}$, and for $n \gg 0$ (see [**H1**], Chapter V, 1.10). I should have included this as a guided exercise in the section on cohomology.

Compare Ex. 3.26, or see [**H1**], III.7 and [**R1**], Theorem 2.12 for a more general treatment of the dualising sheaf $\omega_D$ and the adjunction formula

$$\omega_D = \mathcal{O}_D(K_D) = \mathcal{O}_D(K_Y + D) = \omega_Y(D)_{|D}.$$

**Proof.** This follows formally from duality on $Y$: the short exact sequence

$$0 \to \mathcal{O}_Y(L - D) \to \mathcal{O}_Y(L) \to \mathcal{O}_D(L) \to 0$$

leads to

$$\begin{aligned} \cdots \to H^i(\mathcal{O}_Y(L-D)) \to H^i(\mathcal{O}_Y(L)) \to H^i(\mathcal{O}_D(L)) \\ \to H^{i+1}(\mathcal{O}_Y(L-D)) \to H^{i+1}(\mathcal{O}_Y(L)) \to \cdots \end{aligned}$$

which by duality on $Y$ is dual term by term to

$$\begin{aligned} \cdots \leftarrow H^{2-i}(\mathcal{O}_Y(K_Y + D - L)) \leftarrow H^{2-i}(\mathcal{O}_Y(K_Y - L)) \leftarrow \cdots \\ \cdots \leftarrow H^{1-i}(\mathcal{O}_Y(K_Y + D - L)) \leftarrow H^{1-i}(\mathcal{O}_Y(K_Y - L)) \leftarrow \cdots \end{aligned}$$

But this can be identified with the cohomology long exact sequence of

$$0 \to \mathcal{O}_Y(K_Y - L) \to \mathcal{O}_Y(K_Y + D - L) \to \mathcal{O}_D(K_Y + D - L) \to 0,$$

which identifies $H^{1-i}(\mathcal{O}_D(K_Y + D - L))$ with the dual of $H^i(\mathcal{O}_D(L))$, as required. Q.E.D.

## 4.11. RR on a divisor $D$

In the same notation, define $\deg_D L = DL = \sum_i n_i \deg_{\Gamma_i} L$. Then

$$\chi(\mathcal{O}_D) = -\tfrac{1}{2}(D^2 + DK_Y),$$

and

$$\chi(\mathcal{O}_D(L)) = \chi(\mathcal{O}_D) + \deg_D L.$$

The proof is an exercise in same style as 4.10, using RR on $Y$. Notice that the formulas discussed in 3.6 involving the arithmetic genus $p_aD$,

$$\chi(\mathcal{O}_D) = 1 - p_aD \quad \text{and} \quad 2p_aD - 2 = (D + K_Y)D$$

are a particular case. If $D = C$ is a nonsingular curve, then $g(C) = p_aC = H^1(\mathcal{O}_C)$, and the formula here reduces to the usual RR formula $\chi(\mathcal{O}_C(L)) = 1 - g + \deg L$.

## 4.12. Rational singularities

**Definition.** Let $f\colon Y \to X$ be a resolution of a normal surface singularity $P \in X$. Then $P \in X$ is a *rational singularity* if $R^1f_*\mathcal{O}_Y = 0$. It is an *elliptic singularity* if $R^1f_*\mathcal{O}_Y$ is 1-dimensional. Elliptic Gorenstein singularities are treated in 4.21–23 below. Without the Gorenstein condition, elliptic is too weak a condition to be very interesting.

**Proposition** (Numerical characterisation of rational singularities).

1. *$P \in X$ is rational if and only if*
$$\chi(\mathcal{O}_D) \geq 1 \quad \textit{for every effective divisor } D \textit{ supported on } f^{-1}P.$$
2. *For a rational singularity, the numerical cycle $Z_{\text{num}}$ is the maximal divisor with $\chi(\mathcal{O}_D) = 1$. In other words,*
$$h^0(\mathcal{O}_D) \geq 2 \quad \textit{for every } D > Z_{\text{num}}.$$
3. *$Z_{\text{num}}$ is numerically 1-connected (see Definition 3.9).*

**Remark.** By the adjunction formula and RR,
$$(K_Y + D)D = 2p_aD - 2 = -2\chi(\mathcal{O}_D),$$
so that the condition in (1) is equivalent to $p_aD \leq 0$ or $(K_Y + D)D < 0$.

**Proof.** (1) The implication $\Longrightarrow$ is clear: $R^1f_*\mathcal{O}_Y = 0$ means that $H^1(\mathcal{O}_D) = 0$ for every $D$, and $H^0(\mathcal{O}_D) \neq 0$ (because it contains the constant functions $k \subset H^0(\mathcal{O}_Y) \to H^0(\mathcal{O}_D)$), so $\chi(\mathcal{O}_D) \geq 1$.

For $\Longleftarrow$, if $R^1f_*\mathcal{O}_Y \neq 0$, there exists a divisor $D > 0$ such that $H^1(\mathcal{O}_D) \neq 0$, and therefore there exists a minimal one; that is $H^1(\mathcal{O}_D) \neq 0$, but $H^1(\mathcal{O}_{D-\Gamma}) = 0$ for every $\Gamma \subset D$.

By considering the exact sequence
$$H^1(\mathcal{O}_\Gamma(-D + \Gamma)) \to H^1(\mathcal{O}_D) \to H^1(\mathcal{O}_{D-\Gamma}) = 0,$$
we see that $H^1(\mathcal{O}_\Gamma(-D + \Gamma)) \neq 0$ for every $\Gamma \subset D$. It follows by Serre duality that $H^0(\mathcal{O}_\Gamma(K_Y + \Gamma + D - \Gamma)) \neq 0$, that is, $H^0(\mathcal{O}_\Gamma(K_Y + D)) \neq 0$. Now $\Gamma$ is an irreducible curve, and $\mathcal{O}_\Gamma(K_Y + D)$ a line bundle over it, so that the existence of a section implies deg $\geq 0$, that is $\Gamma(K_Y + D) \geq 0$ for every $\Gamma \subset D$. Summing over all components of $D$ gives $D(K_Y + D) \geq 0$, that is $p_aD \geq 1$.

(2) and (3) are easy exercises (see Ex. 4.11–12). Q.E.D.

**Example.** Consider the bunch of $\mathbb{P}^1$s consisting of a central $-3$-curve $\times$ meeting 3 simple chains of $-2$-curves of length $p-1, q-1, r-1$ respectively:

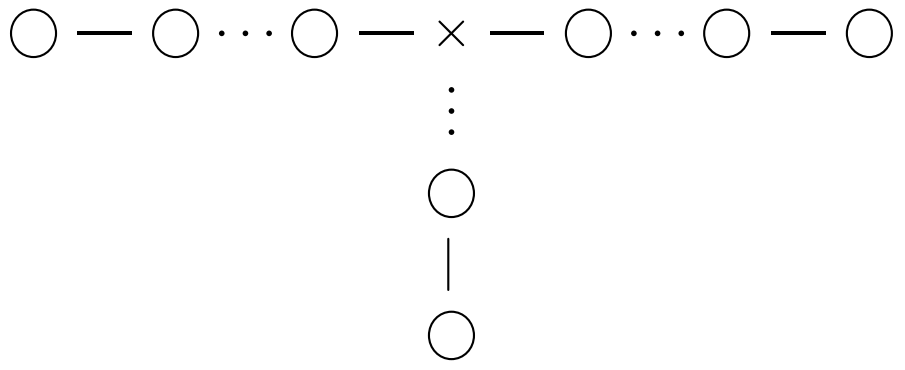

Check that $Z_{\text{num}}$ is the reduced sum of all the curves, and that the configuration is rational. It's a fun calculation to obtain this configuration by resolving the singular surface $X \subset \mathbb{C}^4$ defined by
$$\begin{aligned} xy + xt^q + t^{p+q} &= 0 \\ yz + yt^r + t^{q+r} &= 0 \\ xz + zt^p + t^{p+r} &= 0 \end{aligned} \qquad \text{or} \qquad \operatorname{rank} \begin{pmatrix} z + t^r & -x & -t^q \\ -t^r & x + t^p & -y \\ -z & -t^p & y + t^q \end{pmatrix} \leq 1.$$
[Hint: Check that if $p = q = r = 1$, you get the cone over the rational normal curve of degree 3. Then check that if $p = 0$, $X$ has two Du Val singularities of type $A_{q-1}, A_{r-1}$ on the $x$-axis. Finally, if $p, q, r > 1$, the blowup ($x = x_1t$, etc.) gives 3 coordinate axes as exceptional locus, through a point of the same type with numbers $p-1, q-1, r-1$, and you continue inductively.]

## 4.13. Structure of Pic $D$

Let $\{\Gamma_i\}_{i=1}^k$ be a bunch of curves on $Y$ and $D = \sum n_i\Gamma_i$ an effective divisor with $n_i > 0$ for each $i$. $\mathcal{O}_D^*$ denotes the sheaf of units of $\mathcal{O}_D$, so that $\mathcal{O}_Y^* \to \mathcal{O}_D^*$ is surjective (sheaves for the Zariski topology). Then

$$H^1(\mathcal{O}_D^*) = \big\{\text{isomorphism classes of invertible sheaves over } \mathcal{O}_D\big\} = \operatorname{Pic} D.$$

For any curve $\Gamma = \Gamma_i$, there is a natural group homomorphism

$$\rho_\Gamma\colon H^1(\mathcal{O}_D^*) \to H^1(\mathcal{O}_\Gamma^*) = \operatorname{Pic}\Gamma.$$

Geometrically, this takes an invertible sheaf $L$ over $D$ to $L_{|\Gamma} = L \otimes \mathcal{O}_\Gamma$. Now compose with the degree map $\operatorname{Pic}\Gamma \to \mathbb{Z}$ to get

$$\deg_\Gamma\colon H^1(\mathcal{O}_D^*) = \operatorname{Pic} D \to \mathbb{Z}.$$

Putting together this map for each component $\Gamma_i$ gives a map

$$\deg\colon H^1(\mathcal{O}_D^*) \to \mathbb{Z}^k$$

defined by $L \mapsto \deg_{\Gamma_i} L$. It is easy to see that deg is surjective: choose a general point $Q$ on any component $\Gamma \subset D$, and construct a Cartier divisor with support $Q$ and degree 1 on $\Gamma$, whose local equation $s \in \mathcal{O}_{D,Q}$ restricts to a local equation of $Q$ in $\mathcal{O}_{\Gamma,Q}$. (This operation is important in what follows, and is called a *transverse cut* at $Q$.)

**4.14. Proposition** (Artin). *(i) Write $D' = D_{\text{red}}$. Then the two groups*

$$K = \ker\{H^1(\mathcal{O}_D) \to H^1(\mathcal{O}_{D'})\}$$

*and*

$$K^* = \ker\{H^1(\mathcal{O}_D^*) \to H^1(\mathcal{O}_{D'}^*)\}$$

*have filtrations with isomorphic quotients.*

*(ii)* deg: $H^1(\mathcal{O}_D^*) \to \mathbb{Z}^k$ *is an isomorphism if and only if* $H^1(\mathcal{O}_D) = 0$.

**Proof.** (i) Suppose $C$ is a divisor such that $2C \subset D$; then define the subsheaves $J \subset \mathcal{O}_D$ and $J^* \subset \mathcal{O}_D^*$ to be the kernels in the exact sequences

$$0 \to J \to \mathcal{O}_D \to \mathcal{O}_{D-C} \to 0$$
$$0 \to J^* \to \mathcal{O}_D^* \to \mathcal{O}_{D-C}^* \to 0$$

(here of course $J = \mathcal{O}_C(-D+C)$ by Lemma 3.10); then $J$ and $J^*$ are isomorphic as sheaves of groups, because $J^2 = 0$, so that

$$\varepsilon\colon \Gamma(U,J) \xrightarrow{\sim} \Gamma(U,J^*) \quad \text{defined by } x \mapsto 1 + x$$

gives an isomorphism from the additive to the multiplicative group, a kind of "first order exponential". Also it is clear that $H^0(\mathcal{O}_{D-C}^*) \subset H^0(\mathcal{O}_{D-C})$ (as a set), and contains a Zariski neighbourhood of 1. A trivial cocycle calculation based on this proves that

$$H^0(\mathcal{O}_{D-C}) \to H^1(J) \quad \text{and} \quad H^0(\mathcal{O}^*_{D-C}) \to H^1(J^*) = \varepsilon(H^1(J))$$

have the same images. Hence

$$\ker\{H^1(\mathcal{O}_D) \to H^1(\mathcal{O}_{D-C})\} \cong \ker\{H^1(\mathcal{O}^*_D) \to H^1(\mathcal{O}^*_{D-C})\}.$$

(Of course, the right-hand group is a priori not a $k$-vector space.) The argument continues by induction from $D$ down to $D'$; $\ker\{H^1(\mathcal{O}^*_D) \to H^1(\mathcal{O}^*_{D'})\}$ and $\ker\{H^1(\mathcal{O}_D) \to H^1(\mathcal{O}_{D'})\}$ thus have filtrations by subgroups with isomorphic quotients.

(ii) For a reduced curve $D'$, it can be shown easily that

$$\operatorname{Pic} D' \cong \mathbb{Z}^k \iff H^1(\mathcal{O}_{D'}) = 0.$$

$\Longleftarrow$ By RR, a line bundle $L$ on $D'$ with $\deg_\Gamma L = 0$ for every $\Gamma \subset D'$ has

$$h^0(L) \geq 1 - h^1(\mathcal{O}_{D'}),$$

so $H^1(\mathcal{O}_{D'}) = 0 \implies H^0(L) \neq 0$. Then it is not hard to see that the section $s$ has the property that $\mathcal{O}_{D'} \cdot s = L$.

$\Longrightarrow$ If $L$ is a line bundle with sufficiently high degree on every component of $D'$ then $H^1(L) = 0$, so by RR, $h^0(L) = 1 - H^1(\mathcal{O}_{D'}) + \deg_{D'} L$. If I pick $\deg_{\Gamma_i} L$ "general" points on each component $\Gamma_i$ then obviously the resulting set of points $\{P_i\}$ gives a Cartier divisor $\sum P_i$ on $D'$ with the same degrees as $L$, so $L_0 = L(-\sum P_i)$ is a line bundle of degree 0 on each component. But a simple argument shows that these points can be chosen one at a time to strictly decrease the dimension $h^0(L(-\sum P_i))$ until we get to zero, so that if $H^1(\mathcal{O}_{D'}) \neq 0$, I get a line bundle $L_0$ of degree zero on each component, but $h^0(L_0) \neq 0$, so that $L_0 \cong \mathcal{O}_{D'}$. Q.E.D.

**Remark.** For a projective variety $X$, the group $\operatorname{Pic} X$ is an extension of a discrete group by the group $\operatorname{Pic}^0 X$ of divisors algebraically equivalent to zero modulo divisors linearly equivalent to zero. The elements of $\operatorname{Pic}^0 X$ are naturally parametrised (in the sense of moduli problems) by an algebraic group or algebraic group scheme, the *Picard variety* of $X$. If $X$ is nonsingular in characteristic zero, then $\operatorname{Pic}^0 X$ is an Abelian variety. The Jacobian of a curve is a well-known example. For singular, or nonreduced $X$, the algebraic group $\operatorname{Pic}^0 X$ is in general more complicated; however, it has the structure of an extension of an Abelian variety by an affine group (that is, an extension of additive and multiplicative groups $\mathbb{G}_a$s and $\mathbb{G}_m$s), and the Abelian variety comes entirely from the resolution of reduced components of $X$. Compare the discussion in 3.14.

## 4.15. Contractibility of a rational bunch of curves

Let $Y$ be a nonsingular quasiprojective surface and $\{\Gamma_i\}$ a bunch of curves on $Y$. By Propositions 4.12 and 4.14, we know that the following 3 conditions are equivalent:

1. $H^1(\mathcal{O}_D) = 0$ for every divisor $D$ supported on $\bigcup \Gamma_i$;
2. $p_a D \leq 0$ for every $D$;
3. $\deg\colon H^1(\mathcal{O}^*_D) = \operatorname{Pic} D \xrightarrow{\sim} \mathbb{Z}^k$ for every $D$, where $k$ is the number of components of $D$.

If $f\colon Y \to X$ is a resolution of a normal surface singularity $P \in X$ with exceptional locus $\{\Gamma_i\} = f^{-1}P \subset Y$, these conditions are necessary and sufficient for $P \in X$ to be rational.

**Contraction Theorem** (Artin). *A connected bunch of curves $\{\Gamma_i\}$ on $Y$ is the exceptional locus of a resolution of a rational singularity $P \in X$ if and only if*

*(a) the intersection matrix $(\Gamma_i\Gamma_j)$ is negative definite;*

*(b) $p_a D \le 0$ for every $D$.*

**Examples.** (i) A $-1$-*curve* on a surface $Y$ is a curve $L$ satisfying $L \cong \mathbb{P}^1$, $L^2 = -1$ and $K_Y L = -1$. This satisfies the rationality condition, so can be contracted.

Let's check the rationality condition: obviously if $L = \mathbb{P}^1$ then $H^1(\mathcal{O}_L) = 0$, and for $n \ge 2$,

$$0 \to \mathcal{O}_L(-(n-1)L) \to \mathcal{O}_{nL} \to \mathcal{O}_{(n-1)L} \to 0$$

gives

$$0 = H^1(\mathbb{P}^1, \mathcal{O}(n-1)) \to H^1(\mathcal{O}_{nL}) \to H^1(\mathcal{O}_{(n-1)L}) \to 0.$$

So $H^1(\mathcal{O}_{nL}) = 0$ for all $n$.

This case has the special feature that after contraction, the point $P \in X$ is a nonsingular point, and $Y \to X$ is the usual blowup of $P$ (this is proved in Theorem 4.17). In this case the contraction theorem is Castelnuovo's contractibility condition.

(ii) Similarly, for any $d > 0$, a curve $L \cong \mathbb{P}^1$, with $L^2 = -d$ satisfies the rationality condition. If $d = 1, 2$ we get $-1$-curve, $-2$-curve respectively.

(iii) A bunch of $-2$-curves $\{\Gamma_i\}$ has negative definite intersection matrix if and only if it forms one of the configurations $A_n, D_n, E_6, E_7, E_8$; if this happens then $p_a D \le 0$ for every effective divisor $D$ supported on $\Gamma_i$. This can be seen as follows: $K_Y \Gamma_i = 0$ by the adjunction formula, and so $K_Y D = 0$, so $K_Y D + D^2 < 0$, and $p_a D = \frac{1}{2}(K_Y + D)D + 1 \le 0$. In this case, Theorem 4.17 proves that the image $P \in X$ is a hypersurface double point.

## 4.16. Proof of Theorem

The conditions are clearly necessary (compare Theorem A.7 and Proposition 4.12). I now assume (a) and (b), and prove that the contraction exists. Suppose that $H$ is ample on $Y$. Then $H\Gamma_i > 0$ for every curve $\Gamma_i$. Since the matrix $(\Gamma_i\Gamma_j)$ is nonsingular, there exists rational numbers $q_i$ such that $(H + \sum q_i\Gamma_i)\Gamma_j = 0$ for every $j$. By taking $n$ to be any common denominators of the $q_i$, I get a divisor $D = \sum a_i\Gamma_i$ such that $(nH + D)\Gamma_i = 0$ for every $\Gamma_i$. Take $n$ to be sufficiently large, in a sense to be specified below. This means in particular that $D\Gamma_i < 0$ for every $i$.

I claim that $a_i > 0$: indeed, write $D = A - B$ with $A, B$ effective and with no common components. Then $DB = AB - B^2 \ge 0$. So $B = 0$, or in other words $D > 0$; but then the only way that $D\Gamma_i < 0$ is possible is for $\Gamma_i \subset D$, so $a_i > 0$.

Consider the invertible sheaf $\mathcal{L} = \mathcal{O}_Y(nH + D)$; by construction, $\deg_{\Gamma_i} \mathcal{L} = 0$ for each $i$. So by Proposition 4.14, (ii), $\mathcal{L}_{D'} \cong \mathcal{O}_{D'}$ for any divisor $D'$ supported on $\{\Gamma_i\}$.

Consider the restriction exact sequence

$$0 \to \mathcal{O}_Y(nH) \to \mathcal{O}_Y(nH + D) \to \mathcal{O}_D \to 0.$$

Now $\mathcal{O}_Y(H)$ is ample, so that by Serre vanishing,

$$H^1(\mathcal{O}_Y(nH)) = 0 \quad \text{for every } n \gg 0,$$

and so $H^0(\mathcal{O}_Y(nH+D)) \to H^0(\mathcal{O}_D) \to 0$ is exact.

I now claim that a suitable space of sections of $\mathcal{O}_Y(nH+D)$ defines the required contraction $f\colon Y \to X'$, at least set theoretically. Let $n$ and $s_0, \dots, s_N \in H^0(\mathcal{O}_Y(nH))$ be chosen so that $\varphi\colon Y \to \mathbb{P}^N$ defined by $(s_0 : \cdots : s_N)$ is an embedding; and let $s_{N+1} \in H^0(\mathcal{O}_Y(nH+D))$ be any element which maps to the identity section $1 \in H^0(\mathcal{O}_D)$. Consider the sections $(s_0, \dots, s_{N+1}) \in H^0(Y, \mathcal{L})$, and the map $f'\colon Y \to \mathbb{P}^{N+1}$ by $(s_0 : \cdots : s_{N+1})$.

Then (a) at every point $Q \in Y$, some $s_i \in \mathcal{L}$ is a local generator; (take $s_i$ if $Q \notin D$, and $s_{N+1}$ if $Q \in D$).

(b) $f'(D) = P = (0, \dots, 0, 1)$.

(c) $f'_{|(Y \setminus D)}$ is an embedding (since already the ratios $(s_0 : \cdots : s_N)$ give an embedding). Write $X' = f'(Y)$. Then there is a unique way to factorise the map $f'$ as $Y \xrightarrow{f} X \to X'$ where $X$ is normal, and this is the required contraction. $X$ is projective since $X'$ is, and the normalisation $X \to X'$ is finite. Q.E.D.

## 4.17. Multiplicity and embedding dimension

To understand what singularities arise on making contractions such as those in Theorem 4.15, I need to be able to determine invariants of a singularity $P \in X$ in terms of its resolution $Y$ and the bunch of curves $\{\Gamma_i\}$. Multiplicity and embedding dimension are two simple invariants, which give an easy characterisation of the nonsingular surface point $P \in X$ in Castelnuovo's criterion, and hypersurface double points arising on contracting a bunch of $-2$-curves.

Multiplicity was defined in 4.6 above. The *embedding dimension* of $P \in X$ is the dimension of the Zariski tangent space $\operatorname{emb}\dim(P \in X) = \dim m_P/m_P^2$. This is clearly the dimension of the smallest ambient space in which an analytic neighbourhood $P \in X$ can be embedded.

**Theorem** (Artin). *For a rational singularity $P \in X$, let $Z = Z_{\mathrm{num}}$ be the numerical cycle (see Definition 4.5), and set $d = -Z^2$. Then*

1. $m_P\mathcal{O}_Y = \mathcal{O}_Y(-Z)$, *and therefore* $\operatorname{mult}(P \in X) = d$;
2. $m_P/m_P^2 = \mathcal{O}_Z(-Z)$, *and* $\operatorname{emb}\dim(P \in X) = d+1$.

*In particular, $Z^2 = -1$ implies $P \in X$ is nonsingular, and $Z^2 = -2$ implies $P \in X$ is a hypersurface double point.*

**Remarks.** (a) The proof of (2) gives at once that $m_P^k/m_P^{k+1} = \mathcal{O}_Z(-kZ)$ and $\dim m_P^k/m_P^{k+1} = kd+1$ for every $k \geq 1$, so that the knowledgeable reader who so prefers can use the more professional definition of multiplicity in terms of the Hilbert–Samuel function (compare the footnote to 4.6).

(b) The example of the cone over the rational normal curve helps remember the numbers, but don't be misled into thinking that rational surface singularities are all that simple.

**Lemma.** *Let $f\colon Y \to X$ be a resolution of a rational surface singularity $P \in X$. Suppose that $X$ is affine, and let $\mathcal{L}$ be a line bundle on $Y$ which is nef on $f^{-1}P$ (that is, $\mathcal{L}\Gamma_i \geq 0$ for each $i$).*

*Then* $H^1(Y, \mathcal{L}) = 0$ *and* $|\mathcal{L}|$ *is free near* $f^{-1}P$.

**Proof.** Given any exceptional divisor $D$, I can construct a divisor $A$ on $D$ as a sum of local transverse cuts through general points of the $\Gamma_i$, so that $A$ is an effective Cartier divisor on $D$ with the same degree as $\mathcal{L}$; then by Proposition 4.14, (ii), $\mathcal{L}_{|D} \cong \mathcal{O}_D(A)$. Now the exact sequence of sheaves

$$0 \to \mathcal{O}_D \to \mathcal{O}_D(A) \to \mathcal{O}_A(A) \to 0$$

shows that $H^1(\mathcal{O}_D(A)) = 0$. Now $H^1(Y, \mathcal{L}) = R^1 f_* \mathcal{L}$ can be calculated by the arguments of 4.8 as an inverse limit of the groups $H^1(\mathcal{L}_{|D})$, all of which are zero (compare Ex. 4.15).

The transverse cuts making up $A$ can be chosen freely, so that the linear system $|\mathcal{L}|_D = |A|_D$ is free on any $A$. Apply this to $D = Z$. Then $\mathcal{L}(-Z)$ is nef, so that by the first part, $H^1(\mathcal{L}(-Z)) = 0$, and therefore

$$H^0(Y, \mathcal{L}) \twoheadrightarrow H^0(Z, \mathcal{L}_{|Z}).$$

Thus $|\mathcal{L}|$ on $Y$ cuts out the complete linear system $|\mathcal{L}|_Z$ on $Z$, and I just said that this is free. Q.E.D.

**Proof of Theorem 4.17, (1).** As discussed above, every function $g \in m_P \subset \mathcal{O}_{X,P}$ vanishes on $Z = Z_{\text{num}}$, so that $f_* \mathcal{O}_Y(-Z) = m_P$. By the lemma, $|\mathcal{O}_Y(-Z)| = |-Z|$ is free, which means that $m_P \mathcal{O}_Y = \mathcal{O}_Y(-Z)$ (compare the discussion in 4.6). This means $Z = Z_{\text{f}}$ and the assumption of Lemma 4.6, (i) holds, so $\operatorname{mult}(P \in X) = -Z^2$ follows.

## 4.18. Graded rings over divisors

Before proving (2), I make some general remarks on the graded ring

$$R(Z, \mathcal{L}) = \bigoplus_{k \geq 0} H^0(Z, \mathcal{L}^{\otimes k})$$

corresponding to a divisor $Z = \sum n_i \Gamma_i$ on a nonsingular surface and a nef line bundle $\mathcal{L} \in \operatorname{Pic} Z$ on $Z$ with $d = \deg \mathcal{L} > 0$.

I add a brief word of explanation in case you haven't seen this kind of thing before. There is a graded ring $R(V, \mathcal{L})$ corresponding to any line bundle over a variety (or scheme) $V$: by definition of tensor product of sheaves, there are local multiplication maps $\mathcal{L} \times \mathcal{L} \to \mathcal{L}^{\otimes 2}$ (more generally, $\mathcal{L}^{\otimes a} \times \mathcal{L}^{\otimes b} \to \mathcal{L}^{\otimes a+b}$). The graded ring $R(V, \mathcal{L})$ puts together all global sections of all powers of $\mathcal{L}$, with multiplication between the $H^0$ induced by the local maps of sheaves.

The case to bear in mind is when $\mathcal{L}$ is very ample, so that $V \subset \mathbb{P}^n$ with $H^0(\mathcal{L}) = \langle x_0, \dots, x_n \rangle$, the vector space of homogeneous coordinates of $\mathbb{P}^n$. Then the multiplication $H^0(\mathcal{L}) \times \cdots \times H^0(\mathcal{L}) \to H^0(\mathcal{L}^{\otimes k})$ just views monomials in the $x_i$ as restrictions to $V$ of homogeneous forms in $H^0(\mathbb{P}^n, \mathcal{O}(k))$. In this case, $R(V, \mathcal{L})$ generated in degree 1 means

$$H^0(\mathbb{P}^n, \mathcal{O}(k)) \twoheadrightarrow H^0(V, \mathcal{L}^{\otimes k}) \quad \text{for all } k.$$

If $V$ itself is a normal variety, this condition is called *projective normality*, and is equivalent to saying that the affine cone over $V$ is normal.

In the immediate application, $Z$ is the numerical cycle of a rational singularity or bunch of curves, and $\mathcal{L} = \mathcal{O}_Z(-Z)$. Thus $Z$ satisfies the assumptions of the next theorem.

**Theorem.** *Let $Z = \sum n_i\Gamma_i$ be a divisor on a nonsingular surface, and assume that $Z$ is numerically 1-connected and $H^1(\mathcal{O}_Z) = 0$. Let $\mathcal{L} \in \operatorname{Pic} Z$ be a nef line bundle on $Z$ with $d = \deg \mathcal{L} > 0$. Then the graded ring $R(Z, \mathcal{L}) = \bigoplus_{k\ge0} H^0(Z, \mathcal{L}^{\otimes k})$ is generated in degree 1. In particular, the multiplication map*

$$H^0(\mathcal{L}) \times H^0(\mathcal{L}) \to H^0(\mathcal{L}^{\otimes 2})$$

*is onto.*

**Proof.** The proof is completely elementary: apart from the language, it is just the simplest case of Castelnuovo's free pencil trick.

Write $\deg_{\Gamma_i} \mathcal{L} = d_i$ (so that $d = \sum n_i d_i$), and choose $d_i$ disjoint transverse cuts $A_{ij}$ through $\Gamma_i$ (for $j = 1, \dots, d_i$). As a scheme, $A_{ij} \cong \operatorname{Spec} k[x]/x^{n_i}$. Their sum $A = \sum A_{ij}$ is a Cartier divisor $A$ with the same degree as $\mathcal{L}$ on every component of $Z$, so that Proposition 4.14 gives $\mathcal{L} \cong \mathcal{O}_Z(A)$. In other words, $\mathcal{L}$ has a section $s_0$ with divisor of zeros equal to $A$. As in Lemma 4.17, the cuts making up $A$ can be chosen arbitrarily, so that $|\mathcal{L}| = |A|$ is a free linear system on $Z$.

Consider two general sections $s, s_0 \in H^0(\mathcal{L})$ (the eponymous free pencil), and write $A = \operatorname{div}(s_0)$; I assume that $s$ is a local basis of $\mathcal{L}$ at every point of $\operatorname{Supp} A$. Then $s_0$ defines the short exact sequence

$$0 \to \mathcal{O}_Z \to \mathcal{L} \to \mathcal{O}_A \to 0.$$

Here I use the second section $s$ to identify $\mathcal{L} = \mathcal{O}_Z(A)$ with $\mathcal{O}_Z$ near $A$, which simplifies the notation $\mathcal{L}_{|A} = \mathcal{O}_A(A)$ to $\mathcal{O}_A$.

Now the assumption $H^1(\mathcal{O}_Z) = 0$ obviously gives $H^1(\mathcal{L}) = 0$ and $H^0(\mathcal{L}) = d{+}1$. Let $s_1, \dots, s_d \in H^0(\mathcal{L})$ map to a basis of $H^0(\mathcal{O}_A)$. Then the theorem follows from the more precise claim: for any $k \ge 0$, the vector space $H^0(\mathcal{L}^{\otimes k})$ is based by the $kd + 1$ monomials

$$s_0^k, \quad \text{and} \quad s_0^{k-a} s^{a-1} s_i \quad \text{for } a = 1, \dots, k \text{ and } i = 1, \dots d.$$

Indeed, $H^0(\mathcal{O}_Z)$ is just the constants (because $Z$ is 1-connected), which gives the assertion for $k = 0$; for $k \ge 1$, it follows by induction from the exact sequence

$$H^0(\mathcal{L}^{\otimes k-1}) \hookrightarrow H^0(\mathcal{L}^{\otimes} k) \twoheadrightarrow H^0(\mathcal{O}_A). \quad \text{Q.E.D.}$$

The proof amounts to saying that $\mathcal{L}$ embeds the Artinian scheme $A$ into $\mathbb{P}^{d-1}$ so that $H^0(\mathbb{P}^{d-1}, \mathcal{O}(1)) \to \mathcal{O}_A$ is onto. That is, the image $A \subset \mathbb{P}^{d-1}$ behaves like $d$ points in linearly general position.

**Proof of Theorem 4.17, (2).** Under the identification $m_P = H^0(\mathcal{O}_Y(-Z))$, it is clear that $m_P^2$ maps to $H^0(\mathcal{O}_Y(-2Z))$, which gives a map

$$m_P/m_P^2 \to H^0(\mathcal{O}_Y(-Z))/H^0(\mathcal{O}_Y(-2Z)) = H^0(\mathcal{O}_Z(-Z)) \tag{4}$$

which is well defined and surjective. The point is to prove that it is also injective, so that $m_P/m_P^2 \cong H^0(\mathcal{O}_Z(-Z))$. This completes the proof, since by RR,

$$h^0(\mathcal{O}_Z(-Z)) = 1 - p_a Z + \deg(\mathcal{O}_Z(-Z)) = 1 - Z^2.$$

To see that (4) is injective, it is enough to prove that

$$H^0(Y, \mathcal{O}_Y(-2Z)) = m_P^2;$$

in words, a function vanishing along $2Z$ is a linear combination of products of two functions vanishing along $Z$. I prove this as a corollary of Theorem 4.18. Let $g \in H^0(Y, \mathcal{O}_Y(-2Z))$; its restriction to $Z$ is a section $\overline{g} \in H^0(\mathcal{O}_Z(-2Z))$, so by the theorem, it is in the image of $H^0(\mathcal{O}_Z(-Z)) \otimes H^0(\mathcal{O}_Z(-Z))$; choosing an expression $\overline{g} = \sum \overline{x}_i \overline{y}_i$ and lifting to $x_i, y_i \in H^0(\mathcal{O}_Y(-Z))$ gives

$$g - f_2 \in H^0(\mathcal{O}_Y(-3Z)) \quad \text{with } f_2 = \sum x_i y_i \in m_P^2$$

Proceeding in the same way with $H^0(\mathcal{O}_Z(-3Z))$, etc., gives for any $n$ an expression

$$g - f_2 - f_3 - \cdots - f_{n-1} \in H^0(\mathcal{O}_Y(-nZ)) \quad \text{with } f_i \in m_P^i. \tag{5}$$

On the other hand, I claim that $H^0(\mathcal{O}_Y(-nZ)) \subset m_P^2$ for some $n \gg 0$. Since also $m_P^i \subset m_P^2$, (5) implies that $g \in m_P^2$.

Note that the graded algebra $\bigoplus_{k \geq 0} H^0(\mathcal{O}_Y(-kZ))$ is finitely generated as a $k[X]$-algebra: for $\mathcal{O}_Y(-Z)$ is the pullback of the relatively ample line bundle $\mathcal{O}(1)$ under the morphism $\varphi_{|-Z|}\colon Y \to X \times \mathbb{P}^N$ defined by the free linear system $|-Z|$; (the image is the blowup of $m_P$). The required finite generation follows by projective normalisation. Each $H^0(\mathcal{O}_Y(-kZ))$ is contained in $m_P$, so that for $n$ greater than the degrees of all the generators, the group $H^0(\mathcal{O}_Y(-nZ))$ is contained in a sum of products of at least two of these, therefore in $m_P^2$. Q.E.D.

I'm sorry this proof was so involved. Lots of proofs in algebraic geometry depend on finiteness arguments of this type, and projective normalisation is among the more elementary, provided it's in your vocabulary. A posteriori, the graded algebra is generated in degree 1.

An alternative is to go from (5) to the conclusion $g \in m_P^2$ using some kind of completion; [**A2**] uses Henselisation. If you think of $Y$ as just a thin tubular neighbourhood of $f^{-1}P$ in the complex topology, any general divisor in $|-2Z|$ consists of a number of transverse cuts through general points of components of $f^{-1}P$, with each transverse cut a connected component of the divisor in the complex topology. By considering degrees, it is easy to see that this can be grouped as a sum of two divisor in $|-Z|$, which proves that $H^0(\mathcal{O}_Y(-2Z)) = m_P^2$.

## 4.19. The Gorenstein condition and the canonical cycle $Z_K$

The class of Gorenstein singularities $P \in X$ is a mild generalisation of hypersurface singularities; although I sketch a correct treatment of the definition and simple properties, you could for many practical purposes take the easy way out and assume throughout that $P \in X$ is a hypersurface singularity. Let $P \in X$ be a normal surface singularity; at some future point you need to know that $P \in X$ is automatically Cohen–Macaulay, as discussed in [**YPG**], App. to §3. Let $f\colon Y \to X$ be a resolution of singularities, subsequently usually assumed to be *minimal* (no $-1$-curves in the fibre).

**Proposition–Definition.** *Equivalent conditions:*

1. *The canonical class of $Y$ can be written as an exceptional divisor plus a divisor disjoint from the exceptional locus:*

$$K_Y \overset{\text{lin}}{\sim} -Z_K + M \quad \textit{with} \quad M \cap f^{-1}P = \emptyset.$$

2. *There exists a rational canonical differential $s \in \mathcal{O}_X(K_X)$ on $X$ which is regular and has no zeros on a punctured neighbourhood of $P \in X$.*
3. *The line bundle $\Omega^2_{X^0}$ on the nonsingular locus $X^0 = X \setminus P$ extends over $X$ as a line bundle $\omega_X$.*
4. *The canonical Weil divisor class $K_X$ is Cartier (locally principal) at $P$.*
5. *The Serre–Grothendieck dualising sheaf $\omega_X$ is invertible at $P$.*
6. *The local ring $\mathcal{O}_{X,P}$ is Gorenstein in the sense of commutative algebra textbooks (for example* [**Matsumura**]*, Chapter 6, §18).*

*If these conditions hold, $P \in X$ is* Gorenstein*; the exceptional divisor $Z_K$ is called the* canonical cycle *of the resolution (or the* antidiscrepancy*).*

**Discussion of proof.** The canonical divisor $K_X$ of a normal $n$-fold is the Weil divisor class $\operatorname{div} s$ where $s \in \Omega^n_{k(X)}$ is a rational canonical differential. This was discussed in A.10 under the additional assumption that $X$ is nonsingular, so that $s$ could be compared with a volume form $\mathrm{d}x_1 \wedge \cdots \wedge \mathrm{d}x_n$ at every point of $X$. However, a rational canonical differential $s \in \Omega^n_{k(X)}$ is a birational notion, and to define its divisor $\operatorname{div} s$ as a Weil divisor, we only need to know $X$ at each prime divisor.

It follows from this discussion that, essentially by definition, (1), (3) and (4) are merely restatements of (2). The (pre-) dualising sheaf $\omega_X$ is defined, and its fundamental dualising property proved in [**H1**], III.7. The equivalence of (4) and (5) comes from the fact that the (pre-) dualising sheaf $\omega_X$ equals the divisorial sheaf $\mathcal{O}_X(K_X)$, as discussed in [**C-3f**], App. to §1.

The point of the equivalence between (5) and (6) is to pass between global duality (properties involving $\omega_X$ and the category of quasicoherent sheaves on $X$) and local duality for the local ring $A = \mathcal{O}_{X,P}$. To do this categorically should be an exercise, but one which I can't do convincingly, partly because the literature on local duality leaves me speechless. A practical alternative uses the calculation of $\omega_X$ as an $\mathcal{E}xt$ sheaf: assume that $X \subset \mathbb{P}$ is an embedding into a nonsingular variety; then $X$ Gorenstein at $P$ means that $\mathcal{E}xt^{N-n}_{\mathcal{O}_\mathbb{P}}(\mathcal{O}_X, \omega_\mathbb{P}) = \omega_X$ is invertible at $P$ and $\mathcal{E}xt^i = 0$ for $i \neq N-n$. Localising and playing with injective resolutions, it's not hard to see that this is equivalent to $\operatorname{Ext}^n_A(k, A) = k$ and $\operatorname{Ext}^i = 0$ for $i \neq n$, where $k = A/m_P$. This is one of the many equivalent definitions of a Gorenstein local ring, compare [**Matsumura**], Theorem 18.1.

**Example.** If $P \in X \subset \mathbb{C}^3$ is an isolated hypersurface singularity defined by $f = 0$, a basis of $\mathcal{O}_X(K_X)$ near $P$ is given by

$$s = \operatorname{Res}_X \frac{\mathrm{d}x \wedge \mathrm{d}y \wedge \mathrm{d}z}{f} = \frac{\mathrm{d}x \wedge \mathrm{d}y}{\partial f / \partial z} = \text{etc.}$$

Then $s$ is a rational canonical differential, and its divisor on $Y$ is $\operatorname{div}_Y(s) = -Z_K$. Compare the introductory discussion in [**YPG**], §1.

## 4.20. The canonical and cohomological cycles, $Z_1 = Z_K$

From now on I assume that the resolution $Y$ is minimal, that is, $K_Y$ is nef on the exceptional curves $\{\Gamma_i\}$. Then $-Z_K$ is nef, so that either $Z_K = 0$ or $Z_K \geq Z_{\text{num}}$, where $Z_{\text{num}}$ is the numerical cycle of 4.2. In particular $Z_K \geq 0$, so I can view it as a subscheme $Z_K \subset Y$. It is easy to see that $Z_K = 0$ if and only if $Y \to X$ is a minimal resolution of a Du Val singularity.

**Theorem.** *If $P \in X$ is an isolated Gorenstein surface singularity then $Z_1 = Z_K$. Here $Z_1$ is the cohomological cycle of Theorem 4.8, and $Z_K$ the canonical cycle defined in 4.19.*

**Proof.** First, if $D \geq Z_K$ then $H^1(\mathcal{O}_D) \to H^1(\mathcal{O}_{Z_K})$ is an isomorphism: indeed, it is clearly surjective, with kernel a quotient of

$$H^1(\mathcal{O}_{D'}(-Z_K)) \overset{\text{d}}{=} H^0(\mathcal{O}_{D'}(D')), \quad \text{where } D' = D - Z_K.$$

But it is an easy exercise to see that $H^0(\mathcal{O}_{D'}(D')) = 0$ (see Ex. 4.14). Therefore $Z_K$ achieves the maximum value of $H^1(\mathcal{O}_D)$, and so $Z_1 \leq Z_K$ by definition of $Z_1$ (see Theorem 4.8, (b)).

Next, if $Z_1 < Z_K$, then I claim that $h^1(\mathcal{O}_{Z_1}) < h^1(\mathcal{O}_{Z_K})$, which contradicts the definition of $Z_1$. The duals of these $H^1$s are the groups

$$H^1(\mathcal{O}_{Z_K}) \overset{\text{d}}{=} H^0(\omega_{Z_K}) = H^0(Z_K, \mathcal{O}_{Z_K}(K_Y + Z_K)) = H^0(Z_K, \mathcal{O}_{Z_K})$$
$$H^1(\mathcal{O}_{Z_1}) \overset{\text{d}}{=} H^0(\omega_{Z_1}) = H^0(Z_1, \mathcal{O}_{Z_1}(K_Y + Z_1)) = H^0(Z_1, \mathcal{O}_{Z_1}(Z_1 - Z_K))$$

which fit in the exact sequence

$$0 \to H^0(Z_1, \mathcal{O}_{Z_1}(Z_1 - Z_K)) \to H^0(Z_K, \mathcal{O}_{Z_K}) \to H^0(\mathcal{O}_{D''}), \tag{6}$$

where $D'' = Z_K - Z_1$. But if $D'' \neq 0$, the last map in (6) is certainly nonzero, because it is the identity on the constant functions. This completes the contradiction, and proves the theorem. Q.E.D.

## 4.21. Elliptic Gorenstein surface singularities

**Definition.** A normal surface singularity $P \in X$ is an *elliptic Gorenstein singularity* if it is Gorenstein (Proposition–Definition 4.19) and elliptic (Definition 4.12).

The treatment of elliptic Gorenstein singularities in what follows parallels closely that of rational singularities in 4.12–18.

**Theorem** (Numerical characterisation). *Let $P \in X$ be an elliptic Gorenstein surface singularity, $f\colon Y \to X$ its minimal resolution, and $Z = Z_K$ the canonical cycle. Then $Z$ is the unique effective cycle with $p_a(Z) = 1$, and $p_a(D) \leq 0$ for $D \neq Z$. Moreover, $Z$ is numerically 2-connected and $Z = Z_{\text{num}} = Z_1 = Z_K$.*

*Conversely, suppose that the minimal resolution $Y \to X \ni P$ of a normal surface singularity has an effective cycle $Z$ satisfying $p_a(Z) = 1$, and $p_a(D) \leq 0$ for $D \neq Z$. Then $P \in X$ is elliptic Gorenstein.*

**Proof.** Note that by the adjunction formula 4.10, the canonical cycle $Z_K$ of any Gorenstein singularity automatically satisfies $\omega_Z = \mathcal{O}_Z(K_Y + Z_K) \cong \mathcal{O}_Z$, so that $\chi(\mathcal{O}_{Z_K}) = 0$ and $p_a(Z_K) = 1$.

Suppose first that $P \in X$ is elliptic Gorenstein, and set $Z = Z_1 = Z_K$. Then by Theorem 4.8, (b), any effective exceptional divisor $D$ satisfies $H^1(\mathcal{O}_D) = 0$ if $D \not\geq Z = Z_1$, and so $p_a D \leq 0$. If $D > Z = Z_K$ then

$$2p_a D - 2 = D(D - Z) = (D - Z)^2 + Z(D - Z) \leq (D - Z)^2 < 0, \tag{7}$$

($K_Y = -Z$ is nef on exceptional curves, because $Y$ is minimal), so again $p_a D \leq 0$.

**Claim.** *If an effective divisor $Z$ on a nonsingular surface satisfies $p_a(Z) = 1$ and $p_a(D) \leq 0$ for all $D < Z$ then $D(Z - D) \geq 2$, that is, $Z$ is numerically 2-connected.*

**Proof.** Try it and see. $\square$

This proves the "moreover" part: $-Z$ is nef (because $-Z \overset{\text{num}}{\sim} K_Y$), and for any divisor $D$ with $0 < D < Z$, the claim implies that $D(Z - D) \geq 2$, so $-D$ is not nef. Therefore $Z = Z_{\text{num}}$.

The point of the converse is to prove that $K_Y + Z \overset{\text{lin}}{\sim} 0$ near $f^{-1}P$, so that $P \in X$ is Gorenstein by Proposition 4.19, (1). I divide this up into easy steps.

**Step 1.** *$Z = Z_{\text{num}}$, and in particular $f^{-1}P = \operatorname{Supp} Z$.*

**Proof.** As in (7), applying the assumption $p_a \leq 0$ to $Z + \Gamma_i$ implies that $-Z\Gamma_i \geq 0$ for any exceptional $\Gamma_i$:

$$\begin{aligned} -2 \geq 2p_a(Z + \Gamma) - 2 &= (K_Y + Z + \Gamma)(Z + \Gamma) \\ &= (K_Y + Z)Z + (K_Y + \Gamma)\Gamma + 2Z\Gamma \\ &\geq -2 + 2Z\Gamma. \end{aligned} \tag{8}$$

As in the claim, if $0 < D < Z$ then $-D$ is not nef, so that $Z = Z_{\text{num}}$.

**Step 2.** *$K_Y + Z \overset{\text{num}}{\sim} 0$ on exceptional curves.*

**Proof.** Since $(K_Y + Z)Z = 0$, there is nothing more to prove if $Z$ is irreducible. If $Z = \sum n_i\Gamma_i$ then $\sum n_i(K_Y + Z)\Gamma_i = 0$, and the assumption $p_a(\Gamma_i) = 0$ together with the recent claim gives

$$(K_Y + Z)\Gamma_i = (K_Y + \Gamma_i)\Gamma_i + (Z - \Gamma_i)\Gamma_i = -2 + (Z - \Gamma_i)\Gamma_i \geq 0.$$

Thus $(K_Y + Z)\Gamma_i = 0$ for each $i$, as required.

**Step 3.** *$\mathcal{O}_Z(K_Y + Z) \cong \mathcal{O}_Z$.*

**Proof.** We know that $\mathcal{O}_Z(K_Y + Z)$ is a numerically trivial line bundle on $Z$, and that $Z$ is numerically connected, so that it is enough to observe that it has a section:

$$H^0(\mathcal{O}_Z(K_Y + Z)) \overset{\text{d}}{=} H^1(\mathcal{O}_Z) \neq 0.$$

**Step 4.** *$H^0(\mathcal{O}_Y(K_Y + Z)) \to H^0(\mathcal{O}_Z(K_Y + Z))$ is surjective.*

**Proof.** As in the proof of Theorem 4.8, the cokernel $H^1(\mathcal{O}_Y(-Z))$ can be calculated as an inverse limit $\varprojlim H^1(\mathcal{O}_D(-Z))$, and each of these groups is zero (see Ex. 4.15).

It follows from Steps 3–4 that $\mathcal{O}_Y(K_Y + Z)$ has a section which is a basis at every point of $f^{-1}P$. Therefore $K_Y + Z \overset{\text{lin}}{\sim} 0$ in a neighbourhood of $f^{-1}P$. This completes the proof that $P \in X$ is Gorenstein with canonical cycle $Z_K = Z$. Clearly by Theorems 4.8 and 4.20, $R^1 f_* \mathcal{O}_Y = H^1(\mathcal{O}_Z)$ is 1-dimensional, so that $P \in X$ is elliptic Gorenstein. Q.E.D.

## 4.22. Examples

Just as the affine cone over a rational normal curve of degree $d$ serves as a model example of a rational singularity of multiplicity $d$ (at least for mnemonic purposes), the affine cone over a normally embedded elliptic curve $E_d \subset \mathbb{P}^{d-1}$ of degree $d \geq 3$ is the model example of an elliptic Gorenstein singularity; in singularity theory these cones are called *simple elliptic singularities*. Note that the fact that the elliptic curve $E_d$ has a modulus (the $j$-invariant) means that the analytic type of the singularity also depends on a modulus. The typical example is the ordinary hypersurface triple point

$$0 \in X : (x^3 + y^3 + z^3 + \lambda xyz = 0) \subset \mathbb{C}^3$$

These admit many degenerations (for example, $xyz + x^p + y^q + z^r$ for $p, q, r \geq 4$, known to singularity theorists as $T_{p,q,r}$). This class of singularities includes case that are very complicated from the point of view of classification of singularities.

Elliptic Gorenstein singularities with $d = 1$ or 2 are hypersurface double points, again typified by the simple elliptic singularities

$$x^2 + y^3 + z^6 + \lambda xyz \quad \text{and} \quad x^2 + y^4 + z^4 + \lambda xyz,$$

which this time are cones in a weighted sense: they are the respective affine cones on the curve in weighted projective planes $C_6 \subset \mathbb{P}(3, 2, 1)$ and $C_4 \subset \mathbb{P}(2, 1, 1)$. The statements of the results are slightly different, but the theory is broadly similar to that for $d \geq 3$. They are alternative explicit methods to treat them. Compare Ex. 4.18.

Both [**L2**] and [**R**] contain systematic lists of examples of elliptic Gorenstein singularities with $d = 1, 2, 3$.

## 4.23. The main theorem

Set $d = -Z^2$, defined to be the *degree* of $P \in X$.

**Main Theorem ([R], [L2], [Mori]).**

*(I) If $d \geq 2$ then $m_P \mathcal{O}_Y = \mathcal{O}_Y(-Z)$; in particular $Z = Z_{\mathrm{f}}$ and $|-Z|$ is a free linear system, so that* $\operatorname{mult}_P X = d$.

*(II) If $d \geq 3$, then $f_* \mathcal{O}_Y(-2Z) = m_P^2$, and* $\operatorname{emb\,dim}_P X = d$.

*Moreover, $f_* \mathcal{O}_Y(-kZ) = m_P^k$ for all $k$; that is, the $\mathcal{O}_X$-algebra*

$$\bigoplus_{k \geq 0} f_* \mathcal{O}_Y(-kZ) = \bigoplus_{k \geq 0} m_P^k \tag{9}$$

*is generated by its component $m_P$ in degree* 1.

The results for $d = 1$ or 2 are slightly longer to state, and are deferred to 4.25 for the sake of clarity; however, the proofs are closely related, and are given together when appropriate.

**Discussion.** (II) says that the ordinary blow $X_1 \to X$ of $X$ at $P$ coincides with the relative canonical model of $Y \to X$. Indeed, the blowup is defined to be Proj of the algebra $\bigoplus m_P^k$ on the right-hand side of (9). On the other hand, by definition of the canonical cycle, $\mathcal{O}_Y(-Z) \cong \mathcal{O}_Y(K_Y)$ in a neighbourhood of $f^{-1}P$, so that the left-hand side of (9) is the relative canonical algebra $\mathcal{R}(Y, K_Y)$ of $Y \to X$. It is easy to see that $\operatorname{Proj}\mathcal{R}(Y, K_Y)$ is the normal surface $\overline{Y}$ with only Du Val singularities and ample $K_{\overline{Y}}$ obtained from $Y$ by contracting exactly the $-2$-curves of $f^{-1}P$.

**Proof of (I).** Most of the proof takes place on $Z$. First note that the restriction $H^0(Y, \mathcal{O}_Y(-Z)) \to H^0(Z, \mathcal{O}_Z(-Z))$ is surjective, by the usual argument: the cokernel is $H^1(Y, \mathcal{O}_Y(-2Z))$, which is zero by Ex. 4.15. For $d = -Z^2 \geq 2$, the following result applied to $\mathcal{L} = \mathcal{O}_Z(-Z)$ gives at once that $|-Z|$ is free on $Y$, which proves (I).

**Lemma.** *Let $\mathcal{L}$ be a line bundle on $Z$ which is nef and of degree $d > 0$. Assume if possible that $Q \in Z$ is a closed point at which $|\mathcal{L}|$ is not free. Then $d = 1$, $Q \in Z$ is a nonsingular point (therefore a Cartier divisor), and $\mathcal{L} \cong \mathcal{O}_Z(Q)$. Then $H^0(\mathcal{L})$ consists of the identity inclusion $s\colon \mathcal{O}_Z \hookrightarrow \mathcal{O}_Z(Q)$ and its scalar multiples.*

*If $\deg \mathcal{L} = 1$ then $h^0(\mathcal{L}) = 1$, so that $|\mathcal{L}|$ is not free, and all these assertions hold for some point $Q \in Z$.*

**Proof.** The assumption is that the evaluation map $H^0(Z, \mathcal{L}) \to k_Q$ is not surjective, so its cokernel $H^1(Z, m_Q\mathcal{L})$ is nonzero. I sidestep the difficulty that this $H^1$ might come from restriction to a smaller divisor by taking $D$ to be a divisor with $0 < D \leq Z$ which is minimal with $H^1(D, m_Q\mathcal{L}_{|D}) \neq 0$. Minimal means:

$$H^1(D - \Gamma_i, m_Q\mathcal{L}_{|D - \Gamma_i}) = 0 \quad \text{for every component } \Gamma_i < D,$$

so that the nonzero $H^1$ comes from the sheaf kernel of the restriction map

$$m_Q\mathcal{L} \to m_Q\mathcal{L}_{|D - \Gamma_i}.$$

Outside $Q$, the kernel is $\mathcal{L} \otimes \mathcal{O}_{\Gamma_i}(-D + \Gamma_i)$ by the dévissage exact sequence of 3.10. To understand what happens at $Q$, note there are two cases: either $Q \in D - \Gamma_i$, and the kernel is $\mathcal{O}_{\Gamma_i}(\mathcal{L} - D + \Gamma_i)$, since everything vanishing on $D - \Gamma_i$ is automatically in $m_Q$; or $Q \notin D - \Gamma_i$ (which means that $\Gamma_i$ is a reduced component of $D$, and is the unique component of $D$ containing $Q$, so that the restriction is an isomorphism near $Q$), and the kernel is $m_Q\mathcal{O}_\Gamma(\mathcal{L} - D + \Gamma_i)$.

Therefore

$$\begin{aligned} &H^1(\Gamma_i, \mathcal{O}_{\Gamma_i}(\mathcal{L} - D + \Gamma_i)) \neq 0 \quad \text{if } Q \in D - \Gamma_i; \\ \text{and} \quad &H^1(\Gamma_i, m_Q\mathcal{O}_{\Gamma_i}(\mathcal{L} - D + \Gamma_i)) \neq 0 \quad \text{if } Q \notin D - \Gamma_i. \end{aligned}$$

Now quite generally, for a torsion free sheaf $\mathcal{F}$ of rank 1 on an irreducible curve $\Gamma$,

$$H^1(\Gamma, \mathcal{F}) \neq 0 \implies \deg_\Gamma \mathcal{F} \leq 2p_a\Gamma - 2. \tag{10}$$

Write $D = \sum a_i\Gamma_i$. Since $K_Y = -Z$ locally, $2p_a\Gamma_i - 2 = (-Z + \Gamma_i)\Gamma_i$. Summing (10) over all components $\Gamma_i$ of $D$, and noticing that replacing $\mathcal{L}$ by $m_Q\mathcal{L}$ decreases the degree by 1 for at most one component $\Gamma_i$, I get

$$\sum a_i\Gamma_i(\mathcal{L} - D + \Gamma_i) \le 1 + \sum a_i\Gamma_i(-Z + \Gamma_i).$$

The quadratic terms $\Gamma_i^2$ on each side cancel, giving

$$D\mathcal{L} + D(Z - D) \le 1. \tag{11}$$

But now $Z$ is 2-connected, so that (11) implies that $D = Z$, $\deg \mathcal{L} = 1$. Thus $h^0(\mathcal{L}) = 1$. Moreover, $Q \in \Gamma_i$ holds for a unique irreducible component $\Gamma_i$, and the section $\mathcal{O}_Z \to \mathcal{L}$ vanishes at $Q$ only, so that $Q$ supports an effective Cartier divisor of degree 1, which means it is a nonsingular point of $Z$. This proves the lemma and completes the proof of (I). Q.E.D.

**Remark.** A minor modification of the same argument proves the following result (see Ex. 4.22): if $Z$ is a 2-connected divisor on a nonsingular surface and $\mathcal{L}$ a line bundle with $\mathcal{L} \otimes \omega_Z^{-1}$ nef and of degree $\ge 1$ then $|\mathcal{L}|$ is free unless $\mathcal{L} \cong \omega_Z(Q)$ with $Q \in Z$ a nonsingular point. Catanese and others [**CFHR**] contains a dual version of the argument, using Serre duality more systematically: $0 \ne H^1(m_Q\mathcal{L}) \overset{\mathrm{d}}{=} \operatorname{Hom}(m_Q\mathcal{L}, \omega_Z)$, so choose a nonzero element $s \in \operatorname{Hom}(m_Q\mathcal{L}, \omega_Z)$, and let $Z = A + B$ be the decomposition such that $A$ is the maximal divisor on which $s$ vanishes. Then $s$ defines an inclusion $m_Q\mathcal{L}_{|B} \hookrightarrow \omega_Z(-A) = \omega_B$. A simple numerical argument based on 2-connectedness and $\mathcal{L} > K_Z$ proves that $A = 0$, and $s$ defines an isomorphism $m_Q\mathcal{L} \cong \omega_Z$, so that $Q \in Z$ is a nonsingular point and $\mathcal{L} \cong \omega_Z(Q)$.

## 4.24. Proof of (II)

The main object of study of this section, as in 4.18, is the graded ring

$$R(Z, \mathcal{L}) = \bigoplus_{k \ge 0} H^0(Z, \mathcal{L}^{\otimes k}) \tag{12}$$

corresponding to a nef line bundle of degree $d > 0$ over a divisor $Z$. For $d \ge 3$, the result is as follows.

**Theorem.** *Let $Z$ be a 2-connected divisor on a nonsingular surface $Y$ for which $\mathcal{O}_Z(K_Y + Z) \cong \mathcal{O}_Z$, and let $\mathcal{L}$ be a nef line bundle on $Z$ with $d = \deg \mathcal{L} \ge 3$. Then the graded ring $R(Z, \mathcal{L})$ is generated by elements of degree $k = 1$. In particular,*

$$H^0(\mathcal{L}) \otimes H^0(\mathcal{L}) \twoheadrightarrow H^0(\mathcal{L}^{\otimes 2})$$

*is surjective.*

**Proof.** Because $\mathcal{L}$ is free, a general section $s_1 \in H^0(\mathcal{L})$ does not vanish along any component of $Z$, so is regular at every point (that is, locally defined by a nonzerodivisor); thus it defines an effective Cartier divisor $A$ of degree $d$ on $Z$, such that $\mathcal{L} \cong \mathcal{O}_Z(A)$. As in 4.18, I choose a second general section $s$ and use it to identify $\mathcal{L}$ with $\mathcal{O}_Z$ in a neighbourhood of every point of $A$, and consider the exact sequence

$$0 \to \mathcal{O}_Z \xrightarrow{s_1} \mathcal{L} \to \mathcal{O}_A \to 0. \tag{13}$$

Now $A$ is an Artinian scheme of degree $d$, and since $h^0(\mathcal{O}_Z) = h^1(\mathcal{O}_Z) = 1$ it follows that $H^0(\mathcal{L}) \to \mathcal{O}_A$ has rank $d-1$. (I write $\mathcal{O}_A = H^0(\mathcal{O}_A)$ out of laziness.)

**Claim 1.** *$H^0(\mathcal{L}) \twoheadrightarrow \mathcal{O}_{A'}$ is surjective for every subscheme $A' \subset A$ of length $d-1$. Geometrically, this means that $\varphi_{\mathcal{L}}$ embeds $A$ in $\mathbb{P}^{d-2}$ as a subscheme of degree $d$ with the property that any subscheme of degree $d-1$ spans $\mathbb{P}^{d-2}$, like $d$ points forming a projective frame of reference.*

To prove this, write $\mathcal{I}_{A'} \subset \mathcal{O}_Z$ for the ideal sheaf of $A'$, so that $\mathcal{I}_A \subset \mathcal{I}_{A'}$ and $\mathcal{I}_{A'}/\mathcal{I}_A \cong k_Q$ for some point $Q \in Z$. Then $\mathcal{F} = \mathcal{I}_{A'}\mathcal{O}_Z(A)$ is an intermediate torsion free sheaf

$$\mathcal{O}_Z \subset \mathcal{F} \subset \mathcal{O}_Z(A),$$

with $\mathcal{F}/\mathcal{O}_Z \cong k_Q$. Now by Serre duality $H^1(\mathcal{F}) \stackrel{\mathrm{d}}{=} \mathrm{Hom}(\mathcal{F}, \mathcal{O}_Z)$. But this group is zero; for a nonzero homomorphism $\mathcal{F} \to \mathcal{O}_Z$ must restrict to a nonzero map $\mathcal{O}_Z \to \mathcal{O}_Z$, necessarily an isomorphism by 2-connectedness, and this would split the exact sequence $\mathcal{O}_Z \to \mathcal{F} \to k_Q$, contradicting $\mathcal{F}$ torsion free. Therefore $H^1(\mathcal{I}_{A'}\mathcal{L}) = 0$, which gives the required surjectivity.

**Claim 2.** *A subscheme $A_d \subset \mathbb{P}^{d-2}$ with the property stated in Claim 1 imposes linearly independent conditions on quadrics. (This means that the restriction map $H^0(\mathbb{P}^{d-2}, \mathcal{O}(2)) \to \mathcal{O}_{A_d}$ is surjective.)*

The cases $d = 3$ and $d = 4$ of Claim 2 are very easy, see Ex. 4.23. Although the proof involves various bits of mess, the proof in substance is similar to, and rather easier than [**UAG**], §1.

Indeed, any Artinian scheme $A$ has a filtration $A_1 \subset \cdots \subset A_{d-1} \subset A$ by subschemes $A_i$ of degree $i$, and by assumption, each $A_i$ has scheme theoretic linear span $\langle A_i \rangle = \mathbb{P}^{i-1}$. Clearly each $A_i$ imposes $i$ conditions on quadrics (since it does on hyperplanes), and I prove inductively that $A_i$ is a scheme theoretic intersection of quadrics. It follows from this that $A_{d-1}$ imposes $d-1$ conditions on quadrics, and that there is a quadric through $A_{d-1}$ not through $A_d$, so that $A_d$ imposes one further condition, which proves Claim 2.

The inductive step is as follows: write

$$\langle A_i \rangle = \mathbb{P}^{i-1} : (x_i = 0) \subset \langle A_{i+1} \rangle = \mathbb{P}^i$$

and $\mathcal{I}_i$ for the ideal sheaf of $A_i$. Suppose that $A_i$ and $A_{i+1}$ differ only at a point $Q$. Then it is easy to check that there is an exact sequence

$$0 \to m_Q\mathcal{O}_{\mathbb{P}^i}(1) \xrightarrow{x_i} \mathcal{I}_{i+1}\mathcal{O}_{\mathbb{P}^i}(2) \to \mathcal{I}_i\mathcal{O}_{\mathbb{P}^{i-1}}(2) \to 0. \tag{14}$$

(See Ex. 4.24.) Clearly, to say that $A_i$ is a scheme theoretic intersection of quadrics means exactly that $\mathcal{I}_i\mathcal{O}_{\mathbb{P}^{i-1}}(2)$ is generated by its $H^0$. Now the right-hand end of (14) is generated by its $H^0$ by the inductive hypothesis, and the left-hand end obviously so. Therefore so is the middle term. This proves Claim 2.

The proof Theorem 4.24 and Theorem 4.23, (II) is now straightforward. Let $s_1, \dots, s_d$ be a basis of $H^0(Z, \mathcal{L})$, with $s_1$ defining an effective Cartier divisor $A$ as in (13). Since $s_2, \dots, s_d$ are the homogeneous coordinates in the $\mathbb{P}^{d-2}$ containing $A_d$, by Claim 2, I can choose a set $\{t_1, \dots, t_d\}$ of quadratic monomials in $s_2, \dots, s_d$ that, viewed as sections of $H^0(Z, \mathcal{O}_Z(2A))$, map to a basis of $\mathcal{O}_A$. For $k \geq 2$, consider the exact sequence

$$0 \to H^0(Z, \mathcal{O}_Z((k-1)A) \xrightarrow{s_1} H^0(Z, \mathcal{O}_Z(kA)) \to \mathcal{O}_A \to 0.$$

Recall that I chose a second section $s \in H^0(Z, \mathcal{L})$ that is a local basis of $\mathcal{L}$ at every point of $\operatorname{Supp} A$, and which I use to identify $\mathcal{L}_{|A}$ with $\mathcal{O}_A$. Then the products $s^{k-2} \cdot \{t_1, \dots, t_d\}$ form a complimentary basis to $s_1 \cdot H^0(\mathcal{O}_Z((k-1)A)$. This proves that the ring $R(Z, \mathcal{L})$ is generated in degree 1.

The final part of the proof of (II) is the same as in 4.18. Since

$$0 \to f_*\mathcal{O}_Y(-(k+1)Z) \to f_*\mathcal{O}_Y(-kZ) \to H^0(Z, \mathcal{O}_Z(-kZ)) \to 0$$

is exact, any $g \in f_*\mathcal{O}_Y(-kZ)$ can be written modulo $f_*\mathcal{O}_Y(-(k+1)Z)$ as a polynomial of degree $k$ in elements of $m_P$. By induction,

$$f_*\mathcal{O}_Y(-kZ) \subset m_P^k + f_*\mathcal{O}_Y(-lZ) \quad \text{for an } l > k.$$

But since the $\mathcal{O}_X$-algebra $\mathcal{R}(Y, \mathcal{O}_Y(-Z))$ is finitely generated, the argument at the end of 4.18 proves that $f_*\mathcal{O}_Y(-lZ) \subset m^k$ for any fixed $k$ and $l \gg 0$. Therefore $f_*\mathcal{O}_Y(-kZ) = m_P^k$, as requested. Q.E.D.

## 4.25. The special Cases $d = 2, 1$

The results for $d = \deg \mathcal{L} = 2$ or 1 are similar in spirit, but the statements must involve weighted homogeneous rings because the ring $R(Z, \mathcal{L})$ is not generated in degree 1.

**Theorem.** *For* $d = 2:$

$$R(Z, \mathcal{L}) = k[x_1, x_2, y]/f, \quad \textit{where } \deg x_1, x_2, y = 1, 1, 2,$$

*and $f$ is homogeneous of degree 4, of the form $f = y^2 +$ terms involving $x_1, x_2$.*

*For* $d = 1:$

$$R(Z, \mathcal{L}) = k[x, y, z]/f, \quad \textit{where } \deg x, y, z = 1, 2, 3,$$

*and $f$ is homogeneous of degree 6, of the form $f = z^2 + y^3 +$ terms involving $x$.*

**Proof.** In either case, there exists a section $x \in H^0(Z, \mathcal{L})$ defining a Cartier divisor $A$ of degree $d$, and therefore an exact sequence

$$0 \to \mathcal{O}_Z \to \mathcal{L} \to \mathcal{O}_A \to 0.$$

When $d = 2$, I set $x = x_1$, and I can find a second section $x_2$ which is a local basis of $\mathcal{L}$ near $\operatorname{Supp} A$. Now choose an element $y \in H^0(Z, \mathcal{L}^{\otimes 2})$ so that $y, x_2^2$ map to a basis of $\mathcal{O}_A$. Arguing as usual on the restriction exact sequence

$$0 \to H^0(Z, \mathcal{O}_Z((k-1)A)) \to H^0(Z, \mathcal{O}(kA)) \to \mathcal{O}_A \to 0,$$

we see that $x_2^k, x_2^{k-2}y$ forms a complimentary basis to $x_1 \cdot H^0(Z, \mathcal{O}_Z((k-1)A))$ in every degree $k \geq 2$. It's very easy to see from this that

$$S^k(x_1, x_2), S^{k-2}(x_1, x_2)y \quad \text{bases} \quad H^0(Z, \mathcal{L}^{\otimes k}).$$

In particular, $x_1, x_2, y$ generate the ring, and there is a unique relation expressing $y^2$ in terms of the stated basis.

The case $d = 1$ is similar, and is left as an easy exercise. [Hint: you just have to pick elements $y, z$ in degrees $2, 3$ that map to a basis of $\mathcal{O}_A = k_Q$.] Q.E.D.

**Corollary.** *If* $d = 2$ *then* $\operatorname{mult}_P X = 2$ *and* $\operatorname{emb\,dim}_P X = 3$*, that is,* $P \in X$ *is a hypersurface double point. The linear system* $|-Z|$ *on* $Y$ *is free; that is,* $Z = Z_{\mathrm{f}}$*, and* $m_P\mathcal{O}_Y = \mathcal{O}_Y(-Z)$*. Moreover,* $f_*\mathcal{O}_Y(-2Z) = m_P^2 + (y)$ *for some element* $y \in m_P \setminus m_P^2$*, and the* $\mathcal{O}_X$*-algebra* $\mathcal{R}(Y, K_Y) = \bigoplus_{k\geq 0} f_*\mathcal{O}_Y(-kZ)$ *is generated by* $m_P$ *in degree* 1 *and* $y$ *in degree* 2 *(together with* $m_P^2$*).*

*If* $d = 1$ *then* $\operatorname{mult}_P X = 2$ *and* $\operatorname{emb\,dim}_P X = 3$*, that is,* $P \in X$ *is a hypersurface double point. The linear system* $|-Z|$ *on* $Y$ *has a single base point* $Q \in Z$*, which is a nonsingular point of* $Z$*; that is,* $Z = Z_{\mathrm{f}}$*, and* $m_P\mathcal{O}_Y = m_Q\mathcal{O}_Y(-Z)$*. The algebra* $\mathcal{R}(Y, K_Y) = \bigoplus_{k\geq 0} f_*\mathcal{O}_Y(-kZ)$ *is generated by* $m_P$ *in degree* 1*, and elements* $y \in m_P \setminus m_P^2$ *in degree* 2 *and* $z \in m_P \setminus (m_P^2, y)$ *in degree* 3*.*

**Remark.** In either case, the system of ideals

$$I_k = f_*\mathcal{O}_Y(-kZ) \subset \mathcal{O}_X \quad \text{for } k \geq 0$$

defines a filtration of $\mathcal{O}_X$, and the relative canonical model $Y \to \overline{Y} \to X$, defined as $\overline{Y} = \operatorname{Proj} \mathcal{R}(Y, K_Y)$, is obtained by blowing up $X$ along this filtration. In either case, the ideals $I_k$ can be defined by a weighting of the ambient space $\mathbb{A}^3$. That is, give

$$y, x_1, x_2 \text{ weights } (2, 1, 1) \text{ in case } d = 2,$$
$$z, y, x \text{ weights } (3, 2, 1) \text{ in case } d = 1.$$

then $I_k$ is the ideal generated by monomials of weighted degree $\geq k$. Taking $\operatorname{Proj} \bigoplus_{ge0} I_k$ can also be described in toric geometry as the weighted blowup.

When $d = 2$, the relative canonical model can also be described as the blowup of $X$ at $P$ followed by normalisation.

## Exercises to Chapter 4

**1.** Let $\mathbb{Z}/2$ act on $\mathbb{A}^2$ by $(u, v) \mapsto (-u, -v)$. Prove that the subring of $k[u, v]$ of invariant polynomials is generated by $x = u^2, y = uv, z = v^2$, and that the morphism $\mathbb{A}^2 \to \mathbb{A}^3$ given by $(u, v) \mapsto (x, y, z)$ defines an embedding of the orbit space $\mathbb{A}^2/(\mathbb{Z}/2)$ as the ordinary quadratic cone $X : (xz = y^2) \subset \mathbb{A}^3$ of 4.1.

Calculate the canonical class $K_V$ and its relation with $K_{\mathbb{A}^2}$ and $K_Y$.

**2.** Let $V \to \mathbb{A}^2$ be the blowup of the origin; show that the action of $\mathbb{Z}/2$ of Ex. 4.1 extends to an action on $V$ which fixes the $-1$-curve pointwise, and that the quotient space $V/(\mathbb{Z}/2)$ is isomorphic to the blowup $Y \to X$, and fits into a commutative diagram

$$\begin{array}{ccc} V & \to & \mathbb{A}^2 \\ \downarrow & & \downarrow \\ Y & \to & X. \end{array}$$

**3.** Carry out blowups to resolve the Du Val singularities

$$\begin{aligned} A_n : &\quad x^2 + y^2 + z^{n+1} = 0, \\ D_n : &\quad x^2 + y^2 z + z^{n-1} = 0, \\ E_6 : &\quad x^2 + y^3 + z^4 = 0, \\ E_7 : &\quad x^2 + y^3 + yz^3 = 0, \\ E_8 : &\quad x^2 + y^3 + z^5 = 0 \end{aligned}$$

in the spirit of 4.2. [Hint: In 4.2, I only calculated one affine piece of the blowup and stated that nothing interesting happens on the others. In general, one has to consider more than one affine piece.] Check that the exceptional locus is a bunch of $-2$-curves with intersections given by the Dynkin diagram.

**4.** Calculate the numerical cycle $Z_{\text{num}}$ for each of the Du Val singularities. Check that $Z_{\text{num}}^2 = -2$. [Hint: Start from $Z_0 = \sum \Gamma_i$, and successively increase the coefficients of $\Gamma_i$ only if $Z_j \Gamma_i > 0$.]

**5.** If $\bigcup \Gamma_i$ is a connected bunch of $-2$-curves with negative definite intersection matrix $\{\Gamma_i \Gamma_j\}$, check that the configuration is given by one of the Dynkin diagrams $A_n, D_n, E_6, E_7, E_8$. [Hint: If the configuration is more complicated, you can write down a combination $D = \sum a_i \Gamma_i$ with $D^2 \geq 0$. All the completed Dynkin diagrams $\widetilde{A}_n, \dots, \widetilde{E}_8$ appear as logical ends of the argument. For example, on the graph

$\widetilde{D}_4$ :

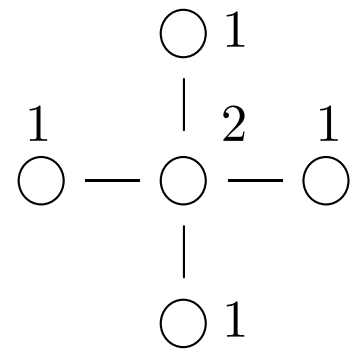

$D = \sum a_i \Gamma_i$ (with $a_i$ as indicated) has $D^2 = 0$, which cannot happen in a negative definite graph $G$, so that $G$ must have valency $\leq 3$. You have to prove that $G$ does not have a loop, or 2 vertexes of valency 3, or 1 vertex of valency 3 and very long arms out of it, etc.]

**6.** Let $\mathbb{Z}/3$ act on $\mathbb{A}^2$ by $(u, v) \mapsto (\varepsilon u, \varepsilon^2 v)$, where $\varepsilon$ is a primitive 3rd root of 1. Study the quotient

$$\mathbb{A}^2/(\mathbb{Z}/3) = X : (xz = y^3) \subset \mathbb{A}^3$$

as in Ex. 4.1 in terms of the quotient map $\mathbb{A}^2 \to \mathbb{A}^3$ defined by

$$x = u^3, y = uv, z = v^3.$$

Let $Y \to X$ be the minimal resolution. Show how to construct a $\mathbb{Z}/3$ cover $V \to Y$ by first blowing up the origin in $\mathbb{A}^2$, then blowing up the two points in which the $-1$-curve meets the coordinate axes.

The point of the question is that the minimal resolution $Y$ does not have a finite $\mathbb{Z}/3$ cover by a nonsingular blowup of $\mathbb{A}^2$.

**7.** Do the same as in Ex. 4.6 for $\mathbb{Z}/4$ acting by $(u, v) \mapsto (iu, i^3 v)$, where $i^2 = -1$.

**8.** Let $\varepsilon$ be a primitive 5th root of 1. Make $\mathbb{Z}/5$ act on $\mathbb{A}^2$ by

$$(u, v) \mapsto (\varepsilon u, \varepsilon^2 v).$$

Show that the subring of $k[u, v]$ of invariant polynomials is generated by 4 monomials, so that the quotient space $\mathbb{A}^2/(\mathbb{Z}/5) \cong X \subset \mathbb{A}^4$. (This is the quotient singularity of type $(5, 2)$ in Hirzebruch's terminology, or of type $\frac{1}{5}(1, 2)$ in that of [**YPG**].) Show that the ideal of relations holding between the invariant monomials is generated by 3 relations that can be written as the maximal minors of a $2 \times 3$ matrix. Show how to resolve the singularity by an explicit blowup, and describe the exceptional curves.

**9.** Let $\varepsilon$ be a primitive $2n$th root of 1 and consider the two matrixes

$$A = \begin{pmatrix} \varepsilon & 0 \\ 0 & \varepsilon^{-1} \end{pmatrix} \qquad \text{and} \qquad B = \begin{pmatrix} 0 & -1 \\ 1 & 0 \end{pmatrix}$$

Show that they generate a binary dihedral subgroup $BD_n \subset \mathrm{SL}(2,\mathbb{C})$ of order $4n$. Find the ring of invariants of $BD_n$ acting on $\mathbb{C}^2$ and prove that $\mathbb{C}^2/BD_n$ is isomorphic to the Du Val singularity $D_{n+2}$. [Hint: Find first the ring of invariants of $\langle A \rangle \cong \mathbb{Z}/2n$, and show that $B$ acts by an involution on it.]

**10.** In Ex. 4.9 you saw that the inclusion $\mathbb{Z}/2n \lhd BD_n$ defines a double cover $A_{2n-1} \to D_{n+2}$ between Du Val singularities. Show also that $BD_n \subset BD_{2n}$ is a normal subgroup of index 2 and that it defines a double cover $D_{n+2} \to D_{2n+2}$ between Du Val singularities.

**11.** Prove Proposition 4.12, (2). [Hint: Several methods are possible. For example, you can prove that $h^0(\mathcal{O}_D(-Z_{\text{num}})) \neq 0$, so that $h^0(\mathcal{O}_{D+Z}) \geq 2$ for every $D > 0$. Or you can calculate $\chi(\mathcal{O}_{D+Z})$ using the numerical games of (1).]

**12.** Prove Proposition 4.12, (3). [Hint: If $Z_{\text{num}} = D_1 + D_2$, write out $p_a(Z_{\text{num}}) = 0$, $p_a(D_1), p_a(D_2) \leq 0$ in terms of the adjunction formula 3.6 or 4.10. Compare the method of 4.21, (7–8)]

**13.** Consider the codimension 2 singularity $P \in X \subset \mathbb{C}^4$ defined by the two equations $x_1^2 = y_1^3 - y_2^3$ and $x_2^2 = y_1^3 + y_2^3$. Prove that $P \in X$ has a resolution $f \colon Y \to X$ such that $f^{-1}P$ is a nonsingular curve $C$ of genus 2 and $C^2 = -1$, but with $H^0(\mathcal{O}_C(-C)) = 0$. Prove that $x_1, x_2$ vanish along $C$ with multiplicity 2 and $y_1, y_2$ with multiplicity 3. Hence $Z_{\mathrm{f}} = 2C$, although obviously $Z_{\text{num}} = C$.

**14.** If $(\Gamma_i \Gamma_j)$ is negative definite, prove that $H^0(\mathcal{O}_D(D)) = 0$ for every effective divisor supported on $\Gamma_i$. More generally, if $\mathcal{L}$ is a line bundle on $D$ such that $\deg_{\Gamma_i} \mathcal{L} \leq 0$ for all $i$ then $H^0(\mathcal{O}_D(D) \otimes \mathcal{L}) = 0$. [Hint: Since $D^2 < 0$, any section $s \in H^0(\mathcal{O}_D(D))$ must vanish on some components of $D$. Now use the argument of Lemma 3.10.]

**15.** Let $Y \to X$ be a resolution of a normal surface singularity $P \in X$, and $\mathcal{L}$ a line bundle on $Y$ such that $\mathcal{L}\Gamma_i \geq K_Y \Gamma_i$ for each exceptional curve $\Gamma_i$. Prove that $H^1(D, \mathcal{L}_{|D}) = 0$ for every exceptional divisor $D$, and deduce that $H^1(Y, \mathcal{L}) = R^1 f_* \mathcal{L} = 0$. This is related to the vanishing theorems of Kodaira and Grauert–Riemenschneider, but works in all characteristic.

**16.** Suppose that $X$ is projective with ample $H$, and $f \colon Y \to X$ as usual. Let $Z = \sum a_i \Gamma_i$ be an effective exceptional divisor with $Z\Gamma_i < 0$. Prove that $n(f^*H) - Z$ is ample on $Y$ for $n \gg 0$. [Hint: Some section on $nH$ on $X$ vanishes on $Z$, so that $n(f^*H) - Z$ is effective for some $n$. Now prove that by taking a larger $n$ if necessary, $(n(f^*H) - Z)C > 0$ for every curve $C \subset Y$. The result follows by the Nakai–Moishezon ampleness criterion, compare footnote to 4.8 and [**H1**], Chapter V, 1.10.]

**17.** Show how to resolve the elliptic Gorenstein singularity

$$T_{p,q,r} \quad \text{defined by} \quad xyz + x^p + y^q + z^r \quad \text{for} \quad \frac{1}{p} + \frac{1}{q} + \frac{1}{r} \leq 1,$$

and describe the configuration of exceptional curves. [Hint: Start with the cases $p, q, r > 3$. Blowing up the maximal ideal gives a surface $X_1$ meeting the exceptional $\mathbb{P}^2$ in the tangent cone, which is the triangle $xyz = 0$, and having Du Val singularities at the 3 corners.]

**18.** Show how to resolve the elliptic Gorenstein singularities of degree 1 given by

$$x^2 + y^3 + z^k = 0 \quad \text{for } k = 6, 7, \dots, 11, \quad \text{and} \quad x^2 + y^3 + yz^l = 0 \quad \text{for } l = 5, 7.$$

[Hint: Use the weighted blowup with weights $3, 2, 1$, follows by the resolution of Du Val singularities.]

**19.** Ex. 4.8 is the special case $p = 2, q = r = 1$ of the rational triple points of Example 4.12. If you enjoyed it and want more of the same, carry out the calculations to verify the assertions made in Example 4.12.

**20 (harder).** The possible configurations of exceptional curves of a resolution $Y \to X$ of a rational triple points $P \in X$ are listed in [**A2**], p. 135. Verify this list following the hint given there. Find equations for some of these following the ideas of Example 4.12.

**21 (harder).** Rational surface singularities include quotient singularities $\mathbb{C}^2/G$ where $G \subset \mathrm{GL}(2, \mathbb{C})$ is a finite group not containing any quasireflections. Prove that $P \in X$ is a quotient singularity if and only if it is a rational singularity, and the minimal resolution $Y \to X$ has $-K_Y \overset{\text{num}}{\sim} \sum a_i \Gamma_i$ with $a_i \in \mathbb{Q}$, $a_i < 1$. List all the cases. [Hint: Because $P \in X$ is rational, some multiple of $K_X$ is a Cartier divisor. Then $P \in X$ has a cyclic cover which is a Du Val singularity. The lists can be found in Alexeev [**Utah2**], Chapter 3.]

**22.** Let $Z$ be a 2-connected divisor on a nonsingular surface and $\mathcal{L}$ a line bundle with $\mathcal{L} \otimes \omega_Z^{-1}$ nef and of degree $\geq 1$. Modify the proof of Lemma 4.23 to prove that if $|\mathcal{L}|$ is not free at $Q$ then $Q \in Z$ is a nonsingular point, and $\mathcal{L} \cong \omega_Z(Q)$.

**23.** State and prove 4.24, Claim 2 for $d = 3$ and $d = 4$. [Hint: For $d = 4$, compare [**UAG**], §1.]

**24.** Check the exact sequence 4.24, (14). [Hint: Step 1. There is a restriction map $\mathcal{I}_{i+1} \to \mathcal{I}_i$. The kernel is obviously contained in $\mathcal{O}_{\mathbb{P}^i}(-1)$. Step 2. Because $A_{i+1} \not\subset \mathbb{P}^{i-1}$, at $Q$, the kernel is equal to $m_Q \mathcal{O}_{\mathbb{P}^i}(-1)$. Step 3. It's tricky to see that $\mathcal{I}_{i+1} \to \mathcal{I}_i$ is surjective by local considerations at $Q$ in the crucial case $Q \in \operatorname{Supp} A_i$; however, there is a very easy proof by calculating $\chi$ of the 3 sheaves in the exact sequence

$$0 \to m_Q \mathcal{O}_{\mathbb{P}^i}(-1) \to \mathcal{I}_{i+1} \to \mathcal{I}_i.]$$

# Chapter D. Minimal models of surfaces via Mori theory

The final two chapters explain the proof of the classification of surfaces. The treatment of the classification of surfaces in these chapters is nonclassical, with emphasis on a new point of view, motivated in part by recent developments in 3-folds. For classical treatments, see [**Beauville**], or [**3 authors**]. Some sections of Chapter E are somewhat preliminary in nature.

## Summary

### Chapter D. Minimal models via Mori theory

### Chapter E. Classification of surfaces with $K$ nef

## D.0. Introduction and preview

Let $X$ be a projective nonsingular surface defined over $k$; it's my responsibility to make everything work over any algebraically closed field $k$, but the student who is so inclined may suppose that $k = \mathbb{C}$ to make things easier. The classification of surfaces breaks up into two parts. (For the definition of $K_X$ and of a nef divisor, see the end of D.1 below, or various places earlier in these notes.)

### Case $K_X$ not nef

This corresponds to studying the choice of birational model of $X$, the theory of minimal models and the characterisation of ruled and rational surfaces. Very roughly, the main result is that every surface is birational either to a surface $X$ with $K_X$ nef or to $\mathbb{P}^2$ or $\mathbb{P}^1 \times C$ with $C$ a curve.

The main point to note is the dichotomy: after choosing a suitable birational model, either $K_X$ is nef, or $-K_X$ is ample, or relatively ample on a nontrivial fibration $X \to C$. In this form, the result remains true over a nonclosed field, and is capable of extension to higher dimensions, although of course everything becomes harder.

### Case $K_X$ nef

This part is concerned with dividing surfaces into classes according to the numerical properties of $K_X$, and then proving that this coincides with the classification by the Kodaira dimension $\kappa(X)$. Very roughly, the main result is that $X$ has "sufficiently many" holomorphic differential forms. More precisely:

(i) If $K_X \overset{\text{num}}{\not\sim} 0$ (numerical equivalence) then $h^0(mK_X) \to \infty$ with $m$;

(ii) if $K_X \overset{\text{num}}{\sim} 0$ then $mK_X \overset{\text{lin}}{\sim} 0$ for some $m$. In other words, the topological condition $c_1(X) = 0 \in H^2(X, \mathbb{R})$ implies that $\mathcal{O}_X(mK_X)$ is analytically isomorphic to $\mathcal{O}_X$.

**Remark.** In both parts, the conditions can be relaxed to allow any field $k$, and $X$ to have various singularities; it can also make sense to study a *log surface*, that is, a surface $X$ together marked with a divisor $D$ (say reduced normal crossing), and classify the set-up $X$ with $D$, replacing $K_X$ by $K_X + D$.

## D.1. Preliminaries: 1-cycles versus codimension 1 cycles

Let $X$ be any projective scheme. Write $\operatorname{Pic} X = H^1(X, \mathcal{O}_X^*)$ for the group of invertible sheaves (or divisor classes) on $X$, and $\{\text{1-cycles}(X)\}$ for the group of 1-cycles on $X$, that is, the free group on all curves $\Gamma \subset X$. Then there is a pairing

$$\operatorname{Pic} X \times \{\text{1-cycles}(X)\} \to \mathbb{Z}$$

defined by $(\mathcal{L}, \Gamma) \mapsto \mathcal{L}\Gamma = \deg_\Gamma \mathcal{L}_{|\Gamma}$, and extended by linearity.

### D.1.1. Definition of numerical equivalence

Two invertible sheaves $\mathcal{L}_1, \mathcal{L}_2 \in \operatorname{Pic} X$ are *numerically equivalent* if $\mathcal{L}_1\Gamma = \mathcal{L}_2\Gamma$ for every curve $\Gamma \subset X$; we write $\mathcal{L}_1 \overset{\text{num}}{\sim} \mathcal{L}_2$. Similarly, two 1-cycles $C_1$, $C_2$ are *numerically equivalent* (written $C_1 \overset{\text{num}}{\sim} C_2$) if $\mathcal{L}C_1 = \mathcal{L}C_2$ for every $\mathcal{L} \in \operatorname{Pic} X$. Define

$$N^1X = (\operatorname{Pic} X / \overset{\text{num}}{\sim}) \otimes \mathbb{R} \quad \text{and} \quad N_1X = (\text{1-cycles} \otimes \mathbb{R}) / \overset{\text{num}}{\sim};$$

obviously, by definition $N^1X$ and $N_1X$ are dual $\mathbb{R}$-vector spaces, and $\overset{\text{num}}{\sim}$ is the smallest equivalence relation for which this holds. Clearly, $N^1X$ contains $\mathbb{Z}$- and $\mathbb{Q}$-forms, denoted by $N^1_{\mathbb{Z}}X$, $N^1_{\mathbb{Q}}X$ respectively. The fact that $\rho = \dim_{\mathbb{R}} N^1X$ is finite is the Néron–Severi theorem.

Notice that for a variety defined over $\mathbb{C}$, it is natural to consider $N_1X$ as a subspace of $H_2(X, \mathbb{R})$ and $N^1X$ as the corresponding quotient of $H^2(X, \mathbb{R})$ (in fact both are direct summands). Over $\mathbb{C}$, the finite dimensionality of $N^1X$ and $N_1X$ follows at once from this (see App. to §1 below); thus the content of the Néron–Severi theorem is a generalisation of this to the purely algebraic set-up. Notice also that is the same thing as $\operatorname{NS}(X) \otimes \mathbb{R}$, where $\operatorname{NS}(X)$ is the Néron–Severi group of $X$ (that is, Cartier divisors modulo algebraic equivalence).

### D.1.2. The case of a surface

I have deliberately avoided the following point: if $X$ is a nonsingular surface then, of course, 1-cycles are the same as codimension 1 cycles, so that $N^1X = N_1X$. This can be seen in two ways:

(1) every curve $\Gamma \subset X$ is a Cartier divisor (that is, locally defined by one equation) and

$$\{\text{1-cycles of } X\} = \operatorname{Div} X \to \operatorname{Pic} X = \operatorname{Div} X / \overset{\text{lin}}{\sim}$$

is defined by

$$\sum n_i C_i \mapsto \mathcal{O}_X(\sum n_i C_i).$$

(2) $\operatorname{Pic} X$ has an intersection pairing defined by

$$\mathcal{L}_1\mathcal{L}_2 = \chi(\mathcal{O}_X) - \chi(\mathcal{L}_1) - \chi(\mathcal{L}_2) + \chi(\mathcal{L}_1 \otimes \mathcal{L}_2),$$

which becomes nondegenerate on $N^1X$.

Despite the canonical equality $N^1X = N_1X$, it is nevertheless convenient to distinguish between the two spaces. This is useful for dealing with the higher dimensional case, and also for singular surfaces, when it's no longer true that all divisors are Cartier.

Although very simple, this is one of the key ideas of Mori theory, and came as a surprise to anyone who knew the theory of surfaces before 1980: the quadratic intersection form for curves on a nonsingular surface can for most purposes be replaced by the bilinear pairing between $N^1$ and $N_1$, and in this form generalises to singular varieties and to higher dimension.

### D.1.3. Definition of nef

An element $D \in N^1X$ is *numerically eventually free* or *nef* if $D\Gamma \geq 0$ for every curve $\Gamma \subset X$. This is related to the notion of a free linear system $|D|$: if $|mD|$ is free for some $m > 0$ ("eventually free") then $D$ is nef.

The relation between nef and ample divisors is the content of the Kleiman criterion (see D.2).

### D.1.4. The canonical class

If $X$ is a nonsingular $n$-fold then $\mathcal{O}_X(K_X) = \omega_X = \Omega^n_X$ is a line bundle, locally generated by $\mathrm{d}x_1 \wedge \cdots \wedge \mathrm{d}x_n$. (Compare A.10.) This is the only nontrivial line bundle which can be intrinsically associated with any variety.

## Appendix to D.1. Cycle class groups and homology

**D.1.5. Proposition.** *(i) If $X$ is a variety over $\mathbb{C}$ then $N_1X \subset H_2(X,\mathbb{R})$ and $H^2(X,\mathbb{R}) \twoheadrightarrow N^1X$ is the dual quotient.*

*(ii) In terms of Hodge theory,*

$$N^1X = H^2(X,\mathbb{R}) \cap H^{1,1} = H^2(X,\mathbb{R})/\big(H^{2,0} \oplus H^{0,2} \cap H^2(X,\mathbb{R})\big).$$

*(iii) Algebraic equivalence of 1-cycles implies numerical equivalence.*

**Sketch Proofs.** (i) The inclusion takes a curve $C$ into $[C]$ the 2-cycle obtained by triangulating $C$. The degree of a line bundle $\mathcal{L}$ on a curve is the same thing as the first Chern class $c_1(\mathcal{L}) \in H^2(X,\mathbb{R}) \twoheadrightarrow N^1X$ evaluated on the 2-cycle $[C] \in H_2(X,\mathbb{R})$.

(ii) Write down the exponential sequence of sheaves.

(iii) This follows easily from the fact that the degree of a line bundle is invariant in a flat family of curves. $\square$

## D.2. Easy consequences of RR

Let $X$ be a nonsingular surface. Then RR takes the form

$$\chi(\mathcal{O}_X(D)) = \chi(\mathcal{O}_X) + \tfrac{1}{2}D(D - K_X)$$

The fact that this is a quadratic function of $D$ has an immediate corollary.

**D.2.1. Proposition.** *$D^2 > 0$ implies either $h^0(nD) \neq 0$ or $h^0(-nD) \neq 0$ for $n \gg 0$. The sign of $HD$ (for any ample $H$) distinguishes the two cases.*

**Proof.** If $D^2 > 0$ then

$$h^0(nD) + h^2(nD) \geq \chi(\mathcal{O}_X(nD)) \sim n^2D^2/2,$$

so that either $h^0(nD)$ or $h^2(nD)$ goes to infinity with $n$. Now by Serre duality,

$$h^2(nD) = h^0(K_X - nD).$$

I can write down the same argument with $D$ replaced by $-D$, and conclude that either $h^0(-nD)$ or $h^0(K_X + nD)$ goes to infinity with $n$

Now obviously, $h^0(K_X - nD)$ and $h^0(K_X + nD)$ cannot both go to infinity with $n$; indeed, if $h^0(K_X - nD) \neq 0$, multiplying by a nonzero section $s \in H^0(K_X - nD)$ gives an inclusion $H^0(K_X + nD) \hookrightarrow H^0(2K_X)$ so that $h^0(K_X + nD) \leq h^0(2K_X)$. Therefore either $h^0(nD)$ or $h^0(-nD)$ grows quadratically with $n$.

Here is an alternative way of concluding the argument: it is easy to see that for any divisor $L$ on $X$,

$$|h^0(-nD) - h^0(L - nD)| < \text{const.} n.$$

(Just write $L = A - B$ with $A, B$ curves and consider restriction maps.) Hence if $h^2(nD) = h^0(K_X - nD)$ grows quadratically with $n$, so does $h^0(-nD)$. Q.E.D.

**D.2.2. Corollary** (the Index Theorem). *If $H$ is ample on $X$ then $HD = 0$ implies $D^2 \leq 0$. Moreover if $D^2 = 0$ then $D \overset{\text{num}}{\sim} 0$.*

In other words, the intersection pairing on $N^1X$ has signature $(+1, -(\rho - 1))$. The cone $D^2 > 0$ looks like the "light cone": the "past" and "future" components of the cone are separated by the hyperplane $H^\perp$. Another way of stating this, which is very useful in calculations, is that if $D_1, D_2$ are divisors and $(\lambda D_1 + \mu D_2)^2 > 0$ for some $\lambda, \mu \in \mathbb{R}$ then the determinant

$$\begin{vmatrix} D_1^2 & D_1D_2 \\ D_1D_2 & D_2^2 \end{vmatrix} \leq 0,$$

with equality if and only if $D_1, D_2$ are $\mathbb{Q}$-linearly dependent in $N^1X$.

**Proof.** If $D^2 > 0$ then either $nD$ or $-nD$ is equivalent to a nonzero effective divisor for $n \gg 0$. Since $H$ is ample, either of these conditions implies $HD \neq 0$. This proves (i).

(ii) is left as an exercise. [Hint: if $DA \neq 0$ for some curve $A$ then replace $A$ by $B = A - H$ so that $HB = 0$; now if $DA \neq 0$ also $DB \neq 0$, and some linear combination of $D$ and $B$ has $(D + \alpha B)^2 > 0$, contradicting (i).] Q.E.D.

**D.2.3. Corollary** (weak form of the Kleiman criterion). *$D$ is nef implies*

1. $D^2 \geq 0$ *and*
2. $D + \varepsilon H$ *is ample for any* $\varepsilon \in \mathbb{Q}$, $\varepsilon > 0$.

**Proof.** Introduce the quadratic polynomial $p(t) = (D + tH)^2$. Then $p(t)$ is a continuous increasing function of $t \in \mathbb{Q}$ for $t > 0$, and $p(t) > 0$ for sufficiently large $t$. The following assertion then obviously implies (1):

**D.2.4. Claim.** *Let $t \in \mathbb{Q}$, $t > 0$; then $p(t) > 0$ implies that also $p(t/2) > 0$.*

**Proof.** For $(D + tH)^2 > 0$, together with $H(D + tH) > 0$ implies that $n(D + tH)$ is equivalent to an effective divisor for suitable $n \gg 0$. By the assumption that $D$ is nef, $D(D + tH) \geq 0$, and hence $(D + (t/2)H)^2 = D(D + tH) + (t/2)^2 H^2 > 0$.

For (2) it is enough to notice that $(D + \varepsilon H)\Gamma > 0$ for every curve $\Gamma \subset X$ and $(D + \varepsilon H)^2 > 0$, so that $D + \varepsilon H$ satisfies the conditions of the Nakai–Moishezon criterion (see for example [**H2**], Chapter I or [**H1**], Chapter V, 1.10). Q.E.D.

### D.2.5. The Kleiman criterion

Corollary D.2.3 can be thought of as saying that nef divisors are limits of ample divisors. It is not hard to get from this to (a weak form of) Kleiman's ampleness criterion, which asserts that conversely, the cone of ample divisors is the interior of the nef cone in $N^1X$. For this make the following definitions: let $\operatorname{NE} X \subset N_1X$ be the cone of effective 1-cycles, that is

$$\operatorname{NE} X = \left\{ C \in N_1X \;\middle|\; C = \sum r_i\Gamma_i \text{ with } r_i \in \mathbb{R}, r_i \geq 0 \right\}.$$

Let $\overline{\operatorname{NE}}\, X$ be the closure of $\operatorname{NE} X$ in the real topology of $N_1X$. This is called the *Kleiman–Mori cone.*

One effect of taking the closure is the following trivial observation, which has many important uses in applications: if $H \in N^1_X$ is positive on $\overline{\operatorname{NE}}\, X \setminus 0$, then the section $(Hz = 1) \cap \overline{\operatorname{NE}}\, X$ is compact. Indeed, the projectivisation of the closed cone $\overline{\operatorname{NE}}\, X$ is a closed subset of $\mathbb{P}^{\rho-1} = \mathbb{P}(N_1X)$, and therefore compact, and the section $(Hz = 1)$ projects homeomorphically to it. The same also holds for any face or closed subcone of $\overline{\operatorname{NE}}\, X$.

**Theorem.** *For $D \in \operatorname{Pic} X$, view the class of $D$ in $N^1X$ as a linear form on $N_1X$. Then*

$$D \textit{ is ample} \iff DC > 0 \textit{ for all } C \in \overline{\operatorname{NE}}\, X \setminus 0.$$

Note that it is not true that $DC > 0$ for every curve $C \subset X$ implies $D$ ample. The necessary and sufficient condition $DC > 0$ for every $C \in \overline{\operatorname{NE}}\, X \setminus 0$ is just a tiny bit stronger than that. You can interpret it as saying that $DC$ is "reliably bigger than $C$".

**Sketch proof.** The implication $\Longleftarrow$ comes at once from (ii), since it is easy to see by a compactness argument that for $H$ ample and sufficiently small $\varepsilon > 0$, $(D - \varepsilon H)C > 0$ for all $C \in \overline{\operatorname{NE}}\, X \setminus 0$, so $D = (D - \varepsilon H) + \varepsilon H$ is ample by Corollary D.2.3, (2). For $\Longrightarrow$ you need to use a compactness argument to prove that if a norm $\| \;\|$ is chosen on $N_1X$ then for an ample divisor $D$ there exists $\varepsilon > 0$

such that

$$DC > \varepsilon ||C|| \quad \text{for all } C \in \operatorname{NE} X$$

(see for example [**H2**], Chapter I).

**Remarks.** (1) This is only a weak form of Kleiman's criterion, since $X$ is a priori assumed to be projective. The full strength of Kleiman's criterion gives a necessary and sufficient condition for ampleness in terms of the geometry of $\overline{\operatorname{NE}}\, X$. In particular, it can in be used to prove that a variety is projective (without assuming it!). See [**Kleiman**].

(2) An analogous statement, the Kleiman criterion (weak form), holds for an arbitrary projective $k$-scheme: if $D$ is a nef Cartier divisor then $D^r Z \geq 0$ for every irreducible $r$-dimensional subvariety $Z$, and $D$ is ample if and only if it is positive on $\overline{\operatorname{NE}}\, X \setminus 0$. The proof is more or less the same as (D.2.3–5), although you have to use cohomology and carry out an induction on the dimension. It goes something like this.

**Step 1.** (A version of the Nakai–Moishezon criterion.) If $D^r Z > 0$ for every $r$ and every irreducible $r$-dimensional subvariety $Z$ then

$$H^j(\mathcal{F} \otimes \mathcal{O}_X(nD)) = 0.$$

for every coherent sheaf $\mathcal{F}$, every $j > 0$, and all $n \gg 0$. (Just use induction on $\dim \operatorname{Supp} \mathcal{F}$.)

**Step 2.** If $D$ is nef and $H$ an ample Cartier divisor, then using RR, $(D+tH)^r Z > 0$ implies that $H^0(\mathcal{O}_Z(N(D+tH))) \neq 0$, so that $(D+tH)_{|Z}$ is an effective $\mathbb{Q}$-divisor. Then as in Claim D.2.4, using a further induction, one proves that for $t > 0$,

$$(D + tH)^r Z > 0 \implies (D + (t/r)H)^r Z > 0.$$

## D.3. The rationality lemma and cone theorem

Another easy consequence of the form of RR is the following:

**D.3.1. Rationality lemma.** *Suppose that $H$ is ample and $K_X$ is not nef. Define the* nef threshold $t_0 = t_0(H)$ *by*

$$t_0 = \sup\bigl\{t \bigm| H + tK_X \text{ is nef}\bigr\},$$

*so that $t_0 \in \mathbb{R}$, $t_0 > 0$. Then (i) $t_0 \in \mathbb{Q}$, and (ii) its denominator is $\leq 3$.*

Note that the condition $K_X$ not nef means exactly that $\overline{\operatorname{NE}}\, X$ meets the halfspace $(K_X z < 0)$ of $N^1 X$. Please draw a picture of the hyperplane $H+tK_X$ rotating from its initial position at $t = 0$ outside $\overline{\operatorname{NE}}\, X$ (since $H$ is ample) to its asymptotic position for $t \gg 0$, when it cuts into $\overline{\operatorname{NE}}\, X$; $t_0$ is the threshold value at which it first hits $\overline{\operatorname{NE}}\, X$.

**Discussion.** The underlying reason behind the cone theorem and the contraction theorem, and all their generalisations to higher dimensional singular and log varieties is *vanishing*. Although the arguments get technically quite complicated, the main mechanism is very simple: when $t = t_0 + \varepsilon$, just after the critical value, $H + tK_X$ is not nef, but vanishing applies to it, because $H + (t-1)K_X$ is ample.

The rationality lemma is proved in much greater generality in [**Kawamata**] and [**Kollár**]. For a good introduction to the problem, see [**Utah1**].

**Step 1.** If $n(H + t_1 K_X) = D_1$ is effective for some $n > 0$, with $t_1 \in \mathbb{Q}$ and $t_1 > t_0$ (so that $D_1$ is not nef), then $t_0$ is determined by

$$t_0 = \min_{\Gamma \subset D_1} \left\{ \frac{H\Gamma}{-K_X\Gamma} \right\},$$

where the minimum runs over the irreducible components $\Gamma$ of $D_1$ such that $K_X\Gamma < 0$. In fact for $0 < t < t_1$ the divisor $H + tK_X$ is a positive combination of $H$ and $D$, so that it fails to be nef if and only if $(H + tK_X)\Gamma < 0$ for a component $\Gamma$ of $D$.

**Step 2.** Now if $t_0 \notin \mathbb{Q}$ then for $n, m \in \mathbb{Z}$ with $n < mt_0 < n + 1$, if follows that $mH + nK_X$ is ample but $mH + (n+1)K_X$ is not nef. Since $mH + nK_X$ is ample, by Kodaira vanishing $h^0(mH + (n+1)K_X)$ is given by RR, and will sometimes be positive; this gives either a proof by Step 1 or a contradiction. More precisely, for $m > 0$, set $mt_0 = n + \alpha$ with $0 \le \alpha < 1$, and write $D_0 = H + t_0K_X \in N^1X$. Then RR gives

$$\begin{aligned} H^0(mH + (n+1)K_X) &\ge \chi(\mathcal{O}_X) + \tfrac{1}{2}(mH + (n+1)K_X)(mH + nK_X) \\ &= \chi(\mathcal{O}_X) + \tfrac{1}{2}(m^2D_0^2 + m(1-2\alpha)D_0K_X - \alpha(1-\alpha)K_X^2). \end{aligned}$$

(Here $\ge 0$ comes from the vanishing of $H^2$, which is trivial in any characteristic; if $\operatorname{char} k = 0$, the vanishing of $H^1$ gives equality.)

Hence if $D_0^2 > 0$ then $H^0 \ne 0$ for all large $m$; if $D_0^2 = 0$ and $D_0 \overset{\text{num}}{\not\sim} 0$ then necessarily $D_0K_X < 0$ (because $D_0(H + t_0K_X) = 0$ and $D_0H > 0$), and $H^0 \ne 0$ if $m$ is large and $1 - 2\alpha$ is bounded away from 0. If $D_0^2 = D_0K_X = 0$ then also $D_0H = 0$, so that $D_0 \overset{\text{num}}{\sim} 0$, and this implies that $K_X \overset{\text{num}}{\sim} -(1/t_0)H$; in particular $-K_X$ is ample.

This proves (i).

The proof of (ii) is an exercise. [Hint: if $D_0^2 > 0$ then there must exist a curve $L$ with $D_0L = 0$, since otherwise $D_0$ would be ample; it is easy to see that $D_0L = 0$ implies that $L$ is a $-1$-curve, so that $t_0 \in \mathbb{Z}$. If $D_0^2 = 0$, one sees that $H + tK_X$ cannot be effective for any $t > t_0$. If $D_0 \ne 0$ and $2t_0 \notin \mathbb{Z}$, the above inequality gives a contradiction by taking $m \gg 0$ with $(1 - 2\alpha)$ negative and bounded away from zero. Finally if $D_0 \overset{\text{num}}{\sim} 0$ then $-K_X$ is ample, and it is easy to see (in any characteristic, see D.4 or especially E.2) that $\chi(\mathcal{O}_X) = 1$. Now since $H + t_0K_X \overset{\text{num}}{\sim} 0$, it follows that the denominator of $t_0$ divides $K_X$ in $\operatorname{Pic} X$. It is now not hard to check that if $-K_X = dL$ with $L$ ample, then $d \le 3$: indeed, a simple argument of projective geometry shows that $h^0(L) \le 2 + L^2$, whereas RR gives $h^0(L) \ge 1 + (1/2)(d-1)L^2$. (ii) also follows easily by the method of proof of the contraction theorem.]

The rationality lemma implies the Mori cone theorem by an argument that is pure convex body theory, taken from [**Utah1**], Lecture 11. (This result is not logically necessary for the proof of classification.)

**D.3.2. Cone Theorem.** *Set*

$$\overline{\mathrm{NE}}_{K_X} X = \left\{ C \in \overline{\mathrm{NE}}\, X \;\middle|\; K_XC \ge 0 \right\}.$$

*Then*

$$\overline{\mathrm{NE}}\, X = \overline{\mathrm{NE}}_{K_X} X + \sum R_i, \qquad (*)$$

*where the $R_i$ are extremal rays of $\overline{\mathrm{NE}}\,X$ contained in $(K_Xz < 0)$. Moreover, for any ample Cartier divisor $H > 0$ and any given $\varepsilon > 0$, there are only finitely many extremal rays $R_i$ with $(K_X + \varepsilon H)R_i \le 0$.*

In other words, away from the half-space $K_X \ge 0$, the cone is locally polyhedral. (Draw a picture for God's sake! I'm imprisoned in a word-processor.) It follows by the result of D.4 that each ray $R_i$ is spanned by the class of an irreducible curve.

**Idea of the proof.** Suppose that $F_L = L^\perp \cap \overline{\mathrm{NE}}\,X$ is a face of the cone $\overline{\mathrm{NE}}\,X$ contained in the region $(K_X < 0)$ of the cone; if by bad luck $F_L$ is positive dimensional, I can wiggle $L$ slightly to make it a ray. This shows that each face of the cone contained in $(K_X < 0)$ is spanned by rays. In turn, I then show that each ray is defined by a set of linear equations with bounded denominators, so that the rays are discrete in $(K_X < 0)$.

To start on the formal proof, fix once and for all a basis of $N^1X$ of the form $K_X, H_1, \dots, H_{\rho-1}$, where the $H_i$ are ample and $\rho = \operatorname{rank} N^1$. As in the rationality lemma, for any nef element $L \in N^1X$, set

$$t_0(L) = \max\{t \mid L + tK_X \text{ is nef}\}.$$

If $L \in \operatorname{Pic} X$ is nef and the corresponding face $F_L = L^\perp \cap \overline{\mathrm{NE}}\,X$ is contained in $(K_Xz < 0)$ (for example, if $L$ is ample) then the rationality lemma gives $rt_0(L) \in \mathbb{Z}$ for the fixed integer $r = 6$.

**D.3.3. Main Claim.** *Let $L \in \operatorname{Pic} X$ be a nef divisor which supports a face $F_L$ of $\overline{\mathrm{NE}}\,X$ contained in $(K_Xz < 0)$. Consider $\nu L + H_i$ for all $i$ and for $\nu \gg 0$.*

1. *$t_0(\nu L + H_i)$ in an increasing function of $\nu$, is bounded above, and attains its bound.*
2. *Let $\nu_0$ be any point after $t_0(\nu L + H_i)$ has attained its upper bound, and suppose $\nu > \nu_0$. Set*
$$L' = 6\big(\nu L + H_i + t_0(\nu L + H_i)K_X\big),$$
*(multiplying by 6 is simply to ensure that $L' \in \operatorname{Pic} X$ by the rationality lemma). Then $L'$ supports a face $F_{L'} \subset F_L$.*
3. *If $\dim F_L \ge 2$ then there exists $i$ and $\nu \gg 0$ such that*
$$L'_i = 6\big(\nu L + H_i + t_0(\nu L + H_i)K_X\big)$$
*supports a strictly smaller face $F_{L'} \subsetneq F_L$.*
4. *In particular, $F_L$ contains an extremal ray $R$ of $\overline{\mathrm{NE}}\,X$.*
5. *If $F_L = R$ is an extremal ray and $z \in R$ is a nonzero element then $6H_iz/K_Xz \in \mathbb{Z}$.*
6. *The extremal rays of $\overline{\mathrm{NE}}\,X$ are discrete in $(K_Xz < 0)$.*

**Proof.** (1) is almost obvious: $t_0(\nu L + H_i)$ is an increasing function of $\nu$ by construction. It is bounded above, because any point $z \in F_L \setminus 0$ obviously gives $Lz = 0$ and $(H_i + tK_X)z < 0$ for $t \gg 0$. It attains its bounds because $t_0$ varies in the discrete set $(1/6)\mathbb{Z}$.

(2) Suppose that $t_0 = t_0(\nu L + H_i)$ does not change with $\nu \ge \nu_0$, and let $F_0$ be the face orthogonal to the nef element $\nu_0 L + H_i + t_0 K_X$. Then for $\nu > \nu_0$, any $z \in F_{L'_i}$ satisfies

$$0 = (\nu L + H_i + t_0K_X)z = (\nu - \nu_0)Lz + (\nu L + H_i + t_0K_X)z,$$

and therefore $Lz = 0$.

For (3), consider $\nu \gg 0$, and

$$L'_i = 6\bigl(\nu L + H_i + t_0(\nu L + H_i)K_X\bigr)$$

for each $i$. Since $\nu \gg 0$, these are "small wiggles" of $L$ in $\rho - 1$ linearly independent directions. Each $F_{L'_i} \subset F_L$ is a face of $F_L$. The intersection of all the $F_{L'_i}$ is contained in the set defined by

$$(\nu L + H_i + t_0(\nu L + H_i)K_X)z = 0,$$

which are $\rho - 1$ linearly independent conditions on $z$. Therefore at least one of the $F_{L'_i}$ is strictly smaller than $F_L$. (4) follows obviously from (3).

(5) follows from (2) and the rationality lemma. Indeed, since $F_L = R$ is a ray, (2) implies that $F_{L'_i} = F_L = R$. That is, $R$ is orthogonal to $H_i + t_0 K_X$, and $6t_0 \in \mathbb{Z}$.

Finally, (6) follows from (5), since in $(K_X z < 0)$, every ray contains a unique element $z$ with $K_X z = -1$ and $H_i z \in (1/6)\mathbb{Z}$. Q.E.D.

**Proof of Theorem D.3.2.** Write $B = \overline{\mathrm{NE}}_{K_X} X + \sum R_i$ for the right-hand side of $(*)$. Note that $B \subset \overline{\mathrm{NE}}\, X$ is a closed convex cone: indeed, by (6) the $R_i$ can only have accumulation points in $(K_X z \geq 0)$, that is, in $\overline{\mathrm{NE}}_{K_X} X$.

Suppose that $B \subsetneq \overline{\mathrm{NE}}\, X$. Then there exists an element $M \in N^1 X$ which is nef, and supports a nonzero face $F_M$ of $\overline{\mathrm{NE}}\, X$ disjoint from $B$, necessarily contained in $(K_X z < 0)$. By the usual compactness argument, for sufficiently small $\varepsilon > 0$, $M - \varepsilon K_X$ is ample, and $M + \varepsilon K_X$ is not nef, but is positive on $B$. These are all open conditions on $M$, and any open neighbourhood of $M$ in $N^1 X$ contains a rational element $M$, so that by passing to this and taking a multiple, I can assume that $M \in \operatorname{Pic} X$. But now the rationality lemma and Claim D.3.3, (4) imply that there is a ray $R$ of $\overline{\mathrm{NE}}\, X$ with $MR < 0$, so that $R \not\subset B$. This contradicts the choice of $B$, and proves $(*)$. Q.E.D.

## D.4. Contraction theorem and minimal models of surfaces

A consequence of the rationality lemma is that if $K_X$ is not nef then there exists a divisor $D$ on $X$ which is nef but not ample, and such that $D - \varepsilon K_X$ is ample for $\varepsilon > 0$ (namely the divisor $D = n(H + t_0 K_X)$ appearing in the rationality lemma). Similarly, it is clear from the statement of the cone theorem that for each extremal ray $R = R_i$ of $\overline{\mathrm{NE}}\, X$ there is a divisor class $D$ on $X$ such that $D \geq 0$ on $\overline{\mathrm{NE}}\, X$ and $\overline{\mathrm{NE}}\, X \cap D^\perp = R$, and this divisor satisfies the same condition.

**D.4.1. Contraction theorem.** *Let $D$ be a divisor which is nef but not ample, and such that $D - \varepsilon K_X$ is ample for some $\varepsilon > 0$. Then $|mD|$ is free for all $m \gg 0$.*

*Equivalently, there exists a morphism $\varphi\colon X \to Y$ to a projective variety $Y$ such that $D = \varphi^*(\text{ample})$. Of course, this means that for a curve $C \subset X$,*

$$\varphi(C) = \text{pt.} \iff DC = 0.$$

By taking normalisation, I can assume that $\mathcal{O}_Y = \varphi_* \mathcal{O}_X$. In the 3 possible cases, this morphism is Case $D^2 > 0$: $\varphi$ is the contraction of a finite number of disjoint $-1$-curves $L_i$ to nonsingular points of a surface $Y$.

Case $D^2 = 0$ but $D \overset{\text{num}}{\not\sim} 0$: $\varphi$ is a conic bundle structure $X \to$ (curve).

Case $D \overset{\text{num}}{\sim} 0$: $X$ is a surface with $-K_X$ ample and $p_g = q = 0$ (del Pezzo surface).

**Proof.** If $D^2 > 0$ then by the Nakai–Moishezon ampleness criterion there must be curves $L$ with $DL = 0$. Now any such curve is a $-1$-curve, that is, $L \cong \mathbb{P}^1$ and $L^2 = -1$; for by the index theorem $L^2 < 0$, and since $D - \varepsilon K_X$ is ample $K_X L < 0$, and therefore by the adjunction formula $K_X L = L^2 = -1$ and $p_a L = 0$. The index theorem implies that any two of these are disjoint. The standard classical proof of Castelnuovo's proof then gives the result.

If $D^2 = 0$ but $D \overset{\text{num}}{\not\sim} 0$ then since $D - \varepsilon K_X$ is ample, $K_X D < 0$, and RR gives $h^0(mD) \sim (-K_X D/2)m$, so in particular $h^0(mD) \to \infty$. The mobile part of $|mD|$ then gives a conic bundle structure $X \to$ (curve), and using this it is not hard to see that $|mD|$ itself is free.

Finally, if $D \overset{\text{num}}{\sim} 0$ then the condition on $D - \varepsilon K_X$ implies that $-K_X$ is ample, and the assertion is just that $D \overset{\text{lin}}{\sim} 0$. The key to this is to prove that $\chi(\mathcal{O}_X) = 1$, therefore $\chi(\mathcal{O}_X(D)) = \chi(\mathcal{O}_X) = 1$, and $H^0(\mathcal{O}_X(D)) \neq 0$, giving $D \overset{\text{lin}}{\sim} 0$ are required. If $\operatorname{char} k = 0$ then by Kodaira vanishing $H^i(\mathcal{O}_X) = 0$ for $i = 1, 2$, and so $\chi(\mathcal{O}_X) = 1$. If $\operatorname{char} k = p$ then $p_g = 0$ and $K_X^2 > 0$, so that the following alternative argument implies that $q = 0$ also in the case $\operatorname{char} k = p$. (See E.2 for how to get rid of the reference to Mumford.)

### D.4.2. Traditional argument

Because $p_g = 0$, the Picard scheme is reduced (see E.2, Remark (c), or compare [**M1**], Lecture 27), and $B_1 = B_3 = 2q$. Hence (in etale cohomology), $c_2 = 2 + B_2 - 4q$, so as $B_2 \geq 1$ and $K_X^2 > 0$ I get

$$12\chi(\mathcal{O}_X) = 12(1 - q) = c_1^2 + c_2 \geq 4 - 4q;$$

thus $q \leq 1$, and $q = 1$ is only possible if $B_2 = 1$. But if $q = 1$ then the fibre and section of the Albanese morphism imply $B_2 \geq 2$, so a contradiction. Therefore $\chi(\mathcal{O}_X) = 1$. Q.E.D.

**D.4.3. Corollary.** *Let $R$ be an extremal ray of $\overline{\mathrm{NE}}\,X$; then there exists a morphism $\varphi\colon X \to Y$ to a projective variety $Y$ such that $-K_X$ is relatively ample for $\varphi$, $\mathcal{O}_Y = \varphi_*\mathcal{O}_X$, and for all curves $C \subset X$,*

$$\varphi(C) = \text{pt.} \iff DC = 0.$$

The corresponding classification of extremal rays gives

1. contraction of a $-1$-curve;
2. $\mathbb{P}^1$-bundle $\varphi\colon X \to C$ over a curve $C$;
3. $\mathbb{P}^2$.

Note that in (2), the assertion that $X = \mathbb{P}(\mathcal{E})$ for a rank 2 vector bundle is nontrivial: you need to know that $X \to C$ has a rational section, either by Tsen's theorem (see Theorem C.4.2 or [**Sh**]), or an argument in topology using Poincaré duality, and the fact that because $p_g = 0$ the cohomology $H^2(X, \mathbb{Z})$ is spanned by curves on $X$ (see for example [**Beauville**], Theorem III.4). Similarly, to know that $X \cong \mathbb{P}^2$ in (3), the main point is to prove that $-K_X$ is divisible by 3 in $\operatorname{Pic} X$, which again follows from $K_X^2 = 9$ (Noether's formula) by Poincaré duality.

**D.4.4. Corollary** (Main Theorem of minimal model theory). *Let $X$ be any surface. Then there exists a chain of contractions of $-1$-curves $X \to X_1 \cdots \to X_N = Y$ such that either $K_Y$ is nef, or $Y$ is $\mathbb{P}^2$ or a $\mathbb{P}^1$-bundle over a curve.*

# Chapter E. Classification of surfaces with $K$ nef

## E.1. Statement of the main result

The first part of the classification of surfaces, the theory of minimal models, was concerned with the case $K_X$ not nef. Thus I assume throughout this chapter that $K_X$ is nef. It is convenient to summarise the main theorem of the classification of surfaces in the following form:

**Table–Main Theorem.** *A nonsingular projective surface $X$ with $K_X$ nef belongs to just one of the following cases.*

| Name of Case: | $\nu = 2$ | $\nu = 1$ | $\nu = 0$ |
|---|---|---|---|
| Definition: | $K_X^2 > 0$ | $K_X^2 = 0, K_X \overset{\text{num}}{\not\sim} 0$ | $K_X \overset{\text{num}}{\sim} 0$ |
| $P_m$ and $\kappa$ : | $P_m \sim m^2, \kappa = 2$ | $P_m \sim m, \kappa = 1$ | $P_m = 1\ {}^{\exists}m, \kappa = 0$ |
| Effective results: | $P_m \geq 2$ for all $m \geq 2$ | $P_m \geq 1$ for some $m \in \{1,2,3,4,6\}$<br>$P_m \geq 2$ for some $m \leq 42$ | $mK_X \overset{\text{lin}}{\sim} 0$ for some $m \in \{1,2,3,4,6\}$ |
| Structure result: | $\exists$ canonical model $\varphi\colon X \to Y$ | $\exists$ elliptic fibre space $\varphi\colon X \to C$ | either Abelian or $K3$ or etale quotient by $\mathbb{Z}/m$ (same $m$) |
| | | (or more special structure in char. 2 or 3) | |

In this table, the first two rows define the invariant $\nu(X) = \nu(K_X) = \kappa_{\text{num}}(X)$, the *numerical Kodaira dimension* of $X$: the three possible cases for $K_X^2$ and $K_X$ enumerated in the second row are exclusive and cover all surfaces with $K_X$ nef, and $\nu$ is defined to be 2, 1 or 0 accordingly; more generally, for any nef divisor on a projective variety $X$, the *characteristic dimension* of $D$ is defined by

$$\nu(X, D) = \max\bigl\{k \bigm| D^k \overset{\text{num}}{\not\sim} 0\bigr\}.$$

The third row states the theorem that $\nu(X)$ is equal to the Kodaira dimension $\kappa(X)$, which, as everyone knows, is defined by

$$\begin{aligned} \kappa = 2 &\iff P_m \sim m^2 \text{ for } m \gg 0, \\ \kappa = 1 &\iff P_m \sim m \text{ for } m \gg 0, \\ \kappa = 0 &\iff P_m \leq 1 \text{ for all } m, \text{ and } = 1 \text{ for some } m. \end{aligned}$$

The fourth row is a slightly more precise version of the same result, giving effective values of $m$ for which $P_m$ grows. The last row summarises the structure results in each case.

The special structures referred to in the final line are as follows: a surface with $\kappa = 1$ has in any case a fibre space structure $\varphi\colon X \to C$ with fibres of arithmetic genus 1, but in characteristic 2 and 3 there is an extra possible case, a "quasielliptic" fibre space: it can happen that every geometric fibre is a singular curve isomorphic to a cuspidal plane cubic. A surface with $\kappa = 0$ is either as stated in the table, or (in characteristic 2 or 3 only) a quotient of an Abelian surface by a nonreduced finite group scheme, or "quasi-bielliptic" (see E.8.4), or an Enriques surface in characteristic 2. I do not have time in these notes to treat these fascinating cases in the detail they deserve.

**Crude corollary.** *If $K_X$ is nef then there exists a morphism $\varphi\colon X \to Y$ to a projective variety $Y$ such that*

$$K_X \overset{\text{num}}{\sim} \varphi^*(\text{ample } \mathbb{Q}\text{-divisor on } Y).$$

In higher dimension this statement is called the "abundance conjecture", and remains at the time of writing a hard problem. The 3-fold case has recently been settled in several papers of Miyaoka and Kawamata (see Kollár [**Utah2**] for the latest information.) For surfaces, the nontrivial assertion is

$$\nu = 1 \implies \exists \text{ an elliptic fibre space.}$$

This is proved in E.6, and is really the main aim of this chapter. The remainder of the proof of the main theorem can be viewed simply as tidying up after this, using basically similar arguments.

Note that there is essentially nothing to prove in the theorem if $p_g \geq 2$: the case $\nu = 0$ is excluded, and $\nu = 1 \iff K_X^2 = 0$, in which case all the conclusions in the $\nu = 1$ column are satisfied rather trivially.

**Corollary** (Enriques' criterion)**.** *Let $X$ be a surface, without the condition that $K_X$ nef. Then equivalent conditions*

1. *$X$ is birationally ruled;*
2. *$\nexists$ birational model of $X$ with $K_X$ nef;*
3. *$P_m = 0$ for all $m \in \{1,2,3,4,6\}$;*
4. *$\kappa(X) = -\infty$;*
5. *adjunction terminates;*

*etc.*

**Corollary** (Castelnuovo's criterion)**.** *$X$ is rational if and only if $p_g = P_2 = 0$ and $q = 0$.*

In the current view of minimal models of surfaces via Mori theory, ruled surfaces are characterised in the first instance by the numerical property of having no model with $K_X$ nef; hence the proof of Enriques' criterion, in the form

$$\exists \text{ model with } K_X \text{ nef} \implies P_m \neq 0 \text{ for some } m \in \{1,2,3,4,6\}$$

becomes a problem in the classification of surfaces with $K_X$ nef. Thus minimal model theory is concerned (as a matter of ideology) purely with the numerical properties of $K_X$, leaving analytic questions such as the dimension of $H^0(mK_X)$ to the second stage of classification theory. The classical arguments of Castelnuovo and Enriques work with $P_2 = 0$ via "adjunction terminates", mixing up analytic constructions such as the Albanese map together with numerical considerations of components of an effective divisor $D$ with $H^0(D + K_X) = 0$. See [**Beauville**], Chapter 6 for a readable account.

## E.2. Traditional numerical game

It turns out to be useful to tabulate the possibilities for the numerical invariants of surfaces with $p_g \le 1$ (compare Bombieri and Mumford [**B–M2**], Introduction). For this, use Max Noether's formula

$$({c_1}^2 + c_2)(X) = 12\chi(\mathcal{O}_X),$$

and interpret the terms on the two sides in terms of other invariants of $X$.

In characteristic 0 this is very easy (refer to Chapter C for a discussion of the invariants): on the right-hand side $\chi(\mathcal{O}_X) = 1 - q + p_g$, and on the left-hand side ${c_1}^2(X) = {K_X}^2$ and $c_2(X) = E(X) = 2 - 4q + 2p_g + h^{1,1}$.

In characteristic $p$, (or for non-Kähler complex surfaces) a little care is needed. Write $q' = h^1(\mathcal{O}_X)$, and $q = \dim \operatorname{Pic}^0 X = \dim \operatorname{Alb} X$.

**Proposition** (Igusa)**.** *(i) $q' \ge q$; (ii) $c_2(X) = 2 - 4q(X) + B_2$ where $B_2 \ge \rho(X)$.*

**Discussion of proof.** (i) The connected group scheme $\operatorname{Pic}^0 X$ contains an Abelian variety $A = (\operatorname{Pic}^0 X)_{\text{red}}$ such that $(\operatorname{Pic}^0 X)/A$ is a finite group scheme with only one point (possibly nonreduced). Now $q' = h^1(\mathcal{O}_X)$ is the dimension of the Zariski tangent space to $\operatorname{Pic}^0 X$, and so $q' \ge \dim \operatorname{Pic}^0 X = \dim A$; the Albanese variety of $X$ is the dual Abelian variety to $A$, so that $q = \dim A$.

(ii) The second Chern class $c_2(X) \in H^4_{\text{et}}(X, \mathbb{Q}_\ell)$ equals the Euler characteristic of $X$ in etale cohomology,

$$c_2(X) = B_0 - B_1 + B_2 - B_3 + B_4,$$

where $B_i = B_i(X) = \dim_{\mathbb{Q}_\ell} H^i_{\text{et}}(X, \mathbb{Q}_\ell)$ and $\ell$ is a prime $\ne p$. Now $B_0 = B_4 = 1$ by Poincaré duality, and

$$B_1 = B_3 = 2q,$$

since $H^1_{\text{et}}$ is related to $\pi_1$ (more precisely, their $\ell$-primary parts), and $\pi_1(X)$ to $\pi_1(\operatorname{Alb} X)$ essentially as in the classical case.

Finally, under the cycle class map $\operatorname{NS}(X) \to H^2_{\text{et}}(X, \mathbb{Q}_\ell)$, the intersection number of divisors goes over to the cup product in etale cohomology, which implies that $\operatorname{NS}(X) \hookrightarrow H^2_{\text{et}}$, and hence $B_2 \ge \rho(X)$. Q.E.D.

Putting everything together gives

$$K_X^2 + 2 - 4q + B_2 = 12 - 12q' + 12p_g;$$

that is,

$$K_X^2 + 12(q' - q) + 8q + B_2 = 10 + 12p_g. \qquad (*)$$

**Famous Table.** *$p_g \le 1$ leads to the following table of cases:*

| | | | | |
|---|---|---|---|---|
| 1 | $p_g = 1$ | $q' = q = 0$ | $K_X^2 + B_2 = 22$ | (if $K_X \overset{\text{num}}{\sim} 0$ then $K3$) |
| 2 | $p_g = 1$ | $q' = q = 1$ | $K_X^2 + B_2 = 14$ | (doesn't exist with $K_X \overset{\text{num}}{\sim} 0$) |
| 3 | $p_g = 0$ | $q' = q = 0$ | $K_X^2 + B_2 = 10$ | (if $K_X \overset{\text{num}}{\sim} 0$ then Enriques) |
| $3'$ | $p_g = 1$ | $q' = 1, q = 0$ | $K_X^2 + B_2 = 10$ | (if $K_X \overset{\text{num}}{\sim} 0$ then Enriques, char 2) |
| 4 | $p_g = 1$ | $q' = q = 2$ | $K_X^2 + B_2 = 6$ | (if $K_X \overset{\text{num}}{\sim} 0$ then Abelian) |
| 5 | $p_g = 0$ | $q' = q = 1$ | $K_X^2 = 0, B_2 = 2$ | (if $K_X \overset{\text{num}}{\sim} 0$ then bielliptic) |
| $5'$ | $p_g = 1$ | $q' = 2, q = 1$ | $K_X^2 = 0, B_2 = 2$ | (if $K_X \overset{\text{num}}{\sim} 0$ then quasi-bielliptic) |

**Proof.** If I restrict to surfaces with $p_g \leq 1$ the right-hand side of $(*)$ is $\leq 22$, and the left a sum of positive terms. It's easy to tabulate the possibilities. For the conclusion $K_X^2 = 0$ in (5) and (5′), note that if $q = 1$ then $B_2(X) \geq \rho(X) \geq 2$; indeed, then $\operatorname{Alb} X$ is a curve, and $\alpha\colon X \to \operatorname{Alb} X$ a nonconstant morphism, so that an ample divisor and the fibre of $\alpha$ provide two linearly independent classes in $\operatorname{NS}(X)$. Q.E.D.

**Remarks.** **(a)** The proof of abundance when $p_g \leq 1$ appeals to Famous Table E.2 several times. However, it's not clear whether the numerical game is really central to the classification of surfaces; if so, then it seems to me that the abundance conjecture in higher dimensions is in serious trouble.

**(b)** The last 3 cases have $\chi(\mathcal{O}_X) = 0$, and otherwise $\chi(\mathcal{O}_X) > 0$. A priori, surfaces with $\chi(\mathcal{O}_X) < 0$ may exist, but not under the assumption $p_g \leq 1$, by what I've just said. Note that in characteristic 0, Bogomolov's inequality ${c_1}^2 \leq 4c_2$ implies that $\chi(\mathcal{O}_X) \geq 0$ and $\chi(\mathcal{O}_X) > 0$ if $K_X^2 > 0$.

**(c)** Notice that one reads directly from (3) and (5) of the table that $p_g = 0$ implies $q = q'$, that is, that $\operatorname{Pic}^0 X$ is reduced. More generally, for any projective scheme $X$, in the same way that $H^1(\mathcal{O}_X)$ is the space of 1st order infinitesimal deformations of a line bundle, the vector space $H^2(\mathcal{O}_X)$ is the obstruction space for deformations, hence $H^2(\mathcal{O}_X) = 0$ implies that $\operatorname{Pic}^0 X$ is reduced, so that $q = q'$. Mumford [**M1**], Lecture 27 shows how the tangent space to $(\operatorname{Pic}^0 X)_{\mathrm{red}}$ as a vector subspace of $H^1(\mathcal{O}_X)$ can be determined more generally. However, none of this fancy stuff is required for the current proof.

## E.3. General type

For the proof of the main theorem, I start by separating off the case $\nu = 2$.

**Theorem.** *Suppose $K_X^2 > 0$. Then the curves (if any) with $K_X\Gamma = 0$ are all $-2$-curves (that is, $\Gamma \cong \mathbb{P}^1$ and $\Gamma^2 = -2$), and there is at most a finite set $\{\Gamma_i\}$ of them; each connected component of $\bigcup \Gamma_i$ supports only rational cycles, and can be contracted by a morphism $\varphi\colon X \to Y$, so that $Y$ has only Du Val singularities, and $K_Y$ is an ample Cartier divisor. Thus $Y$ is the canonical model of $X$. Hence $X$ is of general type, and its canonical ring $R(X, K_X) = R(Y, K_Y)$ is finitely generated over $k$.*

*Moreover $P_m \geq 2$ for all $m \geq 2$.*

**Proof.** The first part is standard and the proof is omitted; see Chapter 4 for Du Val singularities (especially Contraction Theorem 4.15), or [**M2**]. If $p_g \geq 2$ then $|K_X|$ is positive dimensional, hence also $|mK_X|$ for each $m \geq 1$. Thus to prove $P_m \geq 2$, I can assume $p_g \leq 1$. This is the first application of Famous Table E.2: I can deal separately with the different cases.

Cases (5) and (5′) are excluded since $K_X^2 > 0$. Cases (1), (2), (3) and (3′) all satisfy $\chi(\mathcal{O}_X) \geq 1$, so RR gives

$$h^0(mK_X) \geq \chi(\mathcal{O}_X(mK_X)) \geq 1 + \binom{m}{2} K_X^2 \geq 2 \qquad \text{for all } m \geq 2.$$

The only remaining case is (4). Then $\chi(\mathcal{O}_X) = 0$ and $c_2 = -K_X^2 < 0$; in characteristic zero, this is impossible, since surfaces of general type have $c_2 > 0$ (for example by Bogomolov's inequality ${c_1}^2 \le 4c_2$). In any case,

$$h^0(mK_X) \ge \chi(\mathcal{O}_X(mK_X)) \ge \binom{m}{2} K_X^2 \ge 1 \qquad \text{for all } m \ge 2,$$

and $\ge 2$ except possibly for $m = 2$ and $K_X^2 = 1$. Together with the assumptions of case (4), this gives

$$p_g = 1, \quad q = q' = 2 \quad \text{and} \quad K_X^2 = 1. \tag{$**$}$$

**Claim** (Shepherd-Barron). $(**)$ *is impossible.*

**Proof.** Consider the Albanese morphism $X \to \operatorname{Alb} X$ and its Stein factorisation $\alpha\colon X \to A$. By [**Sh-Barron**], Theorem 7, $X$ has a covering family of geometrically rational curves, each of which is contracted by $\alpha$. Since $q = 2$, it follows that $\alpha\colon X \to A \subset \operatorname{Alb} X$ is a morphism to a curve of genus 2.

Let $C \in |K_X|$ be the unique curve. Then since $K_X^2 = 1$, $C$ has a unique irreducible component with $K_X\Gamma = 1$, possibly together with some $-2$-curves. If $C$ is reducible, or irreducible but singular, then $\Gamma$ has geometric genus $\le 1$, so has no nonconstant morphism to $A$. It would follow that $\alpha$ contracts the whole of $C$ to a point of $A$, which contradicts the index theorem since $C^2 = K_X^2 = 1$. Therefore $C$ is irreducible and nonsingular.

The restriction of $\alpha$ to $C$ is either an isomorphism, or has $\deg \alpha > 1$ (then $\alpha$ is necessarily inseparable, since $g(C) = g(A) = 2$). Either case leads to a contradiction: if $C \to A$ is an isomorphism, $C$ meets any fibre $F$ of $\alpha$ transversally in one point, so $CF = 1$; but then $K_XF = CF = 1$ contradicts $F^2 = 0$.

Now since $X \to A$ is its own Stein factorisation, the function field $k(A)$ is algebraically closed in $k(X)$. It follows that the differential of $\alpha$ is a nonzero map $s : \alpha^*(\Omega^1_A) \to \Omega^1_X$. Let $k \ge 0$ be the maximum number such that $s$ vanishes $k$ times along $C$, so that $s\colon \alpha^*(\Omega^1_A) \to \Omega^1_X(-kC)$ has nonzero restriction to $C$. Then the diagram

$$\begin{array}{ccccccccc} & & & & \alpha^*(\Omega^1_A) & & & & \\ & & & & \downarrow & & & & \\ & & & & \Omega^1_X(-kC) & & & & \\ & & & & \downarrow & & & & \\ 0 & \to & \mathcal{O}_C(-(k+1)C) & \to & \Omega^1_X(-kC)_{|C} & \to & \Omega^1_C(-kC) & \to & 0 \end{array}$$

gives a nonzero map of $\alpha^*(\Omega^1_A)$ to one of

$$\Omega^1_C(-kC) \quad \text{or} \quad \mathcal{O}_C(-(k+1)C).$$

This is a contradiction, since $\alpha^*(\Omega^1_A)$ has degree $2\deg\alpha > 2$, whereas

$$\deg \Omega^1_C(-kC) = 2 - k \quad \text{and} \quad \deg \mathcal{O}_C(-(k+1)C) = -(k+1). \quad \text{Q.E.D.}$$

## E.4. The cases with $\chi(\mathcal{O}_X) > 0$

The following result deals with all surfaces with $\chi(\mathcal{O}_X) > 0$ and $p_g \le 1$ (that is, cases (1), (2), (3), (3′) of Famous Table E.2).

**Theorem.** *Assume that* $K_X^2 = 0$ *and* $\chi(\mathcal{O}_X) = 1 - q' + p_g > 0$. *Then*
*either* $K_X \overset{\text{lin}}{\sim} 0$ *and* $q' = 0$ *(K3 surface);*
*or* $2K_X \overset{\text{lin}}{\sim} 0$ *and* $p_g = q' = 0$ *(classical Enriques surface);*
*or* $\operatorname{char} k = 2$, $K_X \overset{\text{lin}}{\sim} 0$ *and* $p_g = q' = 1$, $q = 0$ *(nonclassical Enriques surface);*
*or* $P_6 \geq 2$.

**Proof.** The idea of the proof is to plug the assumption $\chi(\mathcal{O}_X) = 1$ or $2$ into RR, giving $\chi(\mathcal{O}_X(mK_X)) \geq 1$. For $m \geq 2$, note that

$$h^0(mK_X) + h^2(mK_X) \geq \chi(\mathcal{O}_X(mK_X)) \geq 1.$$

But $H^2(mK_X) \overset{\text{d}}{=} H^0(-(m-1)K_X)$, so that $H^2 \neq 0$ implies that $-(m-1)K_X$ is effective, hence $(m-1)K_X \overset{\text{lin}}{\sim} 0$, since $K_X$ is nef. The theorem follows easily from a small case division based on this. In more detail:

**Case (1):** $q' = 0$, $p_g = 1$, so that $\chi = 2$. Then $\chi(\mathcal{O}_X(2K_X)) = 2$, and hence either $P_2 = h^0(2K_X) \geq 2$ or $H^2(2K_X) \neq 0$. The second case gives $H^0(-K_X) \neq 0$, so that $K_X \overset{\text{lin}}{\sim} 0$. Note that this leads directly to the definition of K3s as surfaces with $K_X \overset{\text{lin}}{\sim} 0$ and $q' = 0$ .

**Case (2):** $p_g = 1$, $q' = q = 1$; then $\chi(\mathcal{O}_X) = 1$, so RR gives

$$\chi(\mathcal{O}_X(mK_X + \sigma)) = 1 \quad \text{for every } m \text{ and every } \sigma \in \operatorname{Pic}^0 X.$$

Hence there is an effective divisor $D_\sigma \in |K_X + \sigma|$ for every $\sigma \in \operatorname{Pic}^0 X$; divisor of the form $D_\sigma + D_{-\sigma}$ provide infinitely many distinct elements of $|2K_X|$, and therefore $P_2 \geq 2$. In particular $\kappa(X) = 1$.

**Case (3):** $p_g = q' = 0$; then $\chi(\mathcal{O}_X(2K_X)) = 1$. Now as explained at the start of the proof, $p_g = 0$ implies that $H^2(2K_X) = 0$, so that $H^0(2K_X) \neq 0$. Let $0 \leq D_2 \overset{\text{lin}}{\sim} 2K_X$ be an effective divisor. If $D_2 = 0$ then I have the second conclusion of the theorem (a classical Enriques surface). Otherwise $D_2 > 0$ implies that in turn $H^2(3K_X) = 0$, so again by RR $H^0(3K_X) \neq 0$. Choose an effective divisor $0 < D_3 \overset{\text{lin}}{\sim} 3K_X$; now obviously $3D_2 = 2D_3$ would imply that $D_2 = 2D_1$ and $D_3 = 3D_1$, with $0 < D_1 \overset{\text{lin}}{\sim} K_X$, contradicting $p_g = 0$. Hence $|6K_X|$ contains two distinct elements $3D_2$ and $2D_3$, which gives $P_6 \geq 2$.

**Case (3′):** $p_g = q' = 1$, $q = 0$. Suppose first that $K_X \overset{\text{lin}}{\sim} 0$. The argument of [**B–M2**], p. 39 shows that in this case $\operatorname{char} k = 2$. (Roughly, if $\operatorname{char} k = p$, a simple subgroup scheme $\mu_p$ or $\alpha_p$ of the group scheme $\operatorname{Pic}^0 X$ corresponds to a cover $Y \to X$ which is respectively etale and Galois with group $\mathbb{Z}/p$, or a torsor under $\alpha_p$; then $\chi(\mathcal{O}_Y) = p \cdot \chi(\mathcal{O}_X) = p$, but on the other hand one sees that $Y$ is a reduced irreducible surface with $K_Y \overset{\text{lin}}{\sim} 0$, so that $\chi(\mathcal{O}_Y) \leq 2$. Thus $\operatorname{char} k = p \leq 2$.) In [**B–M2**], p. 26, Enriques surfaces are defined by $K_X \overset{\text{num}}{\sim} 0$ and $B_2 = 10$, so that cases (3) and (3′) of Famous Table E.2 are classical and nonclassical Enriques surfaces by definition.

If $K_X \overset{\text{lin}}{\sim} D > 0$, you can prove that $P_2$ or $P_3 \geq 2$ as a very instructive exercise in the techniques of E.6, using the assumption $q' = 1$. [Hint: Using the terminology and basic results of E.6.1, write $K_X \overset{\text{lin}}{\sim} D = rE + D_1$, with $E$ a 0-curve, $r \geq 1$ and $D_1$ disjoint from $E$. If $D_1 \neq 0$ then $D$ is disconnected, and $h^0(\mathcal{O}_D) \geq 2$; by RR and duality, this implies that $h^0(\omega_D) = h^0(\mathcal{O}_D(K_X + D)) \geq 2$, from which by cohomology $P_2 = h^0(\mathcal{O}_X(K_X + D)) \geq 2$. Otherwise, $D_1 = 0$, so that $K_X \overset{\text{lin}}{\sim} D = rE$, and by adjunction $\mathcal{O}_E((r+1)E) \cong \omega_E \cong \mathcal{O}_E$. Then the exact sequence

$$0 \to \mathcal{O}_E(-(r+1)E) \to \mathcal{O}_{(r+2)E} \to \mathcal{O}_{(r+1)E} \to 0$$

implies that $h^0(\mathcal{O}_{(r+2)E}) \geq 2$. Thus also $h^0(\mathcal{O}_{(r+2)E}(K_X + (r+2)E)) \geq 2$. From this, using the assumption $q' = 1$ and cohomology, $|K_X + (r+2)E| = |(2r+2)E|$ moves, so that $P_m \geq 2$ for all $m \geq 4$ (and for all $m \geq 3$ if $r > 1$).] Q.E.D.

**Example**

In every case except (3) and (3′) I've proved that $P_2 \geq 2$. An elliptic surface $f\colon X \to \mathbb{P}^1$ with a double fibre $F_2 = 2E_2$ and a triple fibre $F_3 = 3E_3$, and whose Jacobian fibration is a rational elliptic surface has $\chi(\mathcal{O}_X) = 1$ and

$$K_X = f^*\mathcal{O}(-1) + E_2 + 2E_3 = -F + \frac{1}{2}F_2 + \frac{2}{3}F_3;$$

it's easy to see that $X$ has $p_g = 0$, $P_2 = \cdots = P_5 = 1$, $P_6 = 2$. (Compare with Champion E.7.5.)

## E.5. Digression: the case $\chi(\mathcal{O}_X) = 0$ in characteristic 0

This section is not needed for the main proof, and is intended just to put the difficulty in its place.

In characteristic 0, $c_2 \geq 0$ (for example by Bogomolov's inequality), so also $\chi(\mathcal{O}_X) = (1/12)c_2 \geq 0$. The arguments given above in the case $\chi(\mathcal{O}_X) > 0$ are rather simple, so treating the remaining cases $K_X^2 = c_2 = 0$ should be considered as the heart of the classification of surfaces. The thing that makes the proof awkward in these cases is that RR gives $\chi(\mathcal{F}) = 0$ for any naturally occurring sheaf on $X$ (compare the discussion in E.9.3). Curiously, although in this case the proof of classification is quite involved, the answer itself is very simple (see E.8.6 for the proof):

**E.5.1. Theorem.** *If $X$ is an algebraic surface over $\mathbb{C}$ such that $K_X^2 = 0$ and $c_2 = 0$ then $X$ is either an Abelian surface, or a surface of the form $X = (C \times E)/G$, where $C$ is a curve of genus $\geq 1$, $E$ a curve of genus 1, and $G$ is a finite Abelian group acting anyhow on $C$, and on $E$ by translations, in such a way that the diagonal action on $C \times E$ is free.*

**E.5.2. Definition.** I propose to call the surfaces $(C \times E)/G$ *sesquielliptic surface.* When $K_X \overset{\text{num}}{\sim} 0$ then $g(C) = 1$, so that the surfaces here are the *bielliptic surfaces* which occupy the position $p_g = 0$, $q = 1$ in the classification of surfaces with $K_X \overset{\text{num}}{\sim} 0$ (see [**Beauville**], VI.19–20 and [**B–M2**], §3). The name is good, since they are characterised by having two elliptic fibrations, to $\mathbb{P}^1$ and $\operatorname{Alb} X$. Sesquielliptic surface also have two projections to $C/G$ and $E/G = \operatorname{Alb} X$, the first of which is an elliptic pencil; (according to the dictionary, the prefix "sesqui-", meaning $1\frac{1}{2}$, is a contraction of Latin *semi* + *que* (and)).

**E.5.3. Theorem.** *A surface $X$ with $p_g = 0, q = 1$ is sesquielliptic. A sesquielliptic surface $(C \times E)/G$ has $\nu = 0$ if $g(C) = 1$, and $\nu = 1$ if $g(C) \geq 2$.*

Note that the sesquielliptic conclusion can be interpreted as saying that $X$ is a quotient of $\mathbb{C} \times \mathbb{C}$ or $\mathbb{C} \times \mathcal{H}$ (where $\mathcal{H}$ is the complex upper half-plane). Thus Theorem E.5.1 says that $c_1^2(X) = c_2(X) = 0$ implies that $X$ has a locally homogeneous differential geometric structure. It would be reasonable to look for a differential geometric proof of this fact. See E.9.4 for a wild discussion.

## E.6. The existence of an elliptic pencil

I now give the proof of abundance: if $\nu(X) = 1$ then $X$ has an elliptic pencil. The proof breaks up into two big steps: (I) $X$ has a 0-curve $E$ (an effective divisor having the numerical properties of a nonmultiple fibre of an elliptic pencil); and (II) if $K_X$ is nef and $E$ a 0-curve then a multiple of $E$ moves in an elliptic pencil.

It's not necessary to assume $\operatorname{char} k = 0$, or anything about the values of $c_2$ or $p_g$, although if $p_g \geq 2$ or $c_2 > 0$ everything has been proved in the preceding section.

### E.6.1. 0-curves

I write down some preliminary material, leaving the easy proofs to the reader. For details, see [**Beauville**], Chapter VIII, or [**3 authors**], or [**M3**], p. 332 (or work them out for yourself).

**(a)** A divisor *of elliptic fibre type* is a nef divisor $D \overset{\text{num}}{\not\sim} 0$ with $D^2 = DK_X = 0$. It follows by the Index Theorem that $K_X \overset{\text{num}}{\sim} aD$ with $a \in \mathbb{Q}$ and $a \geq 0$. If $D = \sum n_i C_i$ is effective then $DC_i = K_X C_i = 0$ for every component $C_i$; it is easy to see from this that the intersection matrix $\{C_i C_j\}$ is negative semidefinite (see Theorem A.7, (2)). In other words, $D$ has the numerical properties of a sum of fibres of an elliptic fibre space.

**(b)** A 0-*curve* is an effective divisor $E = \sum n_i C_i$ of elliptic fibre type for which $\operatorname{Supp} E$ is connected and $\operatorname{hcf}(n_i) = 1$. Then $E$ is a Gorenstein curve with $\omega_E \cong \mathcal{O}_E$ and $H^0(\mathcal{O}_E) = 1$; $E$ behaves throughout in every respect as if it were an irreducible nonsingular elliptic curve; the crucial property is that for any line bundle $\mathcal{L}$ on $E$,

$$\deg_{C_i} \mathcal{L} = 0 \text{ for all } i \text{ and } H^0(\mathcal{L}) \neq 0 \quad \Longrightarrow \quad \mathcal{L} \cong \mathcal{O}_E.$$

**(c)** If $E = \sum n_i C_i$ is a 0-curve and $D$ an effective divisor on $X$ with $DC_i = 0$ for all components of $E$ then

$$D = rE + D',$$

with $r \geq 0$ and $D'$ disjoint from $E$; in particular the line bundle $\mathcal{O}_X(D')$ is trivial in a neighbourhood of $E$, so that $\mathcal{O}_{aE}(D') \cong \mathcal{O}_{aE}$ for every $a > 0$.

**(d)** Any effective divisor $D$ of elliptic fibre type is a sum of 0-curves $D = \sum n_i E_i$.

Mumford's rather regrettable terminology is "canonical type" for "elliptic fibre type" and "indecomposable divisor of canonical type" for "0-curve".

### E.6.2. Remark

Miyaoka's inequality for log surfaces [**Miyaoka**] implies that in characteristic 0, a surface with $K_X$ nef, $K_X^2 = 0$ and $c_2 = 0$ does not contain $-2$-curves or nodal or cuspidal elliptic curves. So as far as the main case $\operatorname{char} k = c_2 = 0$ is concerned, a 0-curve $E$ actually is a nonsingular elliptic curve.

### E.6.3. Step I, the existence of a 0-curve

Since $\nu(X) = 1$, it is enough to prove that there exists an effective divisor numerically equivalent to $mK_X$ with $m > 0$, and then apply (d). If $p_g \neq 0$ this is trivial. If $p_g = 0$ then by (3) and (5) of Famous Table E.2, $q' = q \leq 1$; as above, if $q' = 0$ and $K_X \overset{\text{num}}{\not\sim} 0$ then $P_2 > 0$ by RR. This leaves only the case $p_g = 0$, $q' = q = 1$. Then $A = \operatorname{Pic}^0 X$ is reduced, so an elliptic curve.

**E.6.4. Proposition.** *$H^0(K_X + \sigma) \neq 0$ for some $\sigma \in \operatorname{Pic}^0 X$.*

**Proof.** The Albanese variety is the dual Abelian variety $\operatorname{Alb} X = A^\vee = A$; write $0 \in A$ for the zero of the group law. On the product $A \times A$, consider the divisor

$$L = \Delta_A - A \times 0 - 0 \times A.$$

$L$ can be thought of as the universal family parametrising line bundles of degree 0 on $A$, since for each $a \in A$, restricting $\mathcal{O}_{A\times A}(L)$ to $A \times a$ gives $\mathcal{O}_A(a - 0)$.

Consider the Albanese map $\alpha\colon X \to A$ and the pullback $L' = (\alpha \times \mathrm{id})^* L$ via

$$\alpha \times \mathrm{id}\colon X \times A \to A \times A.$$

This is the universal line bundle on $X \times \operatorname{Pic}^0 X$, essentially by definition of the Albanese morphism. Note that $L^2 = -2$, so that $(L')^2$ on $X \times A$ works out as $-2F$, where $F$ is a fibre of $\alpha$.

### E.6.5. Base change

Let $p$, $q$ denote the projections of $X \times \operatorname{Pic}^0 X$ to its factors. The idea of the proof is to compute $\chi = \chi(X \times A, p^*K_X + L')$ in two different ways: first, by RR on a 3-fold and the following easy calculation

$$\begin{aligned}
\chi(X \times A, p^*K_X + L') &= \operatorname{ch}(p^*K_X + L') \cdot \operatorname{Td}_{X\times A} \\
&= \frac{1}{6}(p^*K_X + L')^3 - \frac{1}{4}K_X(p^*K_X + L')^2 \\
&= \frac{1}{4}K_X(L')^2 = 1 - f
\end{aligned}$$

(the other terms all vanish for simple reasons), one sees that $\chi = 1 - f \leq 0$, where $f = p_a F$ is the genus of a fibre of $\alpha$. Secondly, by the base change theorem, there is a complex

$$M^\cdot\colon 0 \to M^0 \xrightarrow{a} M^1 \xrightarrow{b} M^2 \to 0$$

of vector bundles on $A = \operatorname{Pic}^0 X$ such that $R^i q_*(p^*K_X + L')$ is the homology of $M^\cdot$. (To construct $M^\cdot$, take a Čech complex that computes $R^i q_*(p^*K_X + L')$, then use the fact that any bounded complex with coherent cohomology is quasi-isomorphic to a complex of vector bundles; you can truncate the complex so that $M^i = 0$ whenever $H^i = 0$ just by taking ker and coker of maps of constant rank. See [**M4**], Chapter 2, §5 or [**H1**], Chapter III.)

Now $a$ and $b$ are maps between locally free sheaves, and $a$ is injective since $H^0(K_X+\sigma)=0$ for general $\sigma\in A$. The rank of $a$ drops at $\sigma\in A$ if and only if $H^0(K_X+\sigma)\neq 0$, and that of $b$ if and only if $H^2(K_X+\sigma)\neq 0$; either of these only happen at finitely many points $\sigma\in A$, so that $R^iq_*(p^*K_X+L')$ is of finite length; write $r^i(\sigma)$ for the length of the stalk of the sheaf $R^iq_*(p^*K_X+L')$ at $\sigma$. Then the Leray spectral sequence gives

$$\chi=\sum_{\sigma\in A}-r^1(\sigma)+r^2(\sigma).$$

Now $r^2(0)\neq 0$ (since $h^2(K_X)=h^0(\mathcal{O}_X)=1$), and so $\chi\leq 0$ implies $r^1(\sigma)\neq 0$ for some $\sigma\in A$; then the rank of $a$ drops there, as required. Q.E.D.

### E.6.6. Enriques' argument

If the Albanese map $\alpha\colon X\to \operatorname{Alb} X$ has fibre genus $f=1$ then the existence of an elliptic pencil is established. At this point Enriques (see [**M3**], p. 331) argues that if $f\geq 2$ and $F=\alpha^{-1}P$ is a general fibre then the linear system $|2K_X+F|$ contains an element vanishing on some other fibre $F'$; hence $|2K_X+\sigma|\neq\emptyset$ where $\sigma=F-F'\in\operatorname{Pic}^0 X$. Enriques' argument is more geometric and much more picturesque, but it depends on one further dichotomy, and the result it proves directly is slightly weaker than Proposition E.6.4.

### Step II, a 0-curve moves

**E.6.7. Lemma.** *Let $E=\sum n_iC_i$ be a 0-curve and $D$ a divisor with $DC_i=0$ for all $i$; suppose that $n$, $b>0$. Then the exact sequence*

$$0\to\mathcal{O}_{(n+b-1)E}(D+(b-1)E)\to\mathcal{O}_{(n+b)E}(D+bE)\to\mathcal{O}_E(D+bE)\to 0$$

*gives rise to an inclusion $H^0(\mathcal{O}_{(n+b-1)E}(D+(b-1)E))\subset H^0(\mathcal{O}_{(n+b)E}(D+bE))$, and*

$$\text{this inclusion is strict}\quad\iff\quad \mathcal{O}_{(n+b)E}(D+bE)\cong\mathcal{O}_{(n+b)E}.$$

**Proof.** If there exists a section $s$ having nonzero restriction to $E$, then $s$ is a global basis of $\mathcal{O}_{(n+b)E}(D+bE)$ by E.6.1, (b). Q.E.D.

**E.6.8. Corollary.** *If the inclusion $H^0(\mathcal{O}_{nE}(D))\subsetneq H^0(\mathcal{O}_{(n+b)E}(D+bE))$ is strict then there exists $b'$ with $0<b'\leq b$ such that*

$$\mathcal{O}_{(n+b')E}(D+b'E)\cong\mathcal{O}_{(n+b')E}.$$

**Proof.** At least one of the inclusions in the chain

$$H^0(\mathcal{O}_{nE}(D))\subset\cdots\subset H^0(\mathcal{O}_{(n'-1)E}(D+(n'-n-1)E))\subset H^0(\mathcal{O}_{n'E}(D+(n'-n)E))\subset\cdots$$

must be strict, so that this follows from Lemma E.6.7. Q.E.D.

**E.6.9. Theorem.** *Let $X$ be a surface with $K_X$ nef and $E$ a 0-curve on $X$.*

*(i) For every $\mu \geq 2$,*

$$H^0(\mathcal{O}_X(\mu(K_X+E))) \neq 0.$$

*(ii) There exist effective divisors $D_1$, $D_2$ disjoint from $E$ such that*

$$\mu K_X + \nu E \overset{\text{lin}}{\sim} D_1 - D_2.$$

*with $\mu, \nu \in \mathbb{Z}$ and $\mu > \nu$, $\mu > 0$. In particular*

$$\mathcal{O}_{nE}(\mu K_X + \nu E) \cong \mathcal{O}_{nE} \quad \textit{for all } n \geq 1.$$

*(iii) The main point: (ii) implies that*

$$H^0(\mathcal{O}_{nE}(K_X + nE)) \to \infty \quad \textit{as} \quad n \to \infty;$$

*therefore eventually $h^0(\mathcal{O}_{nE}(K_X+nE)) \geq h^1(\mathcal{O}_X) + 2$, and hence by cohomology*

$$\dim |K_X + nE| \geq 1.$$

**E.6.10. Proof.** (i) Consider the exact sequence

$$0 \to \mathcal{O}_X(\mu K_X + (\mu-1)E) \to \mathcal{O}_X(\mu(K_X+E)) \to \mathcal{O}_E \to 0$$

(where I substitute $\mathcal{O}_E(K_X+E) = \omega_E \cong \mathcal{O}_E$ in the last term). Since $H^1(\mathcal{O}_E) \neq 0$ and $H^2(\mathcal{O}_X(\mu K_X + (\mu-1)E)) = 0$ for $\mu \geq 2$ (by duality, using the fact that $K_X$ is nef), it follows that $h^1(\mathcal{O}_X(\mu(K_X+E))) \geq 1$, and therefore

$$h^0(\mathcal{O}_X(\mu(K_X+E))) = \chi(\mathcal{O}_X) + h^1(\mathcal{O}_X(\mu(K_X+E))) \geq \chi(\mathcal{O}_X) + 1.$$

Now if $p_g \neq 0$ then $\mu(K_X+E) = \mu K_X + \mu E$ is obviously effective for $\mu \geq 1$, and there is no problem. If $p_g = 0$ then by (3) and (5) of Famous Table E.2, $\chi(\mathcal{O}_X) \geq 0$. Thus in any case $|\mu(K_X+E)| \neq \emptyset$. Q.E.D.

(ii) Let $D \in |\mu(K_X+E)|$; then by E.6.1, (c),

$$D = rE + D_1,$$

with $r \geq 0$ and $D_1$ disjoint from $E$. Now I have a dichotomy:
**Case** $r > 0$**.** Then (ii) is satisfied with $\nu = \mu - r$.
**Case** $r = 0$**.** Then $D_1 = D \neq 0$, so that it contains a 0-curve $E_1$ disjoint from $E$; in this case apply (i) to $E_1$ to get $D' \in |\mu(K_X+E_1)|$; if

$$D' = sE + D_2$$

with $s \geq 0$ and $D_2$ disjoint from $E$ then

$$\mu K_X - sE \overset{\text{lin}}{\sim} D_2 - \mu E_1,$$

which proves (ii) in this case. Q.E.D.

**E.6.11. Proof of** (ii) $\implies$ (iii). Let $E$ be a 0-curve. By Lemma E.6.7, in the chain

$$0 \subset H^0(\mathcal{O}_E(K_X+E)) \subset \cdots \subset H^0(\mathcal{O}_{(n-1)E}(K_X+(n-1)E)) \\ \subset H^0(\mathcal{O}_{nE}(K_X+nE)) \subset \cdots,$$

the inclusion from $(n-1)$ up to $n$ is strict if and only if $\mathcal{O}_{nE}(K_X+nE) \cong \mathcal{O}_{nE}$; thus (iii) says exactly that this coincidence occurs infinitely often. It first occurs when $n=1$, since $\omega_E = \mathcal{O}_E(K_X+E) \cong \mathcal{O}_E$. It's clearly enough to prove the following:

**E.6.12. Claim.** *Suppose that $\mathcal{O}_{nE}(K_X+nE) \cong \mathcal{O}_{nE}$ for some $n \geq 1$. Then*

$$\mathcal{O}_{n'E}(K_X+n'E) \cong \mathcal{O}_{n'E}$$

*for some $n' > n$.*

To prove this, note that by (ii), also $\mathcal{O}_{nE}(\mu K_X+\nu E) \cong \mathcal{O}_{nE}$. Thus

$$\mathcal{O}_{nE} \cong \mathcal{O}_{nE}(\mu(K_X+nE)-(\mu K_X+\nu E)) \cong \mathcal{O}_{nE}(bE), \quad \text{where } b=\mu n-\nu>0.$$

(notice that $\mu > \nu$ is needed for the case $n=1$). The claim therefore follows from the next lemma.

**E.6.13. Lemma.** *Let $E$ be a 0-curve and $n, b > 0$. Suppose that $\mathcal{O}_{nE}(bE) \cong \mathcal{O}_{nE}$; then there exists $n'$ with $n < n' \leq n+b$ such that*

$$\omega_{n'E} = \mathcal{O}_{n'E}(K_X+n'E) \cong \mathcal{O}_{n'E}.$$

**Proof.** The restriction map

$$H^0(\mathcal{O}_{(n+b)E}) \to H^0(\mathcal{O}_{bE})$$

is nonzero (because of the constant sections). Therefore I have a strict inclusion

$$H^0(\mathcal{O}_{nE}(-bE)) \subsetneq H^0(\mathcal{O}_{(n+b)E}),$$

and in particular $h^0(\mathcal{O}_{nE}(-bE)) < h^0(\mathcal{O}_{(n+b)E})$. By RR and duality,

$$h^0(\mathcal{O}_{nE}(-bE)) = h^1(\mathcal{O}_{nE}(-bE)) = h^0(\omega_{nE}(bE)) \\ h^0(\mathcal{O}_{(n+b)E}) = h^1(\mathcal{O}_{(n+b)E}) = h^0(\omega_{(n+b)E}),$$

so that the preceding inequality gives $h^0(\omega_{nE}(bE)) < h^0(\omega_{(n+b)E})$. Now using the assumption $\mathcal{O}_{nE}(bE) \cong \mathcal{O}_{nE}$, I conclude that the inclusion

$$H^0(\mathcal{O}_{nE}(K_X+nE)) \subset H^0(\mathcal{O}_{(n+b)E}(K_X+(n+b)E)$$

is also strict. The lemma therefore follows by Corollary E.6.8. Q.E.D.

This completes the proof of Theorem E.6.9.

**E.6.14. Remarks**

(i) Every step of the proof is absolutely trivial, so what's going on? The subtle point is the curious argument in (E.6.13)

strict inclusion of vector spaces
$\implies$ strict numerical inequality
$\implies$ (by RR and duality) another strict numerical inequality
$\implies$ another strict inclusion of vector spaces;

thus the bald statement that $h^0(\mathcal{O}_{nE}(-bE)) = h^0(\omega_{nE}(bE))$ conceals the fact that both $H^0(\mathcal{O}_{nE}(-bE))$ and $H^0(\omega_{nE}(bE))$ are modules over $H^0(\mathcal{O}_{nE})$, and the module structures contain nontrivial information about the nilpotents of the scheme $nE$. Hence the essence of the proof is using the numerical fact given by RR and duality to relate the two different nilpotent structures of the sheaves $\mathcal{O}_{(n+b)E}$ and $\omega_{(n+b)E}$.

(ii) The above argument is adapted from [**M3**], pp. 334–5. Mumford uses only $\mu = 1$ and applies duality 4 or 6 times more than necessary (for example, he starts his argument by saying $H^2(\mathcal{O}_X) \to H^2(\mathcal{O}_X(rE))$ is not injective, which can more simply be read $|K_X| \ni D = rE + D'$ with $D'$ disjoint from $E$).

## E.7. Effective growth of plurigenera

The remaining issues in Main Theorem E.1 are the effective results on the $P_m$, and the finer structure results, especially those for $\nu = 0$. As remarked in E.4, if $p_g \geq 2$ or if $\chi(\mathcal{O}_X) \geq 1$ then everything has already been proved. The structure results are treated in E.8. The remaining assertions divide up as follows:

**E.7.1. Theorem** (see especially [**B–M2**]).

1. *Case* $\nu = 0$. *Then* $mK \overset{\text{lin}}{\sim} 0$ *for some* $m \in \{1, 2, 3, 4, 6\}$.
2. *Case* $\nu = 1$. *Then* $P_m \geq 1$ *for some* $m \in \{1, 2, 3, 4, 6\}$.
3. *Case* $\nu = 1$. $P_m \geq 2$ *for some* $m \leq 42$.

**Discussion of proof of 1 and 2.** There is nothing to prove if $p_g = 1$, so we are in Case (5) of Famous Table E.2, that is, $p_g = 0$, $q = q' = 1$. The proof works in either case by applying the adjunction formula (canonical bundle formula) to a suitably chosen elliptic fibre space $X \to B$. You are almost certainly guessing that I'm going to use the Albanese morphism $\alpha\colon X \to A$, which is a map to an elliptic curve; however, it is a striking observation that *this never works*. Indeed, in the case $K_X \overset{\text{lin}}{\sim} 0$ the adjunction formula gives $K_X = \alpha^* L$ with $\deg L = 0$, but at this stage of the argument there's no reason why $L$ should be torsion, so there is no conclusion to be drawn. In the case $\nu = 1$, we don't know that the Albanese map is an elliptic fibre space (in fact it never is).

In the case $\nu = 1$, I use the elliptic fibre space $X \to B$ provided by abundance (E.6), and in the case $\nu = 0$, the following result, which is the first step in the analysis of bielliptic and quasi-bielliptic surfaces.

**E.7.2. Theorem.** *If* $K_X \overset{\text{num}}{\sim} 0$ *and* $p_g = 0$, $q = q' = 1$ *then* $X$ *has an elliptic fibration* $\varphi\colon X \to \mathbb{P}^1$ different from the Albanese morphism.

**Proof.** Given $B_2 = 2$ and the structure of the inner product on $\mathrm{NS}(X)$, if $F$ is the fibre of $\alpha$, it's not hard to find a class $E$ with $EF > 0$ and $E^2 = 0$. Then, because $\chi(\mathcal{O}_X) = EK_X = E^2$, arguing exactly as in E.6.4–5 gives that

$$\chi(X \times A, p^*E + L') = \frac{1}{6}(p^*E + L')^3 = \frac{1}{2}p^*E(L')^2 = -EF,$$

so that the Leray spectral sequence implies $H^0(X, \mathcal{O}_X(E+\sigma)) \neq 0$ for some element $\sigma \in \operatorname{Pic}^0 X$. Thus $E$ is numerically equivalent to an effective divisor. Now some multiple of $E$ moves in a pencil, by Theorem E.6.9. It necessarily maps to $\mathbb{P}^1$, since otherwise $q \geq 2$. Q.E.D.

### E.7.3. Vertical divisors on $X$ and fractional divisors on $B$

The proof of Theorem E.7.1, 1–2 is based on the adjunction formula for the canonical bundle of an elliptic fibration $\varphi\colon X \to B$. For clarity, I start with the *tame* case, which is defined by the condition that $R^1\varphi_*\mathcal{O}_X$ is torsion free. Then $R^1\varphi_*\mathcal{O}_X = L^{-1}$, where $L$ is a line bundle with $\deg L = \chi(\mathcal{O}_X)$, and the adjunction formula for $\varphi$ is simply Kodaira's canonical bundle formula

$$K_X = \varphi^*(K_B + L) + \sum_i (m_i - 1)E_i, \qquad (*)$$

where $F_i = m_iE_i$ are all the multiple fibres (see for example [**3 authors**], Theorem V.12.1). In all our cases, $\chi(\mathcal{O}_X) = 0$, so that $\deg L = 0$ in $(*)$.

All the calculations in what follows take place at the level of fractional divisors on the curve $B$. To explain this, I write a divisor on $X$ of the form $f^*D + \sum a_iE_i$ as the pullback $f^*\Delta$ of the fractional divisor on $B$ of the form $\Delta = D + \sum(a_i/m_i)Q_i$. If $B = \mathbb{P}^1$ and $\deg D = d$ then $D \overset{\text{lin}}{\sim} dQ$ for a chosen general point $Q \in \mathbb{P}^1$, so that I write $\Delta = dQ + \sum(a_i/m_i)Q_i$. Note that the divisorial sheaf associated with a $\mathbb{Q}$-divisor is defined by $\mathcal{O}_B(m\Delta) = \mathcal{O}_B([m\Delta])$, where $[\ ]$ is the integral part or round-down. On $B$, this says that rational functions cannot use up a fractional allocation of poles, and on $X$, it corresponds exactly to the fact that only multiples of $m_iE_i$ can move in a linear system made up of fibres.

You'll be able to follow the argument much better if you work out the details of the following example. (Compare Katsura and Ueno [**K–U**].)

### E.7.4. Champion!

The worst case, when the estimates $P_6 \geq 1$ and $P_{42} \geq 2$ are best possible, arises as follows: let $X \to \mathbb{P}^1$ be an elliptic surface having trivial Jacobian fibration (that is $J \cong E \times \mathbb{P}^1$ with $E$ an elliptic curve), and having multiple fibres $F_2 = 2E_2$, $F_3 = 3E_3$, $F_7 = 7E_7$. Then just as in the example of E.4,

$$K_X = f^*\mathcal{O}_{\mathbb{P}^1}(-2) + E_2 + 2E_3 + 6E_7 = f^*\Delta,$$

where $\Delta$ is the $\mathbb{Q}$-divisor

$$\Delta = -2Q + (1/2)Q_2 + (2/3)Q_3 + (6/7)Q_7$$

on $\mathbb{P}^1$. Then $H^0(mK_X) = H^0(\mathcal{O}_{\mathbb{P}^1}(m\Delta))$, and the canonical ring of $X$ is the graded ring $R(\mathbb{P}^1, D)$. An elementary (but fairly long) calculation shows that this ring is the graded ring $k[x,y,z]/(f)$, where $x$, $y$, $z$ have weights 6, 14, 21 and $f = x^7 + y^3 + z^2$. In particular $P_i = 0$ for all $i \leq 5$, $P_i \leq 1$ for all $i \leq 41$. (Also $P_{43} = 0$ and $P_{85} = 1$. This is responsible for the slightly obscure statement that $P_m \geq 2$ for all $m \geq 86$ that sometimes appears in the literature.)

**E.7.5. Proof of 1 and 2, the tame case.** For the proof of 1, I apply (∗) to the alternative fibration $\varphi\colon X \to \mathbb{P}^1$ obtained in Theorem E.7.2. Then $K_X$ is the pullback of the fractional divisor

$$\Delta = -2Q + \sum \frac{m_i - 1}{m_i} Q_i \quad \text{on } \mathbb{P}^1,$$

where the $Q_i \in \mathbb{P}^1$ are distinct points, and each $m_i \geq 2$. Clearly $K_X \overset{\text{num}}{\sim} 0$ if and only if the fractional divisor has degree 0, that is,

$$-2 + \sum_i \frac{m_i - 1}{m_i} = 0 \quad \text{or} \quad \sum \frac{1}{m_i} = (\text{number of } i) - 2$$

This obviously has exactly 4 solutions

$$(2,2,2,2), \quad (3,3,3), \quad (2,4,4), \quad (2,3,6).$$

Thus $mK_X \overset{\text{lin}}{\sim} 0$ for $m$ respectively $2,3,4,6$. This proves 1.

To prove 2, I apply (∗) to the elliptic or quasielliptic fibre space $X \to B$ provided by abundance, that is, the results of E.6. Thus $K_X$ is the pullback of $\Delta = K_B + L + \sum((m_i-1)/m_i)Q_i$ and from $\nu = 1$ it follows that $\deg \Delta > 0$. The assumptions of Case (5) imply immediately that $\deg L = \chi(\mathcal{O}_X) = 0$ and $g(B) \leq 1$.

Because $K_X$ nef and nonzero, if $g(B) = 1$ there must be at least one multiple fibre; then every term $2((m_i-1)/m_i)Q_i$ makes a positive contribution to $[2\Delta]$, so that $\deg[2\Delta] > 0$, and RR on $B$ gives $P_2 = H^0(\mathcal{O}_B([2\Delta]) \geq 1$. (Similarly, $P_3 > 0$ and $P_m \geq 2$ for every $m \geq 4$. Here and below, results like this leak out for free, and are used in the proof of Theorem E.7.1, 3).

The other possibility is that $B \cong \mathbb{P}^1$. Then $K_X$ is the pullback of $\Delta = -2Q + \sum(m_i - 1/m_i)Q_i$, and $K_X$ nef and nonzero gives

$$-2 + \sum_i \frac{m_i - 1}{m_i} > 0, \quad \text{that is,} \quad \sum \frac{1}{m_i} < (\text{number of } i) - 2$$

Now every point $Q_i$ makes a positive integral contribution to $2\Delta$, so that $P_2 \geq 2$ if there are at least 5 points $Q_i$. If there are 4 points, at least one must be $\geq 3$, and $P_2 \geq 1$, $P_4 \geq 2$. The final part of the proof is a pleasurable exercise. [Hint: If all $m_i \geq 3$, the smallest case is $(3,3,4)$, giving

$$3K_X \overset{\text{lin}}{\sim} \varphi^*(-6Q + 2Q_1 + 2Q_2 + (9/4)Q_3) \overset{\text{lin}}{\sim} E_3 = (1/4)F_3,$$

so $P_3 = 1$ and $P_{12} = 2$. Otherwise, $m_1 = 2$, and if both $m_2, m_3 \geq 4$ the smallest case is $(2,4,5)$, giving

$$4K_X \overset{\text{lin}}{\sim} \varphi^*(-8Q + 2Q_1 + 3Q_2 + (16/5)Q_3) \overset{\text{lin}}{\sim} E_3 = (1/5)F_3,$$

so that $P_4 = 1$ and $P_{20} = 2$. Otherwise, $m_1 = 2, m_2 = 3$, and $m_3 \geq 7$. This leads to Champion E.7.4.]

**E.7.6. Exercise.** Prove that $P_{12} \leq 1$ happens only in 7 cases:

$$(2,5,5) \quad (2,4,5) \quad \text{and} \quad (2,3,m) \text{ with } m \leq 11.$$

### E.7.7. Wild!

Once the terminology and basic properties are established, the proofs in the wild case are essentially the same combinatorial arguments as in E.7.3–6.

Let $\varphi\colon X \to B$ be an elliptic or quasielliptic fibration and $F = mE = \varphi^*P$ a multiple fibre. Following [**B–M2**], I introduce the following three numerical characteristics of $F$: the *multiplicity* $m$, the *order* $n$ of $\mathcal{O}_E(E)$ as a torsion element of $\operatorname{Pic}^0 E$, and the *local canonical class*, that is, the integer $a < n$ such that $\mathcal{O}_X(K_X) \cong \mathcal{O}_X(aE)$ in a neighbourhood of $F$.

In characteristic zero, automatically $a = m-1$ and $n = m$. Indeed, if $t_P$ is a local parameter at $P$ on $B$, then locally near $E$, I have $\operatorname{div}(t_P) = mE$ and $\operatorname{div}(\mathrm{d}t_P) = (m-1)E$. Also taking the $m$th root of $t_P$ and normalising gives rise to an etale cyclic cover of a neighbourhood of $F$ such that the inverse image of $F$ is connected. On the other hand, the restriction to $E$ is the cyclic cover defined by $\mathcal{O}_E(E)$, hence this has order $m$.

**Remark.** The argument here is essentially topological, and can be viewed as a beautiful illustration of Kollár's philosophy that *vanishing* is when a coherent cohomology group has a topological interpretation: in this case, $H^0(\mathcal{O}_E(iE)) = 0$ for all $0 < i < m$ because the inverse image of $E$ in a cyclic cover is connected.

In characteristic $p$, the same argument applies only to the part of $m$ coprime to $p$: if $m = m'p^a$ with $m'$ not divisible by $p$ then the $m'$th root of $t_P$ still defines an etale cyclic cover in which $F$ remains connected; therefore the torsion order $n$ is given by $n = m'p^b$ for some $b \le a$.

Quite generally, a multiple fibre of $X \to B$ can be reduced to a nonmultiple fibre by taking the normalised pullback by a suitable Galois separable ramified cover $B' \to B$ of the base curve, which induces an etale cover of a neighbourhood of $F \subset X$. However, in characteristic $p$, a ramified cover can have a complicated $p$-group as local Galois group (ramification group), and in any case, the typical cyclic cover is of the form $x^p + x = a$ (an *Artin–Schreier extension*).

**Definition.** The fibre $F$ is *tame* if $n = m$ and *wild* if $n < m$. Obviously since both $\mathcal{O}_E(F) = \mathcal{O}_E(mE)$ and $\omega_E = \mathcal{O}_E((a+1)E)$ are trivial, it follows that $n$ divides $m$ and $a+1$.

If $F$ is a wild fibre, then clearly $h^0(\mathcal{O}_F) \ge 2$, so that also

$$h^1(\mathcal{O}_F) > 1.$$

In particular, the stalk of $R^1\varphi_*\mathcal{O}_X$ at $P$ needs $\ge 2$ generators, and it follows that $R^1\varphi_*\mathcal{O}_X = T \oplus L^{-1}$, where $T$ is a torsion sheaf, nonzero at each wild fibre, and $L$ is a line bundle with $\deg L = \chi(\mathcal{O}_X) + \operatorname{length}(T)$. (This follows simply by RR on $B$.) The adjunction formula in the general case is thus

$$K_X = \varphi^*(K_B + L) + \sum a_i E_i, \qquad (**)$$

where $F_i = m_i E_i$ are the multiple fibres, and $a_i$ their local canonical invariants.

**E.7.8. Lemma.** *If* $h^1(\mathcal{O}_X) = 1$ *then* $a = m-1$ *or* $a = m-n-1$.

**Proof.** Since $\mathcal{O}_E(nE) = \mathcal{O}_E$, it follows from Lemma E.6.13, (ii) that

$$h^0(\mathcal{O}_{n'E}(K + n'E)) \ge 2 \quad \text{for some } n' \text{ with } 1 < n' \le n+1.$$

Then by cohomology, $H^0(K_X) \subsetneq H^0(K_X + n'E)$, so that some part of the divisor $(a+n')E$ must move off $E$, that is, $m \le a+n' \le a+n+1$. Thus $n \ge m-(a+1)$, but on the other hand, $n$ divides $m$ and $a+1$. Therefore, $m-(a+1) = 0$ or $n$, as asserted. Q.E.D.

**E.7.9. Proof of 1 and 2 in the wild case.** As before, to prove 1, I work on the alternative fibration $\varphi\colon X \to \mathbb{P}^1$ provided by Theorem E.7.2. In $(**)$ I have $\deg L = \chi(\mathcal{O}_X) + \operatorname{length}(T)$. If $\deg L \ge 2$ then $(**)$ would give $p_g \ge 1$, a contradiction, so it follows that $\deg L = 1$ and there is exactly one wild fibre with $\operatorname{length}(T) = 1$. Set $F_0 = m_0E_0$ for the wild fibre, $n_0$ and $a_0$ for its period and local canonical invariant, and $F_i = m_iE_i$ for $i = 1, \dots, k$ for the tame multiple fibres. Now argue exactly as in E.7.4: $K_X$ is the pullback of the fractional divisor

$$-Q + \frac{a_0}{m_0}Q_0 + \sum_{i=1}^{k} \frac{m_i - 1}{m_i}Q_i.$$

Moreover $a_0 = m_0 - 1$ or $m_0 - n_0 - 1$. In the first case it is an easy exercise to deduce that $K_X \overset{\text{num}}{\sim} 0$ gives $k = 1$ and $(m_0, m_1) = (2, 2)$ so that $P_2 = 1$.

For the second case, it's useful to note that the wild fibre behaves for numerical purposes exactly like two usual multiple fibres of multiplicity $m_0$ and $m_0/n_0$: to see this, it is enough to write the above $\mathbb{Q}$-divisor throughout in the form

$$-2Q + \frac{m_0 - 1}{m_0}Q_0 + \frac{m_0/n_0 - 1}{m_0/n_0}Q_0 + \sum_{i=1}^{k} \frac{m_i - 1}{m_i}Q_i.$$

Thus $K_X \overset{\text{num}}{\sim} 0$ gives

$$k = \frac{1}{m_0} + \frac{n_0}{m_0} + \sum_{i=1}^{k} \frac{1}{m_i}.$$

The right-hand side is a sum of reciprocal integers, and an obvious calculation shows that the set of integers $(m_0, m_0/n_0, \{m_i\}_{i=1}^k)$ is one of the usual list:

$$(2,2,2,2) \quad (3,3,3) \quad (2,4,4) \quad (2,3,6).$$

Thus $mK_X \overset{\text{lin}}{\sim} 0$ for $m = 2, 3, 4$ or $6$. This proves 1.

For the purposes of this proof, I didn't need to figure out which of the integers is $m_0/n_0$. Note that since $n_0 = p^a m_0$, where $p = \operatorname{char} k$, the argument shows that for a surface with $\nu = 0$ in Case (5), a wild fibre can only happen in characteristic 2 or 3.

As before, for 2, I work on the pluricanonical elliptic fibration $X \to B$ provided by abundance. The argument is almost identical to previous work. In view of $p_g = 0$, $q = q' = 1$, the adjunction formula $(**)$ gives that $B = \mathbb{P}^1$, $\deg L = 1$ and there is exactly one wild fibre to which Lemma E.7.8 applies. Now $K_X$ is the pullback of the fractional divisor

$$-Q + \frac{a_0}{m_0}Q_0 + \sum_{i=1}^{k} \frac{m_i - 1}{m_i}Q_i.$$

where $a_0 = m_0 - 1$ or $m_0 - n_0 - 1$. In the first case I have

$$-1 + \frac{m_0 - 1}{m_0} + \sum_{i=1}^{k} \frac{m_i - 1}{m_i} > 0 \quad \text{that is,} \quad \sum_{i=0}^{k} \frac{1}{m_i} < k$$

Thus $k \geq 1$. If $k \geq 2$ then $P_2 \geq 2$; otherwise $k = 1$ and one of $m_0, m_1 \geq 3$, so that $P_2 \geq 1$, $P_3 \geq 1$ and $P_6 \geq 2$ (compare Example E.4). In the second case I get

$$-1 + \frac{m_0 - n_0 - 1}{m_0} + \sum_{i=1}^{k} \frac{m_i - 1}{m_i} > 0 \quad \text{that is,} \quad \frac{1}{m_0} + \frac{n_0}{m_0} + \sum_{i=1}^{k} \frac{1}{m_i} < k$$

Here again $k \geq 1$, and the pleasurable exercise at the end of E.7.5 can be repeated verbatim to prove that $P_m \geq 1$ for some $m \leq 6$ and $P_m \geq 2$ for some $m \leq 42$.

**E.7.10. Proof of Theorem E.7.1, 3.** If $p_g \geq 2$ there is nothing to prove. If $\chi(\mathcal{O}_X) \geq 1$ and $\nu = 1$ then $P_6 \geq 2$ was proved in Theorem E.4. The remaining case with $p_g = 0$ is Case (5), and everything has been proved in E.7.5 in the tame case, and in E.7.9 in the wild case.

The remaining cases are (4), (5′), which have $p_g = 1$ and $q' = 2$. In this case, the proof consists of using the construction of E.6.9–13 as an *effective* method. Write $K_X = rE + D'$ with $E$ a 0-curve, $r \geq 1$ and $D' > 0$ a divisor disjoint from $E$. The arguments of E.6.9–13 give $H^0(\mathcal{O}_{(r+2)E}) \geq 2$, and $H^0(\mathcal{O}_{(r+2)E+D'}) \geq 3$ if $D' \neq 0$, so that by RR and duality $H^0(\mathcal{O}_{(r+2)E+D'}(K_X + (r+2)E + D')) \geq 3$, therefore, by the cohomology long exact sequence, $|K_X + (r+2)E + D'|$ moves. In the worst case $r = 1$, this gives $4K_X > K_X + (r+2)E + D'$, so that $P_4 \geq 2$. If $D' = 0$, then I take the argument of Lemma E.6.13 one step further: then $H^0(\mathcal{O}_{bE}(r+b)E) \cong \mathcal{O}_{bE}$ for some $b \leq r + 2$, so that $H^0(\mathcal{O}_{(r+2b)E}) \geq 3$ and $|K_X + (r + 2b)E|$ moves. Thus in the worst case $D' = 0$, $r = 1$ and $b = 3$, I get that $|8E|$ moves, so that $P_m \geq 2$ for all $m \geq 8$. This completes the proof of Theorem E.7.1. Q.E.D.

**E.7.11. Remark.** I'm rather disappointed that the arguments of E.6.9–13 don't seem to give the effective bound $P_m \geq 42$ also in the cases $p_g = 0$. The interest of the question is that if you could find a better argument here, you might be able to tidy up the main proof of E.6.9–13.

In this case $q = q' \leq 1$, and you could try writing $mK_X = rE + D'$ for some $m = \{1, 2, 3, 4, 6\}$, with $D'$ effective and disjoint from $E$. Thus

$$\mathcal{O}_E(mK_X - rE) \cong \mathcal{O}_E \cong \mathcal{O}_E(K_X + E) \implies \mathcal{O}_E((m+r)E) \cong \mathcal{O}_E.$$

Now arguing as in Theorem E.6.9, (i), there exists $c \leq m + r + 1$ such that $H^0(\mathcal{O}_{cE}) \geq 2$ and $\mathcal{O}_{cE}(K_X + cE) \cong \mathcal{O}_{cE}$. If $D' \neq 0$, you get as before a curve $cE + D'$ with $H^0(\mathcal{O}_{cE+D'}) \geq 3$, and $|K_X + cE + D'|$ moves by cohomology. In the worst case $m = 6$, $r = 1$, $c = 8$, this gives $P_m \geq 2$ for $m \geq 49$. If $D' = 0$, the next step is to say that $\mathcal{O}_{cE}((mc + r)E) \cong \mathcal{O}_{cE}$, therefore $H^0(\mathcal{O}_{((m+1)c+r)E}) \geq 3$, so that $|K_X + ((m+1)c + r)E|$ moves. In the worst case $m = 6$, $r = 1$, $c = 8$, so that $|K_X + 57E| = |343K_X|$ moves (notice that $343 = 7^3$). Of course, this is far from best possible. The way to improve it would be to show that $\mathcal{O}_{8E}(7E) \cong \mathcal{O}_{8E}$, that is, replace 49 by 7. Probably you can show that the kernel of $\operatorname{Pic}^0(8E) \to \operatorname{Pic}^0 E$ is a vector space by the same kind of arguments as in 4.13, so has no torsion (at least in characteristic $\neq 7$)?

## E.8. Abelian, bielliptic and sesquielliptic surfaces

**E.8.1. Theorem** (characterisation of Abelian surfaces). *Let $X$ be a surface with $K_X \overset{\text{lin}}{\sim} 0$, $p_g = 1$, $q = q' = 2$ (that is, $\nu = 0$ and (4)). Then the Albanese map $\alpha\colon X \to \operatorname{Alb} X$ is an isomorphism, so that $X$ is an Abelian surface.*

**E.8.2. Proof.** Write $\alpha\colon X \to A = \operatorname{Alb} X$ for the Albanese map. Recall that it has the universal mapping property (UMP) for morphisms of $X$ to an Abelian variety. I exclude the possibility that $\alpha$ maps to a curve $C \subset A$ by saying that $C$ has genus 2, and thus has etale covers of large genus; hence, by pullback, $X$ has etale covers with $q > 2$, therefore $p_g > 1$, which contradicts $K_X = 0$.

Consider the differential $\mathrm{d}\alpha\colon T_X \to \alpha^* T_A$. In characteristic zero, it is generically an isomorphism, so gives an injective map

$$T_X \hookrightarrow \alpha^* T_A \cong \mathcal{O}_X^{\oplus 2},$$

necessarily an isomorphism because $K_X \overset{\text{lin}}{\sim} 0$; thus $\alpha$ is etale. But an etale cover of an Abelian variety is itself Abelian, so that $X = A$ follows from the UMP of $\alpha\colon X \to A$.

In characteristic $p$, the same proof goes through unchanged if $\mathrm{d}\alpha$ is generically injective. On the other hand, $\mathrm{d}\alpha$ is not generically zero. Indeed, $\mathrm{d}\alpha = 0$ means that $\alpha^* k(A) \subset k(X)^p$ (the subfield of $p$th powers of elements of $k(X)$), or equivalently, $k(A)^{1/p} \subset k(X)$. This means that $\alpha$ factors through the geometric Frobenius map $A^{(-1)} \to A$, which again contradicts its UMP.[7]

Therefore I can assume that $\mathrm{d}\alpha$ has generic rank 1, and argue on foliations as in Rudakov and Shafarevich [**R**–**S**] and [**Ekedahl**]. The foliation $\ker \mathrm{d}\alpha \subset T_X$ is a rank 1 subsheaf, and is necessarily saturated; write $\ker \mathrm{d}\alpha = \mathcal{O}_X(D)$. Then there is an exact sequence

$$0 \to \mathcal{O}_X(D) \to T_X \to \mathcal{I}_Z \mathcal{O}_X(-D) \to 0,$$

where $\mathcal{I}_Z$ is the ideal sheaf of a zero dimensional subscheme, and $D \geq 0$, since by assumption there is a nonzero homomorphism $\mathcal{I}_Z \mathcal{O}_X(-D) \to \alpha^* T_A \cong \mathcal{O}_X^{\oplus 2}$.

Consider the subfield of the function field $k(X)$ generated by $k(X)^p$ and $k(A)$, that is, the composite field $K = k(A)k(X)^p$. Then since $[k(X) : k(X)^p] = p^2$ and $k(A) \not\subset k(X)^p$, I deduce that $K \subset k(X)$ is inseparable of degree $p$. It follows that $\alpha$ factorises as $X \to Y \to A$, where $Y$ is the surface obtained as the normalisation of $A$ in $K$, the morphism $\pi\colon X \to Y$ is inseparable of degree $p$, and $\ker \mathrm{d}\pi = \mathcal{O}_X(D) \subset T_X$. Therefore $\mathcal{O}_X(D)$ is a $p$-closed foliation and $Y = X/\mathcal{O}_X(D)$ is the corresponding quotient. (See [**R**–**Sh**] for the terminology.)

**Case** $D > 0$**.** By the adjunction formula for $X \to Y$ of [**R**–**Sh**], §2, Proposition 2 (p. 1211 of translation), $\pi^* K_Y = -(p-1)D$, so that $K_{Y'}$ cannot be nef on any birational model $Y'$ of $Y$, and $Y$ is rational or ruled by the main result of minimal model theory (Corollary D.4.4). This contradicts $Y \to A$ generically finite.

[7] Here $A^{(-1)}$ is the normalisation of $A$ in the field $k(A)^{1/p}$; it is isomorphic as a scheme to $A$ via the absolute Frobenius, so is an Abelian variety in its own right. The point of the funny notation is mainly to avoid offending experts: $A^{(-1)}$ and the morphism $A^{(-1)} \to A$ are defined over $k$, and $A^{(-1)}$ is conjugate to $A$ by the Frobenius automorphism of $k$, but is not isomorphic to $A$ over $k$.

**Case $D = 0$.** Then $Z = 0$, since $c_2(X) = 0$. By [**R–Sh**], §1, Corollary of Theorem 1 (p. 1208 of translation), the surface $Y$ is nonsingular, and $K_Y \overset{\text{num}}{\sim} 0$ by the same adjunction formula. Thus $p_g(Y) \le 1$; but since $q(Y) \ge 2$, it follows that $Y$ is also in Case (4) of Famous Table E.2, that is, $Y$ satisfies all the same hypotheses as $X$. I can therefore assume by induction on $\deg(Y/A)$ that $Y = A$, that is, $\alpha\colon X \to A$ is inseparable of degree $p$.

Moreover, by the exact sequence $0 \to \mathcal{O}_X(D) \to T_X \to \mathcal{I}_Z\mathcal{O}_X(-D) \to 0$, the differential $\mathrm{d}\alpha\colon T_X \to \alpha^*T_A \cong \mathcal{O}_X^{\oplus 2}$ has image isomorphic to $\mathcal{O}_X$, which is thus a direct summand. Thus $T_A = L_1 \oplus L_2$ is a direct sum of two trivial line bundles $L_i \cong \mathcal{O}_A$, such that $\mathrm{d}\alpha$ maps surjectively to $\alpha^*L_1$.

Now $X \to A$ is inseparable of degree $p$, so that the inclusion $k(X)^p \subset k(A)$ defines a factorisation $F\colon X \to A \to X^{(1)}$ of the Frobenius. On the other hand, since $\mathrm{d}\alpha$ has image $\alpha^*L_1$, and $\mathrm{d}F = 0$, it follows that $\beta\colon A \to X^{(1)}$ has $\ker \mathrm{d}\beta = L_1$, and $\beta$ is the quotient by $L_1$. Therefore $L_1$ is $p$-closed.

Now also $L_1 \cong \mathcal{O}_A$, and is a direct summand of $T_A$, so that $L_1$ is a $p$-Lie subalgebra of the algebra of tangent fields on $A$. The Frobenius map $F\colon A \to A^{(1)}$ is a homomorphism of algebraic group schemes, and the vector space $H^0(T_A)$ is in a natural way the $p$-Lie algebra corresponding to the finite subgroup scheme $\ker F$. From the fact that $L_1$ is $p$-closed, it follows that $H^0(L_1)$ is a $p$-Lie subalgebra, so that by the Lie correspondence ([**M4**], Chapter III, §14, Theorem), it is the Lie algebra of a subgroup scheme $G \subset \ker F$. Therefore the quotient by $L_1$ is the same thing as the quotient by $G$, and it follows that $X^{(1)}$ is itself an Abelian variety. Therefore $X$ is also an Abelian variety, and the factorisation $X \to A$ contradicts the UMP of $\alpha$ as before. Q.E.D.

**E.8.3. Remarks.** (1) The proof in [**B–M2**], pp. 40–41 is an extremely weird reduction to a finite field. My proof (based largely on suggestions of Nick Shepherd-Barron) is an expanded and simplified version of the argument of [**Ekedahl**], Proposition 4.3.

(2) It's an exercise to generalise the treatment of the final case to prove the following theorem of Igusa and Serre ([**Serre**], §2.6, Theorem 4): if $f\colon X \to A$ is a generically finite surjective morphism of an $n$-fold to an Abelian variety, the two conditions (a) and (b) are equivalent:

**(a)** There exist a 1-form $\omega \in H^0(\Omega^1_A)$ such that $f^*\omega = 0 \in \Omega^1_X$.

**(b)** $f$ factors via $B \to A$, an inseparable morphism of height 1 between Abelian varieties. (This means that $B \to A$ has degree $> 1$ and factors the Frobenius morphism $B \to A \to B^{(1)}$; or, more simply, $k(B)^p \subset k(A) \subsetneq k(B)$.)

It follows that if the Albanese map $\alpha\colon X \to \operatorname{Alb} X$ is surjective for any variety $X$, then its differential is an injective map $\alpha^*\colon H^0(\Omega^1_A) \hookrightarrow H^0(\Omega^1_X)$ .

**E.8.4. Theorem.**

1. *Case* $K_X \overset{\text{num}}{\sim} 0$, $p_g = 0$, $q = q' = 1$ *(that is,* $\nu = 0$ *and (5)). Then* $X$ *is bielliptic (see Definition E.5.2).*
2. *Case* $K_X \overset{\text{lin}}{\sim} 0$, $p_g = 1$, $q = 1$, $q' = 2$ *(that is,* $\nu = 0$ *and (5′)). Then* $\operatorname{char} k = 2$ *or* $3$*, and* $X$ *is quasi-bielliptic.*

Here quasi-bielliptic means that $X$ is an quasielliptic fibre space $X \to A$ over an elliptic curve, which becomes a $\mathbb{P}^1$-bundle on making an inseparable cover of the base curve.

**E.8.5. Sketch proof.** In either case (5–5′), consider the Albanese map $\alpha\colon X \to A$, and the alternative fibration $\varphi\colon X \to \mathbb{P}^1$ of Theorem E.7.2. Any fibre $E$ of $\varphi$ is a curve with $p_a E = 1$ having a finite morphism $\pi\colon E \to A$. Thus $\pi$ is necessarily an isogeny of elliptic curves (after choosing base points), so that $E$ is a nonsingular elliptic curve, and $E \to A = E/G$ is the quotient by the finite subgroup scheme $G = \ker \pi \subset E$.

Suppose first that $\alpha$ is an elliptic fibration, that is, that the general geometric fibre is nonsingular. Then the adjunction formula (∗∗) applied to $\alpha$ implies that $\deg L = 0$, and all $a_i = 0$, and therefore, every fibre $F$ is the same nonsingular elliptic curve. Now pulling back $X \to A$ by the finite cover $E \to A$ gives an elliptic fibre space $X_E = X \times_A E \to E$ having a section $E$ (the diagonal $\Delta_E \subset E \times_A E$), and with all geometric fibres the same curve $F$. Therefore $X_E$ is isomorphic to the product $E \times F$ and $X_E \to X$ is a quotient by the subgroup scheme $G \subset E$, acting on $E$ by translation, and acting somehow on $F$. However, since $X$ has $mK_X \overset{\text{lin}}{\sim} 0$ for some $m \in \{1, 2, 3, 4, 6\}$ it follows that the action of $G$ on $H^0(K_F)$ factors via the group scheme $\mu_m$ of $m$th roots of unity, where $m$ is the order of $K_X$ in $\operatorname{Pic} X$. The kernel of this action is a translation subgroup scheme of $E \times F$, and the corresponding quotient is an Abelian variety. Therefore $X$ is isomorphic to the quotient of an Abelian variety by $\mu_m$.

If $\operatorname{char} k \neq 2$ or 3, it follows that $X$ is an etale quotient of an Abelian variety by a cyclic group of order $m$. In any case, $m = 1$ would give that $X$ is an Abelian surface, and contradict the case assumption (5–5′). But all surfaces with $\nu = 0$ in Case (5′) have $p_g = 1$, therefore $K_X \overset{\text{lin}}{\sim} 0$, so that $\alpha$ cannot be an elliptic fibration for these surfaces.

Next suppose that $\alpha\colon X \to A$ is a quasielliptic fibre space. As before, $\deg L = 0$ and all $a_i = 0$ in the adjunction formula (∗∗), so that every fibre $F$ of $\alpha$ is a reduced irreducible curve with a cusp. Pulling back by the isogeny $E \to A$ gives a ruled surface with a section, so that the normalised pullback $X_E$ is isomorphic to a $\mathbb{P}^1$-bundle over $E$, and $X_E \to X$ is again a quotient by a subgroup scheme $G \subset E$. Q.E.D.

**E.8.6. Sketch proof of Theorem E.5.1 and E.5.3.** Consider a surface $X$ over $\mathbb{C}$ with $K_X$ nef, $\nu = 1$ and $c_2 = 0$. We know that $X$ is an elliptic surface $f\colon X \to B$ with $c_2 = 0$. Hence by the Euler number calculation, the only degenerate fibres are multiple nonsingular fibres, say $E_i = m_i e_i$. In particular the modular invariant of the fibre is bounded, and therefore constant. Hence the Jacobian fibration is trivial, $J = E \times B$. Now the canonical bundle formula gives $K_X = f^* K_B + \sum (m_i - 1) e_i$, and it's easy to see that $p_g(X) = g(B)$, $q(X) = g(B) + 1$. By the complete reducibility property of Abelian varieties, $\operatorname{Alb} X$ has a projection to an elliptic curve $A$ which restricts to a nonconstant morphism from the general fibre $E$ of $X \to B$. In particular, the Albanese morphism $\alpha\colon X \to \operatorname{Alb} X$ never equals the canonical elliptic fibre space $f\colon X \to B$. The normalised pullback $X \times_A E$ is then a product. But $E \to A$ must be Galois with Abelian group, so that $X$ is obtained as stated in the theorem. Q.E.D.

## E.9. Any questions or comments?

This chapter has given a detailed, complete, self-contained proof of the main results of the classification of surfaces, following Enriques' argument as rewritten by Mumford. My reason for writing this in this form is my impression that while hundreds of mathematicians (and nowadays also a growing number of theoretical physicists and specialists in science fiction) need to know and use the results of the classification, very few have ever been through a proof. Although the coincidence $\nu = \kappa$ may seem at first sight such an obvious and coarse result, its proof is logically quite intricate and beautiful. Moreover, the role of the bielliptic and sesquielliptic surfaces with $c_2 = 0$ as a logical bottleneck of the proof seems to me to deserve some emphasis.

### E.9.1. Abundance as a logical bottleneck

As far as I know, there is only one complete proof of abundance in the literature that does not pass via Enriques' argument at the central point, namely that given in [**3 authors**]. They first give Ueno's proof of Iitaka's additivity conjecture $C_{2,1}$, which uses the theory of moduli of curves. Recall that if $X$ has $p_g = 0$, $q = 1$ and $K_X^2 = 0$ then $X$ has an Albanese morphism $\alpha\colon X \to \operatorname{Alb} X = A$ to an elliptic curve $A$, with fibres of genus 1 if $\nu = 0$, or genus $\geq 1$ if $\nu = 1$. Then $C_{2,1}$ proves at once that $P_m \geq 1$ for some $m$, and $P_m > 1$ for some $m$ if the fibres of $\alpha$ have genus $\geq 2$, and with a bit more work also if $\nu = 1$.

### E.9.2. Finding cohomology when $\chi(\mathcal{F}) = 0$

The central problem in proving abundance is to find nonzero cohomology on $X$. The point is that

$$\chi(\mathcal{O}_X(nK_X)) = \chi(\mathcal{O}_X) \quad \text{for all } n,$$

so that RR on its own does not imply that $H^0(nK_X) \to \infty$ as $n \to \infty$. But if $P_m \geq 2$ for some $m$, then $D \in |mK_X|$ will consist of many fibres of the elliptic fibration $\varphi_{mK_X}$, so $H^0(\mathcal{O}_D)$ is large. It's interesting to analyse the argument of E.6–7 to note the key points in Proposition E.6.4, Theorem E.6.9 and Proposition E.7 at which one has to work to squeeze out nonzero cohomology one drop at a time.

### E.9.3. Finding curves $C$ with $KC = 0$

Abundance can be stated as saying that a suitably defined "moduli space" of maps $\varphi\colon C \to X$ with $K_X\varphi(C) = 0$ has the right dimension. Our experience with Mori theory teaches us that there is essentially only one way of proving that a nonsingular projective variety $X$ with $K_X$ not nef contains rational curves, namely Mori's bending-and-breaking argument: very roughly, a curve $C \subset X$ with $K_XC < 0$ moves in a positive dimensional family (for easy reasons), and must break off a rational component at some point (for delicate easy reasons). Put this crudely, this suggests that a similar argument might be capable of locating the curves with $K_XC = 0$ on a $n$-fold with $K_X$ nef. Unfortunately, it seems almost certain that this kind of approach cannot work.

The point is that Mori's bending-and-breaking argument actually works with the deformation theory of morphisms $\varphi\colon C \to X$ with $C$ a fixed curve; if you allow the moduli of $C$ to vary, you might gain some parameters, but you completely lose the ability to predict that the family breaks up. Although, if you believe abundance, the moduli space of curves in $X$ with $K_XC = 0$ has the right dimension, it seems impossible to approximate them by curves with $K_XC$ "small".

In the rest of this chapter I discuss three other approaches to the proof of abundance, none of which has been completely worked out (to my knowledge). Assume $k = \mathbb{C}$ for a partial sense of security.

### E.9.4. Kähler–Einstein

This is the question left at the end of E.5. Consider minimal surfaces with $K_X^2 = c_2 = 0$. By a theorem of Uhlenbeck–Yau and Donaldson, if we knew that $T_X$ were stable, it could be given a Hermitian–Einstein metric. Of course, a posteriori $T_X$ is stable only for Abelian surfaces. Nevertheless, it seems to me (as an outsider in the analytic side of differential geometry) quite likely that one can prove a priori that there exists a metric that makes the Kähler structure locally symmetric.

H. Tsuji seems to have some general results on this kind of question in higher dimensions, which I paraphrase, probably rather wildly. Try to construct a "Kähler pseudometric" on the line bundle $K_X$ whose curvature form vanishes along the fibres of the canonical map $\varphi\colon X \to Y$, and is positive in directions transverse to the fibres. Of course, a priori we don't know $\varphi$ or any of its fibres, but if $K_X^n = 0$ then the curvature of a Kähler pseudometric must vanish on a subspace of the tangent bundle. Then there is a kind of heat equation argument saying that a Kähler pseudometric should exist as a distribution. If you could prove that this exists as an analytic pseudometric then its kernel should give the fibres of $\varphi$, at least as a subbundle of $T_X$ or as a complex foliation. Then to prove that the canonical fibration exists is equivalent to proving the regularity of the Kähler pseudometric distribution. (??)

### E.9.5. Miyaoka's argument

Let $X$ be a 3-fold with $K_X$ nef and $\nu = 1$ (that is $K_X \overset{\text{num}}{\not\sim} 0$, $K_X^2 \overset{\text{num}}{\sim} 0$), and assume $\kappa \geq 0$ to avoid appealing to Miyaoka's previous hard theorems. Let $E$ be a suitable divisor chosen from the components of a divisor $D \in |mK_X|$, with $E^2 \overset{\text{num}}{\sim} 0$ and $K_{X|E} \overset{\text{num}}{\sim} 0$. Miyaoka proves that some multiple of $E$ moves in a pencil on $X$. The idea of the proof is to restrict $X$ to a tubular neighbourhood of $E$; then some easy homotopy theory shows that it is possible to make cyclic covers of $X$ branched along $E$, and then to pass to a minimal model $F \subset Y$ via Kulikov's results. This model has the following good properties: $F$ is a degenerate Abelian or $K3$ surface (that is, global normal crossings, $K_F \overset{\text{lin}}{\sim} 0$, and all components of $F$ are algebraic surfaces), $Y$ is nonsingular, and $K_{Y|F} \overset{\text{lin}}{\sim} F_{|F} \overset{\text{lin}}{\sim} 0$.

Now $F$ has a local deformation space as an abstract complex space; next, Miyaoka's idea is to prove that there are sufficiently many abstract infinitesimal deformations of $F$ that can be mapped to $Y$. (Specialists know that Miyaoka's appeal to Friedman's results in deformation theory in the case of degenerate Abelian surfaces is not valid as written, but I understand that his proof has been fixed up by slightly different arguments. Miyaoka's proof contains in any case a whole string of important new ideas. See [**Utah2**] for the current status of this.)

It seems that in characteristic 0, this argument can be translated back to the surface case without much trouble, and can be used to replace Step II, (E.6.7–13). (??)

### E.9.6. There must be a better way

It seems conceivable to me that we are just overlooking some much coarser and much more fundamental reason why abundance must be true for surfaces, an argument that doesn't resort to miserable cookery with the constant term in RR. Here is one attempt at such an argument.

Let $X$ be a surface with $\nu = 1$. Consider the set of curves $E \subset X$ with $K_X E = 0$. By the Hodge index theorem, the intersection form on these is negative semidefinite. Thus either there are only finitely many, or at most a single 1-dimensional family plus finitely many; if there is a positive dimensional family, it must be an elliptic pencil, and everything's lovely. I say that there must be a way of proving directly that there are infinitely many curves $E$ with $K_X E = 0$, or better, that there is a curve $E$ with $K_X E = 0$ through every sufficiently general point $P \in X$.

Let $H$ be a divisor on $X$ with $d = K_X H > 0$. Note that replacing $H$ by $H + nK_X$ does not change $K_X H$, so that there is no loss of generality in assuming that also $H^2 > 0$.

Now resolutely ignoring the constant term in RR gives

$$h^0(\mathcal{O}_X(H + nK_X)) \sim (n/2)d.$$

**E.9.7. Conjecture.** *There exists integers* $n \gg 0$ *and* $m < \sqrt{(H + nK_X)^2} - 1$*, and a decomposition*

$$H + nK_X = D + E \qquad \text{with } D, E > 0 \text{ and } DE \leq m.$$

*Moreover,* $E$ *can be chosen to pass through a sufficiently general point* $P$ *of* $X$*.*

The statement itself follows from abundance: just take $E$ to be the fibre of $\varphi_{mK}$ through $P$. The conjectural part is that this statement can be proved directly without appealing to abundance.

**E.9.8. Lemma.** *The conjecture implies abundance.*

**Proof.** The form $D^2E^2 \leq (DE)^2$ of the Hodge index theorem implies at once that $D^2$ or $E^2 \leq 0$, since otherwise

$$(H + nK_X)^2 = D^2 + 2DE + E^2 \leq \max_{1 \leq x \leq m^2} \{x + 2m + m^2/x\} = (1 + m)^2$$

contradicts the assumption on $m$. Exchanging $D$ and $E$ allows me to assume $E^2 \leq 0$. Then $D^2 > 0$.

Set $d_1 = K_X D$ and $d_2 = K_X E$, so that $d = d_1 + d_2$ and $d_1$, $d_2 \geq 0$ (since $K_X$ is nef). Since $D^2 > 0$, the Hodge index theorem gives $d_1 > 0$.

If $d_2 \neq 0$ then $d_1 < d$, and I just replace $H \mapsto D$ and repeat the same argument with a smaller value of $d$. Since $d$ can only decrease a finite number of times, eventually every decomposition $D + E$ as in the statement of the conjecture satisfies $K_X E = 0$. Q.E.D.

### E.9.9. Special clusters

The suggestion for proving that $H + nK_X$ has a decomposition of the required type is in terms of Reider's method. I say that $Z \subset X$ is a *special cluster* of degree $m$ for $|H+(n+1)K_X|$ if $Z$ is a zero dimensional subscheme (for example $Z = P_1+\cdots+P_m$ with distinct $P_i \in X$) which are special in the sense that they impose dependent conditions on $|H+(n+1)K_X|$. In other words, the evaluation homomorphism at $Z$

$$H^0(H+(n+1)K_X) \to \mathcal{O}_Z \qquad (*)$$

is not surjective. $Z$ behaves like an analog of a special divisor on a general curve $C \in |H+nK_X|$.

The cokernel of $(*)$ is the group $H^1(I_Z \cdot \mathcal{O}_X(H+(n+1)K_X))$, which has Serre dual corresponding to extension classes

$$0 \to \mathcal{O}_X \to \mathcal{E} \to I_Z \cdot \mathcal{O}_X(H+nK_X) \to 0;$$

in most cases one can arrange that $\mathcal{E}$ is locally free, so a rank 2 vector bundle.

Now ${c_1}^2(\mathcal{E}) = (H+nK_X)^2$ and $c_2(\mathcal{E}) = m$, so if I take $m$ to be fairly small, $E$ has the Bogomolov numerical instability property $({c_1}^2 - 4c_2)(\mathcal{E}) > 0$. It's not hard to get from this to a decomposition as in the conjecture.

The hard problem is the existence of special clusters for suitable values of $n$ and $m$. One can set up a formalism of vector bundles over the Hilbert scheme of clusters, and define locuses $W_m^r$ of special clusters in analogy with the theory of special linear systems on curves. It's easy to do a dimension count in the style of Brill-Noether to prove that suitable locuses $W_m^r$ have dimension $> 0$, but I don't know how to do the intersection theory to prove that they are nonempty. (??)

### E.9.10

Although this proof runs into technical problems, it seems to me to provide at least a possible reason of principle why abundance may be true, not depending on frail numerical considerations: since $K_X$ is not ample, there must be lots of arbitrarily large cohomology groups around; if you can find one large group, e.g. $H^1(I_Z \cdot \mathcal{O}_X(H+nK_X))$, then it should be possible to chase this back to $H^0(nK_X)$.